\documentclass[a4paper,12pt]{article}
\pdfoutput=1
\usepackage{etex}

\usepackage{amssymb}
\usepackage{amsmath}
\usepackage{amsthm}
\usepackage{bbm}

\usepackage{lineno}
\usepackage[cm]{fullpage}
\usepackage{graphicx}
\usepackage{wrapfig}
\usepackage{subfig}
\usepackage{float}
\usepackage{booktabs}
\usepackage{caption}

\usepackage[inline]{asymptote}

\allowdisplaybreaks

\makeatletter
\newsavebox{\sfe@box}
\newenvironment{subfloatenv}[1]{%
\def\sfe@caption{#1}%
\setbox\sfe@box\hbox\bgroup\color@setgroup}%
{\color@endgroup\egroup\subfloat[\sfe@caption]%
{\usebox{\sfe@box}}}
\makeatother

\newcommand{\norm}[1]{\ensuremath{\lVert#1\rVert}}
\newcommand{\grade}[1]{\ensuremath{\langle#1\rangle}}

\newcommand{\reverse}[1]{\ensuremath{\widetilde{#1}}}
\newcommand{\involute}[1]{\ensuremath{\widehat{#1}}}

\newcommand{\R}[1]{\ensuremath{\mathbb{R}^{#1}}}
\newcommand{\T}[1]{\ensuremath{\mathbb{T}^{#1}}}
\newcommand{\E}[1]{\ensuremath{\mathbb{E}^{#1}}}
\newcommand{\M}[1]{\ensuremath{\mathbb{M}^{#1}}}

\newcommand{\El}[1]{\ensuremath{{E}_{#1}}}
\newcommand{\Hy}[1]{\ensuremath{{H}_{#1}}}

\newcommand{\dS}[1]{\ensuremath{{dS}_{#1}}}
\newcommand{\AdS}[1]{\ensuremath{{AdS}_{#1}}}

\newcommand{\tb}[1]{\ensuremath{\textbf{#1}}}
\newcommand{\mb}[1]{\ensuremath{\boldsymbol{#1}}}
\newcommand{\ts}[1]{\ensuremath{\mathsf{#1}}}

\newcommand{\e}{\tb{e}} 
\newcommand{\I}{\tb{I}} 
\newcommand{\Id}{\ensuremath{I}} 
\newcommand{\J}{\ensuremath{J}} 
\newcommand{\Sub}[1]{\ensuremath{\mathrm{Sub}({#1})}} 

\newcommand{\overbar}[1]{\mkern 1.5mu\overline{\mkern-1.5mu#1\mkern-1.5mu}\mkern 1.5mu}
\newcommand{\m}{\makebox[0pt][r]{\(-\)}}

\DeclareMathOperator{\tr}{tr}

\title{Clifford algebra and the projective model of\\ homogeneous metric spaces: Foundations}
\author{Andrey Sokolov}

\begin{document}
\maketitle

\tableofcontents

\newpage
\section{Introduction}
This paper is to serve as a key to the projective (homogeneous) model developed by Charles Gunn.
The goal is to explain the underlying concepts in a simple language and give plenty of examples.
It is targeted to physicists and engineers and the emphasis is on explanation rather than rigorous proof.
All results stated without proof are left as an exercise for the reader.
Gunn's model is described in his PhD thesis \cite{gunn2011geometry} and in his paper
 ``On the Homogeneous Model Of Euclidean Geometry'' published in
\emph {Guide to Geometric Algebra in Practice} ed. by {Dorst, L. and Lasenby, J.}, {2011}, {Springer} \cite{gunn2011homogeneous}.
The projective model supplements and enhances vector and matrix  algebras.
It also subsumes complex numbers and quaternions.

\begin{wrapfigure}[18]{l}{0.5\textwidth}
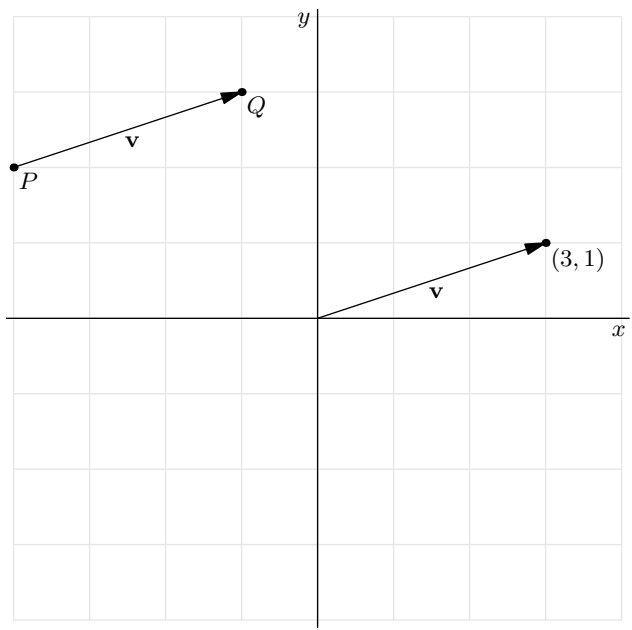

\vspace{-20pt}
\begin{asy}
import Drawing2D;
Drawing2D drawing = Drawing2D(4, 0.1);
drawing.grid();
drawing.target_axes();

pair P = (-4,2);
pair Q = (-1,3);

dot("$P$", P, SE);
dot("$Q$", Q, SE);
draw(Label("$\textbf{v}$",0.5), P--Q, Arrow);
draw(Label("$\textbf{v}$",0.5), (0,0)--(Q-P), Arrow);

dot("$(3,1)$", (3,1), SE);
drawing.crop();
\end{asy}
\caption{Points and vectors in \R{2}.}
\label{vectors in R2}
\vspace{-20pt}
\end{wrapfigure}

Consider points \(P\) and \(Q\) with the coordinates  \((-4,2)\) and \((-1,3)\), respectively (see Figure~\ref{vectors in R2}). 
The points define an oriented line segment \(\overrightarrow{PQ}\) extending from  \(P\) to \(Q\), which is called a vector.
Let \(\tb{v}=\overrightarrow{PQ}\).
In terms of the coordinates, \(\tb{v}\) is obtained by subtracting \(P\) from \(Q\), which gives \(\tb{v}=(3,1)\).
In this representation, there is nothing to distinguish vector \(\tb{v}\) from a point at \((3,1)\).
In fact, vectors are usually attached to the origin, so that 
vector \(\tb{v}\) can also be visualised as connecting the origin~\((0,0)\) and the point~\((3,1)\).

Viewed as oriented line segments that stem from the origin, 
vectors can be used to represent 1-dimensional subspaces of a vector space, i.e.~lines which pass through the origin.
Grassmann algebra is an extension of this idea to arbitrary subspaces, e.g.~planes and 3-dimensional hyperplanes 
passing through the origin. These subspaces are represented by multivectors.
Grassmann algebra uses the outer product, which does not require a metric.
When a metric is available, the inner product of two vectors can be defined, 
which in combination with the outer product gives rise to the geometric product and  Clifford algebra. 
Clifford algebra is the algebra of linear subspaces of a vector space.
Its geometric product captures various geometric relations between the subspaces.

An introduction to Clifford algebra, also known as geometric algebra, and numerous examples of its application
in physics and engineering can be found in the following books:

\emph{Geometric Algebra for Computer Science. An Object Oriented Approach to Geometry} \\
by {Dorst, L. Fontijne, D. and Mann, S.}, {2009}, {Elsevier/Morgan Kaufmann} \cite{dorst2009geometric}.

\emph{Geometric Algebra for Physicists}\\
by {Chris Doran C., Lasenby A.}, {2003}, {Cambridge University Press} \cite{d2003geometric}.

Clifford algebra is tied to linear subspaces, i.e.~lines, planes, etc which pass through the origin,
so it suffers from some of the same problems that afflict vector algebra.
Namely, locations in space are represented by vectors, 
fundamental geometric objects such as lines do not have a direct representation, and so on.
Both vector and Clifford algebras model translations, linear shifts in space,
by addition with a vector representing the shift.
However, this only works for vectors representing locations in space.
It fails for other geometric objects.

These problems are addressed in projective geometry.
The idea is to embed the target \(n\)-dimensional space \R{n}, whose geometry we want to model, 
into a vector space \R{n+1} with one extra dimension (the model space)
and represent geometric objects in the target space with linear subspaces of the model space.
To help distinguish target and model spaces, I will refer to the \(n\)-dimensional target space as \T{n}, rather than \R{n}.
The model space for \T{n} is an \((n+1)\)-dimensional space, which will be denoted by \R{n+1}.
The target space is initially non-metric (later on, when a metric is introduced, 
\T{} will be replaced with another symbol to indicate the metric space, e.g.~\E{} for Euclidean space).

The projective model uses projective geometry and Clifford algebra to model geometric objects and
their transformations in Euclidean and other homogeneous metric spaces.
The key feature of the projective model is its use of the dual space \T{n*}, which provides the top-down view of geometry in the target space.
The dual space \T{n*} is embedded in the dual model space \R{(n+1)*}.
The metric is defined on the dual model space \R{(n+1)*}, so it assumes a prominent role in the projective model.
Most algebraic manipulations are performed in the dual model space where the metric and, therefore, Clifford algebra are defined.

Projective geometry augmented with Clifford algebra provides a unified algebraic framework for describing points, lines, planes, etc,
and their transformations, such as rotations, reflections, projections, and translations.
Geometric relations between various objects are captured by the geometric product, 
facilitating efficient solution of geometric problems.

A modern introduction to projective geometry (without the use of Clifford algebra) can be found in \\
\emph {Perspectives on Projective Geometry: A Guided Tour Through Real and Complex Geometry} \\
by {Richter-Gebert, J.}, {2011}, {Springer} \cite{richtergebert2011perspectives}.

The following classic by Felix Klein provides an elementary introduction to many relevant issues\\
\emph{Elementary mathematics from an advanced standpoint. Geometry} by Felix Klein, 1939 \cite{klein1939elementary}.

\textbf{Brief outline of the projective model}\\
1) Define the dual space \T{n*} and discuss its relationship with the target space \T{n}\\
2) Explain the top-down and bottom-up views of geometry\\
3) Embed   \T{n} into the target model space \R{n+1} and
 \T{n*} into the dual model space \R{(n+1)*}\\
5) Model points, lines, and planes  as linear subspaces of the model spaces\\
6) Use multivectors in \R{n+1} and \R{(n+1)*} to represent linear subspaces \\
7) Establish an isomorphism between the target model space \R{n+1} and the dual model space \R{(n+1)*}\\
8) Define a metric on the dual model space \R{(n+1)*} \\
9) Use Clifford algebra of the dual model space \R{(n+1)*} to study geometry in the target space \T{n}

I will first develop the model in the 2-dimensional space since it already contains most of the essential features of the model
but is relatively simple.
This will be followed by the 3- and 4-dimensional spaces and a brief description of the 1-dimensional space.
To facilitate understanding, all 3D figures in this tutorial are interactive and can be activated within Acrobat Reader.

\newpage
\section{2-dimensional geometry}

The non-metric aspects of the projective model are based on projective geometry and Grassmann algebra.
The most relevant concepts from projective geometry are introduced in the next section.
Of these, the most important is the concept of projective duality and the related concept of the top-down view of geometry.

\subsection{Projective foundations}
\subsubsection{Projective duality}

In projective geometry, duality is defined in a way that does not require a metric.
In the target space \T{2}, a line that does not pass through the origin can be defined by the equation
\begin{equation}
ax+by+1=0,
\label{line in T2}
\end{equation}
where \(a\) and \(b\) are fixed, while \(x\) and \(y\) are variable and \((x,y)\) is in \T{2}. 
So, to define a line in \T{2} one needs to specify a pair of numbers \((a,b)\).
Let \(L\) be a line in \T{2} defined by \(-x-3y+1=0\).
The line is determined by specifying the coefficients \((a,b)=(-1,-3)\).
The dual space \T{2*} is a linear space of the coefficients \((a,b)\) that define lines in \T{2} via Equation~(\ref{line in T2}).
The meaning of \((a,b)=(0,0)\) will be clarified below.
Equation~(\ref{line in T2}) fixes a correspondence between lines in the target space \T{2} and points in the dual space \T{2*}.
For instance, \(L\)  corresponds to a point \(Q\) in \T{2*} with the coordinates \((-1,-3)\).
I will also say that \(Q\) is dual to \(L\) and write \(Q=L^*\).
Equally, \(L\) is dual to \(Q\), i.e.\ \(L=Q^*\).

\begin{figure}[h]
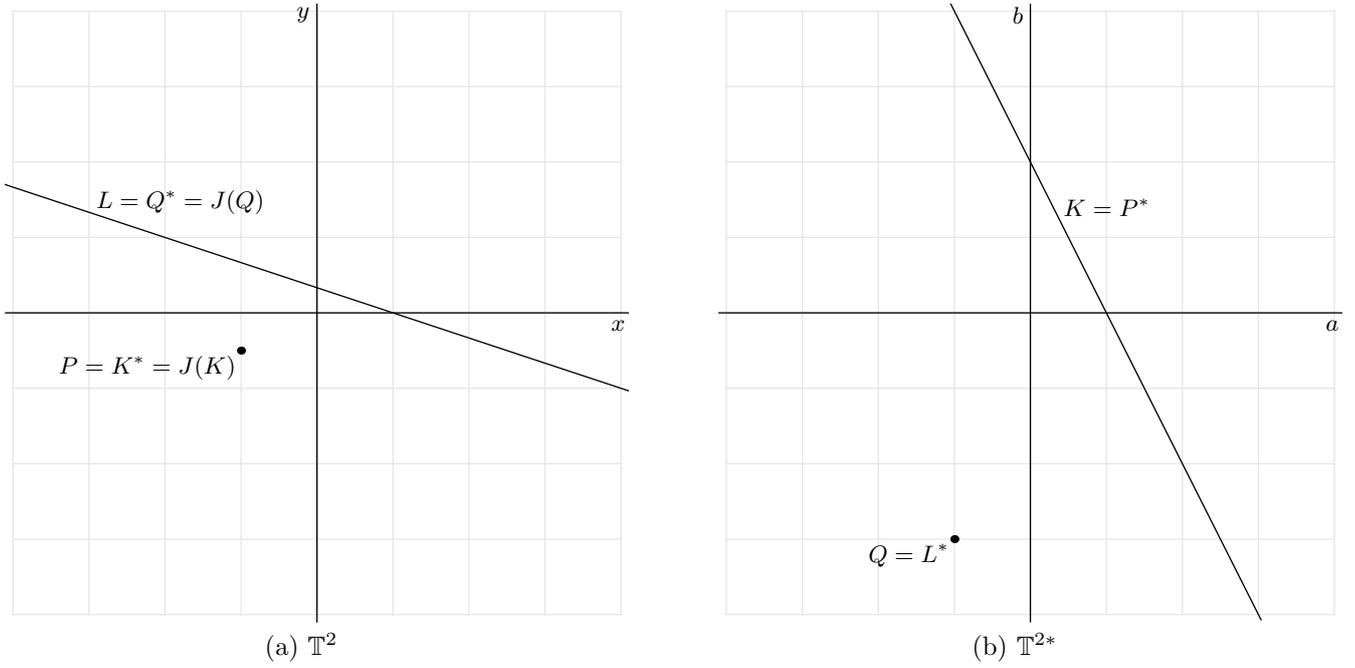

\begin{subfloatenv}{\T{2}}
\begin{asy}
import Drawing2D;
Drawing2D drawing = Drawing2D(4, 0.1);
drawing.grid();
drawing.target_axes();

path l = line(center=(0.1,0.3), direction=(-3,1), extent_from_center=10);
draw(Label("$L=Q^*=J(Q)$",0.6,(0,4)), l);

pair P = (-1,-1/2);
dot(Label("$P=K^*=J(K)$", align=SW), P);
drawing.crop();
\end{asy}
\end{subfloatenv}\hfill%
\begin{subfloatenv}{\T{2*}}
\begin{asy}
import Drawing2D;
Drawing2D drawing = Drawing2D(4, 0.1);
drawing.grid();
drawing.dual_axes();

path l = line(center=(0.8,0.4), direction=(-1,2), extent_from_center=5);
draw(Label("$K=P^*$",0.6), l);

dot(Label("$Q=L^*$", align=SW), (-1,-3));
drawing.crop();
\end{asy}
\end{subfloatenv}
\caption{Target and dual spaces and the mutual correspondence between points and lines}
\label{target and dual spaces}
\end{figure}

Equation~(\ref{line in T2}) is symmetric with respect to \((a,b)\) and \((x,y)\).
One can view \(x\) and \(y\) as fixed  and \(a\) and \(b\) as variable.
This defines a line in the dual space \T{2*}.
Let \(K\) be a line in \T{2*} defined by \(-a-\tfrac{1}{2}b+1=0\).
The line is determined by specifying a point \(P\) in \T{2} with the coordinates \((-1,-\tfrac{1}{2})\).
Consequently, the target space \T{2} can be viewed as 
a linear space of the coefficients \((x,y)\) that define lines in \T{2*} via  Equation~(\ref{line in T2}).
Points in the target space \T{2} correspond to lines in the dual space \T{2*}.
\(K\) is dual to \(P\), i.e\ \(K=P^*\), and equally \(P\) is dual to \(K\), i.e.\ \(P=K^*\)
(see Figure~\ref{target and dual spaces}).

The duality transformation \(\J\) is defined on points and lines in the dual space \T{2*} 
and gives the corresponding dual lines and points in the target space \T{2}.
For instance,
\(
\J(Q)=Q^*=L
\)
for the point \(Q\) in \T{2*} and the line \(L\) in \T{2}, 
\(
\J(K)=K^*=P
\)
for the line \(K\) in \T{2*} and the point \(P\) in \T{2}.
For the inverse transformation, I have \(\J^{-1}(P)=P^*=K\) and \(\J^{-1}(L)=L^*=Q\).

\begin{figure}[t!]
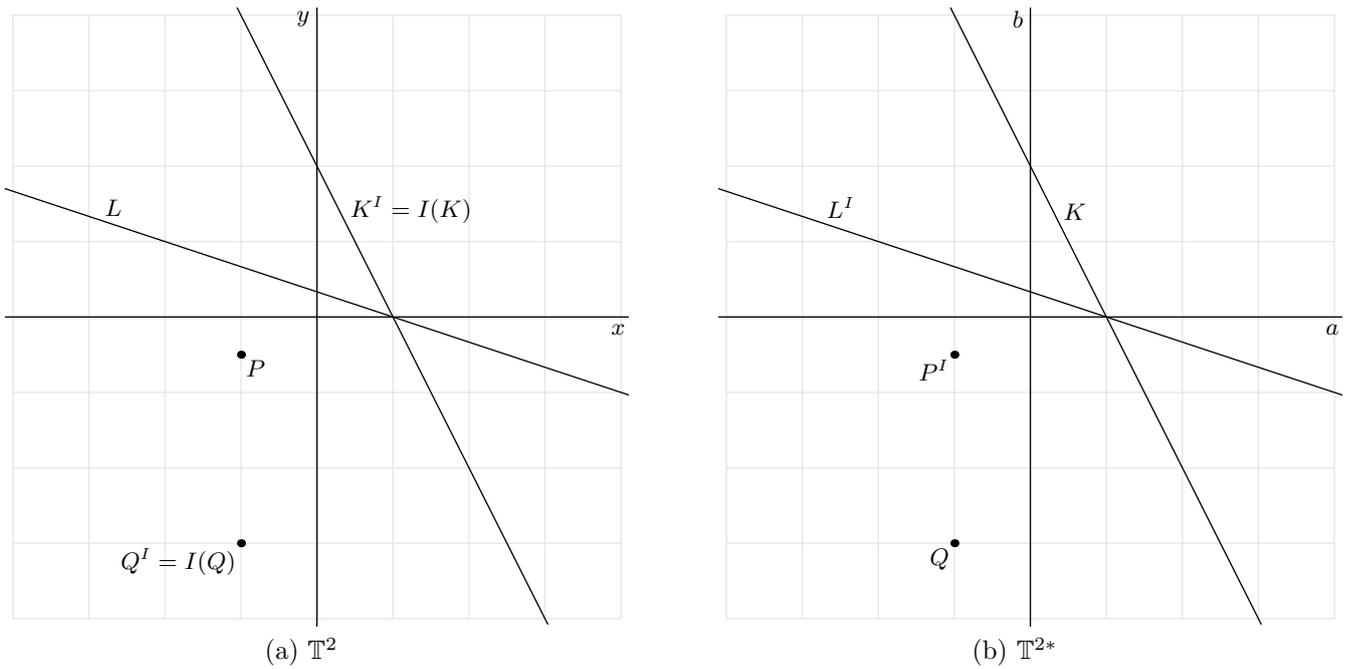

\begin{subfloatenv}{\T{2}}
\begin{asy}
import Drawing2D;
Drawing2D drawing = Drawing2D(4, 0.1);
drawing.grid();
drawing.target_axes();

path l = line(center=(0.1,0.3), direction=(-3,1), extent_from_center=5);
draw(Label("$L$",0.8), l);

pair P = (-1,-1/2);
dot(Label("$P$", align=SE), P);

path l = line(center=(0.8,0.4), direction=(-1,2), extent_from_center=5);
draw(Label("$K^I=I(K)$",0.6), l);

dot(Label("$Q^I=I(Q)$", align=SW), (-1,-3));

drawing.crop();
\end{asy}
\end{subfloatenv}\hfill%
\begin{subfloatenv}{\T{2*}}
\begin{asy}
import Drawing2D;
Drawing2D drawing = Drawing2D(4, 0.1);
drawing.grid();
drawing.dual_axes();

path l = line(center=(0.1,0.3), direction=(-3,1), extent_from_center=5);
draw(Label("$L^I$",0.8,NE), l);

pair P = (-1,-1/2);
dot(Label("$P^I$", align=SW), P);

path l = line(center=(0.8,0.4), direction=(-1,2), extent_from_center=5);
draw(Label("$K$",0.6), l);

dot(Label("$Q$", align=SW), (-1,-3));

drawing.crop();
\end{asy}
\end{subfloatenv}
\caption{Identity transformation}
\label{identity transform}
\end{figure}

For each point or line in \T{2}, there is an identical point or line in \T{2*}.
For instance, the point \(P^I=(a,b)=(-1,-\tfrac{1}{2})\) is identical to \(P=(x,y)=(-1,-\tfrac{1}{2})\) and
the line \(L^I=\{(a,b)|-a-3b+1=0\}\) is identical to \(L=\{(x,y)|-x-3y+1=0\}\).
The same applies to points and lines in \T{2*}.
For instance,  the point \(Q^I=(x,y)=(-1,-3)\)  is identical to \(Q=(a,b)=(-1,-3)\)
 and the line \(K^I=\{(x,y)|-x-\tfrac{1}{2}y+1=0\}\) is identical to \(K=\{(a,b)|-a-\tfrac{1}{2}b+1=0\}\).
The identity transformation \(\Id\) is defined on points and lines in the dual space \T{2*} 
and gives the identical points and lines in the target space \T{2}.
For instance, \(\Id(Q)=Q^I\) and \(\Id(K)=K^I\).
Both duality and identity transformations relate geometric objects of the dual space \T{2*} to
geometric objects of the target space \T{2}.
The duality transformation turns points into lines and lines into points,
whereas the identity transformation gives points for points and lines for lines.
The definition of the identity transformation is illustrated in Figure~\ref{identity transform}.

\begin{wrapfigure}[18]{L}{0.5\textwidth}
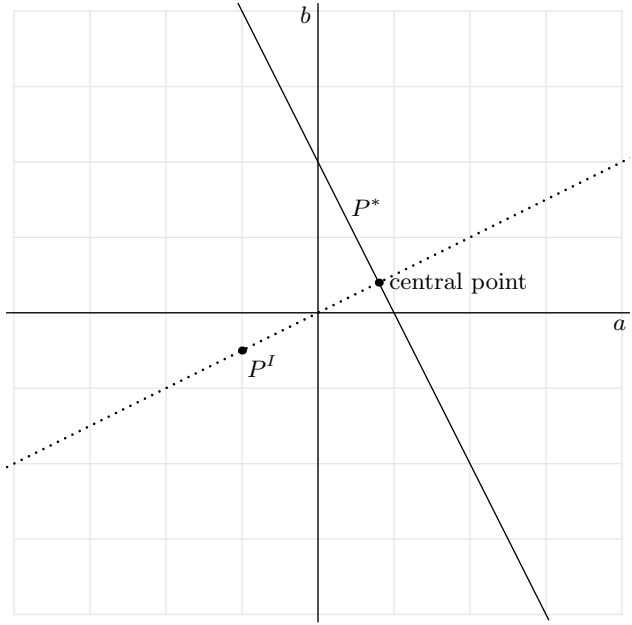

\vspace{-20pt}
\begin{center}
\begin{asy}
import Drawing2D;
Drawing2D drawing = Drawing2D(4, 0.1);
drawing.grid();
drawing.dual_axes();

pair P = (-1,-1/2);
dot(Label("$P^I$", align=SE), P);

path l = line(center=(0.8,0.4), direction=(-1,2), extent_from_center=5);
draw(Label("$P^*$",0.6), l);

path cl = line(center=(0,0), direction=(2,1), extent_from_center=5);
draw(cl, Dotted);

dot(Label("central point", align=(1.5,0)), (0.8,0.4));
drawing.crop();
\end{asy}
\caption{The central point of the line \(P^*\) in \T{2*}}
\label{central point}
\end{center}
\vspace{-20pt}
\end{wrapfigure}

Consider a point \(P\) in \T{2} and the corresponding line \(P^*\) in \T{2*}.
Draw a line in \T{2*} which passes through the point \(P^I\) and the origin of \T{2*}.
Its intersection with the line \(P^*\) will be called the central point of  \(P^*\)
(see Figure~\ref{central point} where \(P=(-1,-\tfrac{1}{2})\)).
For \(P=(x_P,y_P)\), the coordinates of the central point of the dual line \(P^*\) are given by 
\begin{equation}
a_c= \frac{-x_P}{x_P^2+y_P^2},\quad b_c = \frac{-y_P}{x_P^2+y_P^2}.
\end{equation}
These expressions are non-metric.
The central point of a line passing through the origin is at the origin by definition.
For \(P=(-1,-\tfrac{1}{2})\), the central point of \(P^*\) is at \(C(P^*)=(0.8,0.4)\). 
Identical construction can be carried out for points in \T{2*} and their dual lines in \T{2}.
For \(Q=(a_Q,b_Q)\), 
the coordinates of \(C(Q^*)\) are given by
\begin{equation}
x_c= \frac{-a_P}{a_Q^2+b_Q^2},\quad y_c = \frac{-b_Q}{a_Q^2+b_Q^2}.
\end{equation}
For \(Q=(-1,-3)\), the central point of the line \(Q^*\) is at \(C(Q^*)=(0.1,0.3)\). 

In Euclidean space, the line that passes through the point \(Q^I\) and the origin is perpendicular to the line \(Q^*\).
So that the central point is the point of the closest approach to the origin.
The distance from the origin to \(C(Q^*)\) equals
\(d_{C(Q^*)}=\sqrt{0.1^2+0.3^2}=1/\sqrt{10}\).
On the other hand, the distance from the origin to \(Q^I\) equals
\(d_{Q^I}=\sqrt{(-1)^2+(-3)^2}=\sqrt{10}\).
So, these distances are inversely related.
This is true in general in Euclidean space, i.e.~\(d_{C(Q^*)}d_{Q^I}=1\) for any finite point \(Q\), except the origin.
Note that in other metric spaces another point on the line \(Q^*\) could be closer to the origin.

Since Equation~(\ref{line in T2}) is symmetric with respect to \((a,b)\) and \((x,y)\),
the target and dual spaces play equal roles in projective geometry.
In fact, I could choose to refer to \T{2*} as the primary space and call \T{2} the dual space.

\subsubsection{Top-down view of geometry}
\begin{figure}[t!]
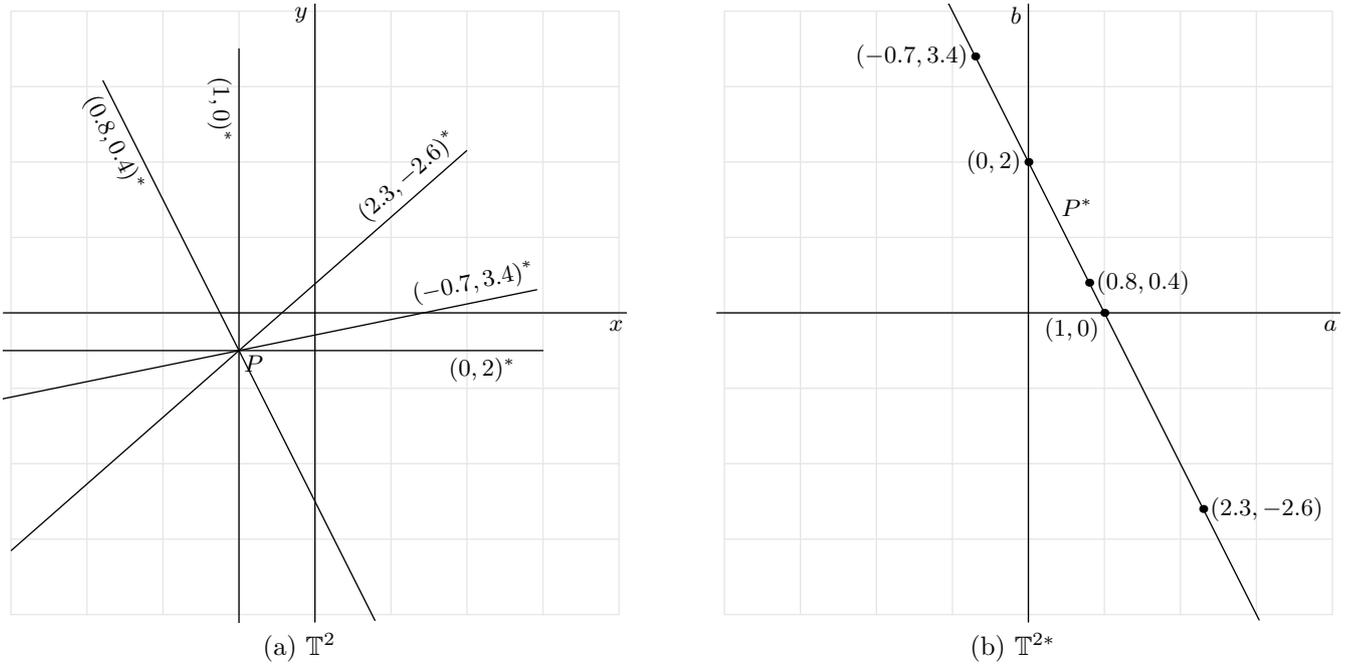

\begin{subfloatenv}{\T{2}}
\begin{asy}
import Drawing2D;
Drawing2D drawing = Drawing2D(4, 0.1);
drawing.grid();
drawing.target_axes();

pair n = normalise((-1,2));
pair c = (0.8,0.4);
labeled_pair[] points = points_on_line(center=c, direction=n, separation=1.5*sqrt(5), N=1);
labeled_pair[] extra_points = {labeled_pair((1,0),'(1,0)'), labeled_pair((0,2), '(0,2)')};
points.append(extra_points);

pair intersection = (-1,-1/2);
real[] relative_position = { 0.1, 0.1, 0.9, 0.1, 0.9};
align[] label_align = {NoAlign, NoAlign, (-0.2,1), NoAlign, NoAlign,};
for(int i=0; i<points.length; ++i) { 
  pair n = perp(points[i].pair);
  path l = line(center=intersection, direction=n, extent_from_center=4); 
  real th = degrees(n);
  if (i==0) { th -= 180; } 
  write(points[i].label, ' ', (string)th);
  draw(rotate(th)*Label('$'+points[i].label+'^*$', relative_position[i], align=label_align[i]), l); 
}
label('$P$', intersection, SE);
drawing.crop();
\end{asy}
\end{subfloatenv}\hfill%
\begin{subfloatenv}{\T{2*}}
\begin{asy}
import Drawing2D;
Drawing2D drawing = Drawing2D(4, 0.1);
drawing.grid();
drawing.dual_axes();

pair n = normalise((-1,2));
pair c = (0.8,0.4);
labeled_pair[] points = points_on_line(center=c, direction=n, separation=1.5*sqrt(5), N=1);
labeled_pair[] extra_points = {labeled_pair((1,0),'(1,0)'), labeled_pair((0,2), '(0,2)')};
points.append(extra_points);
align[] point_label_align = {E,E,W,SW,W};
for(int i=0; i<points.length; ++i) { dot(Label('$'+points[i].label+'$', align=point_label_align[i]), points[i].pair);}
path l = line(center=c, direction=n, extent_from_center=5); 
draw(Label('$P^*$',0.6), l);

drawing.crop();
\end{asy}
\end{subfloatenv}
\caption{Points on line \(P^*\) in \T{2*} and the corresponding lines in \T{2}}
\label{line sheaf}
\end{figure}

Consider a point \(P\) in \T{2} with the coordinates \((-1,-\tfrac{1}{2})\).
Its dual line \(P^*\) is defined by \(-a-\tfrac{1}{2}b+1=0\).
There is a relationship between the points that lie on \(P^*\) and the lines in \T{2} passing through the point \(P\).
To understand this relationship, consider the points 
\((1,0)\),
\((0,2)\),
\((0.8,0.4)\),
\((2.3,-2.6)\),
\((-0.7,3.4)\),
which all lie on the dual line \(P^*\).
Every point in this list corresponds to a line in \T{2},
e.g.\ \((0.8,0.4)\) corresponds to a line defined by \(0.8x+0.4y+1=0\).
Thus, I find five lines in \T{2} corresponding to the five points.
Plotting these lines reveals that they all intersect at the point \(P\) (see Figure~\ref{line sheaf}).
In fact, this is true for any point on \(P^*\), i.e.\ 
the dual line of any point on \(P^*\) passes through \(P\).

The dual lines of the points that lie on  \(P^*\)  form a sheaf\footnote{It is also called a pencil of lines attached to a finite point \(P\).}
 of lines in the target space \T{2}, 
in which every line passes through the point \(P\),
with the exception of the line  that also passes through  the origin.
Thus, a point  can be viewed as a sheaf of lines, which gives rise to the top-down view of geometry.
In the top-down view of \T{2}, lines are considered fundamental and 
lower-dimensional objects, i.e.\ points, are build from lines as their intersection.
So, a point is viewed as a derivative object represented by a sheaf of lines passing through the point.
This contrasts with the bottom-up view of geometry, in which points are 
fundamental and higher-dimensional objects, i.e.\ lines, are build from points as their union or extension.
It takes at least two points to define a line in the bottom-up view.
Similarly, it takes at least two lines to define a point in the top-down view.

In view of the above, it can be said that the dual space \T{2*} enables the top-down view of geometry in the target space \T{2}.
Similar construction can be carried out for points in the dual space \T{2*}.
So, it can equally be said that the target space \T{2} enables the top-down view of the dual space \T{2*}.

\subsubsection{Points at infinity}
A point at infinity\footnote{
This terminology is appropriate for Euclidean space, since any \emph{point at infinity} in \E{2} is at the infinite distance from the origin.
In other metric spaces, \emph{points at infinity} can be at a finite distance from the origin.
However, \emph{points at infinity} correspond to stacks of lines rather than sheaves, regardless of the metric.} 
 in \T{2} is represented by a stack\footnote{It can also be called a pencil of lines attached to a point at infinity.} of lines.
A stack of lines is a set of lines in \T{2}, in which every line is defined by 
\begin{equation}
ax+by+d=0,
\end{equation}
where \(a\) and \(b\) are fixed and \(d\) spans all possible values in \R{}.
For instance, the lines defined by \(2x-y+\tfrac{1}{2}=0\), \(2x-y+1=0\), and \(2x-y-2=0\) belong to the same stack.
Identical definition applies in \T{2*}.
Note that in Euclidean space, a stack of lines consists of parallel lines.

\begin{figure}[h]
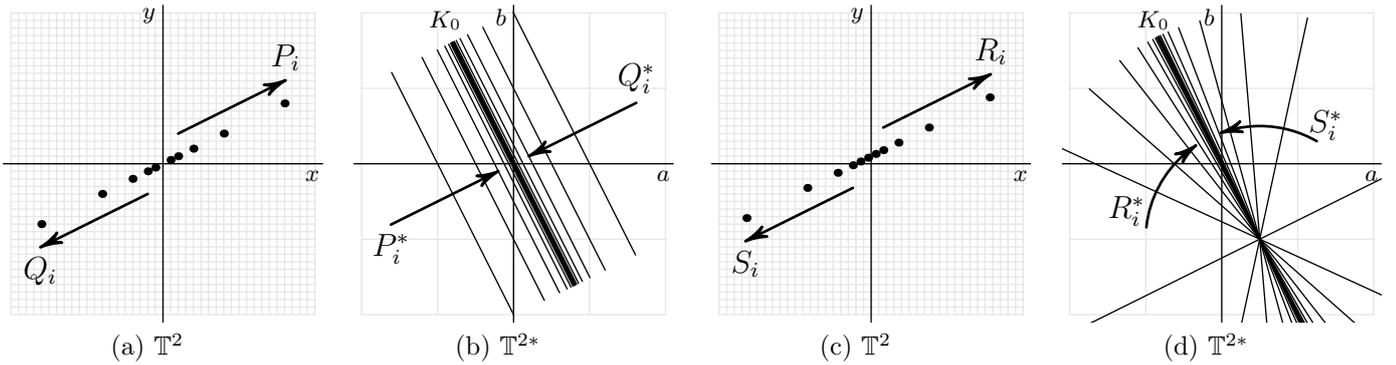

\begin{subfloatenv}{\T{2}}
\begin{asy}
import Drawing2D;
real unit = 0.1;
real norm = 1/unit;

Drawing2D drawing = Drawing2D(2, 0.1, unit_size=unit);
drawing.grid();
drawing.target_axes();

pair n = (1,1/2);
labeled_pair[] points = points_on_line(center=(0,0), direction=n, separation=0.1, N=20);

pair[] points;
int number_of_points = 10;
for(int i=0; i<number_of_points; ++i) { pair lp = (2^i)*n; points.push(lp); points.push(-lp);}

for(pair p: points) { dot(p/norm);}

pair shift = (0.2, 0.4);
path direction_to_infinity = shift--(sqrt(2)*n + shift);
draw(Label("$P_i$", 1, align=N), direction_to_infinity, linewidth(1 bp), Arrow(HookHead, size=5));

path opposite_direction_to_infinity = -shift--(sqrt(2)*(-n) - shift);
draw(Label("$Q_i$", 1, align=S), opposite_direction_to_infinity, linewidth(1 bp), Arrow(HookHead, size=5));

drawing.crop();
\end{asy}
\end{subfloatenv}%
\begin{subfloatenv}{\T{2*}}
\begin{asy}
import Drawing2D;
real unit = 1.0;//0.1;
Drawing2D drawing = Drawing2D(2, 0.1, unit_size=unit);
drawing.grid();
drawing.dual_axes();

pair n = (1,1/2);

pair[] points;
int number_of_points = 10;
for(int i=0; i<number_of_points; ++i) { pair lp = (2^i)*n; points.push(lp); points.push(-lp);}

for(pair p: points) { path l = line(center=-p/length(p)^2, direction=perp(p), extent_from_center=1.8); draw(l); }

pair shift = (0.2, 0.1);
path direction_to_infinity = ((sqrt(2)*n + shift)--shift);
draw(Label("$Q_i^*$",0, align=N), direction_to_infinity, linewidth(1 bp), Arrow(HookHead, size=5));

path opposite_direction_to_infinity = ((sqrt(2)*(-n) + -shift)--(-shift));
draw(Label("$P_i^*$",0, align=S), opposite_direction_to_infinity, linewidth(1 bp), Arrow(HookHead, size=5));

label("$K_0$",(-0.9,1.9));
drawing.crop();
\end{asy}
\end{subfloatenv}
\begin{subfloatenv}{\T{2}}
\begin{asy}
import Drawing2D;
real unit = 0.1;
real norm = 1/unit;

Drawing2D drawing = Drawing2D(2, 0.1, unit_size=unit);
drawing.grid();
drawing.target_axes();

pair n = (1,1/2);

pair[] points;
int number_of_points = 10;
pair center = (-2,4)/5;
points.push(center);
for(int i=0; i<number_of_points; ++i) { pair lp = (2^i)*n; points.push(lp+center); points.push(-lp+center);}

for(pair p: points) { dot(p/norm);}

pair shift = (0.2, 0.4);
path direction_to_infinity = (shift+center/norm)--(sqrt(2)*n + shift+center/norm);
draw(Label("$R_i$", 1, align=N), direction_to_infinity, linewidth(1 bp), Arrow(HookHead, size=5));

path opposite_direction_to_infinity = (-shift+center/norm)--(sqrt(2)*(-n) - shift + center/norm);
draw(Label("$S_i$", 1, align=S), opposite_direction_to_infinity, linewidth(1 bp), Arrow(HookHead, size=5));

drawing.crop();
\end{asy}
\end{subfloatenv}%
\begin{subfloatenv}{\T{2*}}
\begin{asy}
import Drawing2D;
real unit = 1.0;//0.1;
Drawing2D drawing = Drawing2D(2, 0.1, unit_size=unit);
drawing.grid();
drawing.dual_axes();

pair n = (1,1/2);

pair[] points;
int number_of_points = 10;
pair center = (-2,4)/5;
points.push(center);
for(int i=0; i<number_of_points; ++i) { pair lp = (2^i)*n; points.push(lp+center); points.push(-lp+center);}

pair c = (1/2,-1);
for(pair p: points) { path l = line(center=c, direction=perp(p), extent_from_center=3); draw(l); }

path cw1 = arc(c, 1.5, 174, 124, CW);
draw(Label("$R_i^*$",0,align=(-0.15,0.15)), cw1, linewidth(1 bp), Arrow(HookHead, size=5));

path ccw1 = arc(c, 1.5, 60, 110, CCW);
draw(Label("$S_i^*$",0,align=(0.05,0.1)), ccw1, linewidth(1 bp), Arrow(HookHead, size=5));

label("$K_0$",(-0.9,1.9));
drawing.crop();
\end{asy}
\end{subfloatenv}
\caption{Approaching a point at infinity in \T{2} and the corresponding line \(K_0\) in \T{2*}}
\label{approaching point at infinity}
\end{figure}

To gain some intuition about points at infinity, consider the following example.
Let \(L_d\) be a line in \T{2} defined by \(\tfrac{1}{2}x-y+d=0\), where \(d\in\R{}\), and 
\(K_0\) be a line in \T{2*} defined by \(2a+b=0\).
In the target space \T{2}, consider a series of points  
\(\{P_i\}_{i=0}^{\infty}\),
where \(P_i=2^i(1,\tfrac{1}{2})\).
The points of the series lie on the line \(L_0\), i.e.\ \(L_d\) with \(d=0\),
and the series  approaches infinity as \(i\to\infty\).
The dual lines \(P_i^*\) converge to the line \(K_0\)  as \(i\to\infty\) (see Figure~\ref{approaching point at infinity}(ab)).
Another series of points   \(\{Q_i\}_{i=0}^{\infty}\), where \(Q_i=2^i(-1,-\tfrac{1}{2})\),
consists of the points that lies on \(L_0\) as well but it approaches infinity in the opposite direction.
However,  the dual lines \(Q_i^*\) converge to the same line \(K_0\) as \(i\to\infty\).
Now consider series 
\(\{R_i\}_{i=0}^{\infty}\) and \(\{S_i\}_{i=0}^{\infty}\), where \(R_i=2^i(1,\tfrac{1}{2})+(-\tfrac{2}{5},\tfrac{4}{5})\)
and 
\( S_i=2^i(-1,-\tfrac{1}{2})+(-\tfrac{2}{5},\tfrac{4}{5})\).
Points \(R_i\) and \(S_i\) lie on the line 
\(L_1\),  i.e.\ \(L_d\) with \(d=1\),
which belongs to the same stack of lines as \(L_0\).
The dual lines \(R_i^*\) and \(S_i^*\) intersect at the point \((\tfrac{1}{2},-1)\), which is dual to \(L_1\),
and both \(R_i^*\) and \(S_i^*\) converge to \(K_0\) as \(i\to\infty\)  (see Figure~\ref{approaching point at infinity}(cd)).
This result holds in general for any line \(L_d\) from the stack \(\{L_d|d\in\R{}\}\), i.e.\ 
 approaching infinity along any of the lines in the stack yields the same line \(K_0\) in \T{2*}.
It is natural to conclude that the line \(K_0\) is dual a point at infinity in \T{2}
which can be approached by moving along any of the lines in the stack.

\begin{figure}[h]
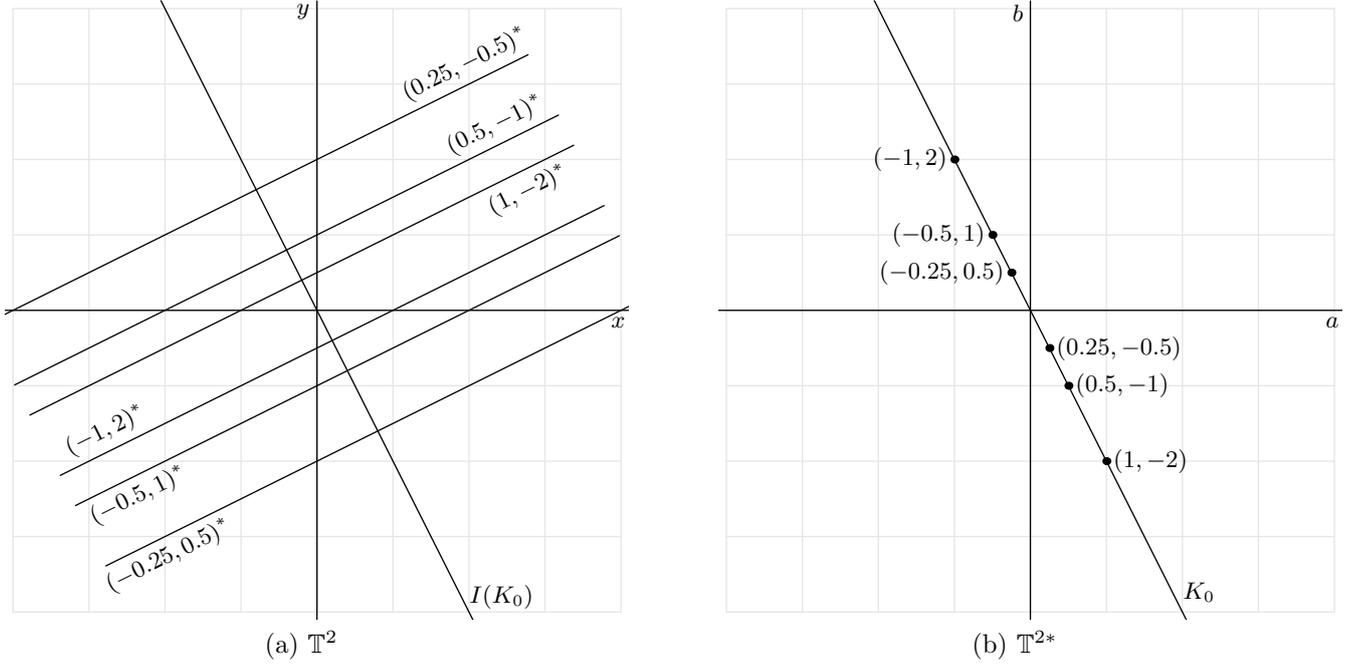

\begin{subfloatenv}{\T{2}}
\begin{asy}
import Drawing2D;
Drawing2D drawing = Drawing2D(4, 0.1);
drawing.grid();
drawing.target_axes();

pair n = normalise((-1,2));
labeled_pair[] points = points_on_line(center=(0,0), direction=n, separation=sqrt(5)/2, N=2);

int index_of_origin = quotient(points.length-1, 2);
points.delete(index_of_origin); // delete (0,0)

labeled_pair[] extra_points = {labeled_pair((-0.25,0.5),'(-0.25,0.5)'), labeled_pair((0.25,-0.5), '(0.25,-0.5)')};
points.append(extra_points);

real[] relative_position = { 0.1, 0.1, 0.1, 0.1, 0.1, 0.1};
align[] label_align = {(0.5,-1), NoAlign, NoAlign, (-0.5,1), NoAlign, NoAlign};
for(int i=0; i<points.length; ++i) { 
  pair n = perp(points[i].pair);
  real a = points[i].pair.x;
  real b = points[i].pair.y;
  pair intersection = -1/5 * (1/a, 4/b);
  path l = line(center=intersection, direction=n, extent_from_center=4); 
  real th = degrees(n);
  if (th > 180) { th -= 180; }
  write(points[i].label, ' ', (string)th);
  draw(rotate(th)*Label('$'+points[i].label+'^*$', relative_position[i], align=label_align[i]), l); 
}

path l = line(center=(0,0), direction=n, extent_from_center=5); 
draw(Label('$\Id(K_0)$',0.07), l);

drawing.crop();
\end{asy}
\end{subfloatenv}\hfill%
\begin{subfloatenv}{\T{2*}}
\begin{asy}
import Drawing2D;
Drawing2D drawing = Drawing2D(4, 0.1);
drawing.grid();
drawing.dual_axes();

pair n = normalise((-1,2));
labeled_pair[] points = points_on_line(center=(0,0), direction=n, separation=sqrt(5)/2, N=2);

int index_of_origin = quotient(points.length-1, 2);
points.delete(index_of_origin); // betel (0,0)

labeled_pair[] extra_points = {labeled_pair((-0.25,0.5),'(-0.25,0.5)'), labeled_pair((0.25,-0.5), '(0.25,-0.5)')};
points.append(extra_points);

align[] point_label_align = {E,E,W,W,W,E};
for(int i=0; i<points.length; ++i) { dot(Label('$'+points[i].label+'$', align=point_label_align[i]), points[i].pair);}

path l = line(center=(0,0), direction=n, extent_from_center=5); 
draw(Label('$K_0$',0.07), l);

drawing.crop();
\end{asy}
\end{subfloatenv}
\caption{Points on line \(K_0\) in \T{2*} and the corresponding stack of lines in \T{2}}
\label{stack of lines}
\end{figure}

There is a relationship between the points that lie on \(K_0\) and the lines in the stack \(\{L_d|d\in\R{}\}\).
Consider the points \((-1,2), (-0.5,1), (-025,0.5), (0.25,-0.5), (0.5,-1)\), and \((1,-2)\), 
which all lie on the line \(K_0\). 
Every point in this list is dual to a line in \T{2}.
Plotting the dual lines shows that they belongs to the stack \(\{L_d|d\in\R{}\}\) (see Figure~\ref{stack of lines}). 
The same hold for any point on \(K_0\), i.e.\
the dual line to any point on \(K_0\) belongs to the same stack.
So, the stack \(\{L_d|d\in\R{}\}\) consists of the lines dual to the points on \(K_0\),
with the exception of the line  that passes through the origin.

The same construction can be carried out for any stack in \T{2}, resulting in different points at infinity in \T{2}.
A stack of lines represents a single points at infinity in \T{2}.
This point can be approached by moving to infinity along any of the lines in the stack, in either direction along each line.
In the top-down view, points at infinity are not fundamentally different from the finite points.
One can say informally that the lines in the stack intersect at a point at infinity.

Exactly the same reasoning applies in the dual space \T{2*}, 
where points at infinity are represented by stacks of lines in \T{2*},
which in turn correspond to lines in \T{2} passing through the origin.

\subsubsection{Extending lines with points at infinity}
\begin{figure}[h]
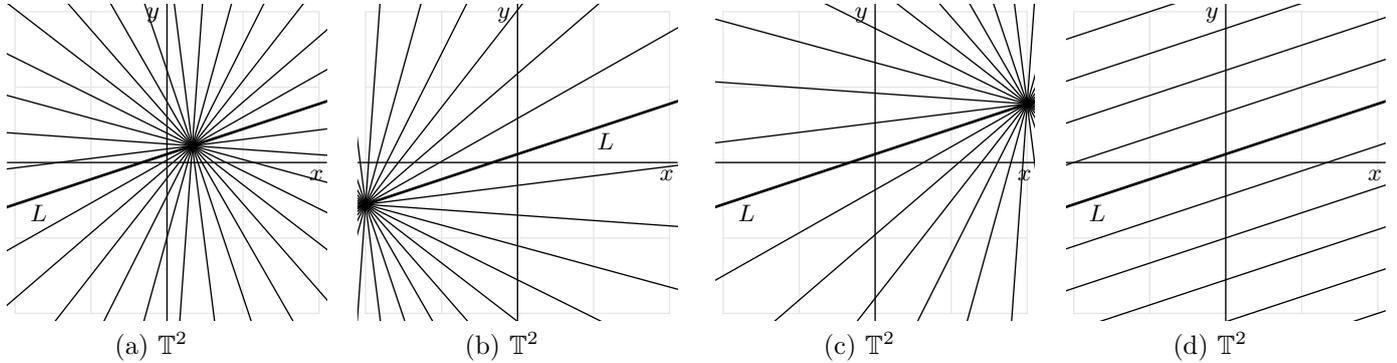

\begin{subfloatenv}{\T{2}}
\begin{asy}
import Drawing2D;
Drawing2D drawing = Drawing2D(2, 0.1);
drawing.grid();
drawing.target_axes();

pair n = (1,1/3);
pair c = (-1/3,0);
path l = line(center=c, direction=n, extent_from_center=5); 
draw(Label("$L$", 0.35), l, currentpen+1);

pair[] ns = {};
int N=16;
real phi0 = angle(n);
for(int i=1; i<N; ++i) { pair n=expi(i*pi/N + phi0); ns.push(n); }

pair P = (1/3,2/9);
for(pair n: ns) { path l = line(center=P, direction=n, extent_from_center=5);  draw(l); }

drawing.crop();
\end{asy}
\end{subfloatenv}%
\begin{subfloatenv}{\T{2}}
\begin{asy}
import Drawing2D;
Drawing2D drawing = Drawing2D(2, 0.1);
drawing.grid();
drawing.target_axes();

pair n = (1,1/3);
pair c = (-1/3,0);
path l = line(center=c, direction=n, extent_from_center=5); 
draw(Label("$L$", 0.65), l, currentpen+1);

pair[] ns = {};
int N=16;
real phi0 = angle(n);
for(int i=1; i<N; ++i) { pair n=expi(i*pi/N + phi0); ns.push(n); }

pair P = (-2,-5/9);
for(pair n: ns) { path l = line(center=P, direction=n, extent_from_center=5);  draw(l); }

drawing.crop();
\end{asy}
\end{subfloatenv}
\begin{subfloatenv}{\T{2}}
\begin{asy}
import Drawing2D;
Drawing2D drawing = Drawing2D(2, 0.1);
drawing.grid();
drawing.target_axes();

pair n = (1,1/3);
pair c = (-1/3,0);
path l = line(center=c, direction=n, extent_from_center=5); 
draw(Label("$L$", 0.35), l, currentpen+1);

pair[] ns = {};
int N=16;
real phi0 = angle(n);
for(int i=1; i<N; ++i) { pair n=expi(i*pi/N + phi0); ns.push(n); }

pair P = (2,7/9);
for(pair n: ns) { path l = line(center=P, direction=n, extent_from_center=5);  draw(l); }

drawing.crop();
\end{asy}
\end{subfloatenv}%
\begin{subfloatenv}{\T{2}}
\begin{asy}
import Drawing2D;
Drawing2D drawing = Drawing2D(2, 0.1);
drawing.grid();
drawing.target_axes();

pair n = (1,1/3);
pair c = (-1/3,0);
path l = line(center=c, direction=n, extent_from_center=5); 
draw(Label("$L$", 0.35), l, currentpen+1);

pair[] ps = {};
int N=10;
for(int i=1; i<N; ++i) { pair P=c+(i/3,-i)/2; ps.push(P); }
for(int i=1; i<N; ++i) { pair P=c+(-i/3,i)/2; ps.push(P); }

for(pair P: ps) { path l = line(center=P, direction=n, extent_from_center=5);  draw(l); }

drawing.crop();
\end{asy}
\end{subfloatenv}
\caption{Line \(L\) and its point at infinity as a stack of lines}
\label{line and its point at infinity}
\end{figure}

Points  can be viewed as sheaves of lines. 
The sheaf of a point which lies on the line \(L\) defined by \(-x-3y+1=0\) is shown in Figure~\ref{line and its point at infinity}(a).
If I move the point along the line and away from the origin, 
the appearance of its sheaf near the origin starts to resemble the stack of lines
that contains \(L\) (see Figure~\ref{line and its point at infinity}(bc)).
So, it is natural to consider the point at infinity
represented by the stack 
(see Figure~\ref{line and its point at infinity}(d)) as belonging to the line \(L\).
I will also say that this point at infinity lies on \(L\).

From now on, every line in \T{2} is extended to infinity by adding the relevant point at infinity
to the set of finite points that constitute the line.
Exactly the same applies in \T{2*}.

\emph{Completing the definition of a sheaf of lines}.\\
According to the description of a sheaf of lines I gave previously,
a sheaf consists of all lines that pass through a single point, except for the line that also passes through the origin.
Indeed, the missing line does not correspond to any finite point on the line dual to the point defined by the sheaf.
However, I have just extended the dual line by a point at infinity. 
This point corresponds to a line in the target space which passes through the point defined by the sheaf and
through the origin.
In other words, the point at infinity on the dual line corresponds to the missing line.
Since, the point at infinity lies on the dual line, it is natural to add the missing line
to the definition of a sheaf.
Thanks to this amendment, a sheaf can be said to consist of all lines passing through a single point, without exception.
This completes the definition of a sheaf  given previously.

\subsubsection{The line at infinity}
\begin{figure}[ht]
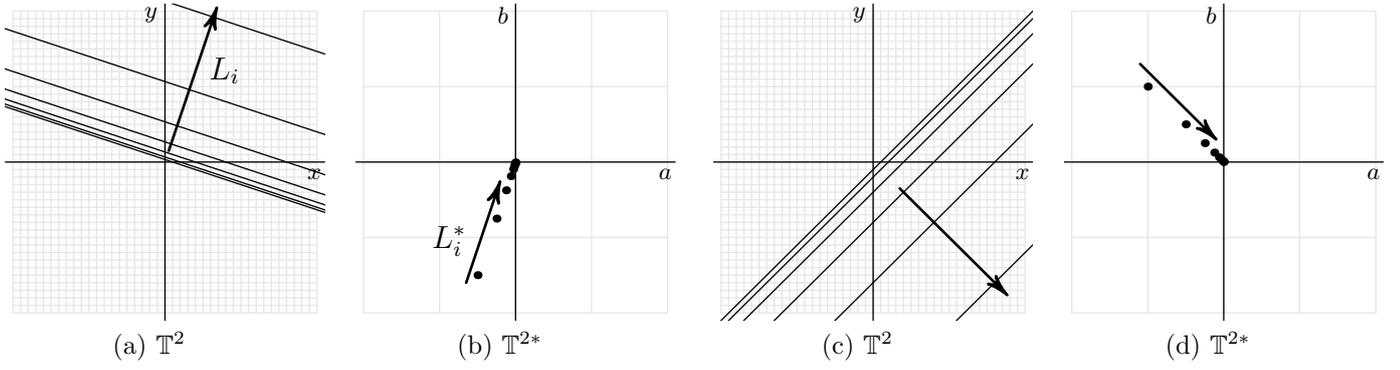

\begin{subfloatenv}{\T{2}}
\begin{asy}
import Drawing2D;
real unit = 0.1;
real norm = 1/unit;

Drawing2D drawing = Drawing2D(2, 0.1, unit_size=unit);
drawing.grid();
drawing.target_axes();

pair n = normalise((1,-1/3));
int number_of_lines = 10;
for(int i=0; i<number_of_lines; ++i) { path l = line(center=(0,2^i/3)/norm, direction=n, extent_from_center=5); draw(l); }

pair shift = (1/2, 3/2)/norm;
path direction_to_infinity = shift--(-2*perp(n) + shift);
draw(Label("$L_i$", 0.7, align=SE), direction_to_infinity, linewidth(1 bp), Arrow(HookHead, size=5));

drawing.crop();
\end{asy}
\end{subfloatenv}%
\begin{subfloatenv}{\T{2*}}
\begin{asy}
import Drawing2D;
real unit = 1.0;//0.1;
Drawing2D drawing = Drawing2D(2, 0.1, unit_size=unit);
drawing.grid();
drawing.dual_axes();

int number_of_points = 10;
for(int i=0; i<number_of_points; ++i) { pair P  = (-1/2^i,-3/2^i); dot(P);}

pair n = normalise((1,-1/3));
pair shift = -(1, 2) + (0.35,0.4);
path direction_to_infinity = shift--(-sqrt(2)*perp(n) + shift);
draw(Label("$L_i^*$", 0.2, align=NW), direction_to_infinity, linewidth(1 bp), Arrow(HookHead, size=5));

drawing.crop();
\end{asy}
\end{subfloatenv}
\begin{subfloatenv}{\T{2}}
\begin{asy}
import Drawing2D;
real unit = 0.1;
real norm = 1/unit;

Drawing2D drawing = Drawing2D(2, 0.1, unit_size=unit);
drawing.grid();
drawing.target_axes();

pair n = normalise((1,1));
int number_of_lines = 10;
for(int i=0; i<number_of_lines; ++i) { path l = line(center=(0,-2^i)/norm, direction=n, extent_from_center=5); draw(l); }

pair shift = (3.5, -3.5)/norm;
path direction_to_infinity = shift--(2*perp(n) + shift);
draw(direction_to_infinity, linewidth(1 bp), Arrow(HookHead, size=5));

drawing.crop();
\end{asy}
\end{subfloatenv}%
\begin{subfloatenv}{\T{2*}}
\begin{asy}
import Drawing2D;
real unit = 1.0;//0.1;
Drawing2D drawing = Drawing2D(2, 0.1, unit_size=unit);
drawing.grid();
drawing.dual_axes();

int number_of_points = 10;
for(int i=0; i<number_of_points; ++i) { pair P  = (-1/2^i, 1/2^i); dot(P);}

pair n = normalise((1,1));
pair shift = (-1, 1) + (-0.1,0.3);
path direction_to_infinity = shift--(sqrt(2)*perp(n) + shift);
draw(direction_to_infinity, linewidth(1 bp), Arrow(HookHead, size=5));

drawing.crop();
\end{asy}
\end{subfloatenv}
\caption{Approaching the line at infinity in \T{2} and the origin in \T{2*}}
\label{line at infinity}
\end{figure}

If \((a,b)=(0,0)\), none of the finite points \((x,y)\in\T{2}\) satisfies Equation~(\ref{line in T2}).
Consider the following sequence of lines in the target space \T{2}:
\begin{equation}
\{L_i\}_{i=0}^{\infty}, \textrm{where } L_i=\{(x,y)|-x-3y+2^i=0\}.
\end{equation}
The sequence approaches infinity as \(i\to\infty\).
In other words, the limit of this sequence is the line at infinity.
\(L_i\) corresponds to the point \(L_i^*\) in \T{2*}.
As \(i\to\infty\), the sequence of points \(\{L_i^*\}_{i=0}^{\infty}\) 
converges to the origin \((0,0)\) of \T{2*} (see Figure~\ref{line at infinity}(ab)).
The same is true for any sequence of lines approaching the line at infinity in \T{2} 
(see Figure~\ref{line at infinity}(cd) for another sequence).

So, the origin of \T{2*} corresponds to the line at infinity in \T{2}.
Similarly, the origin of the target space \T{2} corresponds to the line at infinity in \T{2*}.

\begin{figure}[ht]
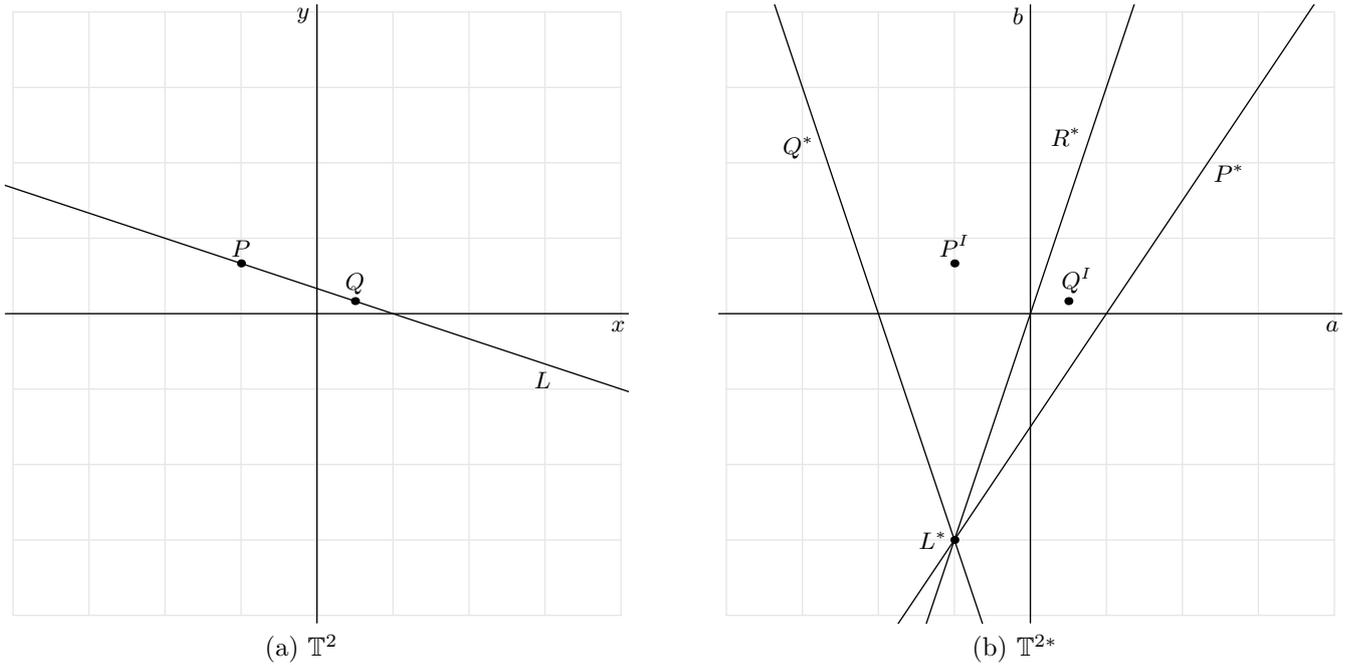

\begin{subfloatenv}{\T{2}}
\begin{asy}
import Drawing2D;
Drawing2D drawing = Drawing2D(4, 0.1);
drawing.grid();
drawing.target_axes();

pair n = normalise((1,-1/3));
pair c = (0,1/3);
path l = line(center=c, direction=n, extent_from_center=8); 
draw(Label("$L$",0.7),l);

pair P1 = (-1,2/3);
dot(Label("$P$",align=N), P1);

pair P2 = (1/2,1/6);
dot(Label("$Q$",align=N), P2);

drawing.crop();
\end{asy}
\end{subfloatenv}\hfill%
\begin{subfloatenv}{\T{2*}}
\begin{asy}
import Drawing2D;
Drawing2D drawing = Drawing2D(4, 0.1);
drawing.grid();
drawing.dual_axes();

pair l_dual = (-1,-3);
dot(Label("$L^*$",align=W), l_dual);

pair P1 = (-1,2/3);
pair P2 = (1/2,1/6);

path P1_dual = line(center=(1,0), direction=perp(P1), extent_from_center=8); 
draw(Label("$P^*$",0.65), P1_dual);

path P2_dual = line(center=(-2,0), direction=perp(P2), extent_from_center=8); 
draw(Label("$Q^*$",0.35), P2_dual);

path R_dual = line(center=(0,0), direction=perp((3,-1)), extent_from_center=8); 
draw(Label("$R^*$",0.35), R_dual);

dot(Label("$P^I$",align=N), P1);
dot(Label("$Q^I$",align=N+(0.4,0.1)), P2);

drawing.crop();
\end{asy}
\end{subfloatenv}
\caption{The incidence relation between target and dual spaces}
\label{incidence relation}
\end{figure}

\emph{Incidence of points and lines in the target and dual spaces}.\\
In the target space \T{2}, consider a finite point \(P\ne(0,0)\) and a line \(L\) that does not pass through the origin.
If \(P\) lies on \(L\), then the point \(L^*\) lies on the line \(P^*\) in \T{2*}.
To understand this, consider the following reasoning.
Since \(L\) passes through \(P\), it belongs to the sheaf of lines around \(P\).
Each line in the sheaf, including \(L\), corresponds to a point in the dual space \T{2*}
which lies on the line \(P^*\).
Therefore, the point \(L^*\) also lies on \(P^*\).
This relation between the target and dual spaces is illustrated in Figure~\ref{incidence relation}
for \(P=(-1,\tfrac{2}{3})\), \(Q=(\tfrac{1}{2},\tfrac{1}{6})\), and \(L=\{(x,y)|-x-3y+1=0\}\).
Similar reasoning applies to points at infinity.
Consider a point \(R\) at infinity in \T{2}, whose stack includes \(L\).
Each line in the stack, including \(L\), corresponds to a point in the dual space \T{2*}
which lies on the line \(R^*\).
Therefore, the point \(L^*\) also lies on \(R^*\).

\emph{Completing the definition of a stack of lines}.\\
My previous description of a stack of lines is missing a line passing through the origin.
This missing line corresponds to a point at infinity which lies on the dual line that defines the stack.
So, it is natural to include it in the description of the stack.
Moreover, the dual line passes through the origin of the dual space, but the origin corresponds to the line at infinity in the target space.
So, I include the line at infinity in the definition of the stack as well.
From now on, I assume that every stack of lines includes the line at infinity and, therefore,
every point at infinity lies on the line at infinity.

\subsubsection{Orientation of finite points (sheaves)}
In the bottom-up view of geometry, points lack an orientation, but higher-dimensional objects, e.g. lines, can be assign an orientation.
A line has two possible orientations, both capturing a direction along the line. 
Consider the points which comprise the line and 
imagine stepping from point to point, moving in one direction along the line. This gives one of the orientations.
Now imagine stepping through the same points in the opposite direction. This will give the other orientation.
This way of orienting a line is point-based and, therefore, belongs to the bottom-up view of geometry.

\begin{figure}[h]
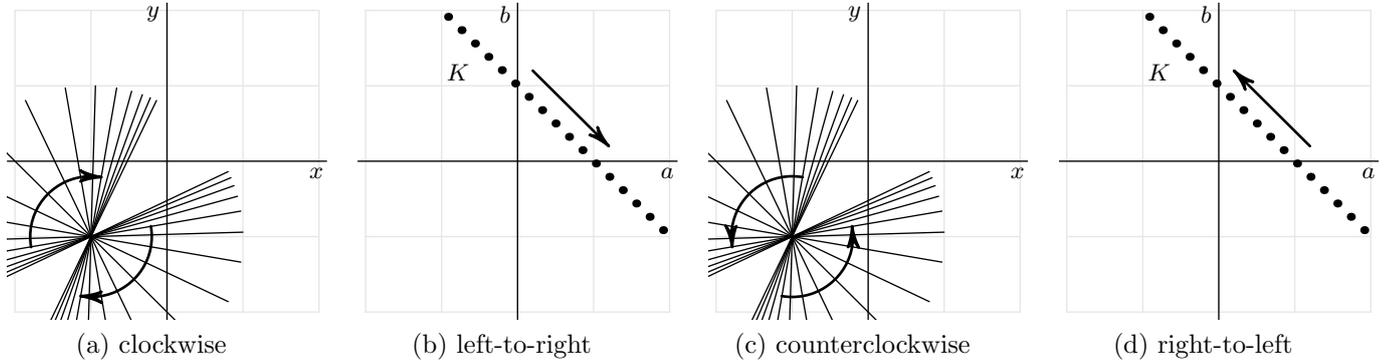

\begin{subfloatenv}{clockwise}
\begin{asy}
import Drawing2D;
Drawing2D drawing = Drawing2D(2, 0.1);
drawing.grid();
drawing.target_axes();

pair n = normalise((1,-1));
labeled_pair[] points = points_on_line(center=(1/2,1/2), direction=n, separation=0.25, N=8);
pair c = (-1,-1);
for(int i=0; i<points.length; ++i) { path l = line(center=c, direction=perp(points[i].pair), extent_from_center=2); draw(l); }

path a1 = arc(c, 0.8, 190, 80, CW);
draw(a1, linewidth(1 bp), Arrow(HookHead, size=5));

path a2 = arc(c, 0.8, 10, -100, CW);
draw(a2, linewidth(1 bp), Arrow(HookHead, size=5));

drawing.crop();
\end{asy}
\end{subfloatenv}%
\begin{subfloatenv}{left-to-right}
\begin{asy}
import Drawing2D;
Drawing2D drawing = Drawing2D(2, 0.1);
drawing.grid();
drawing.dual_axes();

pair n = normalise((1,-1));
labeled_pair[] points = points_on_line(center=(1/2,1/2), direction=n, separation=0.25, N=8);
for(int i=0; i<points.length; ++i) { dot(points[i].pair);}
label("$K$", (-1,1), NE);

pair shift = (0.2, 1.2);
path line_orientation = (shift--(sqrt(2)*n+shift));
draw(line_orientation, linewidth(1 bp), Arrow(HookHead, size=5));

drawing.crop();
\end{asy}
\end{subfloatenv}%
\begin{subfloatenv}{counterclockwise}
\begin{asy}
import Drawing2D;
Drawing2D drawing = Drawing2D(2, 0.1);
drawing.grid();
drawing.target_axes();

pair n = normalise((-1,1));
labeled_pair[] points = points_on_line(center=(1/2,1/2), direction=n, separation=0.25, N=8);
pair c = (-1,-1);
for(int i=0; i<points.length; ++i) { path l = line(center=c, direction=perp(points[i].pair), extent_from_center=2); draw(l); }

path a1 = arc(c, 0.8, 80, 190, CCW);
path a2 = arc(c, 0.8, -100, 10, CCW);
draw(a1, linewidth(1 bp), Arrow(HookHead, size=5));
draw(a2, linewidth(1 bp), Arrow(HookHead, size=5));

drawing.crop();
\end{asy}
\end{subfloatenv}%
\begin{subfloatenv}{right-to-left}
\begin{asy}
import Drawing2D;
Drawing2D drawing = Drawing2D(2, 0.1);
drawing.grid();
drawing.dual_axes();

pair n = normalise((-1,1));
labeled_pair[] points = points_on_line(center=(1/2,1/2), direction=n, separation=0.25, N=8);
for(int i=0; i<points.length; ++i) { dot(points[i].pair);}
label("$K$", (-1,1), NE);

pair shift = (1.2, 0.2);
path line_orientation = (shift--(sqrt(2)*n+shift));
draw(line_orientation, linewidth(1 bp), Arrow(HookHead, size=5));

drawing.crop();
\end{asy}
\end{subfloatenv}
\caption{Line orientation in \T{2*} and the induced point orientation in \T{2}}
\label{orientation}
\end{figure}

\begin{figure}[h]
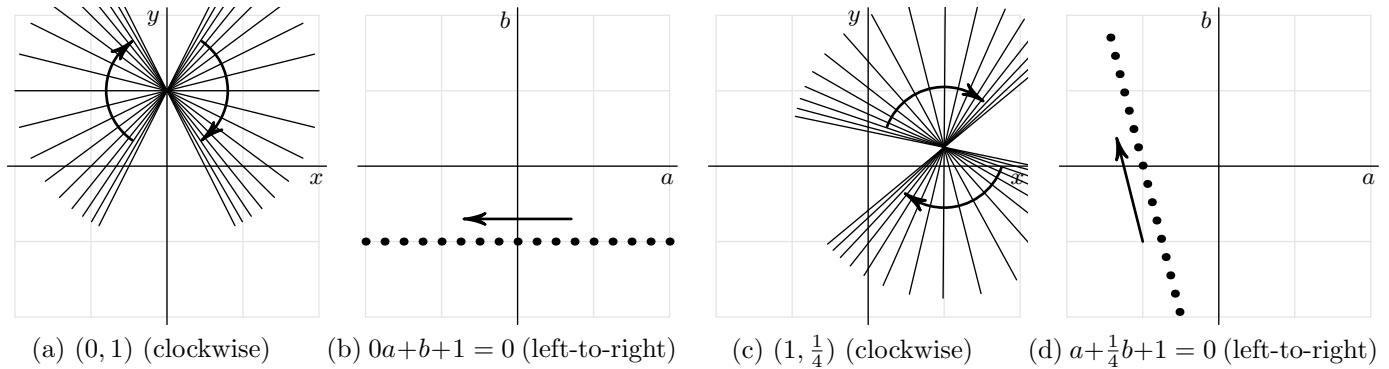

\begin{subfloatenv}{\((0,1)\) (clockwise)}
\begin{asy}
import Drawing2D;
Drawing2D drawing = Drawing2D(2, 0.1);
drawing.grid();
drawing.target_axes();

pair n = normalise((-1,0));
labeled_pair[] points = points_on_line(center=(0,-1), direction=n, separation=0.25, N=8);

pair c = (0,1);
for(int i=0; i<points.length; ++i) { path l = line(center=c, direction=perp(points[i].pair), extent_from_center=2); draw(l); }

path a1 = arc(c, 0.8, 235, 125, CW);
path a2 = arc(c, 0.8, 55,  -55, CW);
draw(a1, linewidth(1 bp), Arrow(HookHead, size=5));
draw(a2, linewidth(1 bp), Arrow(HookHead, size=5));

drawing.crop();
\end{asy}
\end{subfloatenv}%
\begin{subfloatenv}{\(0a+b+1=0\) (left-to-right)}
\begin{asy}
import Drawing2D;
Drawing2D drawing = Drawing2D(2, 0.1);
drawing.grid();
drawing.dual_axes();

pair n = normalise((-1,0));
labeled_pair[] points = points_on_line(center=(0,-1), direction=n, separation=0.25, N=8);
for(int i=0; i<points.length; ++i) { dot(points[i].pair);}

pair shift = (sqrt(2)/2, -0.7);
path line_orientation = (shift--(sqrt(2)*n+shift));
draw(line_orientation, linewidth(1 bp), Arrow(HookHead, size=5));

drawing.crop();
\end{asy}
\end{subfloatenv}%
\begin{subfloatenv}{\((1,\tfrac{1}{4})\) (clockwise)}
\begin{asy}
import Drawing2D;
Drawing2D drawing = Drawing2D(2, 0.1);
drawing.grid();
drawing.target_axes();

pair n = normalise((-1,4));
labeled_pair[] points = points_on_line(center=(-16/17,-4/17), direction=n, separation=0.25, N=8);

pair c = (1,1/4);
for(int i=0; i<points.length; ++i) { path l = line(center=c, direction=perp(points[i].pair), extent_from_center=2); draw(l); }

path a1 = arc(c, 0.8, 160, 50, CW);
path a2 = arc(c, 0.8, -20, -130, CW);
draw(a1, linewidth(1 bp), Arrow(HookHead, size=5));
draw(a2, linewidth(1 bp), Arrow(HookHead, size=5));

drawing.crop();
\end{asy}
\end{subfloatenv}%
\begin{subfloatenv}{\(a+\tfrac{1}{4}b+1=0\) (left-to-right)}
\begin{asy}
import Drawing2D;
Drawing2D drawing = Drawing2D(2, 0.1);
drawing.grid();
drawing.dual_axes();

pair n = normalise((-1,4));
labeled_pair[] points = points_on_line(center=(-16/17,-4/17), direction=n, separation=0.25, N=8);
for(int i=0; i<points.length; ++i) { dot(points[i].pair);}

pair shift = (-1, -1);
path line_orientation = (shift--(sqrt(2)*n+shift));
draw(line_orientation, linewidth(1 bp), Arrow(HookHead, size=5));

drawing.crop();
\end{asy}
\end{subfloatenv}
\caption{Line orientation in \T{2*} and the induced point orientation in \T{2} (additional examples)}
\label{orientation further examples}
\end{figure}

In the top-down view, lines are fundamental and lower-dimensional objects, e.g.\ points, can be assigned an orientation.
A point in the top-down view is a sheaf of intersecting lines. 
As I step from line to line in the sheaf, I perform circular motion around the point where the lines intersect.
I can step through the lines in the clockwise or counterclockwise direction. 
Selecting a sense of circular motion, clockwise or counterclockwise, determines the orientation of the point.

Depending on the stance, i.e.\ bottom-up or top-down, I can orient lines with points or points with lines.
This applies equally to the target space \T{2} and the dual space \T{2*}.
Moreover, the bottom-up orientation of lines in \T{2*} and the top-down orientation of the corresponding points in \T{2} are consistent.
Indeed, if I adopt the bottom-up view of the dual space \T{2*}, I can orient  \(K=\{(a,b)|-a-b+1=0\}\)
in \T{2*} as shown in Figure~\ref{orientation}(bd).
The line \(K\) corresponds to a sheaf of lines in \T{2} intersecting at \((-1,-1)\).
The bottom-up orientation of \(K\), left-to-right as seen from the origin, 
induces the clockwise orientation of the sheaf (see Figure~\ref{orientation}(a)).
The same line with the opposite orientation, i.e.\ right-to-left, induces the counterclockwise orientation of the sheaf (see Figure~\ref{orientation}(c)).
Note that stepping through the line that passes through  the origin \((0,0)\) in \T{2}
corresponds to stepping through a point at infinity in \T{2*}.
Additional examples are shown in Figure~\ref{orientation further examples}.

\subsubsection{Orientation of points at infinity (stacks) and lines}
\begin{figure}[ht!]
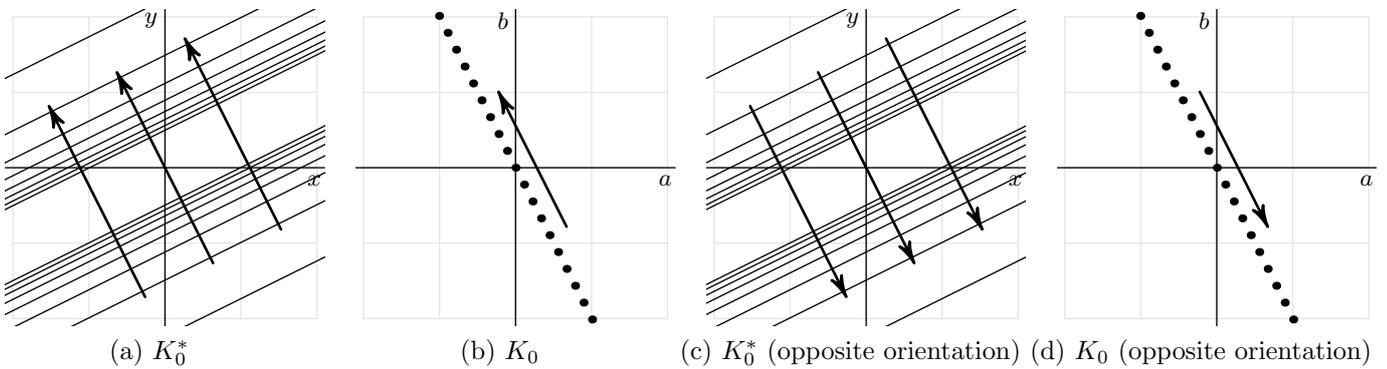

\begin{subfloatenv}{\(K_0^*\)}
\begin{asy}
import Drawing2D;
Drawing2D drawing = Drawing2D(2, 0.1);
drawing.grid();
drawing.target_axes();

pair n = normalise((-1,2));
labeled_pair[] points = points_on_line(center=(0,0), direction=n, separation=0.25, N=9);

int index_of_origin = quotient(points.length-1, 2);
points.delete(index_of_origin);

for(int i=0; i<points.length; ++i) { real a = points[i].pair.x; real b = points[i].pair.y; path l = line(center=-1/5 * (1/a, 4/b), direction=perp(points[i].pair), extent_from_center=4); draw(l); }

pair c = (2/sqrt(5),1/sqrt(5));
pair shift = sqrt(2)*(1/sqrt(5),-2/sqrt(5));
for(int i=-1; i<=1; ++i) { path line_orientation = (i*c+shift)--(i*c-shift); draw(line_orientation, linewidth(1 bp), Arrow(HookHead, size=5)); }

drawing.crop();
\end{asy}
\end{subfloatenv}%
\begin{subfloatenv}{\(K_0\)}
\begin{asy}
import Drawing2D;
Drawing2D drawing = Drawing2D(2, 0.1);
drawing.grid();
drawing.dual_axes();

pair n = normalise((-1,2));
labeled_pair[] points = points_on_line(center=(0,0), direction=n, separation=0.25, N=9);
for(int i=0; i<points.length; ++i) { dot(points[i].pair);}

pair c = 0.25*(2/sqrt(5),1/sqrt(5));
pair shift = (1/sqrt(5),-2/sqrt(5));
path line_orientation = (c+shift)--(c-shift);
draw(line_orientation, linewidth(1 bp), Arrow(HookHead, size=5));

drawing.crop();
\end{asy}
\end{subfloatenv}%
\begin{subfloatenv}{\(K_0^*\) (opposite orientation)}
\begin{asy}
import Drawing2D;
Drawing2D drawing = Drawing2D(2, 0.1);
drawing.grid();
drawing.target_axes();

pair n = normalise((-1,2));
labeled_pair[] points = points_on_line(center=(0,0), direction=n, separation=0.25, N=9);

int index_of_origin = quotient(points.length-1, 2);
points.delete(index_of_origin);

for(int i=0; i<points.length; ++i) { real a = points[i].pair.x; real b = points[i].pair.y; path l = line(center=-1/5 * (1/a, 4/b), direction=perp(points[i].pair), extent_from_center=4); draw(l); }

pair c = (2/sqrt(5),1/sqrt(5));
pair shift = sqrt(2)*(1/sqrt(5),-2/sqrt(5));
for(int i=-1; i<=1; ++i) { path line_orientation = (i*c-shift)--(i*c+shift); draw(line_orientation, linewidth(1 bp), Arrow(HookHead, size=5)); }

drawing.crop();
\end{asy}
\end{subfloatenv}%
\begin{subfloatenv}{\(K_0\) (opposite orientation)}
\begin{asy}
import Drawing2D;
Drawing2D drawing = Drawing2D(2, 0.1);
drawing.grid();
drawing.dual_axes();

pair n = normalise((-1,2));
labeled_pair[] points = points_on_line(center=(0,0), direction=n, separation=0.25, N=9);
for(int i=0; i<points.length; ++i) { dot(points[i].pair);}

pair c = 0.25*(2/sqrt(5),1/sqrt(5));
pair shift = (1/sqrt(5),-2/sqrt(5));
path line_orientation = (c-shift)--(c+shift);
draw(line_orientation, linewidth(1 bp), Arrow(HookHead, size=5));

drawing.crop();
\end{asy}
\end{subfloatenv}
\caption{Orientation of a point at infinity in \T{2}}
\label{orientation parallel lines}
\end{figure}

In the top-down view, a point at infinity is identified with a stack of lines.
I can step from line to line in the stack in two opposite directions.
Selecting one of these directions determines the orientation of the point at infinity.
If I adopt the bottom-up view of \T{2*}, I can assign a specific point-based orientation 
to a line in \T{2*} dual to the stack.
Consider  \(K_0=\{(a,b)|2a+b=0\}\), which corresponds to a stack of lines in \T{2} and, therefore, to a point at infinity in \T{2}.
As I step from point to point along the line \(K_0\), I step from one line in the stack to another
as shown in Figure~\ref{orientation parallel lines}.
Note that stepping through the origin of \T{2*} corresponds to stepping through the line at infinity in \T{2}.
Stepping through the point at infinity on the line \(K_0\) corresponds to stepping through 
a line passing through the origin in \T{2}.

If an oriented finite point can be viewed as a sheaf of rotating lines intersecting at the finite point,
then an oriented point at infinity, i.e.~an oriented stack of lines,
can be viewed informally as a set of lines undergoing a kind of rotation around the point at infinity.

\begin{figure}[h!]
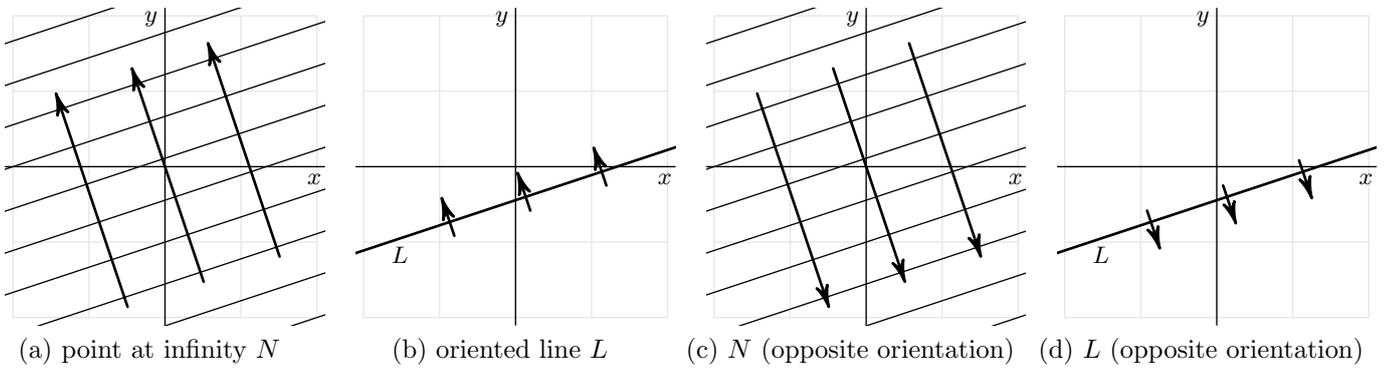

\begin{subfloatenv}{point at infinity \(N\)}
\begin{asy}
import Drawing2D;
Drawing2D drawing = Drawing2D(2, 0.1);
drawing.grid();
drawing.target_axes();

pair n = (1,1/3);
pair c = (-1/3,0);
pair[] ps = {};

int N=10;
for(int i=0; i<N; ++i) { pair P=c+(i/3,-i)/2; ps.push(P); }
for(int i=1; i<N; ++i) { pair P=c+(-i/3,i)/2; ps.push(P); }

for(pair P: ps) { path l = line(center=P, direction=n, extent_from_center=5);  draw(l); }

drawing.crop();

pair c = (1,1/3);
pair s = (1/9,-1/3)/3;
pair shift = sqrt(2)*(1/3,-1);
for(int i=-1; i<=1; ++i) { path line_orientation = (i*c+s+shift)--(i*c+s-shift); draw(line_orientation, linewidth(1 bp), Arrow(HookHead, size=5)); }

drawing.crop();
\end{asy}
\end{subfloatenv}%
\begin{subfloatenv}{oriented line \(L\)}
\begin{asy}
import Drawing2D;
Drawing2D drawing = Drawing2D(2, 0.1);
drawing.grid();
drawing.target_axes();

pair n = (1,1/3);
pair c = (-1/3,0);
path l = line(center=c+(1/3,-1)/2, direction=n, extent_from_center=5); 
draw(Label("$L$", 0.35), l, currentpen+1);

pair c = (1,1/3);
pair s = (1/9,-1/3);
pair shift = 0.25*(1/3,-1);
for(int i=-1; i<=1; ++i) { path line_orientation = (i*c+s+shift)--(i*c+s-shift); draw(line_orientation, linewidth(1 bp), Arrow(HookHead, size=5)); }

drawing.crop();
\end{asy}
\end{subfloatenv}%
\begin{subfloatenv}{\(N\) (opposite orientation)}
\begin{asy}
import Drawing2D;
Drawing2D drawing = Drawing2D(2, 0.1);
drawing.grid();
drawing.target_axes();

pair n = (1,1/3);
pair c = (-1/3,0);
pair[] ps = {};

int N=10;
for(int i=0; i<N; ++i) { pair P=c+(i/3,-i)/2; ps.push(P); }
for(int i=1; i<N; ++i) { pair P=c+(-i/3,i)/2; ps.push(P); }

for(pair P: ps) { path l = line(center=P, direction=n, extent_from_center=5);  draw(l); }

drawing.crop();

pair c = (1,1/3);
pair s = (1/9,-1/3)/3;
pair shift = sqrt(2)*(1/3,-1);
for(int i=-1; i<=1; ++i) { path line_orientation = (i*c+s-shift)--(i*c+s+shift); draw(line_orientation, linewidth(1 bp), Arrow(HookHead, size=5)); }

drawing.crop();
\end{asy}
\end{subfloatenv}%
\begin{subfloatenv}{\(L\) (opposite orientation)}
\begin{asy}
import Drawing2D;
Drawing2D drawing = Drawing2D(2, 0.1);
drawing.grid();
drawing.target_axes();

pair n = (1,1/3);
pair c = (-1/3,0);
path l = line(center=c+(1/3,-1)/2, direction=n, extent_from_center=5); 
draw(Label("$L$", 0.35), l, currentpen+1);

pair c = (1,1/3);
pair s = 1.5*(1/9,-1/3);
pair shift = 0.25*(1/3,-1);
for(int i=-1; i<=1; ++i) { path line_orientation = (i*c+s-shift)--(i*c+s+shift); draw(line_orientation, linewidth(1 bp), Arrow(HookHead, size=5)); }

drawing.crop();
\end{asy}
\end{subfloatenv}
\caption{Orienting a line with a point at infinity}
\label{orientation of lines}
\end{figure}

Lines lack an orientation in the top-down view, but I can give them an orientation by noting
that each finite line \(L\) is closely related to a stack of lines representing the point at infinity that lies on \(L\).
So, I can orient \(L\) by adopting a top-down orientation of the relevant point at infinity.
This is illustrated in Figure~\ref{orientation of lines}, where 
\(L=\{(x,y)|-x+3y+1=0\}\).

This definition can be extended to the line at infinity,
which can be oriented either towards the origin or away from the origin.

\subsubsection{Embedding}
\begin{figure}[h!]
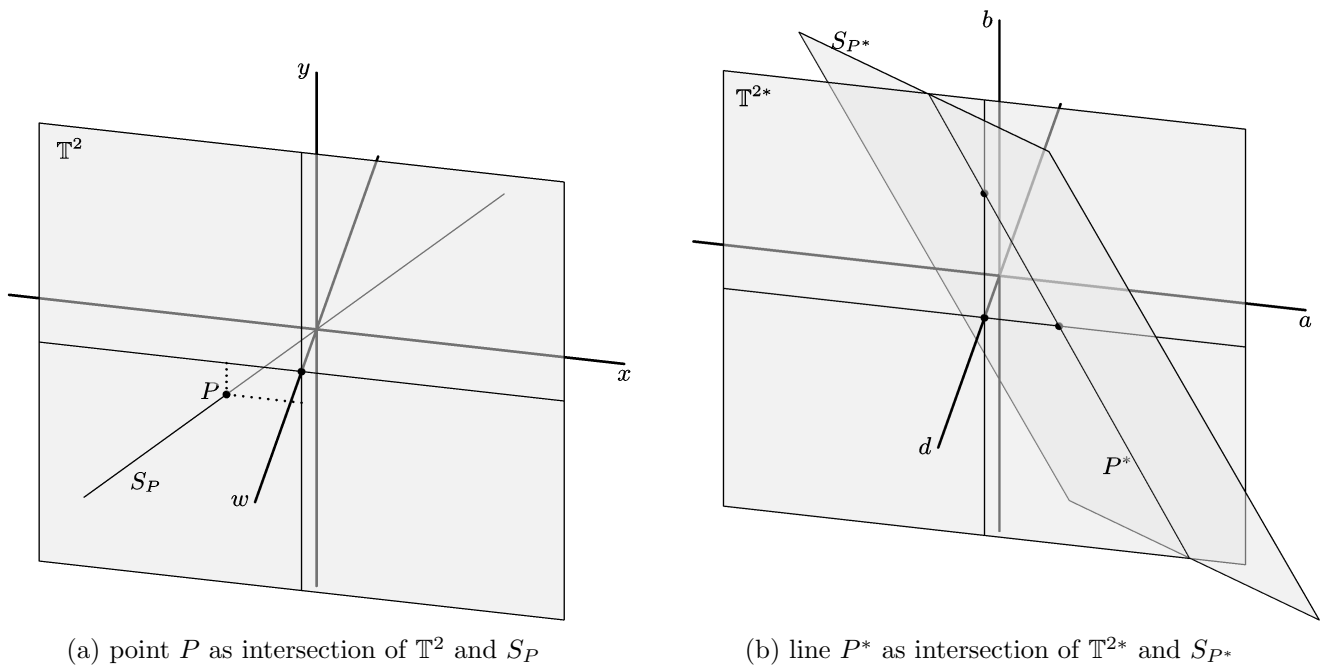

\begin{subfloatenv}{point \(P\) as intersection of \T{2} and \(S_P\)}
\begin{asy}
import Drawing3D;
DrawingR3 drawing = DrawingR3(4,0.1, camera=(10,2,7));	
drawing.target_axes();

path3 p = plane_path3((1,0,0),(1,0,0), v=(0,1,0));
drawing.Drawing3D.plane(p);
dot((1,0,0));

draw((1,-3.5,0)--(1,3.5,0));
draw((1,0,-3.5)--(1,0,3.5));

path3 l = line_path3((1,-1,-1/2)/4, direction=(1,-1,-1/2), extent_from_center=3.5);
draw(Label("$S_P$", 0.9, SE), l);
dot((1,-1,-1/2));

draw((1,-1,0)--(1,-1,-1/2), Dotted);
draw((1,0,-1/2)--(1,-1,-1/2), Dotted);

//label("$(1,-1,-\frac{1}{2})$",(1,-1,-1/2),(1/2,-1,1/2));
label("$P$",(1,-1,-1/2),(1/2,-1,1/2));
label("$\mathbb{T}^{2}$", (1.5,-3,3.5));
\end{asy}
\end{subfloatenv}\hspace{-12pt}%
\begin{subfloatenv}{line \(P^*\) as intersection of \T{2*} and \(S_{P^*}\)}
\begin{asy}
import Drawing3D;
DrawingR3 drawing = DrawingR3(4,0.1, camera=(10,2,7));
drawing.dual_axes();

path3 p = plane_path3((1,0,0),(1,0,0), v=(0,1,0));
drawing.Drawing3D.plane(p);
dot((1,0,0));

draw((1,-3.5,0)--(1,3.5,0));
draw((1,0,-3.5)--(1,0,3.5));

path3 Q = plane_path3((1,-1,-1/2), (1,1,0), theta=-13, scale=7.8);
drawing.Drawing3D.plane(Q);

path3 l = plane_intersection(Q,p);
draw(Label("$P^*$", 0.2,W), l);

label("$S_{P^*}$",(0,-2,3.5));
label("$\mathbb{T}^{2*}$", (1.5,-3,3.5));

//path3 l = line_path3((1,-1,-1/2)/4, direction=(1,-1,-1/2), extent_from_center=3.5);
//draw(Label("$I^{-1}(S_P)$", 0.9, SE), l);
//dot((1,-1,-1/2));

dot((1,1,0));
//label("$(1,1,0)$",(1,1,0),(1,-0.5,-0.5));
dot((1,0,2));
//label("$(1,0,2)$",(1,0,2),(1,-0.8,0));
\end{asy}
\end{subfloatenv}
\caption{Embedding and representation (a point in \T{2} and its dual line)}
\label{embedding T2 in R3}
\end{figure}

The target model space \R{3} consists of triples \((w,x,y)\), where \((1,x,y)\) is identified with \((x,y)\in\T{2}\),
and the dual model space \R{3*} consists of triples \((d,a,b)\), where \((1,a,b)\) is identified with \((a,b)\in\T{2*}\).
In other words, \T{2} is a plane embedded in \R{3} at \(w=1\),
 and \T{2*} is a plane embedded in \R{3*} at \(d=1\).
Neither \T{2} nor \T{2*} are linear subspaces since they do not pass through the origin
of their model spaces.
\T{2} can be obtained from the \(xy\)-plane by shifting it from \(w=0\) to \(w=1\),
and \T{2*} can be obtained from the \(ab\)-plane by shifting it from \(d=0\) to \(d=1\)
(see Figure~\ref{embedding T2 in R3}).

A 1-dimensional linear subspace of \R{3} represents a point in \T{2} as follows.
Consider a finite point \(P=(x,y)\) in \T{2} .
Thanks to the embedding, it can be viewed as a point in \R{3} with the coordinates \((1,x,y)\).
There is a line \(S_P\) which passes through the origin \((0,0,0)\) of \R{3} and 
intersects \T{2} at \((1,x,y)\), i.e. \(S_P\) intersects \T{2} at \(P\)
 (see Figure~\ref{embedding T2 in R3}(a) where \(P=(-1,-\tfrac{1}{2})\)).
\(S_P\) is a linear subspace of \R{3} since it passes through the origin of \R{3} and
it is 1-dimensional since it is a line.
There is a one-to-one correspondence between the point \(P\) and the linear subspace \(S_P\).
Indeed, selecting a 1-dimensional linear subspace which intersects \T{2} determines the intersection point uniquely,
and selecting a finite point in \T{2} determines the linear subspace uniquely.
In this way, \(S_P\) is a 1-dimensional linear subspace of \R{3} that represents the point \(P\) in the target space \T{2}.

Essentially the same construction applies to the 2-dimensional linear subspaces of \R{3}, which represent lines in \T{2}.
Consider a finite line \(L\) in \T{2}.
Thanks to the embedding, it can be viewed as a line in \R{3}.
There is a plane \(S_L\) which passes through the origin of \R{3} and intersects \T{2} at \(L\)
(see Figure~\ref{embedding T2 in R3 another example}(a) where \(L\) is defined by \(\tfrac{1}{2}x-\tfrac{1}{3}y+1=0\)).
\(S_L\) is a linear subspace of \R{3} since it passes through the origin and
it is 2-dimensional since it is a plane.
There is a one-to-one correspondence between the line \(L\) and the linear subspace \(S_L\).
Selecting a 2-dimensional linear subspace which intersects \T{2} determines the intersection line uniquely,
and selecting a finite line in \T{2} determines the linear subspace uniquely.
In this way, \(S_L\) is a 2-dimensional linear subspace of \R{3} that represents the line \(L\) in the target space \T{2}.

\begin{figure}[t!]
\hspace{-1cm}\begin{subfloatenv}{line \(L\) as intersection of \T{2} and \(S_L\)}
\begin{asy}
import Drawing3D;
DrawingR3 drawing = DrawingR3(5.1,0.1, axes_extent=4.0, camera=(10,2,7));	
drawing.target_axes();

path3 p = plane_path3((1,0,0),(1,0,0), v=(0,1,0));
drawing.Drawing3D.plane(p);
dot((1,0,0));
draw((1,-3.5,0)--(1,3.5,0));
draw((1,0,-3.5)--(1,0,3.5));

real x = -1.584;
triple c = (1, x, 1.5*x+3);
path3 Q = plane_path3((1,1/2,-1/3), c, theta=-11, scale=6.92); drawing.Drawing3D.plane(Q);

path3 l = plane_intersection(Q,p); draw(Label("$L$", 0.2,E), l);

dot((1,-2,0));
//label("$(1,-2,0)$",(1,-2,0),(1,1,-0.5));
dot((1,0,3));
//label("$(1,0,3)$",(1,0,3),(1,1.5,-0.5));

label("$S_L$",(2,-5,-1));
label("$\mathbb{T}^{2}$", (1.5,-3,3.5));

//path3 l = line_path3((1,1/2,-1/3)/4, direction=(1,1/2,-1/3), extent_from_center=3.5);
//draw(Label("$I(S_{L^*})$", 0.9, E), l);
//dot((1,1/2,-1/3));
\end{asy}
\end{subfloatenv}\hspace{-12pt}%
\begin{subfloatenv}{point \(L^*\) as intersection of \T{2*} and \(S_{L^*}\)}
\begin{asy}
import Drawing3D;
DrawingR3 drawing = DrawingR3(4,0.1, camera=(10,2,7));
drawing.dual_axes();

path3 p = plane_path3((1,0,0),(1,0,0), v=(0,1,0));
drawing.Drawing3D.plane(p);
dot((1,0,0));

draw((1,-3.5,0)--(1,3.5,0));
draw((1,0,-3.5)--(1,0,3.5));

path3 l = line_path3((1,1/2,-1/3)/4, direction=(1,1/2,-1/3), extent_from_center=3.5);
draw(Label("$S_{L^*}$", 0.9, E), l);
dot((1,1/2,-1/3));

draw((1,1/2,0)--(1,1/2,-1/3), Dotted);
draw((1,0,-1/3)--(1,1/2,-1/3), Dotted);

//label("$(1,\frac{1}{2},-\frac{1}{3})$",(1,1/2,-1/3),E);
label("$L^*$",(1,1/2,-1/3),E);
label("$\mathbb{T}^{2*}$", (1.5,-3,3.5));
\end{asy}
\end{subfloatenv}
\caption{Embedding and representation (a line in \T{2} and its dual point)}
\label{embedding T2 in R3 another example}
\end{figure}%
\begin{figure}[t!]
\begin{subfloatenv}{point \(N\) as a stack of lines represented by \(S_{N}\)}
\begin{asy}
import Drawing3D;
DrawingR3 drawing = DrawingR3(4.,0.1, axes_extent=4.0, camera=(10,2,7));	
drawing.target_axes();

path3 p = plane_path3((1,0,0),(1,0,0), v=(0,1,0));
drawing.Drawing3D.plane(p);
dot((1,0,0));
draw((1,-3.5,0)--(1,3.5,0));
draw((1,0,-3.5)--(1,0,3.5));

path3 l = line_path3((0,0,0), direction=(0,1/2,3), extent_from_center=4);
draw(Label("$S_{N}$", 0.9, E), l);

path3 l = line_path3((1,0,0), (0,1/2,3), 10);
path3 p = drawing.T2_path3();
int N=5;
for(int i=-N; i<=N; ++i) { path3 shifted_l = shift(-i*(0,-1,2)/2)*l; triple[] ps = intersectionpoints(shifted_l,p); draw(ps[0]--ps[1]); }

label("$\mathbb{T}^{2}$", (1.5,-3,3.5));
\end{asy}
\end{subfloatenv}\hspace{-12pt}%
\begin{subfloatenv}{line \(N^*\) as intersection of \T{2*} and \(S_{N^*}\)}
\begin{asy}
import Drawing3D;
DrawingR3 drawing = DrawingR3(4,0.1, camera=(10,2,7));
drawing.dual_axes();

path3 p = plane_path3((1,0,0),(1,0,0), v=(0,1,0));
drawing.Drawing3D.plane(p);
dot((1,0,0));

draw((1,-3.5,0)--(1,3.5,0));
draw((1,0,-3.5)--(1,0,3.5));

path3 Q = plane_path3((0,1/2,3), (0,0,0), theta=50.6, scale=7.12); drawing.Drawing3D.plane(Q);

path3 l = plane_intersection(Q,p); draw(Label("$N^*$", 0.7,N), l);

dot((3.56,0,0));
dot((-3.56,0,0));

label("$S_{N^*}$",(3,2,-0.25));

//path3 l = line_path3((0,0,0), direction=(0,1/2,3), extent_from_center=4);
//draw(Label("$I^{-1}(S_{K_0^*})$", 0.9, E), l);

label("$\mathbb{T}^{2*}$", (1.5,-3,3.5));
\end{asy}
\end{subfloatenv}
\caption{Embedding and representation
(point at infinity in \T{2} and its dual line}
\label{embedding T2 in R3 infinite point}
\end{figure}

Exactly the same construction can be carried out for points and lines in the dual space \T{2*}, which are
represented by the linear subspaces of the dual model space \R{3*}.
See Figure~\ref{embedding T2 in R3}(b), where this construction is applied to
the line \(P^*\) in \T{2*}, which is dual to the point \(P\) in \T{2}.
The line \(P^*\) is represented by the 2-dimensional linear subspace \(S_{P^*}\).
See also Figure~\ref{embedding T2 in R3 another example}(b), where the construction is applied to
the point \(L^*\) in \T{2*}, which is dual to the line \(L\) in \T{2}.
The point \(L^*\) is represented by the 1-dimensional linear subspace \(S_{L^*}\).

Points at infinity in \T{2} are represented by the 1-dimensional linear subspaces of \R{3} that do not intersect \T{2}
(all such subspaces lie in the \(xy\)-plane).
For example, consider a point at infinity \(N\) in \T{2} shown in Figure~\ref{embedding T2 in R3 infinite point}(a) as a stack of lines.
It is represented by the 1-dimensional linear subspace \(S_{N}\).
The point at infinity \(N\) is dual to the line \(N^*\) in \T{2*}, which is represented in \R{3*} by the 2-dimensional linear subspace \(S_{N^*}\).
The line at infinity in \T{2} is represented by the \(xy\)-plane in \R{3}; its dual is the \(d\)-axis in \R{3*}.
Exactly the same construction applies to points at infinity and the line at infinity in \T{2*}.
Namely, the \(ab\)-plane represents the line at infinity in \T{2*}, while the 1-dimensional linear
subspaces in the \(ab\)-plane represent points at infinity in \T{2*}.

I will refer to the subspace \(S_{P^*}\) in \R{3*} as dual to the subspace \(S_P\) in \R{3}, 
since they represent dual geometric objects in \T{2} and \T{2*}.
Similarly, I will refer to \(S_{L^*}\) as dual to \(S_L\), and \(S_{N^*}\) as dual to \(S_{N}\).
The duality transformation \(\J\) is thus extended to the linear subspaces of \R{3} and \R{3*}.
For instance, \(\J(S_{P^*})=S_{P}\), \(\J(S_{L^*})=S_{L}\), and \(\J(S_{N^*})=S_{N}\),
i.e.\ a 1-dimensional linear subspace of \R{3} corresponds to a 2-dimensional linear subspace of \R{3*} and vice versa.
I also extend the identity transformation \(\Id\) to the linear subspaces of \R{3} and \R{3*}.
For instance, \(\Id(S_{L^*})\) is a 1-dimensional linear subspace of \R{3} which passes through \((1,\tfrac{1}{2},-\tfrac{1}{3})\) in \R{3},
whereas \(S_{L^*}\) is a 1-dimensional linear subspace of \R{3*} passing through the point \((1,\tfrac{1}{2},-\tfrac{1}{3})\) in \R{3*}.

\subsection{Grassmann algebra}
Non-metric relations between linear subspaces of \R{3} (and \R{3*}) are described by Grassmann algebra.

I introduce the following notation for the basis of \R{3}: 
\(\e^0=(1,0,0)\), \(\e^1=(0,1,0)\), \(\e^2=(0,0,1)\), so that \((w,x,y)=w\e^0+x\e^1+y\e^2\).
And for \R{3*}: \(\e_0=(1,0,0)\), \(\e_1=(0,1,0)\), \(\e_2=(0,0,1)\), so that \((d,a,b)=d\e_0+a\e_1+b\e_2\).
The superscript indices will be used for the basis vectors in the target space and the subscript indices for the basis vectors in the dual space.

The basis of the Grassmann algebra of the dual model space \R{3*} (denoted by \(\bigwedge\R{3*}\)) 
consists of one scalar, three vectors, three bivectors, and one pseudoscalar (trivector).
General multivectors are formed by linear combinations of scalars, vectors, bivectors, and trivectors.
I will use the symbol \(\wedge\) to denote the outer product in \(\bigwedge\R{3*}\).
For vectors \(\tb{a}, \tb{b}, \tb{c}\in\R{3*}\), the outer product exhibits the following properties:\\
1) \(\tb{a}\wedge\tb{b}=-\tb{b}\wedge\tb{a}\) (anticommutative),\\
2) \((\alpha\tb{a}+\beta\tb{b})\wedge\tb{c}=\alpha\tb{a}\wedge\tb{c}+\beta\tb{b}\wedge\tb{c}\) (distributive),\\
3) \(\tb{a}\wedge(\tb{b}\wedge\tb{c})=(\tb{a}\wedge\tb{b})\wedge\tb{c}\) (associative).\\
It follows that \(\tb{a}\wedge\tb{a}=0\) for any vector \(\tb{a}\in\R{3*}\).
For brevity, I will use a simplified notation for the outer product between the basis vectors.
For instance, \(\e_{10}=\e_1\wedge\e_0\), \(\e_{210}=\e_2\wedge\e_1\wedge\e_0\), etc.
So, the basis of \(\bigwedge\R{3*}\) consists of the following \(2^3=8\) multivectors:\\
1, \\
\(\e_0, \e_1, \e_2\), \\
\(\e_{12}, \e_{20}, \e_{01}\),\\
\(\I_3=\e_{012}\).\\
I will usually drop the subscript of the pseudoscalar \(\I_3\) and write \(\I=\e_{012}\).
I will use capital Roman letters in bold font for bivectors in \(\bigwedge\R{3*}\), e.g.~\(\tb{P}, \tb{Q}\) 
(pseudoscalar \(\I\) is an exception).
So, a general multivector in  \(\bigwedge\R{3*}\) can be written as
\begin{equation}
M=s+\tb{a}+\tb{P}+\alpha\I,
\end{equation}
where \(s, \alpha\in\R{}\).
By definition, \(\alpha\wedge M=\alpha M\) for any scalar \(\alpha\) and multivector \(M\).
The multivector components are graded according to the dimension of the subspace they symbolise, e.g.\ 
scalars are grade-0 and bivectors are grade-2.
Multivectors of each grade form a linear subspace of \(\bigwedge\R{3*}\).
\(\bigwedge^k\R{3*}\) will denote the linear subspace of \(k\)-vectors, e.g.\ \(\bigwedge^2\R{3*}\) for bivectors.
Note that \(\bigwedge^0\R{3*}=\R{}\) and \(\bigwedge^1\R{3*}=\R{3*}\).
A multivector is called even if it consists of the components whose grade is even, e.g.\ scalars, bivectors, etc.
\(\grade{M}_k\) yields the grade-\(k\) component of multivector \(M\), e.g. \(\grade{M}_2=\tb{P}\).
Multivectors consisting of a scalar and a pseudoscalar are called Study numbers.
In \(\bigwedge\R{3*}\), \(A_k\wedge B_l=0\) if \(k+l>3\), where \(A_k\) is a \(k\)-vector and \(B_l\) is an \(l\)-vector.
For instance, the outer product between any two bivectors in \(\bigwedge\R{3*}\) is zero, \(\tb{P}\wedge\tb{Q}=0\).
For even multivectors in general and, therefore, for scalars and bivectors in particular, 
the outer product with an arbitrary multivector \(M\) is commutative: \(\alpha\wedge M=M\wedge\alpha\), \(\tb{P}\wedge M=M\wedge\tb{P}\).
For instance, \(\tb{P}\wedge\tb{a}=\tb{a}\wedge\tb{P}\).

A blade, or a simple \(k\)-vector, is a multivector that can be written as the outer product of \(k\) vectors.
Scalars and vectors are simple by definition.
In  \(\bigwedge\R{3*}\), bivectors and trivectors 
are also simple\footnote{An example of a non-simple bivector is \(\e_{01}+\e_{23}\) in \(\bigwedge\R{4*}\).}.
For instance, \(\e_{12}=\e_1\wedge\e_2\) is the outer product of two vectors and, therefore, it is a simple bivector.
A simple \(k\)-vector symbolises a \(k\)-dimensional linear subspace of \R{3*}.
A multivector that mixes different grades, e.g. \(M=1-2\e_{01}\), is not a \(k\)-vector and does not symbolise a subspace.

The Grassmann algebra of the target model space \R{3} has similar properties.
To avoid confusion, I will use the symbol \(\vee\) to denote the outer product in 
the Grassmann algebra of \R{3}, which will be denoted by \(\bigvee\R{3}\).
The same simplified notation for bivectors and trivectors will be used in \(\bigvee\R{3}\), e.g.~\(\e^{12}=\e^1\vee\e^2\).

A note on notation.
Focusing on \T{2},
it is beneficial to use the symbol \(\vee\), reminiscent of the set-theoretic symbol \(\cup\) (the union),
for the Grassmann algebra of \R{3}, since multivectors of \(\bigvee\R{3}\) enable the bottom-up view of \T{2}
where lines are build from points as their union or extension.
On the other hand, the Grassmann algebra of \R{3*} enables the top-down view of \T{2}
where points are build from lines as their intersection.
So it is natural to use the symbol \(\wedge\), which is reminiscent of the set-theoretic symbol \(\cap\) (the intersection).

A vector \(\tb{u}\in\bigvee\R{3}\) symbolises a 1-dimensional linear subspace of \R{3} consisting of vectors \(\tb{x}=\alpha\tb{u}\),
where \(\alpha\in\R{}\).
All vectors that belong to this subspace can be found by solving 
\begin{equation}
\tb{x}\vee\tb{u}=0.
\label{x from v}
\end{equation}
Indeed, \(\alpha\tb{u}\vee\tb{u}=0\), so vectors \(\alpha\tb{u}\), \(\alpha\in\R{}\), satisfy Equation~(\ref{x from v}).
Such vectors provide the only solution to Equation~(\ref{x from v}), since any vector that is not completely in the
subspace will form a non-zero bivector with \(\tb{u}\).

A bivector \(B=\tb{u}\vee\tb{v}\in\bigvee\R{3}\) symbolises a 2-dimensional linear subspace of \R{3}
consisting of vectors \(\tb{x}=\alpha\tb{u}+\beta\tb{v}\), where \(\alpha,\beta\in\R{}\).
Vectors that belong to this 2-dimensional subspace can be found by solving 
\begin{equation}
\tb{x}\vee B=0.
\label{x from B}
\end{equation}
Indeed, \((\alpha\tb{u}+\beta\tb{v})\vee(\tb{u}\vee\tb{v})=\alpha\tb{u}\vee(\tb{u}\vee\tb{v})+\beta\tb{v}\vee(\tb{u}\vee\tb{v})=
\alpha(\tb{u}\vee\tb{u})\vee\tb{v}-\beta\tb{v}\vee(\tb{v}\vee\tb{u})=
0-\beta(\tb{v}\vee\tb{v})\vee\tb{u}=0\), and
any vector that is not completely in the 2-dimensional subspace will form a non-zero trivector with \(B\).

Consider the 2-dimensional linear subspace \(S_L\) shown in Figure~\ref{embedding T2 in R3 another example}. 
It can be represented by the bivector 
\(
B=(\e^0+3\e^2)\vee(\e^0-2\e^1)=-2\e^0\vee\e^1 +3\e^2\vee\e^0-6\e^2\vee\e^1=
-2\e^{01} +3\e^{20}+6\e^{12}
\),
 since the subspace \(S_L\) passes through  the points \((1,-2,0)\) and \((1,0,3)\).
Substituting \(\tb{x}=w\e^0+x\e^1+y\e^2\) into Equation~(\ref{x from B}) yields
\(
(w\e^0+x\e^1+y\e^2)\vee (-2\e^{01} +3\e^{20}+6\e^{12})=
6w\e^{012}+3x\e^{120}-2y\e^{201}=(6w+3x-2y)\e^{012}=0
\)
and, therefore, \(6w+3x-2y=0\).
To find the intersection of \(S_L\) with \T{2}, I impose \(w=1\) 
and obtain
\(6+3x-2y=0\), which is equivalent to the definition of the line \(L\) shown in Figure~\ref{embedding T2 in R3 another example}.

In the dual model space \R{3*}, in order to find all vectors \(\tb{a}\) that lie in the subspace 
symbolised by a vector \(\tb{s}\) or a bivector \(\tb{P}\),
one needs to solve
\begin{equation}
\tb{a}\wedge\tb{s}=0
\label{a from s}
\end{equation}
or 
\begin{equation}
\tb{a}\wedge\tb{P}=0,
\end{equation}
respectively.

For instance, the linear subspace \(S_{L^*}\) shown in Figure~\ref{embedding T2 in R3 another example}
can be represented by \(\tb{s}=\e_0+\frac{1}{2}\e_1-\frac{1}{3}\e_2\).
Substituting \(\tb{s}\) and \(\tb{a}=d\e_0+a\e_1+b\e_2\)  into Equation~(\ref{a from s})
yields
\[
\begin{split}
&(d\e_0+a\e_1+b\e_2)\wedge (\e_0+\tfrac{1}{2}\e_1-\tfrac{1}{3}\e_2)=\\
&\tfrac{d}{2}\e_{01}-\tfrac{d}{3}\e_{02}
+a\e_{10}-\tfrac{a}{3}\e_{12}+b\e_{20}+\tfrac{b}{2}\e_{21}=\\
&(\tfrac{b}{2}+\tfrac{a}{3})\e_{21}
+(-\tfrac{d}{3}-b)\e_{02}
+(a-\tfrac{d}{2})\e_{10}
=0,
\end{split}
\]
which is equivalent to 
\[
a=\frac{d}{2},\ b=\frac{-d}{3}.
\]
To find the intersection of \(S_{L^*}\) with \T{2*}, I impose \(d=1\) and obtain \(a=\tfrac{1}{2},\ b=-\tfrac{1}{3}\), 
i.e. the point \((\tfrac{1}{2}, -\tfrac{1}{3})\) in \T{2*}.

For a blade \(A_k\) in \(\bigwedge\R{3*}\), I define the linear subspace represented by the blade as follows:
\begin{equation}
\Sub{A_k}=
\{
\tb{a}\in\R{3*}
|
\tb{a}\wedge A_k=0
\}.
\end{equation}
Similar definition applies to  blades in \(\bigvee\R{3}\).

\subsection{Visualising bivectors and trivectors}

For brevity, I consider \(\bigvee\R{3}\) only but exactly the same ideas apply to \(\bigwedge\R{3*}\).

A vector \(\tb{v}\) is completely characterised by its three defining properties:
\vspace{-10pt}\begin{itemize}\setlength{\itemsep}{0mm} 
\item[1)] The subspace that vector \(\tb{v}\) represents in \R{3}.
The subspace consists of vectors \(\tb{x}\in\R{3}\) which satisfy \(\tb{x}\vee \tb{v}=0\).
This defines a line  passing through the origin of \R{3}.
This subspace is called the \textbf{attitude} of \(\tb{v}\).
\item[2)] There are two possible directions along the attitude and
 \(\tb{v}\) points in one of these directions, which is the vector's \textbf{orientation}.
Vectors with the same attitude but pointing in the opposite directions have opposite orientations.
\item[3)] There are many vectors distinct from \(\tb{v}\) that share the same attitude and orientation.
For instance, \(2\tb{v}\) has the same attitude and orientation, but it is obviously different from \(\tb{v}\).
All such vectors are distinguished by a positive scalar factor, which is called the \textbf{weight}.
A vector's weight is a relative concept. It relates two vectors that share the attitude and orientation.
For instance, the weight of \(2\e^0\) with respect to \(\e^0\) equals 2.
\end{itemize}

Given a certain attitude and orientation, it is often convenient to pick a reference vector and define the weight of other equally oriented vectors
in the same attitude  with respect to the chosen reference vector.
It may sometimes be convenient to incorporate the orientation into weight by considering negative weights, 
which gives vectors in the same attitude but with the opposite orientation.
Weight is a pre-metric concept, it does not require  a metric.
Even when a metric is defined, a vector's weight is not related and may not be equal to the length of the vector.
For instance, there are non-zero vectors in Minkowski space whose length is zero,
but their weight is not zero; only the origin has zero weight.
Generally it is not meaningful to compare the weights of vectors with different attitudes.
Vectors' lengths and angles between them, as drawn on paper or screen, do not necessarily refer to their metric lengths and angles.
I use the weight to decide a vector's length and draw vectors in different subspaces as if their weight is comparable.

For scalars, the attitude is the origin of \R{3}, the orientation is the sign, and the weight is the scalar's value with respect to 1 or \(-1\), e.g.\
\(-3\) has the same orientation as \(-1\) and the weight of 3 with respect to \(-1\). 

Vectors and scalars are defined by their attitude, orientation, and weight.
The same characteristics can be used to define and visualise simple bivectors and other blades.

\begin{figure}[t!]
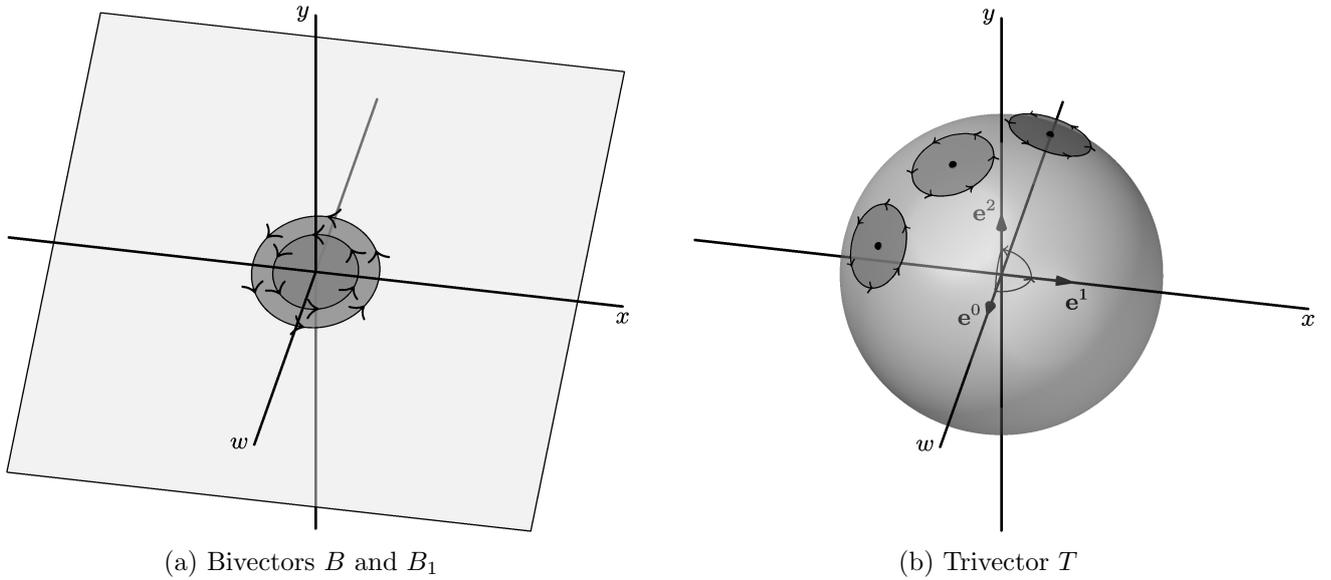

\begin{subfloatenv}{Bivectors \(B\) and \(B_1\)}
\begin{asy}
import Drawing3D;
DrawingR3 drawing = DrawingR3(4,0.1, camera=(10,2,7));	
drawing.target_axes();

triple n = (1,0,2);
triple u = (0,1,0);
triple v = cross(n,u);

path3 p = plane_path3(n, (0,0,0), theta=0, scale=7, v=v);
drawing.Drawing3D.plane(p);

bivector(u,v);
bivector(u,v/sqrt(5),45);

\end{asy}
\end{subfloatenv}\hspace{-12pt}%
\begin{subfloatenv}{Trivector \(T\)}
\begin{asy}
import Drawing3D;
DrawingR3 drawing = DrawingR3(4,0.1, camera=(10,2,7));

import solids;
real r = 2.12156883589411;
currentlight=light(white, viewport=false,(5,-5,10));
draw(surface(sphere((0,0,0),r,100)),lightgray+opacity(0.4));

drawing.target_axes();
drawing.target_basis();

tangent_bivector(r,-50,20);
tangent_bivector(r,-50,60);
tangent_bivector(r,180-50,20);

path3 a= arc((0,0,0),(0.4,0,0),(0,0.4,0));
draw(a,Arrow3(TeXHead2((0,0,1))));

path3 a= arc((0,0,0),(0,0.4,0),(0,0,0.4));
draw(a,Arrow3(TeXHead2((1,0,0))));

path3 a= arc((0,0,0),(0,0,0.4),(0.4,0,0));
draw(a,Arrow3(TeXHead2((0,1,0))));

\end{asy}
\end{subfloatenv}
\caption{Visualising oriented blades}
\label{visualising blades}
\end{figure}

The attitude of a bivector is a 2-dimensional subspace it represents and
its orientation is a specific sense of rotation around the origin.
I visualise bivectors as filled circles, disks, centred on the origin and lying in the bivector's subspace.
The circle's area, understood in Euclidean sense, represents the bivector's weight with respect to a reference bivector in the same attitude.
The orientation is visualised with arrows indicating a certain direction of motion along the circle.
A bivector is thus oriented with vectors attached to the boundary of the disk representing the bivector 
(actually a single vector is enough).
A disk is a 2-dimensional object and its boundary is a circle, which is 1-dimensional, and, therefore, 
it can be oriented with a 1-dimensional object, i.e.\ a vector. 

Consider the  bivector \(B=\e^1\vee(\e^2-2\e^0)\) in \(\bigvee\R{3}\) shown in Figure~\ref{visualising blades}(a).
In accord with the Euclidean conventions adopted for the purpose of visualisation, 
I use the bivector 
\(B_1=\e^1\vee(\tfrac{1}{\sqrt{5}}\e^2-\tfrac{2}{\sqrt{5}}\e^0)\)
as a reference bivector in this subspace.
The weight of \(B\) with respect to \(B_1\) equals \(\sqrt{5}\), i.e. \(B=\sqrt{5}B_1\), 
so I draw the circle representing \(B\) with the area equal to \(\sqrt{5}\).
The area of the smaller circle representing \(B_1\) equals 1.

Trivectors represent 3-dimensional subspaces.
A trivector can be depicted as a filled sphere with its volume representing the trivector's weight.
Its orientation is defined by attaching a specific bivector  to the sphere, which serves as the boundary,
i.e.\ a 3-dimensional object is oriented with a 2-dimensional object.
A trivector is shown in Figure~\ref{visualising blades}(b) (three orienting bivectors are shown, but one is sufficient).
Its orientation is consistent with a rotation from \(\e^0\) to \(\e^1\), 
followed by a rotation from \(\e^1\) to \(\e^2\), and then from \(\e^2\) back to \(\e^0\)
as shown in Figure~\ref{visualising blades}(b).
The volume of the sphere, representing the trivector's weight, equals 40.
So the depicted trivector \(T=40\e^0\vee\e^1\vee\e^2\).

The way bivectors and trivectors are oriented is applicable to vectors and other blades in general.
A vector, which is a line segment, is oriented by attaching a specific 0-dimensional orientation, i.e.\ scalar orientation (\(+1\) or \(-1\)),
to the boundary of the line segment.
A line segment is depicted in Figure~\ref{line segment}.
I orient it by assigning the negative scalar orientation to  the lower left end and the positive scalar orientation to the upper right end.
This yields  \(2\e^0+\e^1\) shown in Figure~\ref{vector}, 
where it is depicted with the negatively oriented end point at the origin.
\begin{wrapfigure}{l}{0.5\textwidth}
\begin{minipage}{0.25\textwidth}
\begin{asy}
import Drawing2D;
DrawingR2 drawing = DrawingR2(2, 0.1);
drawing.Drawing2D.grid();
drawing.target_axes();
drawing.target_basis();

pair p0 = (-1,-0.5);
pair p1 = (+1,0.5);
path l = p0--p1;

draw(l);
label("$-1$",p0,SW);
label("$+1$",p1,NE);

\end{asy}
\caption{\(2\e^0+\e^1\)}
\label{line segment}
\end{minipage}%
\begin{minipage}{0.25\textwidth}
\hspace{5pt}
\begin{asy}
import Drawing2D;
DrawingR2 drawing = DrawingR2(2, 0.1);
drawing.Drawing2D.grid();
drawing.target_axes();
drawing.target_basis();

pair p0 = (-1,-0.5);
pair p1 = (+1,0.5);
path l = p0--p1;

draw(shift(1,0.5)*l,Arrow);
\end{asy}
\caption{same}
\label{vector}
\end{minipage}
\end{wrapfigure}
This follows traditional practice of visualising vectors, 
but the way this vector is drawn in Figure~\ref{line segment} is more consistent with the depiction of other blades.
Indeed, if bivectors and trivectors are visualised as rounds (circles and spheres), 
it is natural to depict vectors as point pairs centred on the origin.

Bivectors and trivectors are frequently drawn as parallelograms and parallelepipeds, 
but this way of depiction ties them to the outer product of specific vectors.
A given bivector or trivector can be represented by a multitude of different vectors.
For instance, \(\e^0\vee\e^1\) and \(\e^0\vee(\e^1-2\e^0)\) represent the same bivector, but the corresponding parallelograms are distinct.
This property may be useful in applications, but when a bivector or trivector is not tied to specific vectors, 
it is preferable to draw them as an oriented circle or sphere.

\subsection{The metric and Clifford algebra}
In the projective model, the inner product is defined in the dual model space \R{3*}.
I assume that the standard basis vectors  are orthogonal, i.e.\ \(\e_i\cdot\e_j=0\) for \(i\ne j\).
The inner product of vectors \(\tb{a}_1\) and \(\tb{a}_2\)  is given by
\begin{equation}
\tb{a}_1\cdot\tb{a}_2=d_1d_2\e_0\cdot\e_0+a_1a_2\e_1\cdot\e_1+b_1b_2\e_2\cdot\e_2,
\end{equation}
where \(\tb{a}_1=d_1\e_0+a_1\e_1+b_1\e_2\) and  \(\tb{a}_2=d_2\e_0+a_2\e_1+b_2\e_2\).
As will be demonstrated below, the inner product provides the measure of angles and distances between geometric objects in the target space.
Thus, the target space becomes a metric space\footnote{Metric spaces defined in the projective model are known as Cayley-Klein geometries.}.

The following list of options for the inner product defines some metric spaces of interest:
\begin{center}
\begin{tabular}{ccccl}
\(\e_0\cdot\e_0\) & \(\e_1\cdot\e_1\) & \(\e_2\cdot\e_2\)  & \T{2} & \\
0 & 1 & 1 &\E{2}& Euclidean \\ 
 1& 1& 1& \El{2} & Elliptic (spherical)\\
\m1 & 1& 1 & \Hy{2} & Hyperbolic \\
0 & 1  & \m1 & \M{2} &Minkowski (pseudo-Euclidean)  \\
1 & 1 & \m1 & \dS{2} & de-Sitter \\
 \m1 & 1 & \m1 & \AdS{2} & Anti de-Sitter 
\end{tabular}
\end{center}
For instance, \(\tb{a}_1\cdot\tb{a}_2=+a_1a_2+b_1b_2\) in Euclidean space and
\(\tb{a}_1\cdot\tb{a}_2=-d_1d_2+a_1a_2-b_1b_2\)  in Anti de-Sitter space.
The metric is called degenerate if \(\e_i\cdot\e_i=0\) for some \(i=0,1,2\); the corresponding metric space is also called degenerate.
The sign of  \(\e_0\cdot\e_0\) is indicative of the curvature, while \(\e_0\cdot\e_0=0\) implies that the space is flat (zero curvature).
Euclidean and Minkowski spaces are flat and degenerate; the other spaces listed above are non-degenerate.
Elliptic space is closely related to spherical geometry 
(without going into details, I note that Elliptic space is equivalent to spherical space where the antipodal points are identified).
Minkowski, de-Sitter, and Anti de-Sitter spaces admit a relativistic kinematic structure.


The Clifford algebra of \R{3*} is an algebra of multivectors in \(\bigwedge\R{3*}\) with the geometric product
that has the following properties on arbitrary multivectors \(A,B\), and \(C\):\\
1) \(A(B+C)=AB+AC\) and \((B+C)A=BA+CA\) (distributive),\\
2) \(A(BC)=(AB)C\) (associative),\\
3) \(\alpha A=A\alpha\) for \(\alpha\in\R{}\) (commutative for scalars),\\
4) \(\tb{a}\tb{b}=\tb{a}\cdot\tb{b}+\tb{a}\wedge\tb{b}\) for vectors \(\tb{a}, \tb{b}\in\R{3*}\).\\
(following conventions, I do not use any symbol for the geometric product;
so, the absence of any symbol between multivectors implies geometric multiplication).

Unlike the inner and outer products, the geometric product is generally neither commutative nor anti-commutative.
In \(\bigwedge\R{3*}\), the geometric product between a \(k\)-vector \(A_k\) and an \(l\)-vector \(B_l\) generally includes multivectors
of grades \(|k-l|\) and \(k+l\) (the latter only if \(k+l<4\)).
For instance, the geometric product of two vectors (grade-1) includes a scalar (grade-0) and a bivector (grade-2).
The outer product of \(A_k\) and \(B_l\) is inherited from the Grassmann algebra and satisfies 
\begin{equation}
A_k\wedge B_k=\grade{A_kB_l}_{k+l}.
\end{equation}
The inner product of \(A_k\) and \(B_l\) is defined by
\begin{equation}
A_k\cdot B_l=\grade{A_kB_l}_{|k-l|}.
\end{equation}
This definition extends by linearity to arbitrary multivectors, e.g. \((\tb{a}+\tb{P})\cdot\tb{Q}=\tb{a}\cdot\tb{Q}+\tb{P}\cdot\tb{Q}\),
where \(\tb{a}\cdot\tb{Q}=\grade{\tb{a}\tb{Q}}_1\) is a vector and \(\tb{P}\cdot\tb{Q}=\grade{\tb{P\tb{Q}}}_0\) is a scalar.
For arbitrary multivectors, the commutator is defined by
\begin{equation}
A\times B=\tfrac{1}{2}(AB-BA).
\end{equation}
The inner and outer products satisfy the following useful equalities independent of dimension:
\begin{equation}
\begin{split}
&\tb{a}\cdot(\tb{b}\wedge\tb{c})=(\tb{a}\cdot\tb{b})\tb{c}-(\tb{a}\cdot\tb{c})\tb{b},\\
&\tb{a}\cdot(\tb{b}\wedge\tb{c}\wedge\tb{d})
=
(\tb{a} \cdot\tb{b} ) \tb{c}\wedge \tb{d}
-(\tb{a} \cdot\tb{c} )\tb{b} \wedge\tb{d}
+(\tb{a} \cdot \tb{d}) \tb{b}\wedge\tb{c},\\
&(\tb{a} \wedge \tb{b})\cdot(\tb{c} \wedge \tb{d})
=
(\tb{a} \cdot\tb{d} )(\tb{b} \cdot \tb{c})
-
(\tb{a} \cdot \tb{c})(\tb{b} \cdot \tb{d}),\\
&(\tb{a}\wedge\tb{b})\cdot(\tb{c}\wedge\tb{d}\wedge\tb{g})
=
\left((\tb{a} \wedge \tb{b})\cdot( \tb{c}\wedge \tb{d})\right)\tb{g}
-\left((\tb{a} \wedge \tb{b})\cdot( \tb{c}\wedge \tb{g})\right)\tb{d}
+\left((\tb{a} \wedge \tb{b})\cdot( \tb{d}\wedge \tb{g})\right)\tb{c},\\
&(\tb{a} \wedge \tb{b})\times(\tb{c} \wedge \tb{d})
=
(\tb{b}\cdot \tb{c}) \tb{a}\wedge \tb{d} 
+ (\tb{a}\cdot \tb{d}) \tb{b}\wedge \tb{c} 
- (\tb{b}\cdot \tb{d}) \tb{a}\wedge \tb{c}
- (\tb{a}\cdot \tb{c}) \tb{b}\wedge \tb{d}.
\end{split}
\end{equation}

Some of the properties of the geometric product depend on dimension.
In \(\bigwedge\R{3*}\), the following applies:
\begin{equation}
\begin{split}
&\tb{a}\tb{b}=\tb{a}\cdot\tb{b}+\tb{a}\wedge\tb{b} \textrm{ (scalar + bivector)},\quad
\tb{a}\cdot\tb{b}=\tb{b}\cdot\tb{a},\quad
\tb{a}\wedge\tb{b}=-\tb{b}\wedge\tb{a},\\
&\tb{a}\tb{P}=\tb{a}\cdot\tb{P}+\tb{a}\wedge\tb{P} \textrm{ (vector + trivector)},\quad
\tb{a}\cdot\tb{P}=-\tb{P}\cdot\tb{a},\quad
\tb{a}\wedge\tb{P}=\tb{P}\wedge\tb{a},\\
&\tb{a}\I=\tb{a}\cdot\I \textrm{ (bivector)},\quad
\tb{a}\cdot\I=\I\cdot\tb{a},\\
&\tb{P}\tb{Q}=\tb{P}\cdot\tb{Q}+\tb{P}\times\tb{Q} \textrm{ (scalar + bivector)},\quad
\tb{P}\cdot\tb{Q}=\tb{Q}\cdot\tb{P},\\
&\tb{P}\I=\tb{P}\cdot\I \textrm{ (vector)},\quad
\tb{P}\cdot\I=\I\cdot\tb{P}.
\end{split}
\end{equation}

It follows that \(M\I=\I M\) for any multivector \(M\) and
\begin{equation}
\begin{split}
&\tb{a}\cdot\tb{b}=\tfrac{1}{2}(\tb{a}\tb{b}+\tb{b}\tb{a}),\quad 
\tb{a}\wedge\tb{b}=\tfrac{1}{2}(\tb{a}\tb{b}-\tb{b}\tb{a}),\quad \tb{a}\times\tb{b}=\tb{a}\wedge\tb{b},  \\
&\tb{a}\cdot\tb{P}=\tfrac{1}{2}(\tb{a}\tb{P}-\tb{P}\tb{a}),\quad
\tb{a}\wedge\tb{P}=\tfrac{1}{2}(\tb{a}\tb{P}+\tb{P}\tb{a}),\quad  \tb{a}\times\tb{P}=\tb{a}\cdot\tb{P}  \\
&\tb{P}\cdot\tb{Q}=\tfrac{1}{2}(\tb{P}\tb{Q}+\tb{Q}\tb{P}),\quad 
\tb{P}\times\tb{Q}=\tfrac{1}{2}(\tb{P}\tb{Q}-\tb{Q}\tb{P}).\\
\end{split}
\end{equation}
The reverse of multivector \(M=s+\tb{a}+\tb{P}+\alpha\I\) is given by \(\reverse{M}=s+\tb{a}-\tb{P}-\alpha\I\).
The reverse of the geometric product of multivectors \(A\) and \(B\) satisfies \(\reverse{AB}=\reverse{B}\reverse{A}\).
The grade involution of multivector \(M\) is given by \(\involute{M}=s-\tb{a}+\tb{P}-\alpha\I\) and satisfies \(\involute{AB}=\involute{A}\involute{B}\)
for any multivectors \(A\) and \(B\). 
Since \(\reverse{\grade{M}}_k=\grade{\reverse{M}}_k\) and \(\involute{\grade{M}}_k=\grade{\involute{M}}_k\),
the outer and inner products satisfy similar relationships, i.e.\ 
\(\reverse{A\wedge B}=\reverse{B}\wedge\reverse{A}\), \(\reverse{A\cdot B}=\reverse{B}\cdot\reverse{A}\), 
\(\involute{A\wedge B}=\involute{A}\wedge\involute{B}\),   and \(\involute{A\cdot B}=\involute{A}\cdot\involute{B}\).
Furthermore, \(\reverse{A\times B}=\reverse{B}\times\reverse{A}\) and \(\involute{A\times B}=\involute{A}\times\involute{B}\) follow from the definition
of the commutator.

The norm of multivector \(M\) is defined by
\begin{equation}
\norm{M}=|M\reverse{M}(M\reverse{M})_1|^{\frac{1}{4}},
\end{equation}
where \((X)_1\) reverses only the grade-1 component of multivector \(X\), e.g. \((1+\e_2)_1=1-\e_2\).
Note, however, it does not satisfy the triangle inequality in general, e.g. \(\norm{(\e_0+\e_1)+(\e_0-\e_1)}=2\), but \(\norm{\e_0+\e_1}+\norm{\e_0-\e_1}=0\) in \Hy{2},
and \(\norm{M}=0\) does not necessarily imply \(M=0\).
Nevertheless, the utility of this nomenclature will become apparent in the following.
If the norm is not zero, \(M\) has the inverse given by
\begin{equation}
M^{-1}=
\frac{\reverse{M}(M\reverse{M})_1}{M\reverse{M}(M\reverse{M})_1}.
\end{equation}
For a blade \(A_k\), the product \(A_k\reverse{A_k}\) yields a scalar, so \(\norm{A_k}=|A_k\reverse{A_k}|^\frac{1}{2}\) and 
\(
A_k^{-1}=\reverse{A_k} / (A_k\reverse{A_k}).
\)
The same applies to even multivectors and Study numbers.

The existence and form of the inverse depends on the metric, e.g., for \(\tb{a}=\e_0+\e_1\),
the inverse \(\tb{a}^{-1}=\tb{a}\) in Euclidean space and \(\tb{a}^{-1}=\tfrac{1}{2}\tb{a}\) in Elliptic space,
since \(\tb{a}\reverse{\tb{a}}=\tb{a}\tb{a}=(\e_0+\e_1)(\e_0+\e_1)=\e_0^2+\e_1^2=\e_0\cdot\e_0+\e_1\cdot\e_1\)
which equals 1 in Euclidean space and 2 in Elliptic space.
In Minkowski space,  \(\tb{a}=\e_1+\e_2\) lacks the inverse, since \(\tb{a}\reverse{\tb{a}}=\e_1^2+\e_2^2=1-1=0\).
Blades that have a non-zero component along \(\e_0\), e.g.\ \(\e_{01}\), 
do not have the inverse in Euclidean and Minkowski spaces.
The pseudoscalar \(\I\) does not have the inverse in degenerate spaces, since
\(\I\reverse{\I}=\e_0\e_1\e_2\e_2\e_1\e_0=\e_0^2\e_1^2\e_2^2=0\).
On the other hand, in non-degenerate spaces, the inverse \(\I^{-1}\) exists, e.g.\ \(\I^{-1}=\reverse{\I}\) in Elliptic space.
The inverse of an even multivector \(E=1+\e_{12}\) is given by \(E^{-1}=\tfrac{1}{2}(1-\e_{12})\)
in Euclidean space.
Indeed, \(EE^{-1}=(1+\e_{12})\tfrac{1}{2}(1-\e_{12})=\tfrac{1}{2}(1-\e_{12}+\e_{12}-\e_{12}\e_{12})=\tfrac{1}{2}(1+\e_{12}\e_{21})=
\tfrac{1}{2}(1+\e_1^2\e_2^2)=\tfrac{1}{2}(1+1)=1\).

The exponential function is defined on an arbitrary multivector \(X\) by
\begin{equation}
e^X=\sum_{n=0}^{\infty}\frac{X^n}{n!}.
\end{equation}
\(e^{X+Y}=e^Xe^Y\) holds for commuting multivectors, i.e.~\(XY=YX\).
The Taylor expansion can be used to define various other functions of multivectors, such as
\begin{equation}
\sin{X}=\sum_{n=0}^{\infty}\frac{(-1)^nX^{2n+1}}{(2n+1)!},
\quad\quad
\cos{X}=\sum_{n=0}^{\infty}\frac{(-1)^nX^{2n}}{(2n)!}.
\end{equation}
These functions depend on the metric.
For instance, \(\sin{\e_0}=\e_0\sin{1}\) in Elliptic space and \(\sin{\e_0}=\e_0\) in Euclidean space.

Even multivectors form a subalgebra, since the sum and the geometric product of two even multivectors is an even multivector.
In a 3-dimensional space, an even multivector can be written as the sum of a scalar and a bivector, e.g.
\begin{equation}
E=s+p_{12}\e_{12}+p_{20}\e_{20}+p_{01}\e_{01}.
\end{equation}
In Elliptic space, the subalgebra of even multivectors is isomorphic to the algebra of quaternions.

\subsection{Orthogonality in \R{3*} with elliptic metric}
In \(\bigwedge\R{3*}\) with the elliptic metric, the pseudoscalar \(\I\) is invertible and the expression
\begin{equation}
A_k^{\perp}=A_k \I^{-1}
\label{orthogonal subspace}
\end{equation}
defines a blade \(A_k^{\perp}\) whose subspace is orthogonal to the subspace of the blade \(A_k\), i.e.
\begin{equation}\label{define orthogonal subspace}
\tb{p}\cdot \tb{a}=0 \textrm{ for any } \tb{p}\in\Sub{A_k^{\perp}} \textrm{ and }  \tb{a}\in\Sub{A_k}.
\end{equation}
In \(\bigwedge\R{3*}\), the grade of blade \(A_k^{\perp}\) equals \(3-k\).

\begin{figure}[t!]
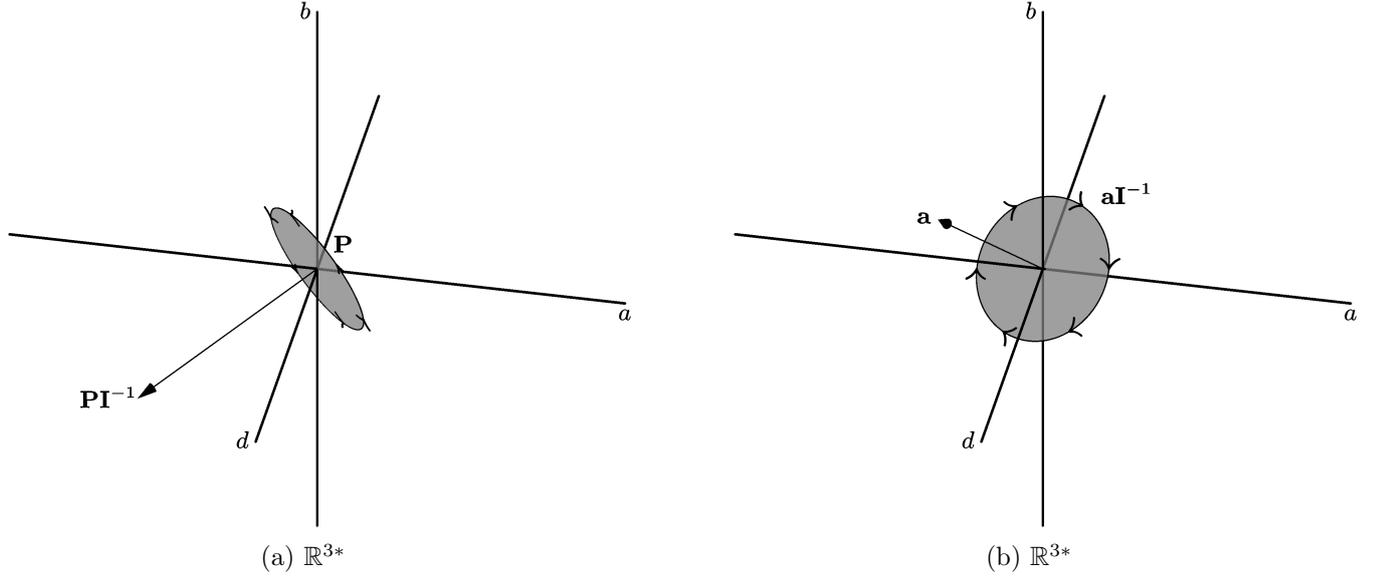

\begin{subfloatenv}{\R{3*}}
\begin{asy}
import Drawing3D;
DrawingR3 drawing = DrawingR3(4,0.1);
drawing.dual_axes();

triple a = (2,-2,-1);
drawing.Drawing3D.vector(a);
label("$\textbf{P}\textbf{I}^{-1}$", a, (0,-0.5,0));

triple u = (1,0,2);
triple v = (1,1,0);
bivector(v,u,40);
label("$\textbf{P}$",(0.8,0.5,1));

\end{asy}
\end{subfloatenv}
\begin{subfloatenv}{\R{3*}}
\begin{asy}
import Drawing3D;
DrawingR3 drawing = DrawingR3(4,0.1);
drawing.dual_axes();

triple a = (2,-1,2);
drawing.Drawing3D.vector(a);
label("$\textbf{a}$", a, (0,-1,0));

triple u = (1,0,-1);
triple v = (1,2,0);

bivector(v,u,40);
label("$\textbf{a}\textbf{I}^{-1}$",(-0.5,1,1));

\end{asy}
\end{subfloatenv}
\caption{Orthogonality in \(\bigwedge\R{3*}\) with elliptic metric}
\label{elliptic orthogonality examples}
\end{figure}

For instance, the bivector \(\e_{01}\) represents the \(da\)-plane, the plane that contains \(\e_0\) and \(\e_1\).
Its orthogonal subspace is represented by 
\(\e_{01}\I^{-1}=\e_{01}\e_{210}=\e_{01}\e_{102}=\e_0^2\e_1^2\e_2=\e_2\).
Similarly, the subspace orthogonal to \(\e_1\) is represented by 
\(\e_{1}\I^{-1}=\e_{1}\e_{210}=-\e_1\e_{120}=-\e_1^2\e_{20}=-\e_{20}\).
A scalar represents a 0-dimensional subspace that consists of the zero vector only, i.e.\ the origin of \R{3*}.
Since the inner product with the zero vector is zero for any vector in \R{3*},
the origin of \R{3*} represented by a scalar is orthogonal to the whole of \R{3*} represented by the pseudoscalar, and vice versa.
For instance, \(2^{\perp}=2\I^{-1}=2\e_{210}=-2\e_{012}\) and \(\I^{\perp}=\I\I^{-1}=1\).

Orthogonality is illustrated for \(\tb{a}=2\e_0-\e_1+2\e_2\) and 
\(\tb{P}=2\e_{12}-2\e_{20}-\e_{01}=(\e_0+\e_1)\wedge(\e_0+2\e_2) \) in Figure~\ref{elliptic orthogonality examples},
where \(\tb{a}\I^{-1}=(2\e_0-\e_1+2\e_2)\e_{210}=-2\e_{12}+\e_{20}-2\e_{01}=(\e_0+2\e_1)\wedge(\e_0-\e_2) \)
and \(\tb{P}\I^{-1}=(2\e_{12}-2\e_{20}-\e_{01})\e_{210}=2\e_0-2\e_1-\e_2\).

The formula~(\ref{orthogonal subspace}) can be generalised to arbitrary multivectors in \(\bigwedge\R{3*}\) by
\begin{equation}
M^{\perp}=M\I^{-1}.
\label{orthogonality for M}
\end{equation}
Since \( (\alpha A+\beta B)^\perp=\alpha A^\perp+\beta B^\perp \) for any multivectors \(A\), \(B\) and \(\alpha,\beta\in\R{}\), 
\(M^{\perp}\) can be computed using linearity and the following table 

\begin{tabular}{l|cccccccc}
\(M\) &         \(1\)       & \(\e_0\)&    \(\e_1\)&         \(\e_2\)&    \(\e_{12}\)& \(\e_{20}\)& \(\e_{01}\)& \(\e_{012}\)\\
\hline
\(M\I^{-1}\) & \(-\e_{012}\) & \(-\e_{12}\)& \(-\e_{20}\)&  \(-\e_{01}\)&    \(\e_0\)&    \(\e_{1}\)& \(\e_{2}\)& \(1\) \\
\end{tabular}

Since~(\ref{orthogonality for M}) relies on the geometric product,  it cannot be applied in a non-metric context.
However, it is easy to define a transformation \(O:\bigwedge\R{3*}\to\bigwedge\R{3*}\)
that acts in the same way as \(M^\perp\), i.e.\ \(O(M)=M^{\perp}\), but does not use the geometric product.
The orthogonality transformation \(O\) is defined on the basis blades as follows:

\begin{tabular}{l|cccccccc}
\(M\) &         \(1\)       & \(\e_0\)&    \(\e_1\)&         \(\e_2\)&    \(\e_{12}\)& \(\e_{20}\)& \(\e_{01}\)& \(\e_{012}\)\\
\hline
\(O(M)\) & \(-\e_{012}\) & \(-\e_{12}\)& \(-\e_{20}\)&  \(-\e_{01}\)&    \(\e_0\)&    \(\e_{1}\)& \(\e_{2}\)& \(1\) \\
\end{tabular}

In addition, \(O\) is assumed to be linear, i.e.
\begin{equation}
O(\alpha A+\beta B)=\alpha O(A) + \beta O(B)
\end{equation}
for arbitrary multivectors \(A,B\) and \(\alpha,\beta\in\R{}\).
So, \(O(M)\) can be computed for any multivector \(M\) by using linearity and applying the table above.
The transformation \(O\) is thus independent of the metric and \(O(M)=M^\perp\) if the metric is elliptic.
Note that applying the transformation twice results in the same multivector but with the opposite orientation, i.e.\ 
\(O(O(M))=-M\).

Table

\begin{tabular}{l|cccccccc}
\(M\) &            \(1\)      & \(\e_0\)&    \(\e_1\)&   \(\e_2\)&    \(\e_{12}\)& \(\e_{20}\)& \(\e_{01}\)& \(\e_{012}\)\\
\hline \\ [-10pt]
\(O^{-1}(M)\) & \(\e_{012}\) & \(\e_{12}\)& \(\e_{20}\)& \(\e_{01}\)& \(-\e_0\)&    \(-\e_{1}\)& \(-\e_{2}\)& \(-1\) \\
\end{tabular}

defines the inverse \(O^{-1}:\bigwedge\R{3*}\to\bigwedge\R{3*}\) on the basis blades,
which extends to arbitrary multivectors by linearity.
\(O^{-1}(O(M))=M\) as expected.

\subsection{Duality transformation on multivectors and the join}
A 2-dimensional subspace of the target model space \R{3} is defined by
\begin{equation}
dw+ax+by=0,
\label{2d subspace in R3}
\end{equation}
where \(d,a,b\) are fixed coefficients and \(w,x,y\) are coordinates of a variable vector \(\tb{x}=w\e^0+x\e^1+y\e^2\).
This equation is not metric, but due to its form it can be written in \R{3*} with the elliptic metric as
\begin{equation}
\tb{a}\cdot\tb{x}^I=0,
\end{equation}
where \(\tb{a}=d\e_0+a\e_1+b\e_2\) is a fixed vector in \R{3*} and 
\(\tb{x}^I=\Id^{-1}(\tb{x})=w\e_0+x\e_1+y\e_2\) is a variable vector in \R{3*} (\(\Id^{-1}\)~is the inverse of the identity transformation).
It follows that the linear subspace represented by \(\tb{x}^I\) is orthogonal to \(\Sub{\tb{a}}\).
Therefore, this subspace can be represented by the blade \(O(\tb{a})\) in \(\bigwedge\R{3*}\).
Moreover, the blade \(\Id(O(\tb{a}))\) in \(\bigvee\R{3}\) represents a linear subspace of \R{3}
consisting of vectors \(\tb{x}=w\e^0+x\e^1+y\e^2\) which satisfy Equation~(\ref{2d subspace in R3}).
For instance, consider  the \(wz\)-plane in \R{3} defined by \(x=0\).
In this case, \(\tb{a}=\e_1\) and \(\Id(O(\e_1))=\Id(-\e_{20})=-\e^{20}\),
but the subspace of \R{3} represented by \(-\e^{20}\) coincides with the \(wz\)-plane.

Imposing additional constraints \(d=1\) and \(w=1\) gives \(ax+by+1=0\), which defines the duality relationship
between the line \(L=\{(x,y)|ax+by+1=0\}\) in \T{2} and the point \((a,b)\) in \T{2*}. 
Exactly the same relationship is captured by \(\tb{a}\in\R{3*}\), which represents \((a,b)\),
and \(\Id(O(\tb{a}))\) in \(\bigvee\R{3}\),
which represents \(L\).

A 1-dimensional subspace of the target model space \R{3} can be defined as the intersection of two 2-dimensional subspaces:
\begin{equation}
\left\{
\begin{aligned}
&d_1w+a_1x+b_1y=0,\\
&d_2w+a_2x+b_2y=0,\\
\end{aligned}
\right.
\label{1d subspace in R3}
\end{equation}
where \(d_1,a_1,b_1\) and \(d_2,a_2,b_2\) are fixed coefficients and \(w,x,y\) are coordinates of a variable vector \(\tb{x}=w\e^0+x\e^1+y\e^2\).
This 1-dimensional subspace represents a point in the target space \T{2}
and, if the metric is elliptic, consists of vectors \(\tb{x}\) which satisfy both 
\begin{equation}
\tb{a}_1\cdot\tb{x}^I=0 \textrm{ and } \tb{a}_2\cdot\tb{x}^I=0,
\end{equation}
where \(\tb{a}_1=d_1\e_0+a_1\e_1+b_1\e_2\) and  \(\tb{a}_2=d_2\e_0+a_2\e_1+b_2\e_2\).
These two equations are equivalent to
\begin{equation}
(\tb{a}_1\wedge\tb{a}_2)\cdot\tb{x}^I=0.
\end{equation}
So, the subspace represented by \(\tb{x}^I\) is orthogonal to \(\Sub{\tb{a}_1\wedge\tb{a}_2}\).
It can be represented by \(O(\tb{a}_1\wedge\tb{a}_2)\), which is a vector in \R{3*}.
Then, the vector \(\Id(O(\tb{a}_1\wedge\tb{a}_2))\) represents a 1-dimensional subspace of \R{3}
defined by~(\ref{1d subspace in R3}).
For instance, \(x=0\) and \(y=0\) define the \(w\)-axis in \R{3}.
In this case, \(\tb{a}_1=\e_1\), \(\tb{a}_2=\e_2\), and \(\Id(O(\e_1\wedge\e_2))=
\Id(O(\e_{12}))=
\Id(\e_0)=\e^0\),
but \(\e^0\) represents a linear subspace of \R{3} which coincides with the \(w\)-axis.

The subspaces of dimension 0 and 3 are special cases, but they conform to the same pattern observed above.
The 0-dimensional subspace, i.e.\ the origin, can be defined as the intersection of three 2-dimensional
subspaces:
\begin{equation}
\left\{
\begin{aligned}
&d_1w+a_1x+b_1y=0,\\
&d_2w+a_2x+b_2y=0,\\
&d_3w+a_3x+b_3y=0.\\
\end{aligned}
\right.
\label{0d subspace in R3}
\end{equation}
This 0-dimensional subspace can be represented by \(\Id(O(\tb{a}_1\wedge\tb{a}_2\wedge\tb{a}_3))\), which is a scalar.
So the origin of \R{3} can be represented dually by the trivector \(T=\tb{a}_1\wedge\tb{a}_2\wedge\tb{a}_3\)
and directly by the scalar~\(\Id(O(T))\).
The~3-dimensional subspace of \R{3} coincides with the whole space.
It can be defined by \(0w+0x+0y=0\), which corresponds to the origin of \R{3*} represented by a scalar \(\alpha\) in \(\bigwedge\R{3*}\).
So the 3-dimensional subspace can be represented dually by \(\alpha\)
and directly by \(\Id(O(\alpha))\), which is a trivector in \(\bigvee\R{3}\).

To summarise, \(\Id(O(A_k))\) directly represents the same linear subspace of \R{3} as that represented dually by 
a blade \(A_k\in\bigwedge\R{3*}\). 
For vectors and bivectors, these subspaces represent points and lines in \T{2} and \T{2*} related by duality.
Therefore, the duality transformation can be defined on blades in \(\bigwedge\R{3*}\) by \(\J(A_k)=\Id(O(A_k))\).
Both \(O\) and \(\Id\) are linear and admit general multivectors, so I can define the duality transformation
\(\J:\bigwedge\R{3*}\to\bigvee\R{3}\) on an arbitrary multivector \(X\in\bigwedge\R{3*}\) by
\begin{equation}
\J(X)=\Id(O(X)).
\end{equation}
On the basis blades, the duality transformation is given by

\begin{tabular}{l|cccccccc}
\(X\) &         \(1\)       & \(\e_0\)&    \(\e_1\)&         \(\e_2\)&    \(\e_{12}\)& \(\e_{20}\)& \(\e_{01}\)& \(\e_{012}\)\\
\hline \\ [-10pt]
\(\J(X)\phantom{^{-1}}\) & \(-\e^{012}\) & \(-\e^{12}\)& \(-\e^{20}\)&  \(-\e^{01}\)&    \(\phantom{-}\e^0\)&    \(\phantom{-}\e^{1}\)& \(\phantom{-}\e^{2}\)& \(\phantom{-}1\) \\
\end{tabular}

This table and the property of linearity can be used to compute \(J\) for any arbitrary multivector.
The inverse \(\J^{-1}:\bigvee\R{3}\to\bigwedge\R{3*}\) for an arbitrary multivector \(Y\in\bigvee\R{3}\) is given by
\begin{equation}
\J^{-1}(Y)=O^{-1}(\Id^{-1}(Y)).
\end{equation}
On the basis of \(\bigvee\R{3}\), I get

\begin{tabular}{l|cccccccc}
\(Y\) &    \(1\)                & \(\e^0\)&    \(\e^1\)&   \(\e^2\)&    \(\e^{12}\)& \(\e^{20}\)& \(\e^{01}\)& \(\e^{012}\)\\
\hline \\ [-10pt]
\(\J^{-1}(Y)\) & \(\phantom{-}\e_{012}\) & \(\phantom{-}\e_{12}\)& \(\phantom{-}\e_{20}\)& \(\phantom{-}\e_{01}\)& \(-\e_0\)&    \(-\e_{1}\)& \(-\e_{2}\)& \(-1\) \\
\end{tabular}

\(\J^{-1}(\J(X))=X\) as expected.

\emph{The join.}\\
The duality transformation makes it possible to define the join of two multivectors \(A\) and \(B\) in the dual model space \(\bigwedge\R{3*}\)
as follows
\begin{equation}
A\vee B = \J^{-1}(\J(A)\vee\J(B)),
\end{equation}
where I use the same symbol \(\vee\) for the join in \(\bigwedge\R{3*}\) as the one used for the outer product in \(\bigvee\R{3}\).
The join is metric-independent, since it only depends on the duality transformation and the outer product in \(\bigvee\R{3}\).

For a vector \(\tb{a}=d\e_0+a\e_1+b\e_2\) and a bivector \(\tb{P}=w\e_{12}+x\e_{20}+y\e_{01}\) in 
\(\bigwedge\R{3*}\) the following equality holds
\begin{equation}
\tb{a}\wedge\tb{P}=(\tb{a}\vee\tb{P})\I, \textrm{ where } \tb{a}\vee\tb{P}=dw+ax+by.
\end{equation}
Note that \(\tb{a}\vee\tb{P}=0\) if the point dually represented by \(\tb{P}\) lies on the line dually represented by \(\tb{a}\).

The join of bivectors \(\tb{P}=w_P\e_{12}+x_P\e_{20}+y_P\e_{01}\) and \(\tb{Q}=w_Q\e_{12}+x_Q\e_{20}+y_Q\e_{01}\) is given by
\begin{equation}
\tb{P}\vee\tb{Q}
=
(y_Px_Q-x_Py_Q)\e_0
+(w_Py_Q-w_Qy_P)\e_1
+(w_Qx_P-w_Px_Q)\e_2.
\end{equation}
It is a vector representing the line that passes through the points dually represented by \(\tb{P}\) and \(\tb{Q}\).
For \(w_P=1\) and \(w_Q=1\), the equation that defines the line becomes
\begin{equation}
(x_Q-x_P)y=(y_Q-y_P)x+x_Qy_P-x_Py_Q.
\end{equation}
The join of a multivector with the pseudoscalar yields the same multivector, i.e.~\(M\vee\I=M\).
In \(\bigwedge\R{3*}\), \(A_k\vee B_l=0\) if \(k+l<3\), where \(A_k\) is a \(k\)-vector and \(B_l\) is an \(l\)-vector.
For instance, \(\tb{a}\vee\tb{b}=0\) and \(\alpha\vee\tb{P}=0\).
Furthermore, \(\tb{a}\vee\tb{P}=\tb{P}\vee\tb{a}\) and \(\tb{P}\vee\tb{Q}=-\tb{Q}\vee\tb{P}\).

The equality \(\tb{P}\times\tb{Q}=(\tb{P}\vee\tb{Q})\I\) holds for any bivectors \(\tb{P}\) and \(\tb{Q}\) in \(\bigwedge\R{3*}\).

The inner product of a \(k\)-blade \(A_k\) and an \(l\)-blade \(B_l\), where \(k\le l\),  is related to the join and the outer product  as follows:
\begin{equation}
\begin{aligned}
&A_k\cdot B_l = (A_k\I)\vee B_l,\\
&(A_k\cdot B_l)\I = A_k\wedge (B_l\I).
\label{dot to vee}
\end{aligned}
\end{equation}
For instance, \(\tb{a}\cdot\tb{b}=(\tb{a}\I)\vee\tb{b}\) and \((\tb{a}\cdot\tb{P})\I=\tb{a}\wedge(\tb{P}\I)\).
The equalities~(\ref{dot to vee}) are applicable in \R{4*} and \R{5*} as well.

\subsection{Generalising duality transformation to \R{(n+1)*}}
The duality transformation defined above for \(n=2\) can be easily generalised to the case of arbitrary \(n\).
Since it is linear, it is sufficient to define  \( \J_{n} : \bigwedge\R{(n+1)*} \to \bigvee\R{n+1}  \)  on the basis blades
of \(\bigwedge\R{(n+1)*}\).
For a basis blade \(\e_{S}\), where \(S\) is a list of indices for this blade, I define
\begin{equation} 
\J_{n}(\e_S)=\e^{\reverse{S^\perp}}\label{J},
\end{equation}
where \({S}^\perp\) is a list of indices complementary to \(S\) such that \({S}^\perp S\) is an even permutation of \(01\dots n\),
and \(\reverse{S^{\perp}}\) is the reversion of the indices in \(S^{\perp}\), 
e.g.\ \(\reverse{0123}=3210\).
It is assumed that \(\e_{\varnothing}=1\) and \(\e^{\varnothing}=1\) 
(note that I use the same notation for scalars in \(\bigwedge\R{3*}\) and \(\bigvee\R{3}\) for simplicity).
The reversion has no effect on an empty list or a list consisting of one index, 
e.g.\ \(\reverse{\varnothing}= \varnothing\) and \(\reverse{1}=1\).
For example, \(S=1\)  yields \({S}^\perp=20\) for \(n=2\), since \(S^{\perp}S=201\) is an even permutation of \(012\), so that
\(\J_2(\e_1)=\e^{\reverse{20}}=\e^{02}=-\e^{20}\).
The transformation \(\J_{n}\) extends to arbitrary multivectors by linearity, e.g. \(\J_3(2\e_0-\e_{23})=2\J_3(\e_0)-\J_3(\e_{23})\).
The duality transformation thus defined is independent of the metric.

To understand this definition, observe that \(\e_{S^{\perp}S}=\e_{01\dots n}=\I_{n+1}\) 
for any basis blade \(\e_S\) in \(\bigwedge\R{(n+1)*}\).
So, in \(\bigwedge\R{(n+1)*}\) with the elliptic metric,
\[
\e_S\I_{n+1}^{-1}=
\e_S\reverse{\I}_{n+1}=
\e_S\reverse{\e}_{S^{\perp}S}=
\e_S\reverse{\e_{S^{\perp}}\e_S}=
\e_S\reverse{\e}_S \reverse{\e}_{S^{\perp}}=
\reverse{\e}_{S^{\perp}}=
\e_{\reverse{S^{\perp}}}
\]
and 
\begin{equation}
\J_n(\e_S)=\Id(O(\e_S))=\e^{\reverse{S^{\perp}}},
\end{equation}
where \(\Id\) and \(O\) are the appropriate identity and orthogonality transformations.

The inverse transformation is defined on an arbitrary multivector \(Y\in\bigvee\R{n+1}\) by
\begin{equation}
\J_n^{-1}(Y)=O^{-1}(I^{-1}(Y)),
\label{inverse of J}
\end{equation}
which gives \(J_n^{-1}(\e^S)=\e_S\I_{n+1}=\e_S \I_{n+1}^{-1}\I_{n+1}^{2}\) for the elliptic metric.
Since \(\I_{n+1}^2=-1\) if \(n=1,2\) and \(+1\) if \(n=3,4\), 
the following applies
\begin{equation}
\J_{n}^{-1}(\e^S)
=
\left\{
\begin{aligned}
& -\e_{\reverse{S^\perp}},\quad\textrm{if }n=1,2\\
& +\e_{\reverse{S^\perp}},\quad\textrm{if }n=3,4\\
\end{aligned}
\right.
\end{equation}
regardless of the metric.
Therefore, in \(\bigwedge\R{4*}\) and \(\bigwedge\R{5*}\), \(\J_{n}^{-1}\) can be obtained from \(\J_{n}\) by lowering and raising the indices, 
e.g.\ \(\J_4(\e_{230})=\e^{41}\) gives \(\J_4^{-1}(\e^{230})=\e_{41}\).
In \(\bigwedge\R{2*}\) and \(\bigwedge\R{3*}\), in addition to lowering and raising the indices it is also necessary to change the sign,
e.g.\ \(\J_2(\e_{12})=\e^0\) gives \(\J_2^{-1}(\e^{12})=-\e_0\).

\begin{table}[H]
\caption{\T{3}, \R{4}, and \R{4*}\label{Jthree}}
\end{table}
\vspace{-1cm}
Bivectors\\[5pt]
\begin{tabular*}{\textwidth}{lccccccr}
\(S\)                 & 10   & 20   & 30   &  23   & 31   &  12  & \\
\cline{1-7} \\[-10pt]
\(S^\perp S\)         & 3210 & 1320 & 2130  & 0123 & 0231 & 0312& {\small even permutations of} 0123\\ 
\cline{1-7} \\[-10pt]
\(S^\perp\)           & 32   & 13   & 21    & 01   & 02   & 03  & \\
\cline{1-7} \\[-10pt]
\(\reverse{S^\perp}\) & 23   & 31   & 12    & 10   & 20   & 30  & \\
\end{tabular*}

Other multivectors\\[5pt]
\begin{tabular*}{\textwidth}{lcccccccccc}
\(S\)                 & \(\varnothing\) & 0    & 1    & 2    & 3   &  123  & 320  & 130  &   210 & 0123            \\
\cline{1-11} \\[-10pt]
\(S^\perp S\)         & 3210            & 3210 & 0231 & 0312 & 0123&  0123 & 1320 & 2130 &   3210& 0123            \\ 
\cline{1-11} \\[-10pt]  
\(S^\perp\)           & 3210            & 321  & 023  & 031  & 012 &  0    & 1    & 2    &   3   & \(\varnothing\) \\
\cline{1-11} \\[-10pt]
\(\reverse{S^\perp}\) & 0123            & 123  & 320  & 130  & 210 &  0    & 1    & 2    &   3   & \(\varnothing\) \\
\end{tabular*}
\vspace{10pt}

\begin{table}[H]
\caption{\T{4}, \R{5}, and \R{5*}\label{Jfour bivectors}}
\end{table}
\vspace{-1cm}
Bivectors\\[5pt]
\begin{tabular*}{\textwidth}{lccccccccccr}
\(S\)                  & 10   & 20   & 30  & 40    & 23   & 31   & 12    &  41   & 42   & 43   & \\
\cline{1-11} \\[-10pt]
\(S^\perp S\)          & 43210& 41320& 42130& 12340 &01423& 02431& 03412&  03241& 01342& 02143 & \\
\cline{1-11} \\[-10pt]
\(S^\perp\)            & 432   & 413    & 421    & 123    & 014  & 024   & 034    &  032   & 013   & 021    & \\
\cline{1-11} \\[-10pt]
\(\reverse{S^\perp}\)  & 234  & 314  & 124   & 321    & 410  & 420   & 430    &  230   & 310   & 120   & \\
\end{tabular*}

Trivectors\\[5pt]
\begin{tabular*}{\textwidth}{lcccccccccc}
\(S\)                      & 234   & 314   & 124   & 321      & 410  &420   & 430    &  230   & 310   & 120    \\
\cline{1-11} \\[-10pt]
\(S^\perp S\)          & 01234& 02314& 03124  & 04321  & 32410& 13420& 21430& 14230 & 24310& 34120  \\
\cline{1-11} \\[-10pt]
\(S^\perp\)            & 01   & 02   & 03       & 04 & 32   & 13   & 21    & 14    & 24   & 34     \\
\cline{1-11} \\[-10pt]
\(\reverse{S^\perp}\)  & 10   & 20   & 30   & 40 & 23   & 31   & 12    & 41    & 42   & 43     \\
\end{tabular*}

Other multivectors\\[5pt]
\begin{tabular*}{\textwidth}{lcccccccccccc}
\(S\)                 & \(\varnothing\) & 0    & 1    & 2    & 3   &  4   & 1234  & 2340  & 3140  & 1240  & 3210  & 01234          \\
\cline{1-13} \\[-10pt]
\(S^\perp S\)         & 43210           & 43210& 04321& 04132&04213&01234 & 01234 & 12340 & 23140 & 31240 & 43210 & 01234          \\ 
\cline{1-13} \\[-10pt]  
\(S^\perp\)           & 43210           & 4321 & 0432 & 0413 & 0421& 0123 & 0     & 1     & 2     & 3     & 4     & \(\varnothing\) \\
\cline{1-13} \\[-10pt]
\(\reverse{S^\perp}\) & 01234           & 1234 & 2340 & 3140 & 1240& 3210 & 0     & 1     & 2     & 3     & 4     & \(\varnothing\) \\
\end{tabular*}
\vspace{10pt}

The formula~(\ref{J}) alone is not sufficient in the case of \(\bigwedge\R{2*}\), because 01 is the only possible even permutation of 01, which gives 
\(\J_1(1)=-\e^{01}\),
\(\J_1(\e_1)=\e^0\),
\(\J_1(\e_{01})=1\).
Using the orthogonality transformation directly gives \(\J_1(\e_0)=-\e^1\).
I also get the following for the inverse:
\(\J_1^{-1}(1)=\e_{01}\),
\(\J_1^{-1}(\e^0)=\e_1\),
\(\J_1^{-1}(\e^1)=-\e_0\),
\(\J_1^{-1}(\e^{01})=-1\).

The transformation \(\J_{n}\) is an isomorphism in the sense that a blade  \(A_k\in\bigwedge\R{(n+1)*}\) represents dually 
the same geometric object in \T{n} as that represented directly by \(\J_n(A_k) \in\bigvee\R{n+1}\).
The target and dual Grassmann algebras are not isomorphic with respect to the outer product, 
i.e.\ \(\J_n(A\wedge B)\ne\J_n(A)\vee\J_n(B)\).

\subsection{Representing orientation and weight with blades}
Since blades in  \(\bigvee\R{3}\) and \(\bigwedge\R{3*}\) are oriented, 
their orientation can be used to represent the orientation of points and lines in \T{2}.
The orientation in \T{2} inherited from \(\bigvee\R{3}\) is called the bottom-up orientation,
while \(\bigwedge\R{3*}\) provides  the top-down orientation.

\begin{figure}[t!]
\begin{subfloatenv}{\R{3}}
\begin{asy}
import Drawing3D;
DrawingR3 drawing = DrawingR3(4,0.1);	
drawing.target_axes();
drawing.T2();

triple P = (1,0,0);

drawing.Drawing3D.vector(P);
draw_point_orientation(P);

\end{asy}
\end{subfloatenv}\hspace{-12pt}%
\begin{subfloatenv}{\R{3*}}
\begin{asy}
import Drawing3D;
DrawingR3 drawing = DrawingR3(4,0.1);
drawing.dual_axes();
drawing.dualT2();

triple u = (0,0,1);
triple v = (0,1,0);
bivector(v,u);

\end{asy}
\end{subfloatenv}
\caption{The origin of \T{2} with the counterclockwise orientation and \(\e_{12}\) in \R{3*}}
\label{counterclockwise origin in T2 and bivector in R3*}
\begin{subfloatenv}{\R{3}}
\begin{asy}
import Drawing3D;
DrawingR3 drawing = DrawingR3(4,0.1);
drawing.target_axes();
drawing.T2();

triple v = (2,-2,-1);
drawing.Drawing3D.vector(v);

triple P = v/2;
dot(P);
draw((P.x,P.y,0)--P, Dotted);
draw((P.x,0,P.z)--P, Dotted);

draw_point_orientation(P);

label("$J(\textbf{P})$", v, (0,-0.5,0));

\end{asy}
\end{subfloatenv}\hspace{-12pt}%
\begin{subfloatenv}{\R{3*}}
\begin{asy}
import Drawing3D;
DrawingR3 drawing = DrawingR3(4,0.1);
drawing.dual_axes();
drawing.dualT2();

path3 p = drawing.T2_path3();
path3 Q = plane_path3((1,-1,-1/2), (1,1,0), theta=-13, scale=7.8);
path3 l = plane_intersection(Q,p);
draw(Label("$L$", 0.2,W), l);

triple u = (1,0,2);
triple v = (1,1,0);
bivector(v,u);
label("$\textbf{P}$",(0.8,0.5,1));
\end{asy}
\end{subfloatenv}
\caption{The counterclockwise point \((-1,-\tfrac{1}{2})\) in \T{2} represented directly by \(\J(\tb{P})=2\e^0-2\e^1-\e^2\)
and dually by \(\tb{P}=(\e_0+\e_1)\wedge(\e_0+2\e_2)=2\e_{12}-2\e_{20}-\e_{01}\)
}
\label{counterclockwise point in T2 and bivector in R3*}
\end{figure}

The point at the origin of \T{2} is represented by the bivector \(\e_{12}\) as shown in Figure~\ref{counterclockwise origin in T2 and bivector in R3*}.
It is natural to assign the orientation of  \(\e_{12}\) as a bivector in \R{3*} to represents the counterclockwise orientation of the point at the origin of \T{2}.
In general, any bivector in \(\bigwedge\R{3*}\) can be written as \(w\e_{12}+\e_0\wedge\tb{a}\), 
where \(\alpha\in\R{}\) and \(\tb{a}\in\R{3*}\).
Its dual is a vector \(w\e^0+\J(\e_0\wedge\tb{a})\) in \R{3}.
I assign the counterclockwise orientation to the points represented dually by \(w\e_{12}+\e_0\wedge\tb{a}\) with \(w>0\),
and the clockwise orientation to those with \(w<0\).
This is illustrated in Figure~\ref{counterclockwise point in T2 and bivector in R3*} for a counterclockwise point at \((-1,-\tfrac{1}{2})\).

In the bottom-up view of \T{2}, points lack an orientation but the representing vectors in \R{3} do have an orientation.
The finite points are represented by vectors \(\tb{x}=w\e^0+x\e^1+y\e^2\), where \(w\ne0\).
The finite points represented by vectors with \(w>0\) can be thought of as having the positive bottom-up orientation,
and those with \(w<0\) as having the negative bottom-up orientation.
The positive bottom-up orientation corresponds to the counterclockwise top-down orientation, while
the negative bottom-up orientation corresponds to the clockwise top-down orientation.
Hence, points shown in Figures~\ref{counterclockwise origin in T2 and bivector in R3*}(a) 
and \ref{counterclockwise point in T2 and bivector in R3*}(a)
have the positive orientation.

Any bivector that can be written as \(\e_0\wedge\tb{a}\), where \(\tb{a}\in\R{3*}\), represents a point at infinity in \T{2}.
Since \(\e_0\wedge\e_0=0\), it is sufficient to consider vectors \(\tb{a}\) that lie in the \(ab\)-plane, i.e.\ \(\tb{a}=a\e_1+b\e_2\). 
The relationship between such bivectors and the top-down orientation of the points at infinity they dually represent is
illustrated in Figures~\ref{infinite point in T2 and bivector in R3*} and \ref{infinite point and bivector opposite orientation},
where the top-down orientation is shown with a solid vector in \T{2}.
The top-down orientation of a point at infinity dually represented by \(\e_0\wedge\tb{a}\) 
is consistent with the direction of \(\Id(\tb{a})=a\e^1+b\e^2\) in \R{3}
and \((a,b)\) in \T{2}.
As there is a natural isomorphism between vectors in the \(ab\)-plane of \R{3} and vectors in \T{2},
I will tacitly identify such vectors in \R{3} with the corresponding vectors in \T{2}.
So, the top-down orientation of a point at infinity dually represented by \(\e_0\wedge\tb{a}\) is provided by 
\(\Id(\tb{a})\) seen as a vector in \T{2}.

Points at infinity in \T{2} are directly represented by vectors \(\tb{x}=x\e^1+y\e^2\), which lie in the plane \(w=0\).
The bottom-up orientation of such points is consistent with the orientation of the representing vector in \R{3},
as illustrated in Figures ~\ref{infinite point in T2 and bivector in R3*}(a) and \ref{infinite point and bivector opposite orientation}(a),
where the bottom-up orientation is shown with a dashed vector.
So, the bottom-up orientation of a point at infinity dually represented by \(\e_0\wedge\tb{a}\) is provided by 
\(\J(\e_0\wedge\tb{a})=-b\e^1+a\e^2\) seen as a vector in \T{2}.
The bottom-up orientation of points at infinity is somewhat misleading as it suggests that a point at infinity is located in the direction of 
the orienting vector, but the same point lies at infinity in the opposite direction. 

\begin{figure}[t!]
\begin{subfloatenv}{\R{3}}
\begin{asy}
import Drawing3D;
DrawingR3 drawing = DrawingR3(4,0.1);	
drawing.target_axes();
drawing.T2();

path3 l = line_path3((1,0,0), (0,2,1), 10);
path3 p = drawing.T2_path3();
int N=2;
for(int i=-N; i<=N; ++i) { path3 shifted_l = shift(-i*(0,-1,2)/3)*l; triple[] ps = intersectionpoints(shifted_l,p); draw(ps[0]--ps[1]); }

triple c = (1,0,0);
triple s = 0.7*(0,1,-2);
int No=2;
path3 lo = (c--(c-s)); 
draw(lo, linewidth(1 bp), EndArrow3(HookHead2((1,0,0)), size=5));
//for(int i=-No; i<=No; ++i) { path3 lo = shift(i*(0,2,1)/2)*((c+s)--(c-s)); draw(lo, linewidth(1 bp), EndArrow3(HookHead2((1,0,0)), size=5)); }

drawing.Drawing3D.vector((0,-2,-1), "$J(\textbf{a})$", (0,-0.5,-1));
triple s = 0.7*(0,-2,-1);
path3 lo = (c--(c+s)); 
draw(lo, linetype(new real[] {4, 4},adjust=false)+linewidth(1 bp), EndArrow3(HookHead2((1,0,0)), size=5));

\end{asy}
\end{subfloatenv}\hspace{-12pt}%
\begin{subfloatenv}{\R{3*}}
\begin{asy}
import Drawing3D;
DrawingR3 drawing = DrawingR3(4,0.1);
drawing.dual_axes();
drawing.dualT2();

triple u = (0,-1,2);
triple v = (1,0,0);
bivector(v,u,20);

path3 l = line_path3((1,0,0), (0,-1,2), 10);
path3 p = drawing.T2_path3();
triple[] ps = intersectionpoints(l,p); 
draw(ps[0]--ps[1]);

drawing.Drawing3D.vector((0,-1,2), "$\textbf{a}$", (0,1,0));

\end{asy}
\end{subfloatenv}
\caption{A point at infinity in \T{2} represented directly by \(J(\e_0\wedge\tb{a})=-2\e^1-\e^2\)
and dually by \(\e_0\wedge\tb{a}\), where \(\tb{a}=-\e_1+2\e_2\) (the point lies on a line defined by \(-x+2y=0\))}
\label{infinite point in T2 and bivector in R3*}
\begin{subfloatenv}{\R{3}}
\begin{asy}
import Drawing3D;
DrawingR3 drawing = DrawingR3(4,0.1);	
drawing.target_axes();
drawing.T2();

path3 l = line_path3((1,0,0), (0,2,1), 10);
path3 p = drawing.T2_path3();
int N=2;
for(int i=-N; i<=N; ++i) { path3 shifted_l = shift(-i*(0,-1,2)/3)*l; triple[] ps = intersectionpoints(shifted_l,p); draw(ps[0]--ps[1]); }

triple c = (1,0,0);
triple s = 0.7*(0,1,-2);
int No=2;
//for(int i=-No; i<=No; ++i) { path3 lo = shift(i*(0,2,1)/2)*((c-s)--(c+s)); draw(lo, linewidth(1 bp), EndArrow3(HookHead2((1,0,0)), size=5)); }
path3 lo = (c--(c+s)); 
draw(lo, linewidth(1 bp), EndArrow3(HookHead2((1,0,0)), size=5));

drawing.Drawing3D.vector((0,2,1), "$J(\textbf{a})$", (0,0.5,-0.5));
triple s = 0.7*(0,2,1);
path3 lo = (c--(c+s)); 
draw(lo, linetype(new real[] {4, 4},adjust=false)+linewidth(1 bp), EndArrow3(HookHead2((1,0,0)), size=5));

\end{asy}
\end{subfloatenv}\hspace{-12pt}%
\begin{subfloatenv}{\R{3*}}
\begin{asy}
import Drawing3D;
DrawingR3 drawing = DrawingR3(4,0.1);
drawing.dual_axes();
drawing.dualT2();

triple u = (0,-1,2);
triple v = (1,0,0);
bivector(u,v,20);

path3 l = line_path3((1,0,0), (0,-1,2), 10);
path3 p = drawing.T2_path3();
triple[] ps = intersectionpoints(l,p); 
draw(ps[0]--ps[1]);

drawing.Drawing3D.vector((0,1,-2), "$\textbf{a}$", (0,1,0));

\end{asy}
\end{subfloatenv}
\caption{The same point at infinity as in Figure~\ref{infinite point in T2 and bivector in R3*} but with the opposite orientation,
i.e.~\(\tb{a}=\e_1-2\e_2\)}
\label{infinite point and bivector opposite orientation}
\end{figure}

\begin{figure}[t!]

\begin{subfloatenv}{\R{3}}
\begin{asy}
import Drawing3D;
DrawingR3 drawing = DrawingR3(4,0.1);	
drawing.target_axes();
drawing.T2();

path3 l = line_path3((1,1,-1/2), (0,2,1), 10);
path3 p = drawing.T2_path3();
triple[] ps = intersectionpoints(l,p); 
draw(ps[0]--ps[1]);

triple u = (1,0,-1);
triple v = (1,2,0);

bivector(v,u);
label("$J(\textbf{a})$",(-0.5,1,1),(0,0,0));

triple c = (1,0.4,-0.8);
triple s = 0.3*(0,-1,2);
int No=2;
//for(int i=-No; i<=No; ++i) { path3 lo = shift(i*(0,2,1)/2)*((c-s)--(c+s)); draw(lo, linewidth(1 bp), EndArrow3(HookHead2((1,0,0)), size=5)); }
path3 lo = (c--(c+s)); 
draw(lo, linewidth(1 bp), EndArrow3(HookHead2((1,0,0)), size=5));

triple s = 0.3*(0,-2,-1);
path3 lo = (c--(c+s)); 
draw(lo, linetype(new real[] {4, 4},adjust=false)+linewidth(1 bp), EndArrow3(HookHead2((1,0,0)), size=5));

\end{asy}
\end{subfloatenv}\hspace{-12pt}%
\begin{subfloatenv}{\R{3*}}
\begin{asy}
import Drawing3D;
DrawingR3 drawing = DrawingR3(4,0.1);
drawing.dual_axes();
drawing.dualT2();

drawing.Drawing3D.vector((2,-1,2),label="$\textbf{a}$",(0,0,0.5));

triple P = (1,-1/2,1);
dot(P);
draw((P.x,P.y,0)--P, Dotted);
draw((P.x,0,P.z)--P, Dotted);

label("$(1,-\frac{1}{2},1)$",(1,-1/2,1),(1/2,-1,-1));

\end{asy}
\end{subfloatenv}
\caption{Line \(L\) defined by \(-\tfrac{1}{2}x+y+1=0\) and dually represented by \(\tb{a}=2\e_0-\e_1+2\e_2\) in \R{3*} and
directly represented by \(\J(\tb{a})
=-2\e^{12}+\e^{20}-2\e^{01}
=
(\e^0+2\e^1)\vee(\e^0-\e^2)
\)}
\label{orienting a line in T2}
%
\begin{subfloatenv}{\R{3}}
\begin{asy}
import Drawing3D;
DrawingR3 drawing = DrawingR3(4,0.1);	
drawing.target_axes();
drawing.T2();

path3 l = line_path3((1,1,-1/2), (0,2,1), 10);
path3 p = drawing.T2_path3();
triple[] ps = intersectionpoints(l,p); 
draw(ps[0]--ps[1]);

triple u = (1,0,-1);
triple v = (1,2,0);

bivector(u,v);
label("$J(\textbf{a})$",(-0.5,1,1),(0,0,0));

triple c = (1,0.4,-0.8);
triple s = -0.3*(0,-1,2);
int No=2;
//for(int i=-No; i<=No; ++i) { path3 lo = shift(i*(0,2,1)/2)*((c-s)--(c+s)); draw(lo, linewidth(1 bp), EndArrow3(HookHead2((1,0,0)), size=5)); }

path3 lo = (c--(c+s)); 
draw(lo, linewidth(1 bp), EndArrow3(HookHead2((1,0,0)), size=5));

triple s = 0.3*(0,2,1);
path3 lo = (c--(c+s)); 
draw(lo, linetype(new real[] {4, 4},adjust=false)+linewidth(1 bp), EndArrow3(HookHead2((1,0,0)), size=5));

\end{asy}
\end{subfloatenv}\hspace{-12pt}%
\begin{subfloatenv}{\R{3*}}
\begin{asy}
import Drawing3D;
DrawingR3 drawing = DrawingR3(4,0.1);
drawing.dual_axes();
drawing.dualT2();

drawing.Drawing3D.vector(-(2,-1,2),label="$\textbf{a}$",(0,0.5,1));

triple P = (1,-1/2,1);
dot(P);
draw((P.x,P.y,0)--P, Dotted);
draw((P.x,0,P.z)--P, Dotted);

label("$(1,-\frac{1}{2},1)$",(1,-1/2,1),(1/2,-1,-1));

\end{asy}
\end{subfloatenv}
\caption{Same line as in Figure~\ref{orienting a line in T2} but with the opposite orientation, i.e.\ \(\tb{a}=-2\e_0+\e_1-2\e_2\)}
\label{orienting a line in T2 opposite}
\end{figure}

A similar relationship exists between vectors in \R{3*} and oriented lines in \T{2}.
Consider a line \(L\) defined by \(ax+by+d=0\).
This line can be represented dually by \(\tb{a}=d\e_0+a\e_1+b\e_2\) and directly by \(\J(\tb{a})\).
The relationship between the line's top-down orientation and \(\tb{a}\) is illustrated in Figure~\ref{orienting a line in T2},
where the top-down orientation is shown with a solid vector in \T{2}.
Since the value of \(d\) has no effect on the line's orientation and is only responsible for the line's shift from the origin,
it can be ignored.
To obtain the top-down orientation, take \(\Id(\tb{a})=d\e^0+a\e^1+b\e^2\) and remove the component along \(\e^0\),
which yields  \(\tb{x}=a\e^1+b\e^2\) in \R{3} and \((a,b)\) in \T{2}.
The line is oriented towards the origin if \(d>0\) and away from the origin if \(d<0\)
(compare Figures~\ref{orienting a line in T2} and \ref{orienting a line in T2 opposite}).
The line at infinity in \T{2} corresponds to the origin of \T{2*} and, therefore,  can be represented by \(d\e_0\), where \(d\ne0\).
The line at infinity is oriented towards the origin if \(d>0\) and away from the origin if \(d<0\).

The line \(L\) is represented directly by \(\J(\tb{a})\), which can be used to get the bottom-up orientation.
To get the orienting vector along the line, obtain \(\J(\tb{a})\) and eliminate the component along \(\e^{12}\),
since it only causes a shift of the line from the origin without affecting its orientation.
If \(\tb{a}=d\e_0+a\e_1+b\e_2\), then \(\J(\tb{a})=-d\e^{12}-a\e^{20}-b\e^{01}\) and setting \(d=0\) yields
\(B=-a\e^{20}-b\e^{01}\), which can be written as \(B=\e^0\vee\tb{x}\), where
\(\tb{x}=-b\e^1+a\e^2\) can be identified with \((-b,a)\) in \T{2} and provides the bottom-up orientation of \(L\).
This is illustrated in Figures~\ref{orienting a line in T2}(a) and \ref{orienting a line in T2 opposite}(b), where
the bottom-up orientation is shown with a dashed vector in \T{2}.

The bottom-up orientation of the line at infinity is either clockwise or counterclockwise, 
depending whether the points at infinity, which constitute the line at infinity, are traversed in the clockwise or
counterclockwise direction.
So, for instance, \(\e^{12}\) is counterclockwise and \(-\e^{12}\) is clockwise.

A point in \T{2} does not have a weight.
However, since the blade it is dually represented by does have a weight, it is natural to use the 
blade's weight as a characteristic of the point itself, i.e.\ the point's weight is identified with the blade's weight.
The same applies to a line, whose weight can be identified with the weight of the blade that dually represents the line.


Summary (from \(\bigwedge\R{3*}\) to \T{2} )\\[12pt]
Bivectors \(\rightarrow\) points in \T{2}:\\
\(\tb{P}=\e_{12}+x\e_{20}+y\e_{01}\) \(\rightarrow\) a counterclockwise (positive) point at \((x,y)\). \\
\(\alpha\tb{P}\), \(\alpha>0\) \(\rightarrow\) same as above but with a different weight.\\
\(-\tb{P}\) \(\rightarrow\) same as \(\tb{P}\) but with the opposite orientation (clockwise, or negative).\\
\(\tb{N}=\e_0\wedge(a\e_1+b\e_2)\) \(\rightarrow\) a point at infinity which lies on the line \(ax+by=0\);
its top-down and bottom-up orientations are given by \((a,b)\) and \((-b,a)\) respectively.\\
\(\alpha\tb{N}\), \(\alpha>0\) \(\rightarrow\) same as above but with a different  weight.\\
\(-\tb{N}\) \(\rightarrow\) same as \(\tb{N}\) but with the opposite orientation.\\[12pt]
Vectors \(\rightarrow\) finite lines in \T{2} and the line at infinity:\\
\(\tb{a}=d\e_0+a\e_1+b\e_2\) \(\rightarrow\) a line in \T{2} defined by \(ax+by+d=0\);
the line's top-down and bottom-up orientations are given by \((a,b)\) and \((-b,a)\) respectively.\\
\(\alpha\tb{a}\), \(\alpha>0\) \(\rightarrow\) same line with a different weight.\\
\(-\tb{a}\) \(\rightarrow\) same line with the opposite orientation\\ 
\(d\e_0\) \(\rightarrow\) the line at infinity in \T{2} 
(in the top-down view, it is oriented towards the origin if \(d>0\) or away from the origin if \(d<0\)).\\[12pt]

In the following, I identify blades in \R{3*} with the geometric objects in \T{2} that these blades dually represent.
For instance, when referring to \(\tb{a}\) as a line and \(\tb{P}\) as a point, 
I mean a line in \T{2} dually represented by \(\tb{a}\) and a point in \T{2} dually represented by \(\tb{P}\), respectively.
This terminology will apply to the dual representation only.
I will use blades in \R{3*} to label geometric objects in \T{2}.

\subsection{Euclidean plane \E{2}}
For a line \(\tb{a}=d\e_0+a\e_1+b\e_2\) and a  point \(\tb{P}=w\e_{12}+x\e_{20}+y\e_{01}\),
I have \(\tb{a}^2=a^2+b^2\), \(\tb{P}^2=-w^2\),
\(\norm{\tb{a}}=(a^2+b^2)^{\tfrac{1}{2}}\), and \(\norm{\tb{P}}=|w|\).
Points at infinity and the line at infinity have zero norm.
A finite line \(\tb{a}\) and a finite point \(\tb{P}\) are normalised if \(\norm{\tb{a}}=1\) and \(\norm{\tb{P}}=1\).
Note that a normalised point can be either clockwise or counterclockwise.

\begin{figure}[h!]
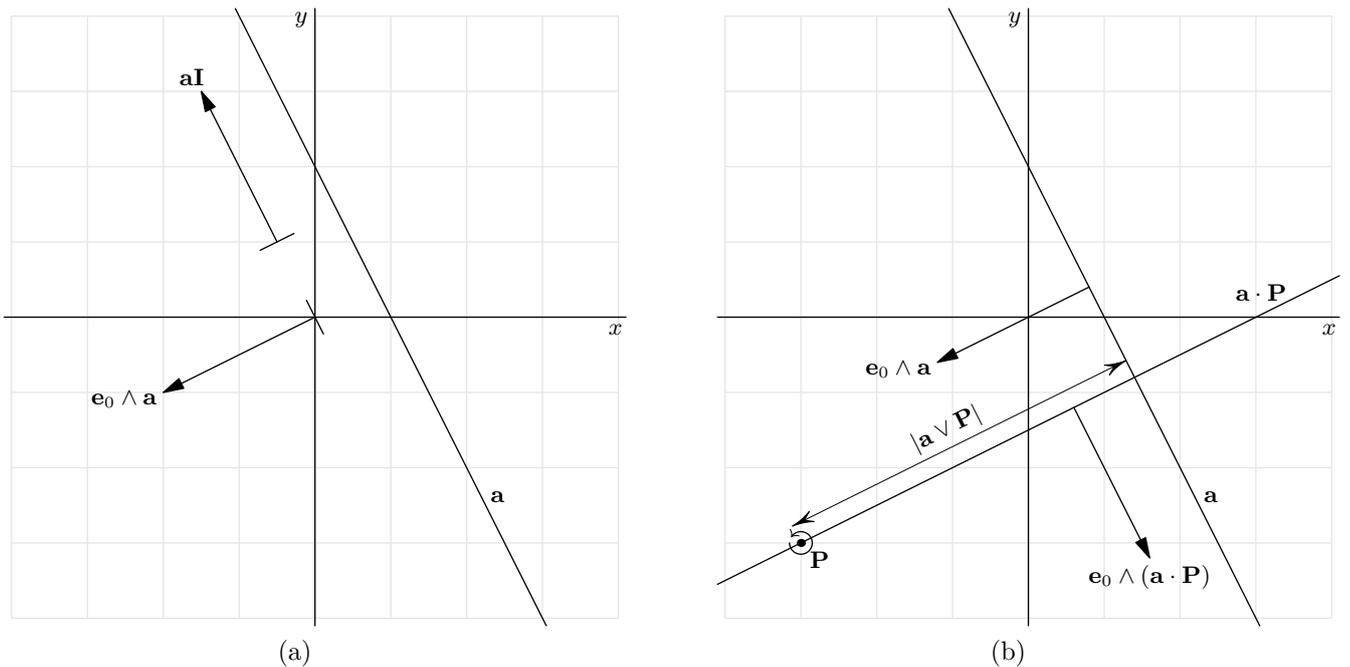

\begin{subfloatenv}{ }
\begin{asy}
import Projective2D;
import Graph2D;
Drawing d = Drawing();

var a = Line(2,-2,-1);
var N = wedge(e_0, a);

d.line(a, "$\textbf{a}$", draw_orientation=false);
d.point_at_infinity(N, "$\textbf{e}_0\wedge\textbf{a}$");
d.point_at_infinity(a*I, "$\textbf{aI}$", O=(-0.5,1));
\end{asy}
\end{subfloatenv}\hfill%
\begin{subfloatenv}{ }
\begin{asy}
import Projective2D;
import Graph2D;
Drawing d = Drawing();

var a = Line(2,-2,-1);
var P = Point(1,-3,-3);

d.line(a, "$\textbf{a}$", "$\textbf{e}_0\wedge\textbf{a}$");
d.point(P, "$\textbf{P}$");

var b = dot(a,P);
d.line(b, "$\textbf{a}\cdot\textbf{P}$", position=0.1, "$\textbf{e}_0\wedge(\textbf{a}\cdot\textbf{P})$", orientation_align=S);
d.stretch(P, dot(P,a)/a, "$|\textbf{a}\vee\textbf{P}|$");
\end{asy}
\end{subfloatenv}
\caption{Basic properties of points and lines in \E{2}}
\label{basic line and point}
\end{figure}

The point at infinity \(\e_0\wedge\tb{a}\) lies on a line \(\tb{a}\) and gives the orientation of the line.
In figures, I depict it as an arrow with a small base. The base can be thought of as a small segment of one of the lines
in the stack representing the point at infinity. The arrow indicates the top-down orientation of the stack. 
Speaking informally, it selects one of the two possible directions of rotation around the point at infinity.
When a point at infinity is depicted in this fashion, it can be placed anywhere in the figure since that has no effect on the stack or its orientation.
For the line \(\tb{a}=2\e_0-2\e_1-\e_2\) shown in Figure~\ref{basic line and point}(a),
I have \(\e_0\wedge\tb{a}=\e_0\wedge(-2\e_1-\e_2)=\e_{20}-2\e_{01}\).
The length of the arrow and its direction are consistent with the vector \((-2,-1)\).
On the other hand, \(\tb{a}\I=(2\e_0-2\e_1-\e_2)\e_{012}=-2\e_{20}-\e_{01}=\e_0\wedge(-\e_1+2\e_2)\) 
is a point at infinity that
lies in the direction perpendicular to the line \(\tb{a}\).
The point at infinity \(\tb{a}\I\)   is called the polar point of \(\tb{a}\).

\begin{figure}[h!]
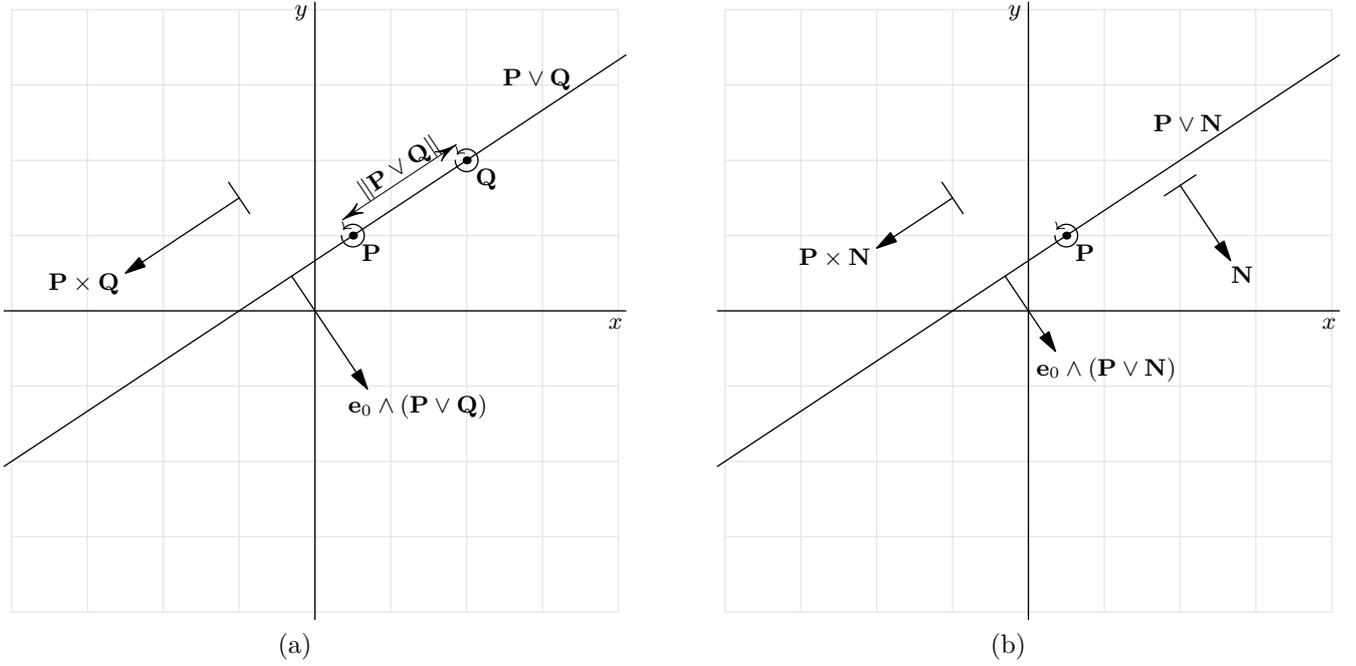

\begin{subfloatenv}{ }
\begin{asy}
import Projective2D;
import Graph2D;
Drawing d = Drawing();

MV P = Point(1, 0.5, 1);
MV Q = Point(1, 2, 2);

d.point(P, "\textbf{P}");
d.point(Q, "\textbf{Q}");
d.line(join(P,Q), "$\textbf{P}\vee\textbf{Q}$", position=0.1, "$\textbf{e}_0\wedge(\textbf{P}\vee\textbf{Q})$");
d.point_at_infinity(cross(P,Q),"$\textbf{P}\times\textbf{Q}$", O=(-1,1.5));
d.stretch(P, Q, "$\norm{\textbf{P}\vee\textbf{Q}}$");
\end{asy}
\end{subfloatenv}\hfill%
\begin{subfloatenv}{ }
\begin{asy}
import Projective2D;
import Graph2D;
Drawing d = Drawing();

MV P = Point(1, 0.5, 1);
MV N = wedge(e_0, Line(0, 2/3, -1));

d.point(P, "\textbf{P}");
d.point_at_infinity(N, "\textbf{N}", O=(2,2-1/3));
d.line(join(P,N), "$\textbf{P}\vee\textbf{N}$", "$\textbf{e}_0\wedge(\textbf{P}\vee\textbf{N})$");
d.point_at_infinity(cross(P,N),"$\textbf{P}\times\textbf{N}$", O=(-1,1.5));
\end{asy}
\end{subfloatenv}
\caption{The join of two points in \E{2}}
\label{join of points in E2}
\end{figure}

The inner product \(\tb{a}\cdot\tb{P}\) gives a line passing through \(\tb{P}\) and is perpendicular to  \(\tb{a}\).
An example is shown in Figure~\ref{basic line and point}(b), where \(\tb{a}=2\e_0-2\e_1-\e_2\) and \(\tb{P}=\e_{12}-3\e_{20}-3\e_{01}\),
which gives \(\tb{a}\cdot\tb{P}=-3\e_0+\e_1-2\e_2\).
The line \(\tb{a}\cdot\tb{P}\) is normalised if both \(\tb{a}\) and \(\tb{P}\)  are normalised.

For a normalised counterclockwise point \(\tb{P}\) and a normalised line \(\tb{a}\), I have \(\tb{a}\wedge\tb{P}=\pm r\I\), 
where \(r=|d+ax+by|\) is the distance between the point and the line (the plus sign obtains if the line is oriented towards the point).
Since \(\tb{a}\wedge\tb{P}=(\tb{a}\vee\tb{P})\I\), I get \(r = |\tb{a}\vee\tb{P}|\).
Furthermore, \(|d|\) gives the distance between the line and the origin 
(the sign of \(d\) is positive if the line is oriented towards the origin, otherwise it is negative).

For a finite point \(\tb{P}=w\e_{12}+x\e_{20}+y\e_{01}\), I get
\(\e_0\wedge\tb{P}=w\I\) and \(\tb{P}\I=-w\e_0\), the latter being the line at infinity weighted by \(w\).
For a point at infinity \(\tb{N}=\e_0\wedge(a\e_1+b\e_2)\), I get \(\e_0\wedge\tb{N}=0\), 
i.e.\ \(\tb {N}\) lies on the line \(\e_0\), and \(\tb{N}\I=0\).

\begin{figure}[h!]
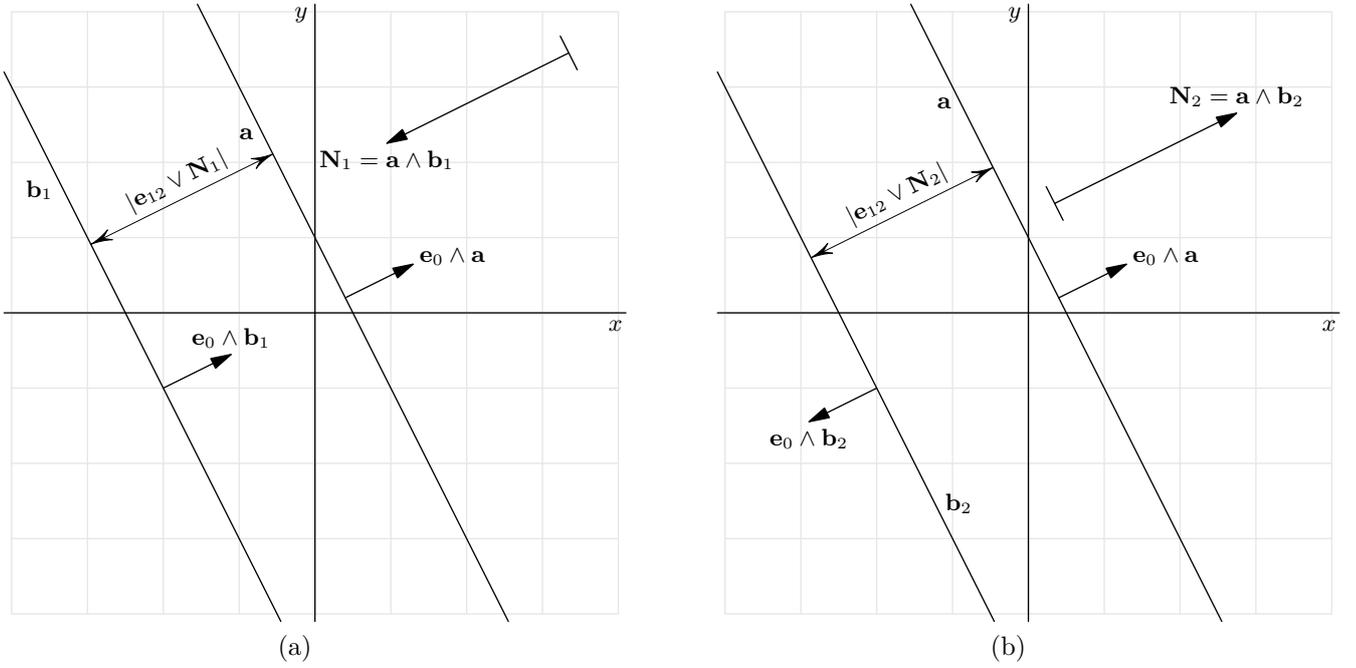

\begin{subfloatenv}{ }
\begin{asy}
import Projective2D;
import Graph2D;
Drawing d = Drawing();

MV a = Line(-1,2,1);
a = a/norm(a);
MV Na = wedge(e_0,a);

MV b = Line(5,2,1);
b = b/norm(b);
MV Nb = wedge(e_0,b);

MV N = wedge(a,b);

d.point_at_infinity(N, O=(3.34655,3.4535), "$\textbf{N}_1=\textbf{a}\wedge\textbf{b}_1$", align=S);
d.line(a, "$\textbf{a}$", "$\textbf{e}_0\wedge\textbf{a}$");
d.line(b, "$\textbf{b}_1$", "$\textbf{e}_0\wedge\textbf{b}_1$", orientation_align=(0,1));

MV a_ = join(Point(1,-3,1), a*I);
d.stretch(wedge(b,a_), wedge(a,a_), "$|\textbf{e}_{12}\vee\textbf{N}_1|$", shift_distance=-0.1);

\end{asy}
\end{subfloatenv}\hfill%
\begin{subfloatenv}{ }
\begin{asy}
import Projective2D;
import Graph2D;
Drawing d = Drawing();

MV a = Line(-1,2,1);
a = a/norm(a);
MV Na = wedge(e_0,a);

MV b = -Line(5,2,1);
b = b/norm(b);
MV Nb = wedge(e_0,b);

MV N = wedge(a,b);

d.point_at_infinity(N, O=(0.34655,1.4535), "$\textbf{N}_2=\textbf{a}\wedge\textbf{b}_2$", align=(0,1));
d.line(a, "$\textbf{a}$", position=0.15, "$\textbf{e}_0\wedge\textbf{a}$");
d.line(b, "$\textbf{b}_2$", "$\textbf{e}_0\wedge\textbf{b}_2$", orientation_align=(0,-1));

MV a_ = join(Point(1,-3,1), a*I);
d.stretch(wedge(b,a_), wedge(a,a_), "$|\textbf{e}_{12}\vee\textbf{N}_2|$", shift_distance=-0.3);

\end{asy}
\end{subfloatenv}
\caption{Parallel lines in \E{2}}
\label{depicting points at infinity}
\end{figure}

The join \(\tb{P}\vee\tb{Q}\) of two finite points gives a line passing through the points.
For example, the join of \(\tb{P}=\e_{12}+\tfrac{1}{2}\e_{20}+\e_{01}\) and 
\(\tb{Q}=\e_{12}+2\e_{20}+2\e_{01}\) yields \(\tb{P}\vee\tb{Q}=\e_0+\e_1-\frac{3}{2}\e_2\) (see Figure~\ref{join of points in E2}(a)).
The commutator \(\tb{P}\times\tb{Q}\) is a point at infinity that lies in the direction perpendicular to
the line passing through \(\tb{P}\) and \(\tb{Q}\);  
the commutator's orientation is  from \(\tb{Q}\) to \(\tb{P}\) if both points have the same orientation.
This is consistent with the equality \(\tb{P}\times\tb{Q}=(\tb{P}\vee\tb{Q})\I\).
For normalised points \(\tb{P}\) and \(\tb{Q}\), the distance \(r\) between them is given by 
\begin{equation}
r=\norm{\tb{P}\vee\tb{Q}}.
\end{equation}
So, for the points used in the example, I get \(\norm{\tb{P}\vee\tb{Q}}=\sqrt{1^2+\left(-\frac{3}{2}\right)^2}=\frac{\sqrt{13}}{2}\).

Joining \(\tb{P}=\e_{12}+\tfrac{1}{2}\e_{20}+\e_{01}\) and the point at infinity \(\tb{N}=\e_0\wedge(\tfrac{2}{3}\e_1 - \e_2)\)
results in \(\tb{P}\vee\tb{N}=\tfrac{2}{3}\e_0+\tfrac{2}{3}\e_1-\e_2\), which is a line that contains both \(\tb{P}\) and \(\tb{N}\) 
(see Figure~\ref{join of points in E2}(b)).
The orientation of \(\tb{P}\vee\tb{N}\) coincides with the orientation of \(\tb{N}\) 
if  \(\tb{P}\) is counterclockwise.
The join of two points at infinity gives the line at infinity.

For normalised parallel (and anti-parallel) lines,
\begin{equation}
\tb{a}\tb{b}=\pm1+\tb{a}\wedge\tb{b},
\end{equation}
where \(\tb{a}\wedge\tb{b}\) is a point at infinity,
oriented from \(\tb{a}\) to \(\tb{b}\) if both lines have the same orientation
(or from \(\tb{b}\) to \(\tb{a}\) if the orientations are opposite, i.e. the lines are anti-parallel).
The distance between the lines can be found by computing \(\norm{\e_{12}\vee(\tb{a}\wedge\tb{b})}\),
provided that the lines are normalised.
Examples are shown in Figure~\ref{depicting points at infinity}.

\begin{figure}[t!]
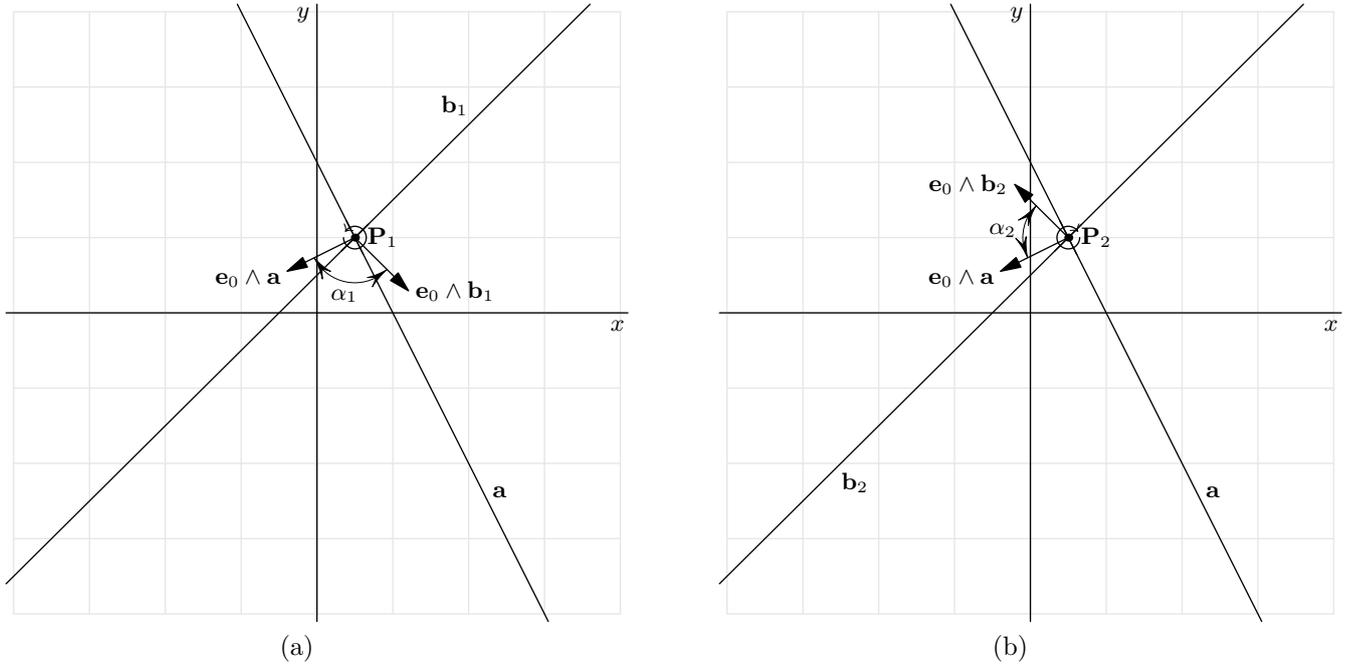

\begin{subfloatenv}{ }
\begin{asy}
import Projective2D;
import Graph2D;
Drawing d = Drawing();

MV a = Line(2,-2,-1);
a = a/norm(a);
MV Na = wedge(e_0,a);

MV b = Line(1,2,-2);
b = b/norm(b);
MV Nb = wedge(e_0,b);

MV P = wedge(a,b);
pair O = topair(P);

d.point(P, "$\textbf{P}_1$", align=E);
d.line(a, "$\textbf{a}$", draw_orientation=false);
d.point_at_infinity(Na, O=O, "$\textbf{e}_0\wedge\textbf{a}$");
d.line(b, "$\textbf{b}_1$", draw_orientation=false);
d.point_at_infinity(Nb, O=O, "$\textbf{e}_0\wedge\textbf{b}_1$", align=E);

d.arc(P, Na,Nb,"$\alpha_1$", CCW);
\end{asy}
\end{subfloatenv}\hfill%
\begin{subfloatenv}{  }
\begin{asy}
import Projective2D;
import Graph2D;
Drawing d = Drawing();

MV a = Line(2,-2,-1);
a = a/norm(a);
MV Na = wedge(e_0,a);

MV b = -Line(1,2,-2);
b = b/norm(b);
MV Nb = wedge(e_0,b);

MV P = wedge(a,b);
pair O = topair(P);

d.point(P, "$\textbf{P}_2$", align=E);
d.line(a, "$\textbf{a}$", draw_orientation=false);
d.point_at_infinity(Na, O=O, "$\textbf{e}_0\wedge\textbf{a}$");
d.line(b, "$\textbf{b}_2$", draw_orientation=false);
d.point_at_infinity(Nb, O=O, "$\textbf{e}_0\wedge\textbf{b}_2$", align=W);

d.arc(P, Na,Nb,"$\alpha_2$", W, CW);
\end{asy}
\end{subfloatenv}
\caption{Intersecting  lines in \E{2}}
\label{intersecting lines in E2}
\end{figure}

For  normalised intersecting lines \(\tb{a}\) and \(\tb{b}\),
\begin{equation}
\tb{a}\tb{b}=\cos{\alpha}+\tb{P}\sin{\alpha},
\label{geometry_ab}
\end{equation}
where \(\alpha\in(0,\pi)\) is the angle between 
the orientation arrows of the lines \(\tb{a}\) and \(\tb{b}\)
and \(\tb{P}=\tb{a}\wedge\tb{b}/\sin{\alpha}\) is a normalised point where the lines intersect.
The resulting point \(\tb{P}\) can be clockwise or counterclockwise depending on the orientation of the lines.
So, the angle \(\alpha\) between the normalised lines \(\tb{a}\) and \(\tb{b}\) is defined by
\begin{equation}
\cos\alpha=\tb{a}\cdot\tb{b}.
\end{equation}
For example, consider lines \(\tb{a}=\tfrac{1}{\sqrt{5}}(2\e_0-2\e_1-\e_2)\)  and \(\tb{b}_1=\tfrac{1}{\sqrt{8}}(\e_0+2\e_1-2\e_2)\),
which have been normalised.
The lines and their orientation are shown in Figure~\ref{intersecting lines in E2}(a).
Computing the inner and outer products gives
\begin{equation*}
\tb{a}\tb{b}_1
=
\tb{a}\cdot\tb{b}_1 + \tb{a}\wedge\tb{b}_1
=
-\tfrac{1}{\sqrt{10}}
+
\tfrac{3}{\sqrt{10}}(\e_{12}+\tfrac{1}{2}\e_{20}+\e_{01}),
\end{equation*}
so that the point \(\tb{P}_1\) where the lines intersect
is counterclockwise and is located at \((\tfrac{1}{2},1)\), while \(\alpha_1=\arccos\left(-\tfrac{1}{\sqrt{10}}\right)\approx 108^\circ\).
On the other hand, for \(\tb{b}_2=-\tb{b}_1\), I get
\begin{equation*}
\tb{a}\tb{b}_2
=
\tb{a}\cdot\tb{b}_2 + \tb{a}\wedge\tb{b}_2
=
\tfrac{1}{\sqrt{10}}
-
\tfrac{3}{\sqrt{10}}(\e_{12}+\tfrac{1}{2}\e_{20}+\e_{01}),
\end{equation*}
so that the intersection point \(\tb{P}_2\) is clockwise and \(\alpha_2=\arccos\left(\tfrac{1}{\sqrt{10}}\right)\approx 72^\circ\)
(see Figure~\ref{intersecting lines in E2}(b)).

\begin{figure}[t!]
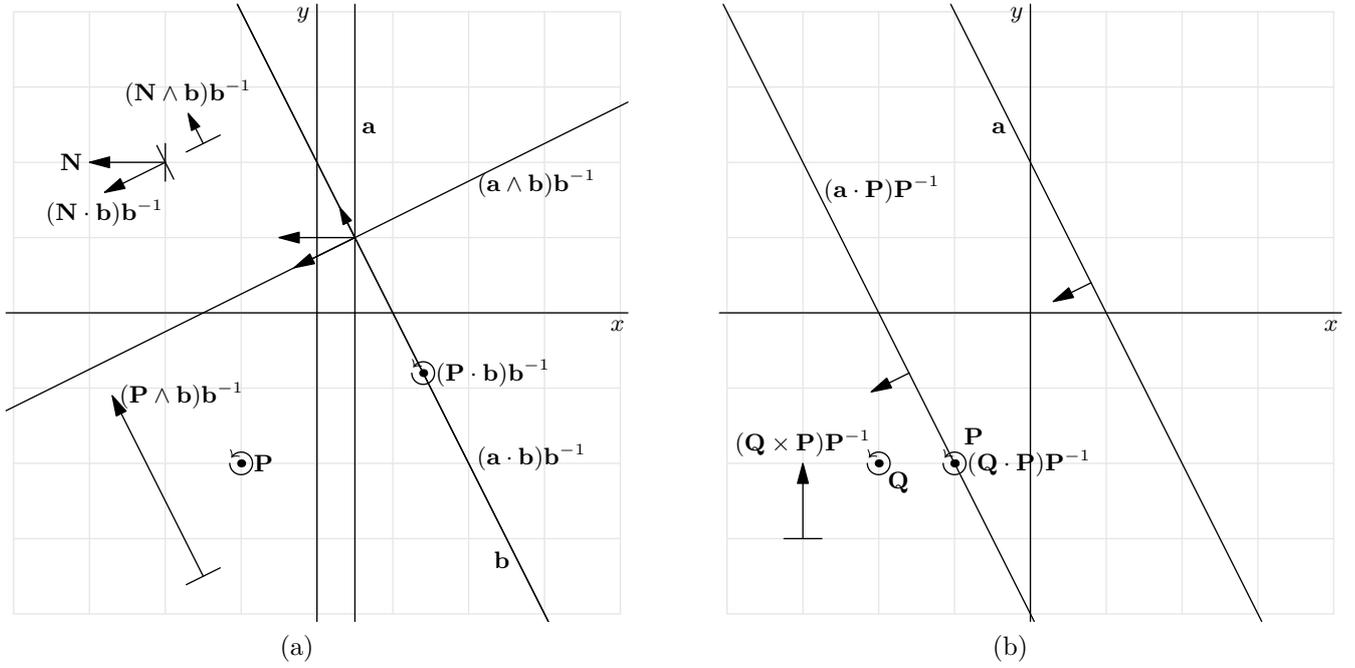

\begin{subfloatenv}{ }
\begin{asy}
import Projective2D;
import Graph2D;
Drawing d = Drawing();

var b = Line(2,-2,-1)/4;
var P = Point(1,-1,-2);

d.line(b, "$\textbf{b}$", position=0.1,align=W, draw_orientation=false);
d.point(P, "$\textbf{P}$", align=E);

var Q = dot(P,b)/b;
d.point(Q, "$(\textbf{P}\cdot\textbf{b})\textbf{b}^{-1}$", align=E);
d.point_at_infinity(wedge(P,b)/b,"$(\textbf{P}\wedge\textbf{b})\textbf{b}^{-1}$",O=(-1.5,-3.5),align=E);

var N = wedge(e_0,Line(0,-1,0));
d.point_at_infinity(N, "$\textbf{N}$", O=(-2,2));
d.point_at_infinity(dot(N,b)/b, "$(\textbf{N}\cdot\textbf{b})\textbf{b}^{-1}$", O=(-2,2),align=S);
d.point_at_infinity(wedge(N,b)/b, "$(\textbf{N}\wedge\textbf{b})\textbf{b}^{-1}$", O=(-1.5,2+1/4), align=(0,1));

var a = Line(1/2,-1,0);
d.line(a,"$\textbf{a}$", position=0.8, draw_orientation=true, O=(1/2,1));
d.line(dot(a,b)/b, "$(\textbf{a}\cdot\textbf{b})\textbf{b}^{-1}$", position=0.25, draw_orientation=true, O=(1/2,1));
d.line(wedge(a,b)/b, "$(\textbf{a}\wedge\textbf{b})\textbf{b}^{-1}$", position=0.8, draw_orientation=true, O=(1/2,1));

\end{asy}
\end{subfloatenv}\hfill%
\begin{subfloatenv}{ }
\begin{asy}
import Projective2D;
import Graph2D;
Drawing d = Drawing();

var a = Line(2,-2,-1)/4;
var P = Point(1,-1,-2);

d.line(a, "$\textbf{a}$", position=0.8,align=W);
d.point(P, "$\textbf{P}$", align=(1,1.5), draw_orientation=false);

var b = dot(a,P)/P;
d.line(b, "$(\textbf{a}\cdot\textbf{P})\textbf{P}^{-1}$", position=0.7, align=E);

var Q = Point(1,-2,-2);
d.point(Q,"$\textbf{Q}$");
d.point(dot(Q,P)/P,"$(\textbf{Q}\cdot\textbf{P})\textbf{P}^{-1}$",align=E);
d.point_at_infinity(cross(Q,P)/P,"$(\textbf{Q}\times\textbf{P})\textbf{P}^{-1}$",O=(-3,-3));

\end{asy}
\end{subfloatenv}
\caption{Projection and rejection in \E{2}}
\label{projection and rejection in E2}
\end{figure}

A geometric object, i.e.\ a line or a point, can be split into projection and rejection with respect to another
geometric object as long as it is invertible. This split is captured naturally by the two components of the geometric product.
For instance, \(\tb{a}=(\tb{a}\tb{b})\tb{b}^{-1}=(\tb{a}\cdot\tb{b}+\tb{a}\wedge\tb{b})\tb{b}^{-1}=
(\tb{a}\cdot\tb{b})\tb{b}^{-1}+(\tb{a}\wedge\tb{b})\tb{b}^{-1}\), where
\(\ts{proj}(\tb{a};\tb{b})=(\tb{a}\cdot\tb{b})\tb{b}^{-1}\) is the projection of a line \(\tb{a}\) on a finite line \(\tb{b}\)
and \(\ts{rej}(\tb{a};\tb{b})=(\tb{a}\wedge\tb{b})\tb{b}^{-1}\) is the rejection of  \(\tb{a}\) by \(\tb{b}\).
Note that the attitude of \(\ts{proj}(\tb{a};\tb{b})\) coincides with that of \(\tb{b}\).
\((\tb{P}\cdot\tb{b})\tb{b}^{-1}\) is the projection of a point \(\tb{P}\) on a finite line \(\tb{b}\)
and \((\tb{P}\wedge\tb{b})\tb{b}^{-1}\) is the rejection
(the same applies to a point at infinity \(\tb{N}\)).
This is illustrated in Figure~\ref{projection and rejection in E2}(a), where
\(\tb{a}=\tfrac{1}{2}\e_0-\e_1\), 
\(\tb{P}=\e_{12}-\e_{20}-2\e_{01}\), \(\tb{N}=-\e_{01}\),
and
\(\tb{b}=\tfrac{1}{4}(2\e_0-2\e_1-\e_2)\).
Observe that both \((\tb{P}\cdot\tb{b})\tb{b}^{-1}\) and \((\tb{N}\cdot\tb{b})\tb{b}^{-1}\) lie on the line \(\tb{b}\).

The projection of a line \(\tb{a}\) on a finite point \(\tb{P}\) is given by \((\tb{a}\cdot\tb{P})\tb{P}^{-1}\);
the rejection \((\tb{a}\wedge\tb{P})\tb{P}^{-1}\)  yields the line at infinity.
The projection of a finite point \(\tb{Q}\) on \(\tb{P}\) is given by \((\tb{Q}\cdot\tb{P})\tb{P}^{-1}\)
and the rejection of \(\tb{Q}\) by \(\tb{P}\) is a point at infinity \((\tb{Q}\times\tb{P})\tb{P}^{-1}\).
Since \(\tb{N}\cdot\tb{P}=0\), the rejection \((\tb{N}\times\tb{P})\tb{P}^{-1}\)
  of a point at infinity \(\tb{N}\) by \(\tb{P}\) coincides with  \(\tb{N}\).
See Figure~\ref{projection and rejection in E2}(b).

There are two ways to define a reflection in \E{2} depending on whether one takes the top-down or bottom-up view.
Both reflections yield the same result in terms of attitude and weight, but the orientation may be different.

\begin{figure}[h!]
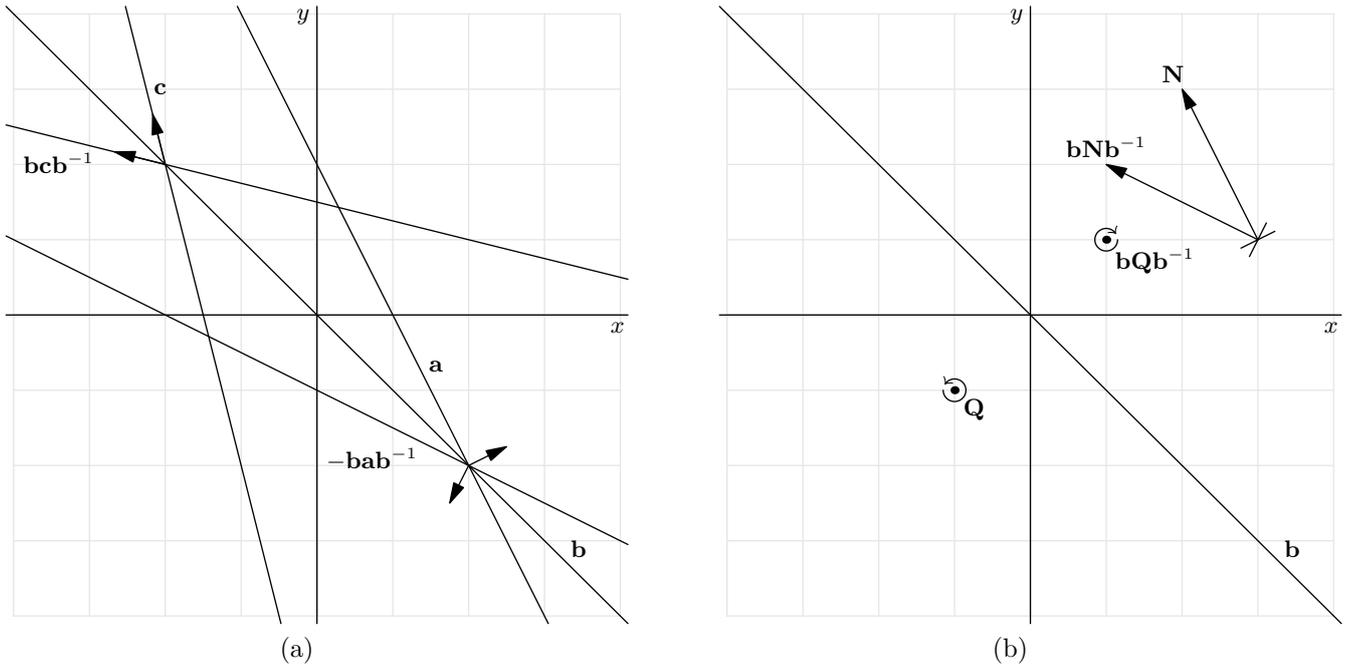

\begin{subfloatenv}{ }
\begin{asy}
import Projective2D;
import Graph2D;
Drawing d = Drawing();

var a = -Line(2,-2,-1)/4;
var b = Line(0,1,1)/4;
write(b);
d.line(a, "$\textbf{a}$", position=0.6, align=NE, draw_orientation=true,O=(2,-2));
d.line(b, "$\textbf{b}$", position=0.9, align=NE, draw_orientation=false);

d.line(-b*a/b, "$-\textbf{b}\textbf{a}\textbf{b}^{-1}$", position=0.33, align=SW, draw_orientation=true,O=(2,-2));

var c = -join(Point(1,-2,2),Point(1,-1.5,0));
c/=3;
d.line(c, "$\textbf{c}$", position=0.15, align=NE, bottomup_orientation=true, draw_orientation=true,O=(-2,2));
d.line(b*c/b, "$\textbf{b}\textbf{c}\textbf{b}^{-1}$", position=0.15, align=SW, bottomup_orientation=true, draw_orientation=true,O=(-2,2));

\end{asy}
\end{subfloatenv}\hfill%
\begin{subfloatenv}{ }
\begin{asy}
import Projective2D;
import Graph2D;
Drawing d = Drawing();

var b = Line(0,1,1)/4;
write(b);
d.line(b, "$\textbf{b}$", position=0.9, align=NE, draw_orientation=false);

var N = wedge(e_0, Line(0,-1,2));
d.point_at_infinity(N, "$\textbf{N}$", O=(3,1));
d.point_at_infinity(b*N/b, "$\textbf{b}\textbf{N}\textbf{b}^{-1}$", align=(0,1), O=(3,1));

var Q = Point(1,-1,-1);
d.point(Q,"$\textbf{Q}$");
d.point(b*Q/b,"$\textbf{b}\textbf{Q}\textbf{b}^{-1}$");

\end{asy}
\end{subfloatenv}
\caption{Reflection in a line in \E{2}}
\label{reflections in a line in E2}
\end{figure}

The reflection of a line \(\tb{a}\) in a finite line \(\tb{b}\) 
defined by \(-\tb{b}\tb{a}\tb{b}^{-1}\) is consistent with the top-down orientation of \(\tb{a}\).
It is illustrated in Figure~\ref{reflections in a line in E2}(a), where
\(\tb{a}= \tfrac{1}{4}(-2\e_0+2\e_1+\e_2)\) and \(\tb{b}=\e_1+\e_2\).
Observe that the top-down reflection of a perpendicular line yields the same line.
Indeed, if \(\tb{a}\cdot\tb{b}=0\), then
\(-\tb{b}\tb{a}\tb{b}^{-1}=-(\tb{b}\wedge\tb{a})\tb{b}^{-1}=
(\tb{a}\wedge\tb{b})\tb{b}^{-1}=\tb{a}\tb{b}\tb{b}^{-1}=\tb{a}\).
On the other hand, the top-down reflection of a parallel line reverses its orientation.
For the line at infinity, \(-\tb{b}\e_0\tb{b}^{-1}=\e_0\).
Note also that the orientation and weight of the line in which the reflection is carried out
have no effect on the reflection.

The top-down reflection of a point in a line can be understood by applying suitable line reflections.
For instance, consider the top-down reflection of a point at infinity \(\tb{N}\) in  \(\tb{b}\), which is given by  \(\tb{b}\tb{N}\tb{b}^{-1}\). 
See Figure~\ref{reflections in a line in E2}(b) where \(\tb{N}=\e_0\wedge(-\e_1+2\e_2)\).
The point at infinity \(\tb{N}\) is a stack of lines and its orientation is the direction of stepping
through the lines in the stack.
Reflecting each line in the stack yields a new stack, 
which corresponds to   \(\tb{b}\tb{N}\tb{b}^{-1}\).
The same formula applies to the reflection of finite points, i.e. \(\tb{b}\tb{Q}\tb{b}^{-1}\) for a point \(\tb{Q}\).
An example is shown in Figure~\ref{reflections in a line in E2}(b)
where \(\tb{Q}=\e_{12}-\e_{20}-\e_{01}\). 
Note that the reflected point has the opposite orientation, i.e.\ clockwise.
Indeed, since the sheaf of lines representing \(\tb{Q}\) is counterclockwise,
reflecting each line in the sheaf yields a clockwise sheaf of lines intersecting 
at the point \(\tb{b}\tb{Q}\tb{b}^{-1}\).
Since a point can be written as the intersection of two lines, the top-down reflection of the point \(\tb{Q}=\tb{q}_1\wedge\tb{q}_2\) in a line \(\tb{b}\)
can be found by reflecting \(\tb{q}_1\) and \(\tb{q}_2\) in \(\tb{b}\) and intersecting the reflected lines.
Namely, \(\tb{b}\tb{Q}\tb{b}^{-1}= \tb{b}(\tb{q}_1\wedge\tb{q}_2)\tb{b}^{-1}= (-\tb{b}\tb{q}_1\tb{b}^{-1})\wedge(-\tb{b}\tb{q}_2\tb{b}^{-1})\).
The minus signs cancel out.

\begin{figure}[t!]
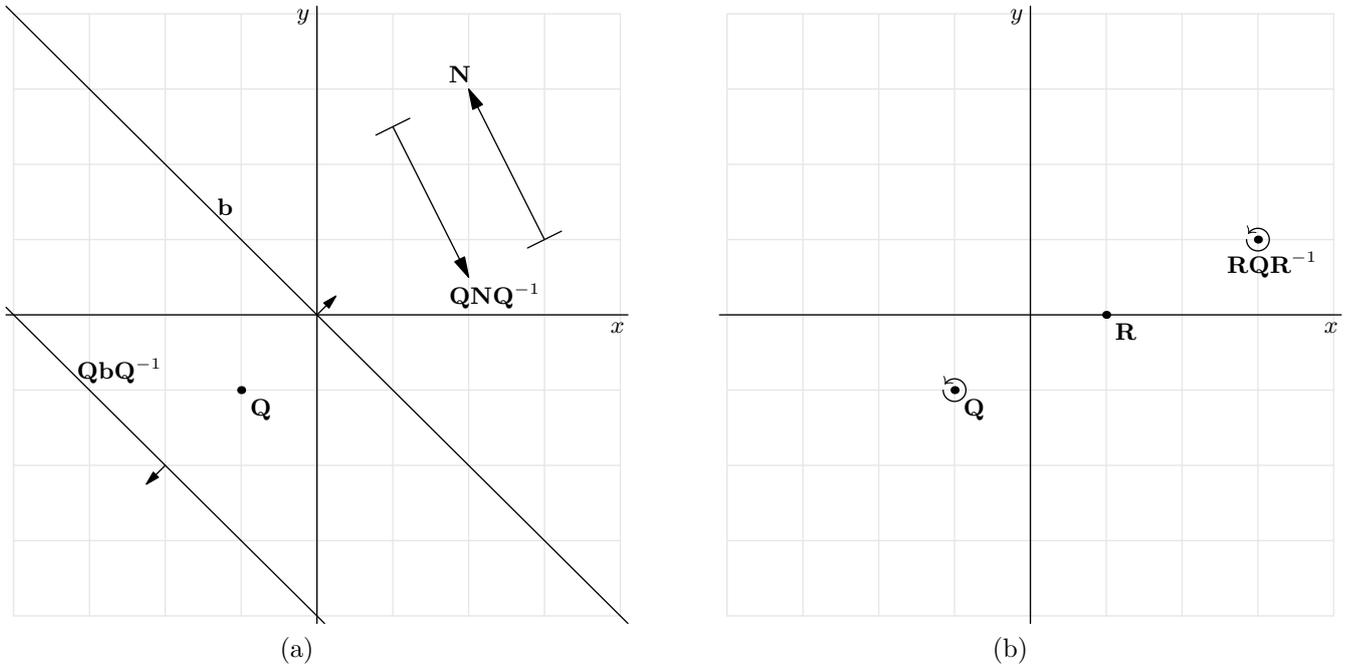

\begin{subfloatenv}{ }
\begin{asy}
import Projective2D;
import Graph2D;
Drawing d = Drawing();

var b = Line(0,1,1)/4;
var Q = Point(1,-1,-1);

d.line(b, "$\textbf{b}$", position=0.35, align=(0.1,1));
d.point(Q, "$\textbf{Q}$", draw_orientation=false);

//d.point(b*Q/b, "$\textbf{b}\textbf{Q}\textbf{b}^{-1}$", align=(0.25,-1));
d.line(Q*b/Q, "$\textbf{Q}\textbf{b}\textbf{Q}^{-1}$", position=0.8, align=E);

var N = wedge(e_0, -e_1+2*e_2);
d.point_at_infinity(N, "$\textbf{N}$", O=(3,1));
d.point_at_infinity(Q*N/Q, "$\textbf{Q}\textbf{N}\textbf{Q}^{-1}$", O=(1,2.5));

\end{asy}
\end{subfloatenv}\hfill%
\begin{subfloatenv}{ }
\begin{asy}
import Projective2D;
import Graph2D;
Drawing d = Drawing();

var Q = Point(1,-1,-1);
var R = Point(1,1,0);

d.point(R, "$\textbf{R}$", draw_orientation=false);
d.point(Q, "$\textbf{Q}$");

d.point(R*Q/R, "$\textbf{R}\textbf{Q}\textbf{R}^{-1}$",align=(0.25,-1));

\end{asy}
\end{subfloatenv}
\caption{Reflection in a point in \E{2}}
\label{reflection in a point in E2}
\end{figure}
\begin{figure}[h!]
\begin{subfloatenv}{ }
\begin{asy}
import Figure2D;

Figure d = Figure();
var b = Line(0,1,1);
d.line(b,"$\textbf{b}$", draw_orientation=false);

var Q = Point(1,0,1);
d.point(Q,"$\textbf{Q}$",align=E);

var a = join(Point(1,3,-3), Q);
a = a/norm(a);
d.line(a,"$\textbf{a}$", position=0.1, align=E);

MV f(real g, MV m) { return m + (g-1)*wedge(m,b)/b; }
MV fI(real g, MV m) { return m + (g-1)*join(m,b)*I/b; }

real g = 3;
d.point(fI(g,Q), "$\mathsf{proj}(\textbf{Q};\textbf{b})+"+string(g)+"\mathsf{rej}(\textbf{Q};\textbf{b})$", align=(0.8,1));
d.line(f(g,a), position=0.1, "$\mathsf{proj}(\textbf{a};\textbf{b})+"+string(g)+"\mathsf{rej}(\textbf{a};\textbf{b})$", align=E);

\end{asy}
\end{subfloatenv}\hfill%
\begin{subfloatenv}{ }
\begin{asy}
import Figure2D;

Figure d = Figure();
var P = Point(1,0,2);
d.point(P,"$\textbf{P}$", draw_orientation=false, align=E);

var Q = Point(1,0,1);
d.point(Q,"$\textbf{Q}$",align=SW);

var a = join(Point(1,3,-3), Q);
a = a/norm(a);
d.line(a,"$\textbf{a}$", position=0.26, align=W);

MV f(real g, MV m) { return m + (g-1)*join(m,P)*I/P; }

real g = 3;
d.point(f(g,Q), "$\mathsf{proj}(\textbf{Q};\textbf{P})+"+string(g)+"\mathsf{rej}(\textbf{Q};\textbf{P})$", align=SW);
d.line(f(g,a), position=0.55, "$\mathsf{proj}(\textbf{a};\textbf{P})+"+string(g)+"\mathsf{rej}(\textbf{a};\textbf{P})$", align=W);

\end{asy}
\end{subfloatenv}
\caption{Scaling in \E{2}}
\label{scaling}
\end{figure}

The top-down reflection of a line \(\tb{a}\) in a finite point \(\tb{Q}\) is given by \(\tb{Q}\tb{a}\tb{Q}^{-1}\).
The reflection of \(\tb{a}\) in \(\tb{Q}\) is equivalent to the reflection of \(\tb{a}\)
in any two perpendicular lines intersecting at the point \(\tb{Q}\).
In other words, the reflection in a point is equivalent to the consecutive reflections in two perpendicular lines.
A finite point \(\tb{Q}\) can be represented as the intersection of two lines,
one of which is parallel to  \(\tb{a}\) and the other is perpendicular.
Since the top-down reflection of \(\tb{a}\) in a perpendicular line leaves \(\tb{a}\) unchanged,
the top-down reflection of \(\tb{a}\) in \(\tb{Q}\) is  equivalent to the reflection in a line that passes 
through \(\tb{Q}\) and is parallel to  \(\tb{a}\).
This is illustrated in Figure~\ref{reflection in a point in E2}(a).

\(\tb{R}\tb{Q}\tb{R}^{-1}\) gives the top-down reflection of  \(\tb{Q}\) in a finite point \(\tb{R}\).
An example for \(\tb{R}=\e_{12}+\e_{20}\) and \(\tb{Q}=\e_{12}-\e_{20}-\e_{01}\)
is shown in Figure~\ref{reflection in a point in E2}(b).
If \(\tb{Q}=\tb{q}_1\wedge\tb{q}_2\), then 
\(\tb{R}\tb{Q}\tb{R}^{-1}=\tb{R}(\tb{q}_1\wedge\tb{q}_2)\tb{R}^{-1}=(\tb{R}\tb{q}_1\tb{R}^{-1})\wedge(\tb{R}\tb{q}_2\tb{R}^{-1})\).
So the reflection of \(\tb{Q}\) in  \(\tb{R}\) can be constructed by reflecting the lines which constitute \(\tb{Q}\)
and taking the intersection of the reflected lines.
The resulting orientation of the reflected point is the same as that of the original point.
Indeed, 
reflecting each line in the sheaf dually represented by \(\tb{Q}\) gives a sheaf of 
lines intersecting at \(\tb{R}\tb{Q}\tb{R}^{-1}\) with the same top-down orientation as \(\tb{Q}\).
In general, \((-1)^{kl}A_k{B}_l{A^{-1}_k}\) gives the top-down reflection of a blade \(B_l\) in an invertible blade \(A_k\).

I assume that the bottom-up reflection of a point  \(\tb{Q}\) in a finite point \(\tb{R}\) is given by \(\tb{R}\tb{Q}\tb{R}^{-1}\),
which implies that the bottom-up orientation of the reflected point is the same as that of the original point, e.g.\ if \(\tb{Q}\) is positively oriented,
then \(\tb{R}\tb{Q}\tb{R}^{-1}\) is also positively oriented.
The bottom-up reflections of objects other than points are constructed from point reflections.
For instance, the reflection of a line \(\tb{a}=\tb{P}_1\vee\tb{P}_2\), which passes through the points \(\tb{P}_1\) and \(\tb{P}_2\), in a point \(\tb{Q}\)
is obtained by reflecting \(\tb{P}_1\) and \(\tb{P}_2\) in \(\tb{Q}\) and then taking the join of the reflected points.
Namely, \(\tb{Q}\tb{a}\tb{Q}^{-1}=\tb{Q}(\tb{P}_1\vee\tb{P}_2)\tb{Q}^{-1}=(\tb{Q}\tb{P}_1\tb{Q}^{-1})\vee(\tb{Q}\tb{P}_2\tb{Q}^{-1})\).

I further assume that the bottom-up reflection of a point \(\tb{Q}\) in a line \(\tb{b}\) is given by \(\tb{b}\tb{Q}\tb{b}^{-1}\),
which implies that the bottom-up orientation of the reflected point is the opposite of the original point, e.g.\ if \(\tb{Q}\) is positively oriented,
then \(\tb{b}\tb{Q}\tb{b}^{-1}\) is negatively oriented.
The bottom-up reflection of a line \(\tb{a}=\tb{P}_1\vee\tb{P}_2\) in an invertible line \(\tb{b}\) can be constructed by reflecting
\(\tb{P}_1\) and \(\tb{P}_2\) in \(\tb{b}\) and taking the join of the reflected points. 
Namely, \(\tb{b}\tb{a}\tb{b}^{-1} = \tb{b}(\tb{P}_1\vee\tb{P}_2)\tb{b}^{-1}=(\tb{b}\tb{P}_1\tb{b}^{-1})\vee(\tb{b}\tb{P}_2\tb{b}^{-1}) \).
Note that, unlike the top-down reflection, the bottom-up reflection of a line in another line does not require the minus sign.
The bottom-up reflection of a line \(\tb{c}\) in \(\tb{b}\) 
is shown in Figure~\ref{reflections in a line in E2}(a) for \(\tb{c}=\e_0+\tfrac{2}{3}\e_1+\tfrac{1}{6}\e_2 \).
The bottom-up reflection of a perpendicular line reverses its orientation and that of the parallel line
maintains its orientation.

In general, the bottom-up reflection is given by \((-1)^{nk(l-1)}A_k B_l{A^{-1}_k}\), where \(n\) is the dimension of the metric space.
Since \(n=2\) in \E{2}, the above expression for the bottom-up reflection simplifies to \(A_k B_l{A^{-1}_k}\) in \E{2}.

Scaling can be obtained by considering the following combination of projection and rejection:
\begin{equation}
\ts{scale}(\tb{a};\tb{b},\gamma)=\ts{proj}(\tb{a};\tb{b})+\gamma\ts{rej}(\tb{a};\tb{b}),
\label{scaling definition}
\end{equation}
which scales \(\tb{a}\) by a factor of \(\gamma\) with respect to \(\tb{b}\);
\(\gamma=1\) gives \(\tb{a}\).
Equation~(\ref{scaling definition}) can also be written as \(\ts{scale}(\tb{a};\tb{b},\gamma)=\tb{a}+(\gamma-1)(\tb{a}\wedge\tb{b})\tb{b}^{-1}\)
or  \(\ts{scale}(\tb{a};\tb{b},\gamma)=\gamma\tb{a}+(1-\gamma)(\tb{a}\cdot\tb{b})\tb{b}^{-1}\).
In general, \(\ts{scale}(A_k;B_l,\gamma)=\ts{proj}(A_k;B_l)+\gamma\ts{rej}(A_k;B_l)\)
scales \(A_k\) with respect to an invertible blade \(B_l\) as long as 
the geometric product of \(A_k\) and \(B_l\) can be split  into two components for projection and rejection.
Scaling is illustrated in Figure~\ref{scaling}, where \(\tb{b}=\e_1+\e_2\), \(\tb{P}=\e_{12}+2\e_{01}\),
\(\tb{a}=\tfrac{1}{5}(-3\e_0+4\e_1+3\e_2)\), \(\tb{Q}=\e_{12}+\e_{01}\), and the scaling factor \(\gamma=3\).
Scaling a line with respect to  another line changes its norm and scaling elements at infinity changes their weight.

Rotations and translations in \E{2} can be obtained via two consecutive reflections
in lines \(\tb{a}\) and \(\tb{b}\), which satisfy \(\tb{a}^2=1\), \(\tb{b}^2=1\).
For the intersecting lines, \(\tb{a}\tb{b}\tb{P}\tb{b}^{-1}\tb{a}^{-1}\) gives
the rotation of a point \(\tb{P}\) around the point \(\tb{a}\wedge\tb{b}\) by twice the angle between the lines
\(\tb{a}\) and \(\tb{b}\).
If \(\tb{a}\) and \(\tb{b}\) are parallel,
\(\tb{a}\tb{b}\tb{P}\tb{b}^{-1}\tb{a}^{-1}\) gives
the translation of \(\tb{P}\) in the direction perpendicular to \(\tb{a}\) and \(\tb{b}\)
by twice the distance between the lines.
Note that \((\tb{a}\tb{b})^{-1}=\tb{b}^{-1}\tb{a}^{-1}\) and so the transformations
can be written as \(S\tb{P}S^{-1}\), where \(S=\tb{a}\tb{b}\) and \(\tb{a}^2=1\), \(\tb{b}^2=1\).

\begin{figure}[t!]
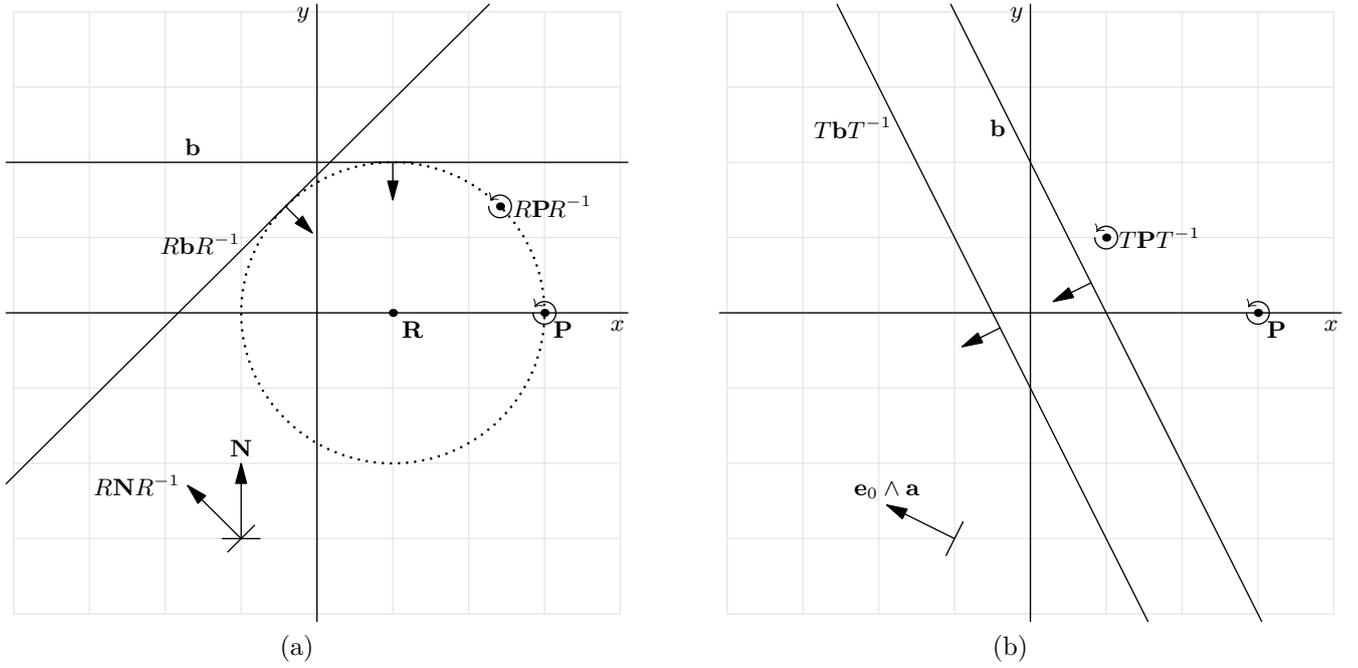

\begin{subfloatenv}{ }
\begin{asy}
import Projective2D;
import Graph2D;
Drawing d = Drawing();

var b = Line(2,0,-1)/2;
var P = Point(1,3,0);
var R = Point(1,1,0);

path arc = arc((1,0), 2, 0, 360);
draw(arc, currentpen+linewidth(0.2 bp)+Dotted);  

d.point(P, "$\textbf{P}$");
d.point(R, "$\textbf{R}$", draw_orientation=false);

real alpha = pi/4;
var Rot = cos(alpha/2)-sin(alpha/2)*R;

d.point(Rot*P/Rot, "$R\textbf{P}\!R^{-1}$",align=E);

d.line(b, "$\textbf{b}$", position=0.7, align=(0,1), draw_orientation=false);
d.point_at_infinity(wedge(e_0, b), O=(1,2));

d.line(Rot*b/Rot, "$R\textbf{b}R^{-1}$", position=0.5, align=(-1,0), draw_orientation=false);
d.point_at_infinity(wedge(e_0, Rot*b/Rot), O=rotate(45, topair(R))*(1,2));

var N = wedge(e_0, Line(0,0,1));
d.point_at_infinity(N, "$\textbf{N}$", O=(-1,-3));
d.point_at_infinity(Rot*N/Rot, "$R\textbf{N}R^{-1}$", align=(-1,0), O=(-1,-3));

\end{asy}
\end{subfloatenv}\hfill%
\begin{subfloatenv}{ }
\begin{asy}
import Projective2D;
import Graph2D;
Drawing d = Drawing();

var a = Line(2,-2,-1)/4;
var P = Point(1,3,0);

d.line(a, "$\textbf{b}$", position=0.8, align=W);
d.point(P, "$\textbf{P}$");

var line = Line(0,-2,1);
var N = wedge(e_0, line/norm(line));
d.point_at_infinity(N, "$\textbf{e}_0\wedge\textbf{a}$", O=(-1,-3), align=(0,1));

var Tran = 1 - sqrt(5)/2 * N;
d.line(Tran*a/Tran, "$T\textbf{b}T^{-1}$", position=0.8, align=W);
d.point(Tran*P/Tran, "$T\textbf{P}T^{-1}$", align=E);

\end{asy}
\end{subfloatenv}
\caption{Rotation and translation in \E{2}}
\label{rotation and translation}
\end{figure}

Multivectors that can be written as the product of an even number of lines,
all of which square to unity, 
form a group called the Spin group.
I will refer to them as spinors.
Spinors are thus responsible for rotations and translations in \E{2}.
Spinors are even multivectors and each spinor satisfies \(S\reverse{S}=1\), so that \(S^{-1}=\reverse{S}\).
Moreover, any proper motion\footnote{A proper motion is equivalent to a solid body motion, 
whereas an improper motion also includes a reflection.} in \E{2}, e.g.\ rotation or translation, can be
generated by a spinor that has the form \(S=e^{A}\) for some bivector \(A\).

Spinors form a group with respect to the geometric product, i.e.\
the product of two spinors is a spinor and the inverse of a spinor is a spinor.
It is not an algebra though, since the sum of two spinors is generally not a spinor.
The Spin group is a smooth manifold, in which
the product and the inverse are smooth functions, so it is a Lie group.
In Euclidean space, this group consists of multivectors \(e^{A}\), where \(A\) is an arbitrary bivector.
Lie algebra of the Spin group consists of bivectors,
with the commutator used for multiplication in the algebra.
Indeed, the commutator of two bivectors is a bivector and, therefore, belongs in the algebra.
The cross product satisfies Jacobi identity
\begin{equation}
(A\times B)\times C + (B\times C)\times A + (C \times A)\times B=0,
\end{equation}
where \(A,B,C\) are bivectors,
so the Lie algebra is not associative, i.e.\
\((A\times B)\times C\ne A\times(B\times C) \)
, but it is anti-commutative since \(A\times B=-B\times A\).

Consider a spinor
\begin{equation}
R
=
e^{-\tfrac{1}{2}\alpha\tb{R}}
=
\cos{\tfrac{\alpha}{2}}-\tb{R}\sin{\tfrac{\alpha}{2}},
\end{equation}
where a finite point \(\tb{R}=\e_{12}+x\e_{20}+y\e_{01}\) is normalised and counterclockwise,
and
\begin{equation}
T=
e^{-\tfrac{1}{2}\lambda\e_0\wedge\tb{a}}
=
1-\tfrac{1}{2}\lambda\e_0\wedge\tb{a},
\end{equation}
where \(\tb{a}=a\e_1+b\e_2\) is normalised.

\(R\tb{P}\!R^{-1}\) is a point that results from the counterclockwise rotation of \(\tb{P}\) around \(\tb{R}\)
by the angle \(\alpha>0\); the rotation is clockwise if \(\alpha<0\).
\(T\tb{P}T^{-1}\) is a point that results from the translation of \(\tb{P}\) by \(\lambda\) in the direction of \((a,b)\),
which is perpendicular to \(\tb{a}\).
\(R\tb{b}R^{-1}\) and \(T\tb{b}T^{-1}\) give the same as above for a line \(\tb{b}\).
Figure~\ref{rotation and translation}(a) shows 
an example of the counterclockwise rotation of \(\tb{P}=\e_{12}+3\e_{20}\),
\(\tb{b}=\e_0-\tfrac{1}{2}\e_2\), and \(\tb{N}=\e_0\wedge\e_2\)
by a spinor  \(R=e^{-\tfrac{1}{2}\tfrac{\pi}{4}\tb{R}}\), 
where \(\tb{R}=\e_{12}+\e_{20}\) is normalised and counterclockwise.
Figure~\ref{rotation and translation}(b) shows the translation of \(\tb{P}=\e_{12}+3\e_{20}\)
and \(\tb{b}=\tfrac{1}{4}(2\e_0-2\e_1-\e_2)\) by the spinor \(T=1-\tfrac{\sqrt{5}}{2}\e_0\wedge\tb{a}\), where
\(\tb{a}=\tfrac{1}{\sqrt{5}}(-2\textbf{e}_1 +\textbf{e}_2)\) is normalised.
Translation has no effect on points at infinity.
If the line to be rotated passes through the point of rotation, i.e.\ \(\tb{b}\wedge\tb{R}=0\) and so
\(\tb{b}\tb{R}=-\tb{R}\tb{b}\), the expression for rotation simplifies: 
\(R\tb{b}R^{-1}=\tb{b}(\cos{\alpha} + \tb{R}\sin{\alpha})=(\cos{\alpha} - \tb{R}\sin{\alpha})\tb{b}\).

In general, the rotation and translation of multivector \(M\) is performed with \(RMR^{-1}\)
and \(TMT^{-1}\), respectively.
For instance, the rotation generated by the spinor \(R\) can be applied to the spinor \(T\), which gives
 \(T'=RTR^{-1}=1-\tfrac{\sqrt{5}}{2}(\e_0\wedge\tb{a}')\).
\(T'\) translates by the same amount as \(T\) but in the direction perpendicular to \(\tb{a}'=R\tb{a}R^{-1}\).

\subsection{Linear functions}

A linear function \(\ts{f}:\R{3*}\to\R{3*}\) is defined on vectors in the dual model space \R{3*} or lines in \T{2}
and satisfies the following properties:
\begin{equation}
\begin{split}
&\ts{f}(\tb{a}+\tb{b})=\ts{f}(\tb{a})+\ts{f}(\tb{b}),\\
&\ts{f}(\alpha\tb{a})=\alpha\ts{f}(\tb{a}).
\end{split}
\end{equation}

\begin{figure}[t!]
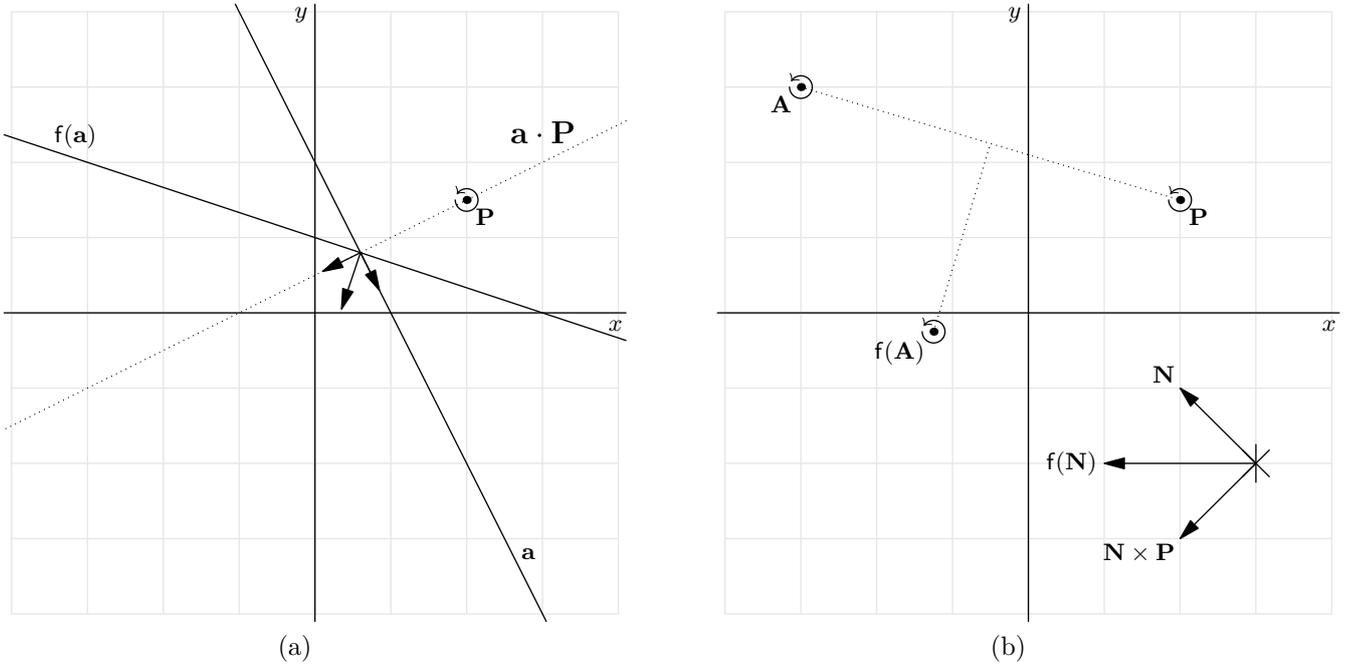

\begin{subfloatenv}{ }
\begin{asy}
import Projective2D;
import Graph2D;
Drawing d = Drawing();

var P = Point(1,2,1.5);

d.point(P, "$\textbf{P}$");
var e1 = Line(0,1,0);
var e2 = Line(0,0,1);

var a = Line(2,-2,-1)/4;

MV f(MV a) {return a + dot(a,P);}

d.line(a, "$\textbf{a}$", position=0.1, draw_orientation=true, O=topair(dot(P,a)/a));
d.line(dot(a,P), "$\textbf{a}\cdot\textbf{P}$", position=0.1, draw_orientation=true, O=topair(dot(P,a)/a),pen=dotted);
d.line(f(a), "$\mathsf{f}(\textbf{a})$", position=0.9, draw_orientation=true, O=topair(dot(P,a)/a));

\end{asy}
\end{subfloatenv}\hfill%
\begin{subfloatenv}{ }
\begin{asy}
import Projective2D;
import Graph2D;
Drawing d = Drawing();

var P = Point(1,2,1.5);
d.point(P, "$\textbf{P}$");

MV f(MV A) { return A + cross(A,P) + dot(A,P)/P; }

var B = Point(1,0,3.5);
//d.point(B, "$\textbf{B}$");
//d.point(f(B), "$\mathsf{f}(\textbf{B})$");

var A = Point(1,-3,3);
d.point(A, "$\textbf{A}$", align=SW);
d.point(f(A), "$\mathsf{f}(\textbf{A})$", align=SW);
draw(topair(A)--topair(P),dotted);
draw(topair(f(A))--topair(A+(1/2-1)*cross(A,P)/P),dotted);

var N = wedge(e_0,Line(0,-1,1));
d.point_at_infinity(N, "$\textbf{N}$", O=(3,-2));
d.point_at_infinity(f(N), "$\mathsf{f}(\textbf{N})$", O=(3,-2));
d.point_at_infinity(cross(N,P), "$\textbf{N}\times\textbf{P}$", O=(3,-2));

\end{asy}
\end{subfloatenv}
\caption{Linear function \ts{f}}
\label{linear function f}
\end{figure}

In this section, for the sake of example I will consider the function 
\begin{equation}
\ts{f}(\tb{a})=\tb{a}+\tb{a}\cdot\tb{P},
\label{example definition}
\end{equation}
where \(\tb{P}=\e_{12}+2\e_{20}+\tfrac{3}{2}\e_{10}\) is fixed and the metric is Euclidean.
The action of \(\ts{f}\) on a line \(\tb{a}\) is shown in Figure~\ref{linear function f}(a).
It is equivalent to 1) a counterclockwise rotation by \(45^\circ\) 
around the projection of  \(\tb{P}\) on the line \(\tb{a}\) and 2) an increase of the weight by a factor of \(\sqrt{2}\).

The definition of a linear function \(\ts{f}\) can be extended to bivectors by
\begin{equation}
\ts{f}(\tb{a}\wedge\tb{b})=\ts{f}(\tb{a})\wedge\ts{f}(\tb{b}).
\label{f on bivector}
\end{equation} 
This implies that the action of  \(\ts{f}\) on  \(\tb{a}\wedge\tb{b}\),
which results from intersecting lines \(\tb{a}\) and \(\tb{b}\), is equivalent to  the intersection of
the lines \(\ts{f}(\tb{a})\) and \(\ts{f}(\tb{b})\).
In other words, the action of \(\ts{f}\) commutes with the outer product.

For the linear function defined by (\ref{example definition}), 
the right-hand side of Equation~(\ref{f on bivector}) can be expressed in terms of a bivector \(\tb{A}=\tb{a}\wedge\tb{b}\)
as 
\begin{equation}
\ts{f}(\tb{A})=\tb{A}+\tb{A}\times\tb{P} - (\tb{A}\cdot\tb{P})\tb{P}.
\label{f on bivector A}
\end{equation} 
Indeed, substituting the definition of \(\ts{f}\) into (\ref{f on bivector}) and expressing
the relevant inner and outer products via the geometric product,
I get 
\[
\begin{split}
\ts{f}(\tb{a}\wedge\tb{b})&=\ts{f}(\tb{a})\wedge\ts{f}(\tb{b})
=
(\tb{a}+\tb{a}\cdot\tb{P})\wedge(\tb{b}+\tb{b}\cdot\tb{P})\\
&=
(\tb{a} +\tfrac{1}{2}\tb{a}\tb{P}  -\tfrac{1}{2}\tb{P}\tb{a})
\wedge
(\tb{b} +\tfrac{1}{2}\tb{b}\tb{P}  -\tfrac{1}{2}\tb{P}\tb{b})\\
&=
\tfrac{1}{2}(\tb{a} +\tfrac{1}{2}\tb{a}\tb{P}  -\tfrac{1}{2}\tb{P}\tb{a})
(\tb{b} +\tfrac{1}{2}\tb{b}\tb{P}  -\tfrac{1}{2}\tb{P}\tb{b})
-
\tfrac{1}{2}
(\tb{b} +\tfrac{1}{2}\tb{b}\tb{P}  -\tfrac{1}{2}\tb{P}\tb{b})
(\tb{a} +\tfrac{1}{2}\tb{a}\tb{P}  -\tfrac{1}{2}\tb{P}\tb{a})\\
&=
\tfrac{1}{2}(\tb{a}\tb{b} - \tb{b}\tb{a})+
\tfrac{1}{4}(\tb{a}\tb{b} - \tb{b}\tb{a})\tb{P}-\tfrac{1}{4}\tb{P}(\tb{a}\tb{b} - \tb{b}\tb{a})\\
&\quad -\tfrac{1}{8}(\tb{a}\tb{b} - \tb{b}\tb{a})\tb{P}^2
-\tfrac{1}{8}\tb{P}(\tb{a}\tb{b} - \tb{b}\tb{a})\tb{P}
+\tfrac{1}{4}(\tb{a}\tb{P}\tb{b}-\tb{b}\tb{P}\tb{a})\tb{P}\\
&=\tb{a}\wedge\tb{b} + (\tb{a}\wedge\tb{b})\times\tb{P}-((\tb{a}\wedge\tb{b})\cdot \tb{P})\tb{P},
\end{split}
\]
where I used \(\tb{P}(\tb{a}\wedge\tb{b})\tb{P}=2((\tb{a}\wedge\tb{b})\cdot\tb{P})\tb{P}-(\tb{a}\wedge\tb{b})\tb{P}^2\)
and  \(\tb{a}\tb{P}\tb{b}-\tb{b}\tb{P}\tb{a}=-2(\tb{a}\wedge\tb{b})\cdot\tb{P}\), which is a scalar.

The action of \(\ts{f}\) on a finite point \(\tb{A}\) and a point at infinity \(\tb{N}\) is illustrated in Figure \ref{linear function f}(b).
Equation~(\ref{f on bivector A}) can be rewritten as
\begin{equation}
\ts{f}(\tb{A})=\tb{A}(1-\tb{P}^2) + (\tb{A}\times\tb{P})\tb{P} + \tb{A}\times\tb{P},
\end{equation}
which becomes
\begin{equation}
\ts{f}(\tb{A})=2\left(\tb{A}+(\tfrac{1}{2}-1)(\tb{A}\times\tb{P})\tb{P}^{-1} + \tfrac{1}{2}\tb{A}\times\tb{P}\right)
\end{equation}
if the metric is Euclidean (\(\tb{P}^2=-1\) and \(\tb{P}^{-1}=-\tb{P}\)).
So, the action of \(\ts{f}\) on a point \(\tb{A}\) consists of 1) scaling \(\tb{A}\) with respect to \(\tb{P}\) by \(\tfrac{1}{2}\),
which gives a point in the middle between \(\tb{A}\) and \(\tb{P}\), and 2) 
shifting the result in the direction perpendicular to the line connecting \(\tb{A}\) and \(\tb{P}\) by half the distance 
between  \(\tb{A}\) and \(\tb{P}\), followed by 3) increasing the weight by 2.

Since the geometric action of \(\ts{f}\) is quite simple, it can be  inverted as follows.
For \(\tb{b}=\ts{f}(\tb{a})\), consider the line \(\tb{b}_1=\tb{b}\cdot\tb{P}\) and the line \(\tb{b}_2\) that results from
scaling  \(\tb{b}\) with respect to \(\tb{P}\) by a factor of 2,
i.e.\ \(\tb{b}_2=\tb{b}+(2-1)(\tb{b}\wedge\tb{P})\tb{P}^{-1}\).
These two lines intersect at a point which lies on \(\tb{a}\),
so \(\tb{a}\) can be found by adding \(\tb{b}_1\) and \(\tb{b}_2\) with some suitable weights.
Guided by this geometric insight, it is easy to find that
\begin{equation}
\ts{f}^{-1}(\tb{b})=\tfrac{1}{2}\left(\tb{b}+(\tb{b}\wedge\tb{P})\tb{P}^{-1} -\tb{b}\cdot\tb{P}\right)
\end{equation}
is indeed the inverse of  \(\ts{f}\).
Simple geometric reasoning also leads to the following expression for the inverse of \(\ts{f}\) on points:
\begin{equation}
\ts{f}^{-1}(\tb{B})=\tfrac{1}{2} (\tb{B}  -\tb{B}\times\tb{P}).
\end{equation}
Indeed, \(\tb{A}\) can be obtained from \(\tb{B}=\ts{f}(\tb{A})\) by shifting \(\tb{B}\)
in the direction perpendicular to the line connecting \(\tb{P}\) and \(\tb{B}\)
by the amount equal to the distance between \(\tb{P}\) and \(\tb{B}\).

A linear function \(\ts{f}\) can be defined on trivectors (pseudoscalars) as follows
\begin{equation}
\ts{f}(\tb{a}\wedge\tb{b}\wedge\tb{c})=\ts{f}(\tb{a})\wedge\ts{f}(\tb{b})\wedge\ts{f}(\tb{c}).
\end{equation}
For the function defined by (\ref{example definition}), I have
\(\ts{f}(\I)
=
\ts{f}(\e_0)\wedge\ts{f}(\e_1)\wedge\ts{f}(\e_2)
=\e_0\wedge(-\tfrac{3}{2}\e_0+\e_1+\e_2)\wedge(2\e_0-\e_1+\e_2)
=2\e_0\wedge\e_1\wedge\e_2
\), i.e.\ 
\(\ts{f}(\I)=2\I\).
In general, for any linear function \(\ts{f}\), there is a unique scalar \(\alpha\) such that \(\ts{f}(\I)=\alpha\I\).
The scalar \(\alpha\) is called the determinant of  \(\ts{f}\) and denoted by \(\det{(\ts{f})}\), so that
\begin{equation}
\ts{f}(\I)=\det{(\ts{f})}\I.
\end{equation}
For a scalar \(a\), I define \(\ts{f}(a)=a\).
This definition is motivated by the fact that \(a\wedge b=ab\) for \(a,b\in\R{}\), so that
\(\ts{f}(1)=\ts{f}(1\wedge 1)=\ts{f}(1)\wedge\ts{f}(1)=\ts{f}(1)^2\), which suggests \(\ts{f}(1)=1\).
Thus, \(\ts{f}\) defined originally on lines is extended to multivectors as follows:
\begin{equation}
\ts{f}(a+\tb{a}+\tb{A}+\alpha\I)
=a 
+ \tb{a} +\tb{a}\cdot\tb{P}
+\tb{A}+\tb{A}\times\tb{P}-(\tb{A}\cdot\tb{P})\tb{P}
+2\alpha\I.
\end{equation}
In Euclidean space, its inverse is given by
\begin{equation}
\ts{f}^{-1}(b+\tb{b}+\tb{B}+\beta\I)
=b
+\tfrac{1}{2}\left(
 \tb{b} +(\tb{b}\wedge\tb{P})\tb{P}^{-1} - \tb{b}\cdot\tb{P}
+\tb{B}-\tb{B}\times\tb{P}
+\beta\I
\right).
\label{example inverse}
\end{equation}
Furthermore, \(\ts{f}(A\wedge B)=\ts{f}(A)\wedge\ts{f}(B)\)
hold for any multivectors \(A,B\in\bigwedge\R{3*}\); the same applies to the inverse \(\ts{f}^{-1}\) and \(\det{(\ts{f}^{-1})}=(\det{\ts{f}})^{-1}\).
Note that the above expression for \(\ts{f}^{-1}\) is not applicable in spaces other than Euclidean.
In fact, the function \(\ts{f}\) defined by (\ref{example definition}) is not invertible in Minkowski space
since \(\det{\ts{f}}=0\) if the metric is Minkowski.

For any line \(\tb{a}_p\) that passes through the point \(\tb{P}\), I have
\(\tb{a}_p\wedge\tb{P}=0\) and, therefore,  
\(\ts{f}(\tb{a}_p)=\tb{a}_p+\tb{a}_p\tb{P}=\tb{a}_p(1+\tb{P})\)
and \(\ts{f}(\tb{a}_p)=(1-\tb{P})\tb{a}_p\).
Vectors dually representing lines that pass through \(\tb{P}\) form a 2-dimensional linear subspace of
\R{3*}. The action of \(\ts{f}\) on this subspace 
is characterised by two eigenvalues, the left-eigenvalue \(\lambda_1=1-\tb{P}\)
and the right-eigenvalue \(\lambda_2=1+\tb{P}\), so that
\(\ts{f}(\tb{a}_p)=\lambda_1\tb{a}_p\) and \(\ts{f}(\tb{a}_p)=\tb{a}_p\lambda_2\).
The action of either eigenvalue on \(\tb{a}_p\)
consists of a counterclockwise rotation around \(\tb{P}\) by \(45^\circ\)
followed by increasing the weight by \(\sqrt{2}\), e.g.\
\(\ts{f}(\tb{a}_p)=\lambda_1\tb{a}_p
=
\sqrt{2}
\left(
\tfrac{1}{\sqrt{2}}-\tfrac{1}{\sqrt{2}}\tb{P}
\right)
\tb{a}_p
\)
and
\(\ts{f}(\tb{a}_p)=\tb{a}_p\lambda_2
=
\sqrt{2}
\tb{a}_p
\left(
\tfrac{1}{\sqrt{2}}+\tfrac{1}{\sqrt{2}}\tb{P}
\right)
\).
These two distinct representations of the same action can be combined in 
\(\ts{f}(\tb{a}_p)
=
\sqrt{2}
S
\tb{a}_p
S^{-1}
\),
where \(S=e^{-\tfrac{1}{2}\tfrac{\pi}{4}\tb{P}}\).
For the line at infinity \(d\e_0\) weighted by \(d\), I have \(\ts{f}(d\e_0)=d\e_0\).
So, a 1-dimensional subspace of \R{3*} consisting of vectors \(d\e_0\), where \(d\in\R{}\), 
is characterised by the eigenvalue \(\lambda_3=1\).
Since any line \(\tb{a}\) can be decomposed into a line that passes through \(\tb{P}\) and the weighted line at infinity
by means of projection and rejection,
the function \(\ts{f}\) defined by (\ref{example definition}) can be written as
\begin{equation}
\ts{f}(\tb{a})=
\sqrt{2}
e^{-\tfrac{1}{2}\tfrac{\pi}{4}\tb{P}}
(\tb{a}\cdot\tb{P})\tb{P}^{-1}
e^{\tfrac{1}{2}\tfrac{\pi}{4}\tb{P}}
+(\tb{a}\wedge\tb{P})\tb{P}^{-1}.
\end{equation}
Moreover, \(\det{(\ts{f})}=\lambda_1\lambda_2\lambda_3=2\) and \(\tr{(\ts{f})}=\lambda_1+\lambda_2+\lambda_3=3\),
where \(\tr{(\ts{f})}\) is the trace of  \(\ts{f}\).

These properties can be extended to points as follows.
For two lines \(\tb{a}_p\) and \(\tb{b}_p\),  which both pass through the point \(\tb{P}\),
\(\ts{f}(\tb{a}_p\wedge\tb{b}_p)=
\ts{f}(\tb{a}_p)\wedge\ts{f}(\tb{b}_p)
=
\sqrt{2}\sqrt{2}\tb{a}_p\wedge\tb{b}_p
=2\tb{a}_p\wedge\tb{b}_p\).
Indeed,
rotating the lines \(\tb{a}_p\) and \(\tb{b}_p\) around  \(\tb{P}\) by the same angle does not change
 \(\tb{a}_p\wedge\tb{b}_p\); only the weight of the bivector increases since the weight of each line increases.
In other words, \(\ts{f}(\tb{P})=2\tb{P}\).
On the other hand, 
\(\ts{f}(\e_0\wedge\tb{a}_p)=
\ts{f}(\e_0)\wedge\ts{f}(\tb{a}_p)
=\e_0\wedge(\tb{a}_p+\tb{a}_p\tb{P})
=\e_0\wedge[\tb{a}_p(1+\tb{P})]=
\e_0[\tb{a}_p(1+\tb{P})]
=(\e_0\tb{a}_p)(1+\tb{P})
=(\e_0\wedge\tb{a}_p)(1+\tb{P})\)
and \(\ts{f}(\e_0\wedge\tb{a}_p)=(1-\tb{P})(\e_0\wedge\tb{a}_p)\).
Since any point at infinity \(\tb{N}\) can be written as \(\e_0\wedge\tb{a}_p\)
for some \(\tb{a}_p\), I have  \(\ts{f}(\tb{N})=\tb{N}(1+\tb{P})\) and
\(\ts{f}(\tb{N})=(1-\tb{P})\tb{N}\).
So, 
\(\bigwedge^2\R{3*}\) splits into two subspaces:
a 1-dimensional subspace consisting of \(\alpha\tb{P}\), where \(\alpha\in\R{}\),
with the eigenvalue \(\mu_1=2\) and a 2-dimensional subspace 
consisting of points at infinity with the left-eigenvalue \(\mu_2=1-\tb{P}\) and the right-eigenvalue \(\mu_3=1+\tb{P}\).
Therefore,
\begin{equation}
\ts{f}(\tb{A})=
2(\tb{A}\cdot\tb{P})\tb{P}^{-1}
+
\sqrt{2}
e^{-\tfrac{1}{2}\tfrac{\pi}{4}\tb{P}}
(\tb{A}\times\tb{P})\tb{P}^{-1}
e^{\tfrac{1}{2}\tfrac{\pi}{4}\tb{P}}.
\end{equation}
The above analysis is explicitly geometric in nature and independent of coordinate, 
in contrast with the standard analysis based on matrices and coordinates, which is summarised below for comparison.

In the standard basis, I have 
\[\ts{f}(d\e_0+a\e_1+b\e_2)=
d\e_0+a(-\tfrac{3}{2}\e_0+\e_1+\e_2)+b(2\e_0-\e_1+\e_2)=
(d-\tfrac{3}{2}a+2b)\e_0+(a-b)\e_1+(a+b)\e_2,
\]
and, therefore, the action of the function \(\ts{f}\) on lines can be represented by a matrix multiplication as follows
\begin{equation}
\begin{pmatrix}
d\\
a\\
b\\
\end{pmatrix}
\rightarrow
\begin{pmatrix}
1&-\tfrac{3}{2}&2\\
0&1&-1\\
0&1&1\\
\end{pmatrix}
\begin{pmatrix}
d\\
a\\
b\\
\end{pmatrix}.
\end{equation}
The determinant of the matrix equals 2, its eigenvalues consist of one real number \(1\) and 
two complex numbers \(1\pm i\), where \(i\) is the imaginary unit.
There is one real eigenvector corresponding to the subspace of the line at infinity with different weights
and two complex eigenvectors corresponding to the subspace of finite lines passing through the point \(\tb{P}\).

Vectors of the dual model space \(\R{3*}\) can be interpreted as linear functionals that act on vectors of the target model space \(\R{3}\)
and yield real numbers, i.e. \(\tb{a}:\R{3}\to\R{}\) for any \(\tb{a}\in\R{3*}\).
In the standard basis, the action of the functional \(\tb{a}=d\e_0+a\e_1+b\e_2\) is defined by
\(\tb{a}[\tb{x}]=dw+ax+by\) for \(\tb{x}=w\e^0+x\e^1+y\e^2\).
This definition is not metric, but in Elliptic space, I have \(\tb{a}[\tb{x}]=\tb{a}\cdot\tb{x}^I\), where \(\tb{x}^I=\Id^{-1}(\tb{x})\) and
\(\Id^{-1}\) is the inverse of the identity transformation.
For a linear function \(\ts{f}\), the adjoint \(\bar{\ts{f}}:\R{3}\to\R{3}\) is defined implicitly by
\begin{equation}
\ts{f}(\tb{a})[\tb{x}]=\tb{a}[\bar{\ts{f}}(\tb{x})].
\label{adjoint}
\end{equation}
The action of  \(\ts{f}(\tb{a})\) on  \(\tb{x}\) is equivalent to the action
of \(\tb{a}\) on  \(\bar{\ts{f}}(\tb{x})\).
In other words, instead of acting on \(\tb{x}\) with \(\ts{f}(\tb{a})\), it is possible to get the same result by acting on \(\bar{\ts{f}}(\tb{x})\) with \(\tb{a}\).
The adjoint is a linear function too, but it is defined on vectors of the target model space \R{3},
whereas \(\ts{f}\) is defined on vectors of \R{3*}.
A linear function \(\ts{f}\) acts on the dual representation of lines and
the adjoint \(\bar{\ts{f}}\) acts on the direct representation of points.
The adjoint is defined on vectors, but it can be extended to multivectors of \(\bigvee\R{3}\)
by using the outer product \(\vee\) 
in the same fashion  \(\ts{f}\) is extended to multivectors of \(\bigwedge\R{3*}\). 
For instance, 
\begin{equation}
\begin{split}
&\bar{\ts{f}}(s)=s,\\
&\bar{\ts{f}}(\tb{x}\vee\tb{y})=\bar{\ts{f}}(\tb{x})\vee\bar{\ts{f}}(\tb{y}),\\
&\bar{\ts{f}}(\tb{x}\vee\tb{y}\vee\tb{z})=\bar{\ts{f}}(\tb{x})\vee\bar{\ts{f}}(\tb{y})\vee\bar{\ts{f}}(\tb{z}),
\end{split}
\end{equation}
where \(s\) is a scalar and \(\tb{x},\tb{y},\tb{z}\) are vectors in \(\R{3}\).
The determinant of the adjoint, defined by \(\bar{\ts{f}}(\e^{012})=\det{(\bar{\ts{f}})}\e^{012}\),
equals that of the function \(\ts{f}\), i.e.\ \(\det{\ts{f}}=\det{\bar{\ts{f}}}\).
If \(\ts{f}\) has the inverse \(\ts{f}^{-1}\), the following equalities hold:
\begin{equation}
\begin{split}
&\ts{f}^{-1}(M) = \J^{-1}\bar{\ts{f}}(\J M)/\det{\bar{\ts{f}}} \quad\textrm{for }{\textstyle M\in\bigwedge\R{3*}},\\
&\bar{\ts{f}}(W) = \J\ts{f}^{-1}(\J^{-1}W)/\det{\ts{f}^{-1}} \quad\textrm{for }{\textstyle W\in\bigvee\R{3}},
\label{adjoint and inverse}
\end{split}
\end{equation}
where \(\J\) is the duality transformation (I use a simplified notation, e.g.\ \(JM=J(M)\), for readability).
Moreover, the adjoint \(\bar{\ts{f}}\) is also invertible and \(\bar{\ts{f}}^{-1}=\overline{\ts{f}^{-1}}\).
So, the action of the adjoint \(\bar{\ts{f}}\) is related to the action of the inverse  \(\ts{f}^{-1}\) if it exists.
However, the adjoint \(\bar{\ts{f}}\) is more general than the inverse,
since it is defined uniquely even if \(\ts{f}\) is not invertible.

A linear function \(\ts{f}\) can be expressed in terms of the inverse of its adjoint, if the inverse exists, as follows:
\begin{equation}
\ts{f}(M)=\J^{-1}\bar{\ts{f}}^{-1}(\J M)/\det{\bar{\ts{f}}^{-1}}.
\end{equation}
Since \(\bar{\ts{f}}^{-1}(W_1\vee W_2)=\bar{\ts{f}}^{-1}(W_1)\vee\bar{\ts{f}}^{-1}(W_2)\)
for any multivectors \(W_1,W_2\in\bigvee\R{3}\), it follows  that
\[
\begin{aligned}
\ts{f}(A\vee B)
&=
\J^{-1}\bar{\ts{f}}^{-1}(\J(A\vee B))/\det{\bar{\ts{f}}^{-1}}=
\J^{-1}\bar{\ts{f}}^{-1}(\J(A)\vee \J(B))/\det{\bar{\ts{f}}^{-1}}\\
&=
\J^{-1}[\bar{\ts{f}}^{-1}(\J A )\vee \bar{\ts{f}}^{-1}(\J B )]/\det{\bar{\ts{f}}^{-1}}=
\J^{-1}[\J\J^{-1}\bar{\ts{f}}^{-1}(\J A)\vee \J\J^{-1}\bar{\ts{f}}^{-1}(\J B)]/\det{\bar{\ts{f}}^{-1}}\\
&=
\J^{-1}\bar{\ts{f}}^{-1}(\J A )\vee \J^{-1}\bar{\ts{f}}^{-1}(\J B)/\det{\bar{\ts{f}}^{-1}}
=
\ts{f}(A)\vee\ts{f}(B)\det{\bar{\ts{f}}^{-1}}
=
\ts{f}(A)\vee\ts{f}(B)/\det{\ts{f}}
\end{aligned}
\]
and, therefore, 
\begin{equation}
\ts{f}(A\vee B)=\frac{\ts{f}(A)\vee\ts{f}(B)}{\det{\ts{f}}},
\end{equation}
which implies
\(
\ts{f}(A\vee B\vee C)=\ts{f}(A)\vee\ts{f}(B)\vee\ts{f}(C)/(\det{\ts{f}})^2\)
for the join of three multivectors.
I substitute \(\tb{P}\vee \tb{Q}\vee \tb{R}\) and use the fact that the join
of three points \(\tb{P}, \tb{Q}, \tb{R}\) is a scalar, so that
\(\ts{f}(\tb{P}\vee \tb{Q}\vee \tb{R})=\tb{P}\vee \tb{Q}\vee \tb{R}\),
to obtain
\begin{equation}
(\det{\ts{f}})^2
=\frac{\ts{f}(\tb{P})\vee\ts{f}(\tb{Q})\vee\ts{f}(\tb{R})}{\tb{P}\vee \tb{Q}\vee \tb{R}}.
\end{equation}
In Euclidean space \E{2}, \(\tfrac{1}{2}|\tb{P}\vee \tb{Q}\vee \tb{R}|\)
gives the area of a triangle defined by the normalised points  \(\tb{P}, \tb{Q}, \tb{R}\),
so
\begin{equation}
(\det{\ts{f}})^2
=\frac{\norm{\ts{f}(\tb{P})}}{\norm{\tb{P}}}
\frac{\norm{\ts{f}(\tb{Q})}}{\norm{\tb{Q}}}
\frac{\norm{\ts{f}(\tb{R})}}{\norm{\tb{R}}}
\frac{AREA[\ts{f}(\tb{P}),\ts{f}(\tb{Q}),\ts{f}(\tb{R})]}
{AREA[\tb{P},\tb{Q},\tb{R}]}
\end{equation}
and, therefore, \((\det{\ts{f}})^2\) measures the amount by which the area of a triangle changes under the action of
the linear function \(\ts{f}\), multiplied by the change in the weight of the points that define the corners of the triangle.

Substituting \(\tb{a}=\e_i\) and \(\tb{x}=\e^j\) into Equation~(\ref{adjoint}), I get
\(
\ts{f}(\e_i)[\e^j]=
(\ts{f}_{0i}\e_0+\ts{f}_{1i}\e_1+\ts{f}_{2i}\e_2)[\e^j]
=\ts{f}_{ji}
\)
and
\(
\e_i[\bar{\ts{f}}(\e^j)]=
\e_i[\bar{\ts{f}}_{0j}\e^0+\bar{\ts{f}}_{1j}\e^1+\bar{\ts{f}}_{2j}\e^2]
=\bar{\ts{f}}_{ij},
\)
so that
\(\bar{\ts{f}}_{ij}=\ts{f}_{ji}\),
where the first and second indices refer to the rows and columns of the matrix, respectively.
In other words, in the standard basis, the matrix representation of 
the adjoint \(\bar{\ts{f}}\) equals the transpose of the matrix representation of the 
linear function \(\ts{f}\).
For \(\ts{f}\) defined by Equation~(\ref{example definition}), I have
\[
\begin{split}
&\bar{\ts{f}}(\e^0)=\bar{\ts{f}}_{00}\e^0+\bar{\ts{f}}_{10}\e^1+\bar{\ts{f}}_{20}\e^2
=\ts{f}_{00}\e^0+\ts{f}_{01}\e^1+\ts{f}_{02}\e^2
=\e^0-\tfrac{3}{2}\e^1+2\e^2,\\
&\bar{\ts{f}}(\e^1)=\bar{\ts{f}}_{01}\e^0+\bar{\ts{f}}_{11}\e^1+\bar{\ts{f}}_{21}\e^2
=\ts{f}_{10}\e^0+\ts{f}_{11}\e^1+\ts{f}_{12}\e^2
=\e^1-\e^2,\\
&\bar{\ts{f}}(\e^2)=\bar{\ts{f}}_{02}\e^0+\bar{\ts{f}}_{12}\e^1+\bar{\ts{f}}_{22}\e^2
=\ts{f}_{20}\e^0+\ts{f}_{21}\e^1+\ts{f}_{22}\e^2
=\e^1+\e^2,
\end{split}
\]
and, therefore,
\[
\bar{\ts{f}}(w\e^0+x\e^1+y\e^2)
=w(\e^0-\tfrac{3}{2}\e^1+2\e^2)+x(\e^1-\e^2)+y(\e^1+\e^2)=
w\e^0+(-\tfrac{3}{2}w+x+y)\e^1+(2w-x+y)\e^2.
\]
So, the action of \(\bar{\ts{f}}\) in the standard basis is indeed given by
\begin{equation}
\begin{pmatrix}
w\\
x\\
y\\
\end{pmatrix}
\rightarrow
\begin{pmatrix}
1&0&0\\
-\tfrac{3}{2}&1&1\\
2&-1&1\\
\end{pmatrix}
\begin{pmatrix}
w\\
x\\
y\\
\end{pmatrix}.
\end{equation}
To illustrate how Equations~(\ref{adjoint and inverse}) work, I compute \(\ts{f}^{-1}(\e_0)\) in terms of the adjoint.
I have
 \[
\begin{aligned}
\bar{\ts{f}}(\e^{12})
&=\bar{\ts{f}}(\e^{1}) \vee\bar{\ts{f}}(\e^{2})=
(\bar{\ts{f}}_{01}\e^0+\bar{\ts{f}}_{11}\e^1+\bar{\ts{f}}_{21}\e^2)
\vee(\bar{\ts{f}}_{02}\e^0+\bar{\ts{f}}_{12}\e^1+\bar{\ts{f}}_{22}\e^2)\\
&=
(\bar{\ts{f}}_{11}\bar{\ts{f}}_{22} -  \bar{\ts{f}}_{21}\bar{\ts{f}}_{12})\e^{12}
+(\bar{\ts{f}}_{21}\bar{\ts{f}}_{02} - \bar{\ts{f}}_{01}\bar{\ts{f}}_{22} )\e^{20}
+(\bar{\ts{f}}_{01}\bar{\ts{f}}_{12} - \bar{\ts{f}}_{11} \bar{\ts{f}}_{02} )\e^{01}
\end{aligned}
\] 
and, therefore,
\[
\begin{aligned}
\J^{-1}\bar{\ts{f}}(\J\e_0)
&=
(\bar{\ts{f}}_{11}\bar{\ts{f}}_{22} -  \bar{\ts{f}}_{21}\bar{\ts{f}}_{12})\e_0
+(\bar{\ts{f}}_{21}\bar{\ts{f}}_{02} - \bar{\ts{f}}_{01}\bar{\ts{f}}_{22} )\e_1
+(\bar{\ts{f}}_{01}\bar{\ts{f}}_{12} - \bar{\ts{f}}_{11} \bar{\ts{f}}_{02} )\e_2\\
&=
(\ts{f}_{11}\ts{f}_{22} -  \ts{f}_{12}\ts{f}_{21})\e_0
+(\ts{f}_{12}\ts{f}_{20} - \ts{f}_{10}\ts{f}_{22} )\e_1
+(\ts{f}_{10}\ts{f}_{21} - \ts{f}_{11}\ts{f}_{20} )\e_2,
\end{aligned}
\] 
which gives
\[
\ts{f}^{-1}(\e_0)=
[(\ts{f}_{11}\ts{f}_{22} -  \ts{f}_{12}\ts{f}_{21})\e_0
+(\ts{f}_{12}\ts{f}_{20} - \ts{f}_{10}\ts{f}_{22} )\e_1
+(\ts{f}_{10}\ts{f}_{21} - \ts{f}_{11}\ts{f}_{20} )\e_2]/\det{\ts{f}}.
\]
This expression determines the first column of the matrix of \(\ts{f}^{-1}\).
Similar expressions can be derived for \(\ts{f}^{-1}(\e_1)\) and \(\ts{f}^{-1}(\e_2)\), which gives
\begin{equation}
\begin{pmatrix}
d\\
a\\
b\\
\end{pmatrix}
\rightarrow
\tfrac{1}{2}
\begin{pmatrix}
2&\tfrac{7}{2}&-\tfrac{1}{2}\\
0&1&1\\
0&-1&\,1\
\end{pmatrix}
\begin{pmatrix}
d\\
a\\
b\\
\end{pmatrix}.
\end{equation}
for the matrix of \(\ts{f}^{-1}\) in the standard basis.
The matrix  depends on the basis used to express the coordinates, whereas
Equation~(\ref{example inverse}) is independent of coordinates.

Consider a linear function that satisfies
\begin{equation}
\ts{f}(M\I)=\ts{f}(M)\I
\end{equation}
for any multivector \(M\).
Substituting \(M=\tb{a}\) gives \(\ts{f}(\tb{a}\I)=\ts{f}(\tb{a})\I \), 
so the action of  \(\ts{f}\) on the polar point of \(\tb{a}\) is consistent with its action on the line itself.
Substituting \(M=1\) gives \(\ts{f}(\I)=\ts{f}(1)\I=\I\) and, therefore, \(\det{\ts{f}}=1\) and \(\ts{f}\) is invertible.
Since \(\tb{a}\cdot\tb{b}\) is a scalar and \(\tb{a}\cdot\tb{b}=(\tb{a}\I)\vee\tb{b}\) , I have
\[
\tb{a}\cdot\tb{b}=\ts{f}(\tb{a}\cdot\tb{b})=
\ts{f}((\tb{a}\I)\vee\tb{b})=\ts{f}(\tb{a}\I)\vee\ts{f}(\tb{b})
=(\ts{f}(\tb{a})\I)\vee\ts{f}(\tb{b})=
\ts{f}(\tb{a})\cdot\ts{f}(\tb{b})
\]
and, therefore, 
\begin{equation}
\label{preserving angles between lines}
\tb{a}\cdot\tb{b}=\ts{f}(\tb{a})\cdot\ts{f}(\tb{b}).
\end{equation}
For \(\tb{a}=\tb{b}\), Equation~(\ref{preserving angles between lines}) gives \(\norm{\tb{a}}=\norm{\ts{f}(\tb{a})}\)
and, therefore,
the function \(\ts{f}\) preserves the angles between lines.
Since \(\tb{P}\vee\tb{Q}\) is a line, I have \(\norm{\tb{P}\vee\tb{Q}}=\norm{\ts{f}(\tb{P}\vee\tb{Q})}\) and, therefore, 
\begin{equation}
\norm{\tb{P}\vee\tb{Q}}=\norm{\ts{f}(\tb{P})\vee\ts{f}(\tb{Q})}.
\end{equation}
\(\tb{P}\cdot\tb{Q}\) is a scalar and \(\tb{P}\cdot\tb{Q}=(\tb{P}\I)\vee\tb{Q}\), so 
\(\tb{P}\cdot\tb{Q}=\ts{f}(\tb{P})\cdot\ts{f}(\tb{Q})\), which for \(\tb{P}=\tb{Q}\) gives \(\norm{\tb{P}}=\norm{\ts{f}(\tb{P})}\).
Therefore, \(\ts{f}\) preserves the distances between points as well.

Defining a linear function on lines in \E{2} and extending it to multivectors as was done above
is one way to obtain a linear function defined on multivectors.
Obviously, this procedure does not exhaust all possible linear functions defined on multivectors.
These functions satisfy
\begin{equation}
\begin{split}
&\ts{f}(A+B)=\ts{f}(A+\ts{f}(B),\\
&\ts{f}(\alpha A)=\alpha \ts{f}(A),
\end{split}
\end{equation}
for \(\alpha\in\R{}\) and \(A,B\in\bigwedge\R{3*}\).
For instance \(\ts{g}(\tb{a})=\e_0\wedge\tb{a}\) defines a linear function \(\ts{g}:\R{3*}\to\bigwedge^2\R{3*}\),\
acting on vectors and yielding bivectors.
Furthermore, one can consider linear functions acting across  spaces.
For instance \(\ts{h}(\tb{P})=\e^1\vee\J(\tb{P})\) defines \(\ts{h}:\bigwedge^2\R{3*}\to\bigvee^2\R{3}\),
acting on bivectors in the dual model space and yielding bivectors in the target model space.
In the standard basis, linear functions such as \(\ts{g}:\bigwedge\R{3*}\to\bigwedge\R{3*}\) and
\(\ts{h}:\bigwedge\R{3*}\to\bigvee\R{3}\) can be represented by \(8\times8\) matrices.

\newpage
\section{3-dimensional geometry}

\subsection{Points and planes}

In the target space \T{3}, a plane that does not pass through the origin can be defined by 
\begin{equation}\label{plane in 3D}
ax+by+cz+1=0,
\end{equation}
where \(a\), \(b\), and \(c\) are fixed and \(x\), \(y\), and \(z\) are variable; \((x,y,z)\) is in \T{3}.
To define a plane, one needs to specify three numbers \((a,b,c)\).
The dual space \T{3*} is a linear space of the coefficients \((a,b,c)\), where every point defines a plane in the target space \T{3}
via Equation (\ref{plane in 3D}). 

\begin{figure}[h]
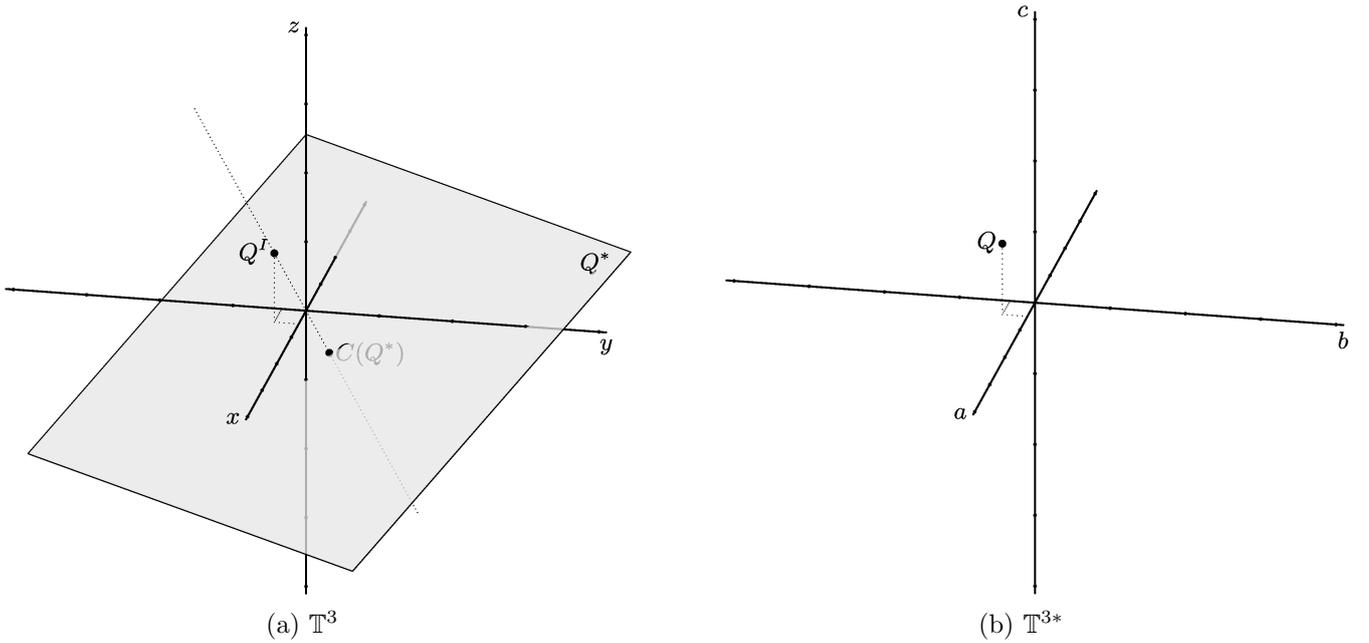
\hspace{-1cm}
\begin{subfloatenv}{\T{3}}
\begin{asy}
import Figure3D;
Figure3D f = Figure3D();

var a = Plane(1,1/2,-1/3,1);
f.plane(a, "$Q^*$", align=(1,-3,-1), draw_orientation=false, draw_central=true);
label("$C(Q^*)$", totriple(a.centre()), align=(0,1,0));

var Q = Point(1,1/2,-1/3,1);
f.point(Q,"$Q^I$", align=(0,-1,0),draw_orientation=false);

var L = join(Q, a.centre());
f.line(L, dotted,draw_orientation=false);

\end{asy}
\end{subfloatenv}\hfill%
\begin{subfloatenv}{\T{3*}}
\begin{asy}
import Figure3D;
Figure3D f = Figure3D(space=DUAL);

var Q = Point(1,1/2,-1/3,1);
f.point(Q,"$Q$", align=(0,-1,0),draw_orientation=false);

\end{asy}
\end{subfloatenv}
\caption{Mutual correspondence between points and planes}
\label{points and planes}
\end{figure}

An example is shown in Figure~\ref{points and planes}, where a point \(Q\) in \T{3*} 
with the coordinates \((a,b,c)=(\tfrac{1}{2},-\tfrac{1}{3},1)\)
is dual to the plane \(Q^*\) in \T{3} defined by \(\tfrac{1}{2}x-\tfrac{1}{3}y+z+1=0\).
The duality and identity transformations are defined in the same way as in the 2-dimensional case,
except that in \T{3} and \T{3*}, points correspond to planes and vice versa.
For instance, \(Q^I=I(Q)=(x,y,z)=(\tfrac{1}{2},-\tfrac{1}{3},1)\) is
a point in \T{3} identical to \(Q\), and \(Q^*=J(Q)\) is a plane in \T{3} dual to the point \(Q\) in \T{3*}.
The central point \(C(Q^*)\) of the plane \(Q^*\) is defined as a point where 
the line in \T{3} that passes through the point \(Q^I\) and the origin of \T{3} 
intersects \(Q^*\).
It is given by 
\begin{equation}
C(Q^*) = - \frac{(a,b,c)}{a^2+b^2+c^2}.
\end{equation}
In Euclidean space, 
\begin{equation}
d_{Q^I}d_{C(Q^*)}=1,
\end{equation}
where \(d_{Q^I}\) is the distance from the origin to point \(Q^I\) and \(d_{C(Q^*)}\) is the distance from the origin to 
the central point of the plane \(Q^*\).

\begin{figure}[t!]
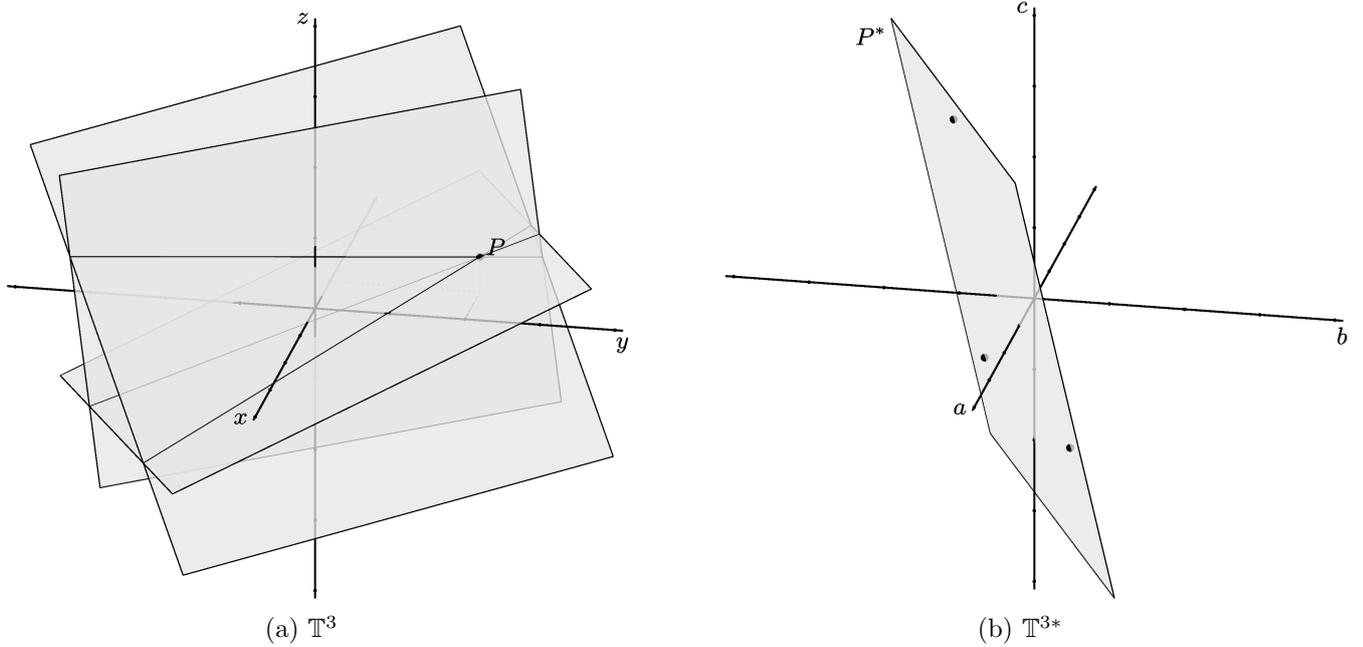

\begin{subfloatenv}{\T{3}}
\begin{asy}
import Figure3D;
Figure3D f = Figure3D();

triple n = (-1,2,1/2);
var a = Plane(1,n.x,n.y,n.z);
triple C = totriple(a.centre());

triple t1 = C+3*cross((1,0,0),n/length(n));
triple[] ts = {};
int tsn = 3;
for(int i=0; i<tsn; ++i) { triple t = rotate(360*i/tsn, C, C+n)*t1; ts.push(t); }

for(triple t: ts) { var p = Plane(1,t.x,t.y,t.z); f.plane(p, u=n, pen=lightgrey, draw_orientation=false); }

var P = Point(1,n.x,n.y,n.z);
f.point(P, "$P$", align=(2,2,2),draw_orientation=false);

\end{asy}
\end{subfloatenv}
\begin{subfloatenv}{\T{3*}}
\begin{asy}
import Figure3D;
Figure3D f = Figure3D(space=DUAL);

triple n = (-1,2,1/2);
var a = Plane(1,n.x,n.y,n.z);
triple C = totriple(a.centre());

f.plane(a, "$P^*$", align=(1,-1,-1), draw_orientation=false, draw_central=false);

triple t1 = C+3*cross((1,0,0),n/length(n));
triple[] ts = {};
int tsn = 3;
for(int i=0; i<tsn; ++i) { triple t = rotate(360*i/tsn, C, C+n)*t1; ts.push(t); }

for(triple t: ts) { dot(t); }

\end{asy}
\end{subfloatenv}
\caption{Top-down view of a point in \T{3}}
\label{top-down point in 3D}
\end{figure}

Similarly, a point \(P\) in \T{3} with the coordinates \((x,y,z)\)  
is dual to the plane \(P^*\) in \T{3*} defined by Equation~(\ref{plane in 3D}),
where \(x\), \(y\), and \(z\) are seen as fixed coefficients and \(a\), \(b\), \(c\) are variable.
The dual plane \(P^*\) enables the top-down view of the point \(P\).
Every point on the plane \(P^*\) corresponds to a plane in \T{3} passing through the point \(P\).
So, a point in the top-down view  is a set of all planes intersecting at the point, called a bundle of planes.
An example is shown in Figure~\ref{top-down point in 3D} where \(P=(x,y,z)=(-1,2,\tfrac{1}{2})\) and the dual plane \(P^*\) in \T{3*}
is defined by \(-a+2b+\tfrac{1}{2}c+1=0\).
Three planes from the bundle of planes around \(P\) are shown in Figure~\ref{top-down point in 3D}(a), 
the corresponding points \T{3*}, which lie on the plane \(P^*\), are shown in Figure~\ref{top-down point in 3D}(b).

\subsection{Lines}

Let \(L\) be a line in \T{3} defined by 
\begin{equation}
\left\{
\begin{aligned}
&p_{12}y-p_{31}z+p_{10}=0,\\
&p_{23}z-p_{12}x+p_{20}=0,\\
&p_{31}x-p_{23}y+p_{30}=0,
\end{aligned}
\right.
\label{equation for line in T3}
\end{equation}
where the parameters satisfy 
\begin{equation}
p_{10}p_{23}+p_{20}p_{31}+p_{30}p_{12}=0.
\label{line constraint T3}
\end{equation}
In this system, only two equations are linearly independent.
For instance, the third equation follows from the first two and the constraint~(\ref{line constraint T3}).
Given that there are three variables involved in (\ref{equation for line in T3}), the system define a line.
Furthermore, (\ref{equation for line in T3}) implies
\begin{equation}
p_{10}x+p_{20}y+p_{30}z=0
\end{equation}
for any point \((x,y,z)\) on the line.

For a line shown in Figure~\ref{lines in 3D}(a),
the parameters have the following values: \((p_{23},p_{31},p_{12})=(1,-1,-\tfrac{1}{3})\) and
\((p_{10},p_{20},p_{30})=(0,-\tfrac{1}{3},1)\).
The line passes through the points \((1,0,0)\) and \((0,1,\tfrac{1}{3})\).
Each point on the line \(L\) corresponds to a plane in \T{3*}.
Seven such planes are shown in Figure~\ref{lines in 3D}(b);
the corresponding points on the line \(L\) are shown in Figure~\ref{lines in 3D}(a).
The planes intersect along the line \(L^*\) defined by
\begin{equation}
\left\{
\begin{aligned}
&p_{30}b-p_{20}c+p_{23}=0,\\
&p_{10}c-p_{30}a+p_{31}=0,\\
&p_{20}a-p_{10}b+p_{12}=0,
\label{definition of line in dual 3D}
\end{aligned}
\right.
\end{equation}
where the parameters have the same values as above.
The equality
\begin{equation}
p_{23}a+p_{31}b+p_{12}c=0
\end{equation}
holds for any point \((a,b,c)\) on \(L^*\).
In Euclidean space, \(L\) is parallel to the vector \((p_{23},p_{31},p_{12})\) and
perpendicular to the vector \((p_{10},p_{20},p_{30})\).

\begin{figure}[h!]
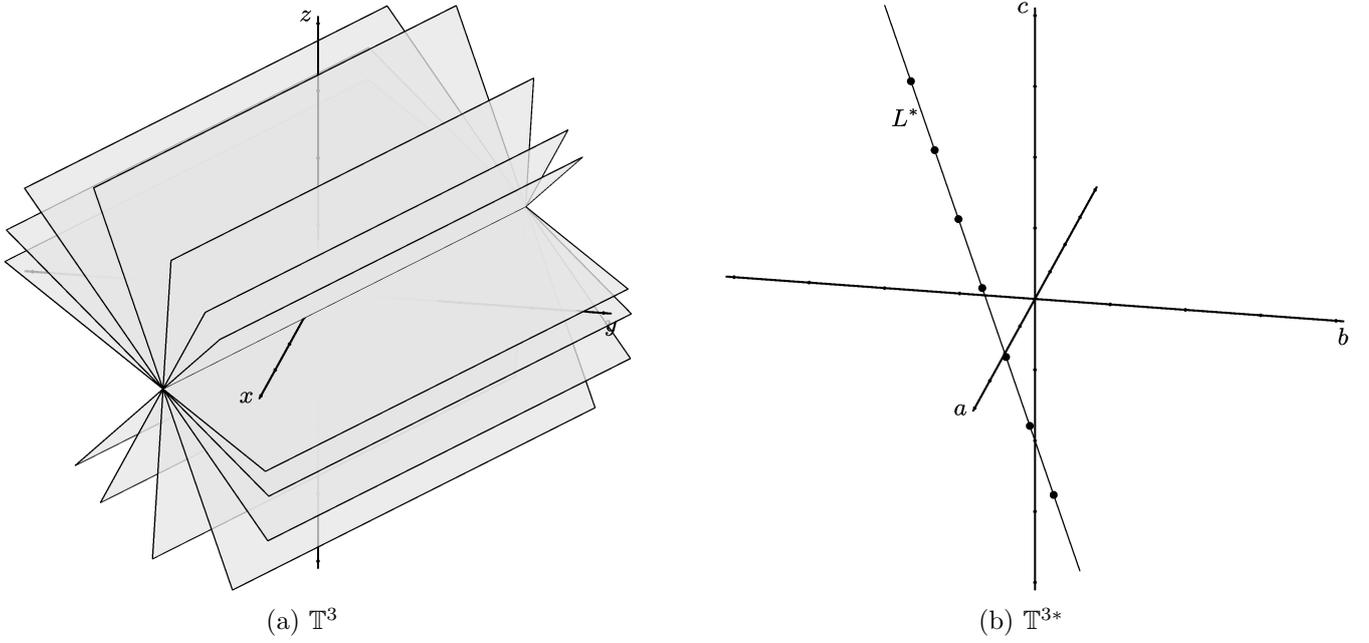
\hspace{-1cm}
\begin{subfloatenv}{\T{3}}
\begin{asy}
import Figure3D;
Figure3D f = Figure3D();

var P = Point(1,1,0,0);
var Q = Point(1,0,1,1/3);
var L = join(P,Q);

Line dualise(Line line) { unravel line; return Line(p23,p31,p12,p10,p20,p30); }
L = dualise(toline(L));

triple C = totriple(toline(L).centre());
triple n = (L.p23,L.p31,L.p12);
n = normalise(n);

triple[] ts = {};
int tsn = 3;
for(int i=-tsn; i<=tsn; ++i) { ts.push(C+i*n);  }

MV[] ps = {};
for(triple t: ts) { var p = Plane(1,t.x,t.y,t.z); ps.push(p); }

for(var p: ps) { f.plane(p,u=n, draw_orientation=false); }

\end{asy}
\end{subfloatenv}\hfill%
\begin{subfloatenv}{\T{3*}}
\begin{asy}
import Figure3D;
Figure3D f = Figure3D(space=DUAL);

var P = Point(1,1,0,0);
var Q = Point(1,0,1,1/3);
var L = join(P,Q);

Line dualise(Line line) { unravel line; return Line(p23,p31,p12,p10,p20,p30); }
L = dualise(toline(L));
f.line(L, "$L^*$", position=0.8, align=(0,-1,0),draw_orientation=false);

triple C = totriple(toline(L).centre());
triple n = (L.p23,L.p31,L.p12);
n = normalise(n);

triple[] ts = {};
int tsn = 3;
for(int i=-tsn; i<=tsn; ++i) { ts.push(C+i*n);  }

for(var t: ts) {dot(t);}

\end{asy}
\end{subfloatenv}
\caption{The top-down view of a line in \T{3}}
\label{top-down lines in 3D}
\end{figure}

In the bottom-up view of \T{3*}, the line \(L^*\) is a set of points, each of which corresponds to a plane in \T{3}
passing through the line \(L\) (see Figure~\ref{top-down lines in 3D}). 
The planes that pass through  \(L\) form a sheaf\footnote{It is also called a pencil of planes.} of planes around \(L\).
So, in the top-down view, a line is a derivative geometric object represented by a sheaf of planes.
The top-down view of the line \(L\) is illustrated in Figure~\ref{top-down lines in 3D}(a).
Its dual line \(L^*\) is shown in Figure~\ref{top-down lines in 3D}(b).

\begin{figure}[t]\hspace{-1cm}
\begin{subfloatenv}{\T{3}}
\begin{asy}
import Figure3D;
Figure3D f = Figure3D();

var P = Point(1,1,0,0);
var Q = Point(1,0,1,1/3);
var L = join(P,Q);
write(L);

//f.point(P,"$P$", align=(0,-2,0));
//f.point(Q,"$Q$", align=(0,0,1));
triple q = totriple(Q);
draw((q.x,q.y,0)--q, currentpen+dotted);
draw((0,0,q.z)--q, currentpen+dotted);

f.line(L, "$L$", align=(0,0,-1),draw_orientation=false);

triple C = totriple(toline(L).centre());
triple n = (L.p23,L.p31,L.p12);
n = normalise(n);

triple[] ts = {};
int tsn = 3;
for(int i=-tsn; i<=tsn; ++i) { ts.push(C+i*n);  }

for(var t: ts) {dot(t);}

\end{asy}
\end{subfloatenv}\hfill%
\begin{subfloatenv}{\T{3*}}
\begin{asy}
import Figure3D;
Figure3D f = Figure3D(space=DUAL);

var P = Point(1,1,0,0);
var Q = Point(1,0,1,1/3);
var L = join(P,Q);

Line dualise(Line line) { unravel line; return Line(p23,p31,p12,p10,p20,p30); }

triple C = totriple(toline(L).centre());
triple n = (L.p23,L.p31,L.p12);
n = normalise(n);

triple[] ts = {};
int tsn = 3;
for(int i=-tsn; i<=tsn; ++i) { ts.push(C+i*n);  }

MV[] ps = {};
for(triple t: ts) { var p = Plane(1,t.x,t.y,t.z); ps.push(p); }

for(var p: ps) { f.plane(p,u=n, draw_orientation=false); }

\end{asy}
\end{subfloatenv}
\caption{A line in \T{3} and its dual line in \T{3*}}
\label{lines in 3D}
\end{figure}

Let \(p_0(L)\) be the plane passing through the line \(L\) and the origin of \T{3}.
The intersection point of the dual line \(L^*\) and the plane \(p_0(L)^I\) in \T{3*}, which is identical to \(p_0(L)\),  
defines the central point of the line \(L^*\).
Its coordinates are given by
\begin{equation}
a_c = -\frac{p_{20}p_{12}-p_{30}p_{31}}{p_{10}^2+p_{20}^2+p_{30}^2},\quad
b_c = -\frac{p_{30}p_{23}-p_{10}p_{12}}{p_{10}^2+p_{20}^2+p_{30}^2},\quad
c_c = -\frac{p_{10}p_{31}-p_{20}p_{23}}{p_{10}^2+p_{20}^2+p_{30}^2}.
\end{equation}
For a line that passes through the origin, the central point is at the origin by definition.
To find the central point of the line \(L\), consider the plane \(q_0(L^*)\) passing through \(L^*\) and the origin of \T{3*}
and its identical plane \(q_0(L^*)^I\) in \T{3}.
The line \(L\) intersects \(q_0(L^*)^I\) at the central point of the line.
The coordinates of the central point are as follows:
\begin{equation}
x_c = \frac{p_{20}p_{12}-p_{30}p_{31}}{p_{23}^2+p_{31}^2+p_{12}^2},\quad
y_c = \frac{p_{30}p_{23}-p_{10}p_{12}}{p_{23}^2+p_{31}^2+p_{12}^2},\quad
z_c = \frac{p_{10}p_{31}-p_{20}p_{23}}{p_{23}^2+p_{31}^2+p_{12}^2}.
\end{equation}

In Euclidean space, 
\begin{equation}
d_{C(L)}d_{C(L^*)^I}=1,
\end{equation} 
where
\(d_{C(L)}\) is the distance from the origin of \T{3} to the central point \(C(L)\) of \(L\)
 and \(d_{C(L^*)^I}\) is the distance from the origin
to  \(C(L^*)^I\), which is the identical counterpart of the central point \(C(L^*)\).
Note that the central point is not affected by the identity transformation, i.e.\ 
if \(C(L)\) is the central point of \(L\), then \(C(L)^I\) is the central point of \(L^I\).

The sheaf of planes representing \(L\) intersects the plane \(q_0(L^*)^I\)
along a structure equivalent to a sheaf of lines residing inside  \(q_0(L^*)^I\).
In fact, the relationship between this sheaf of lines and the line \(L^*\),
which lies in the plane \(q_0(L^*)\), is identical to that observed in the 2-dimensional space.
A sheaf of planes is a simple extension of a sheaf of lines into one extra dimension.

\subsection{Points, lines, and the plane at infinity}
Consider any finite plane in \T{3} and move it to infinity in any direction.
This gives the plane at infinity in \T{3}.
The dual point in  \T{3*} converge to the origin of \T{3*},
so the origin of \T{3*} is dual to the plane at infinity in \T{3}.
Similarly, the origin in the target space \T{3} is dual to the plane at infinity in \T{3*}.

A line at infinity in \T{3} is represented by a stack of planes, which consists of all planes given by
\begin{equation}
ax+by+cz+d=0,
\end{equation}
where the coefficients \((a,b,c)\) are fixed and \(d\) spans all possible values in \R{}.
Taking any finite line in \T{3} and moving it to infinity results in a specific line at infinity.
The corresponding line in \T{3*} moves towards the origin.
In the limit, it arrives at the origin of \T{3*}.
So, any line at infinity in \T{3} corresponds to a line in \T{3*} that passes through the origin of \T{3*}.

\begin{figure}[t]
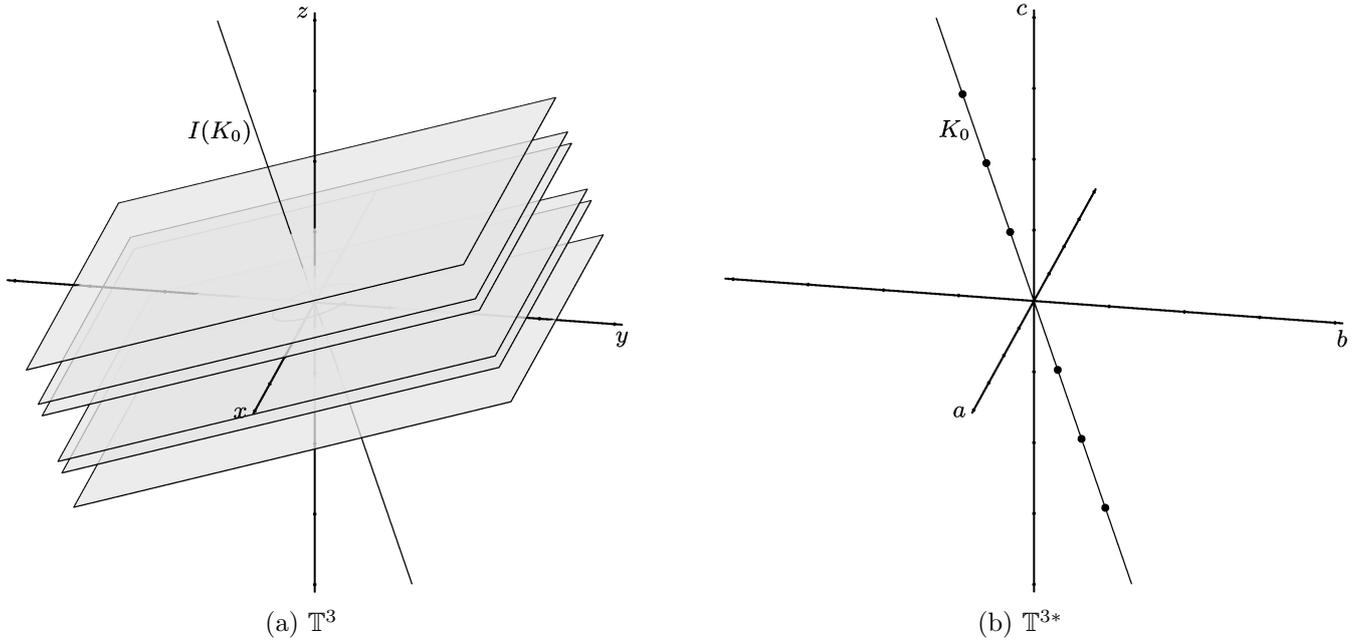
\hspace{-1cm}
\begin{subfloatenv}{\T{3}}
\begin{asy}
import Figure3D;
Figure3D f = Figure3D();

var P = Point(1,1,0,0);
var Q = Point(1,0,1,1/3);
var L = join(P,Q);

Line dualise_through_origin(Line line) { unravel line; return Line(0,0,0,p10,p20,p30); }
L = dualise_through_origin(toline(L));

f.line(L, "$I(K_0)$", position=0.8, align=(0,-1,0));

triple C = totriple(toline(L).centre());
triple n = (L.p23,L.p31,L.p12);
n = normalise(n);

triple[] ts = {};
int tsn = 3;
for(int i=1; i<=tsn; ++i) { ts.push(C+i*n);  ts.push(C-i*n);}

MV[] ps = {};
for(triple t: ts) { var p = Plane(1,t.x,t.y,t.z); ps.push(p); }

for(var p: ps) { f.plane(p, draw_orientation=false); }

\end{asy}
\end{subfloatenv}\hfill%
\begin{subfloatenv}{\T{3*}}
\begin{asy}
import Figure3D;
Figure3D f = Figure3D(space=DUAL);

var P = Point(1,1,0,0);
var Q = Point(1,0,1,1/3);
var L = join(P,Q);

Line dualise_through_origin(Line line) { unravel line; return Line(0,0,0,p10,p20,p30); }
L = dualise_through_origin(toline(L));
f.line(L, "$K_0$", position=0.8, align=(0,-1,0),draw_orientation=false);

triple C = totriple(toline(L).centre());
triple n = (L.p23,L.p31,L.p12);
n = normalise(n);

triple[] ts = {};
int tsn = 3;
for(int i=1; i<=tsn; ++i) { ts.push(C+i*n);  ts.push(C-i*n);}

for(var t: ts) {dot(t);}

\end{asy}
\end{subfloatenv}
\caption{A stack of planes in \T{3} and its dual line in \T{3*}}
\label{line at infinity in T3}
\end{figure}

For example, let  \(K_0\) be a line in \T{3*} defined by~(\ref{definition of line in dual 3D})
with \((p_{10},p_{20},p_{30})=(0,-\tfrac{1}{3},1)\) and \((p_{23},p_{31},p_{12})=(0,0,0)\).
This line is shown in Figure~\ref{line at infinity in T3}(b).
The points on \(K_0\) correspond to the planes in \T{3} which comprise a stack of planes, i.e.\ a line at infinity 
dual to \(K_0\) and denoted by \(K_0^*\).
Some planes from the stack are shown in Figure~\ref{line at infinity in T3}(a).
One can think of the line \(K_0^*\) as the intersection of the planes in the stack.
The concept of a stack of planes is a simple extension of the concept of a stack of lines in \T{2}.

Moving a finite point in \T{3} to infinity in a given direction results in a specific point at infinity.
It corresponds to a plane in \T{3*} passing through the origin of \T{3*}.
Moving the same point to infinity in the opposite direction yields the same plane in \T{3*}.
So, similar to the 2-dimensional case, the same point at infinity can be obtained by approaching infinity in two opposite directions.

\begin{figure}[h!]
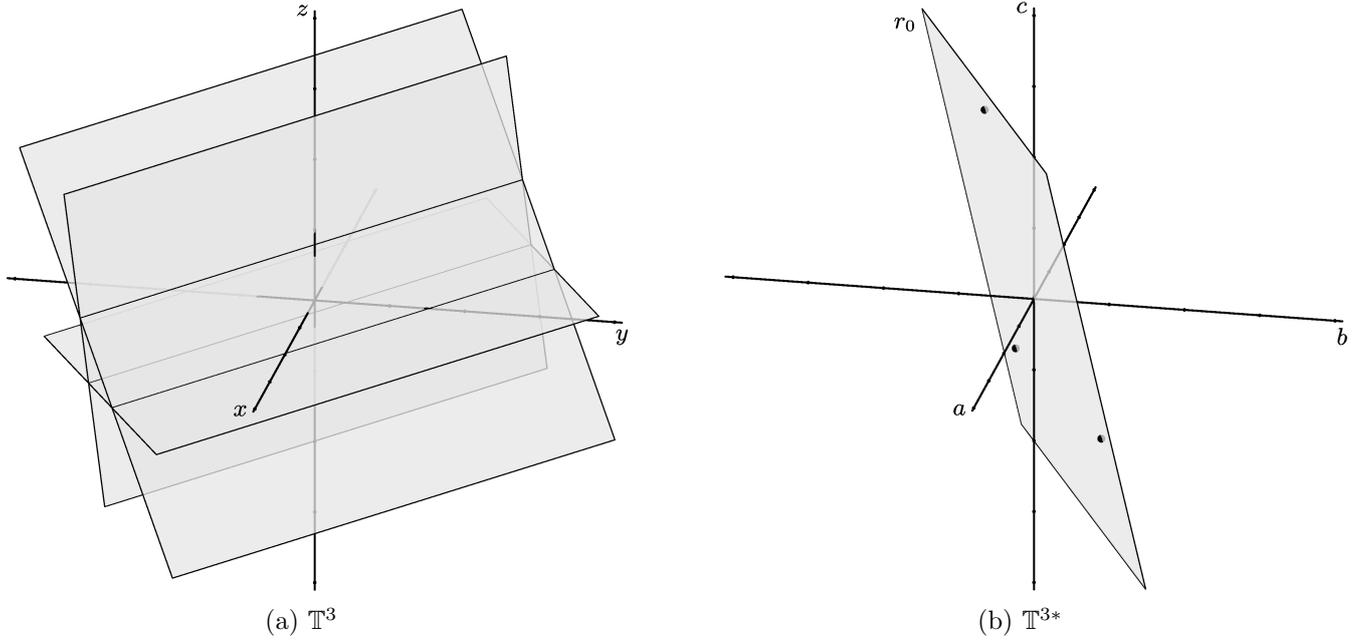
\hspace{-1cm}
\begin{subfloatenv}{\T{3}}
\begin{asy}
import Figure3D;
Figure3D f = Figure3D();

triple n = (-1,2,1/2);
var a = Plane(0,n.x,n.y,n.z);
triple C = totriple(a.centre());

triple t1 = C+3*cross((1,0,0),n/length(n));
triple[] ts = {};
int tsn = 3;
for(int i=0; i<tsn; ++i) { triple t = rotate(360*i/tsn, C, C+n)*t1; ts.push(t); }

for(triple t: ts) { var p = Plane(1,t.x,t.y,t.z); f.plane(p, u=n, pen=lightgrey, draw_orientation=false); }

\end{asy}
\end{subfloatenv}\hfill%
\begin{subfloatenv}{\T{3*}}
\begin{asy}
import Figure3D;
Figure3D f = Figure3D(space=DUAL);

triple n = (-1,2,1/2);
var a = Plane(0,n.x,n.y,n.z);
triple C = totriple(a.centre());

f.plane(a, "$r_0$", align=(1,-1,-1), draw_orientation=false, draw_central=false);

triple t1 = C+3*cross((1,0,0),n/length(n));
triple[] ts = {};
int tsn = 3;
for(int i=0; i<tsn; ++i) { triple t = rotate(360*i/tsn, C, C+n)*t1; ts.push(t); }

for(triple t: ts) { dot(t); }

\end{asy}
\end{subfloatenv}
\caption{A point at infinity in \T{3} and its dual plane in \T{3*}}
\label{point at infinity in 3D}
\end{figure}

For example, let \(r_0\) be a plane in \T{3*} defined by \(-a+2b+\tfrac{1}{2}c=0\).
This plane passes through the origin of \T{3*} (see Figure~\ref{point at infinity in 3D}(b)).
Points in \T{3*} which lie on the plane \(r_0\) corresponds to a set of planes in \T{3} defined by \(ax+by+cz+1=0\),
where the coefficients \((a,b,c)\) satisfy \(-a+2b+\tfrac{1}{2}c=0\).
In Euclidean space, this set consists of the planes that are parallel to vector \((-1,2,\tfrac{1}{2})\).
Three planes from this set are shown in Figure~\ref{point at infinity in 3D}(a).
These planes form a structure in \T{3}
related to the one shown in Figure~\ref{top-down point in 3D}, which depicts the top-down view of a finite point.

Once the infinite points and lines are defined, I can extend finite lines by points at infinity and finite planes by lines at infinity
in the same way I extended finite lines in \T{2} by points at infinity in \T{2}.
With this extension, I augment the description of a sheaf of planes by the relevant plane passing through the origin.
I also augment a bundle of planes by the plane at infinity and the relevant plane passing through the origin.
Likewise, bundles of planes representing points are augmented.

\subsection{Embedding and Grassmann algebra}
The target space \T{3} is embedded in the target model space \R{4} and the dual space \T{3*} is embedded in the dual model space \R{4*}
in the same way as in the 2-dimensional case.
Namely, \((1,x,y,z)\) in \R{4} is identified with  \((x,y,z)\) in \T{3}
and  \((1,a,b,c)\) in \R{4*} is identified with \((a,b,c)\) in \T{3*}.
So, \T{3} is a 3-dimensional hyperplane embedded in the target model space \R{4} at \(w=1\),
and \T{3*} is a hyperplane at \(d=1\) in the dual model space \R{4*}.

The intersection of linear subspaces of \R{4} with \T{3} gives points, lines, and planes in \T{3}.
In particular, 1-dimensional linear subspaces of \R{4} represent points,
2-dimensional subspaces represent lines, and 3-dimensional subspaces represent planes in \T{3}.
The subspaces that do not intersect \T{3} represent geometric objects at infinity.
The same applies to linear subspaces in \R{4*}.
The duality and identity transformations are extended to linear subspaces of the model spaces \R{4} and \R{4*}.

I will use the following notation for the basis in \R{4}: \(\e^0=(1,0,0,0)\), \(\e^1=(0,1,0,0)\), \(\e^2=(0,0,1,0)\), \(\e^3=(0,0,0,1)\),
so that \((w,x,y,z)=w\e^0+x\e^1+y\e^2+z\e^3\).
And for \R{4*}: \(\e_0=(1,0,0,0)\), \(\e_1=(0,1,0,0)\), \(\e_2=(0,0,1,0)\), \(\e_3=(0,0,0,1)\), so that \((d,a,b,c)=d\e_0+a\e_1+b\e_2+c\e_3\).
The subscript indices refer to the basis of the dual model space \R{4*}
and the superscript indices refer to the target model space \R{4}.

The basis of Grassmann algebra \(\bigwedge\R{4*}\) consists of \(2^4=16\) multivectors:\\
1, \\
\(\e_0, \e_1, \e_2, \e_3\), \\
\(\e_{10}, \e_{20}, \e_{30}, \e_{23}, \e_{31}, \e_{12}\),\\
\(\e_{123}, \e_{320}, \e_{130}, \e_{210}\), \\
\(\I_4=\e_{0123}\).\\
I will use \(\I\) instead of \(\I_4\) when there is no ambiguity.
Capital Roman letters in bold font will be used for trivectors in \(\bigwedge\R{4*}\), which dually represent points in \T{3}.
Capital Greek letters in bold font will be used for bivectors in \(\bigwedge\R{4*}\), which dually represent lines in \T{3}.
So, a general multivector in \R{4*} can be written as
\begin{equation}
M=s+\tb{a}+\mb{\Lambda}+\tb{P}+\alpha\I,
\end{equation}
where \(s,\alpha\in\R{}\).

In \(\bigwedge\R{4*}\) most bivectors are not simple, i.e.\ they cannot be written as the outer product of two vectors.
For instance, the bivector \(\e_{01}+\e_{23}\) is not simple.
A trivial calculation reveals that \((\e_{01}+\e_{23})\wedge(\e_{01}+\e_{23})=2\I\).
If a bivector \(\mb{\Lambda}\) is simple, then \(\mb{\Lambda}\wedge\mb{\Lambda}=0\).
The reverse is true as well, i.e.\ if \(\mb{\Lambda}\wedge\mb{\Lambda}=0\), then \(\mb{\Lambda}\) is simple.

In general, a bivector \(\mb{\Lambda}\) can be written as
\begin{equation}
\mb{\Lambda} = p_{10}\e_{10}+p_{20}\e_{20}+ p_{30}\e_{30}+ p_{23}\e_{23}+ p_{31}\e_{31}+ p_{12}\e_{12}.
\label{bivector on basis in R4}
\end{equation}
Computing the outer product of \(\mb{\Lambda}\) with itself yields
\begin{equation}
\mb{\Lambda}\wedge\mb{\Lambda} = -2(p_{10}p_{23}+p_{20}p_{31}+ p_{30}p_{12})\I.
\end{equation}
So, \(\mb{\Lambda}\) is simple if 
\begin{equation}
p_{10}p_{23}+p_{20}p_{31}+ p_{30}p_{12}=0.
\end{equation}
If \(\mb{\Lambda}\) is simple, it represents a linear subspace in \R{4*} which can be obtained by solving
\begin{equation}
\tb{a}\wedge\mb{\Lambda}=0.
\end{equation}
Substituting  \(\tb{a}=d\e_0+a\e_1+b\e_2+c\e_3\) and \(\mb{\Lambda}\) given by (\ref{bivector on basis in R4}) yields
\[
\begin{aligned}
&(d\e_0+a\e_1+b\e_2+c\e_3)\wedge(p_{10}\e_{10}+p_{20}\e_{20}+ p_{30}\e_{30}+ p_{23}\e_{23}+ p_{31}\e_{31}+ p_{12}\e_{12})\\
&=(ap_{23} +bp_{31}+cp_{12})\e_{123} + (dp_{23} + bp_{30}-cp_{20})\e_{023} +\\
&\quad+(dp_{31}+cp_{10}-ap_{30})\e_{031} + (dp_{12}+ap_{20}-bp_{10})\e_{012}=0,
\end{aligned}
\]
which is equivalent to
\begin{equation}
\left\{
\begin{aligned}
&ap_{23} +bp_{31}+cp_{12}=0,\\
&dp_{23}+bp_{30}-cp_{20}=0,\\
&dp_{31}+cp_{10}-ap_{30}=0,\\
&dp_{12}+ap_{20}-bp_{10}=0.
\end{aligned}
\right.
\end{equation}
Of these equations only two are linearly independent, 
so the system defines a plane in \R{4*} whose intersection with \T{3*} is found by setting \(d=1\),
which gives a line in \T{3*}.
In other words, simple bivectors in \R{4*} represent lines in \T{3*} and,
therefore, they dually represent lines in the target space \T{3}.

The subspace represented by a vector \(\tb{s}\) in \R{4*} is found by solving
\begin{equation}
\tb{a}\wedge\tb{s}=0.
\end{equation}
For example,  \(\tb{s}=(2,\tfrac{1}{2},-\tfrac{1}{3},1)=2(1,\tfrac{1}{4},-\tfrac{1}{6},\tfrac{1}{2})\) represents a 1-dimensional subspace given by
\[
\begin{aligned}
&\tb{a}\wedge\tb{s}=(d\e_0+a\e_1+b\e_2+c\e_3)\wedge(2\e_0+\tfrac{1}{2}\e_1-\tfrac{1}{3}\e_2+\e_3)=\\
&
(2a-\tfrac{d}{2})\e_{10} +
(2b+\tfrac{d}{3})\e_{20} +
(2c-d)\e_{30} +
(b+\tfrac{c}{3})\e_{23} +
(\tfrac{c}{2}-a)\e_{31} +
(-\tfrac{a}{3}-\tfrac{b}{2})\e_{12}=0,
\end{aligned}
\]
which is equivalent to 
\[
\begin{aligned}
2a=\frac{d}{2},\quad
2b=\frac{-d}{3},\quad
2c=d,\quad
b=-\frac{c}{3},\quad
\frac{c}{2}=a,\quad
\frac{a}{3}=\frac{-b}{2}.
\end{aligned}
\]
For \(d=1\), I get a point in \R{4*} with the coordinates \(a=\tfrac{1}{4}, b=-\tfrac{1}{6}, c=\tfrac{1}{2}\),
which dually represents a plane in \T{3} defined by \(\tfrac{1}{4}x-\tfrac{1}{6}y+\tfrac{1}{2}z+1=0\).

A trivector \(\tb{P}\) represents a 3-dimensional hyperplane in \R{4*}, which can be obtained by solving
\begin{equation}
\tb{a}\wedge\tb{P}=0.
\end{equation}
For example, for 
\(\tb{P}=2\e_{123}+\e_{320}+3\e_{130}+4\e_{210}=2(\e_{123}+\tfrac{1}{2}\e_{320}+\tfrac{3}{2}\e_{130}+2\e_{210})\), I get
\[
(d\e_0+a\e_1+b\e_2+c\e_3)\wedge(2\e_{123}+\e_{320}+3\e_{130}+4\e_{210})=(2d+a+3b+4c)\e_{0123}=0,
\]
and
\[
2d+a+3b+4c=0,
\]
which is a hyperplane in \R{4*}.
Its intersection with \T{3*} is obtained by setting \(d=1\), which gives a plane
defined by \(2+a+3b+4c=0\) or \(\tfrac{1}{2}a+\tfrac{3}{2}b+2c+1=0\).
This plane corresponds to the point \((\tfrac{1}{2},\tfrac{3}{2},2)\) in the target space \T{3}.

Grassmann algebra of the target model space \(\R{4}\) is defined in a similar way and exhibits similar properties.
As before, I will use the symbol \(\vee\) for the outer product in \(\bigvee\R{4}\).

The duality transformation \(\J:\bigwedge\R{4*}\to\bigvee\R{4}\)  is defined on the standard basis by the following table

\begin{tabular*}{\textwidth}{lcccccccccccccccc}
\(X\)      & 1   & \(\e_{0}\) & \(\e_{1}\) & \(\e_{2}\) & \(\e_{3}\) &  \(\e_{10}\) & \(\e_{20}\) & \(\e_{30}\) & \(\e_{23}\) & \(\e_{31}\) & \(\e_{12}\) & \(\e_{123}\) & \(\e_{320}\) & \(\e_{130}\) & \(\e_{210}\) & \(\e_{0123}\) \\
\cline{1-17} \\[-10pt]
\(\J(X)\) & \(\e^{0123}\) & \(\e^{123}\) & \(\e^{320}\) & \(\e^{130}\) & \(\e^{210}\) & \(\e^{23}\) & \(\e^{31}\) & \(\e^{12}\) & \(\e^{10}\) & \(\e^{20}\) & \(\e^{30}\) & \(\e^{0}\) & \(\e^{1}\) & \(\e^{2}\) & \(\e^{3}\) & 1 \\
\end{tabular*}

It  extends to general multivectors by linearity.
The inverse transformation is obtained by lowering and raising the indices.

The join of multivectors \(A\) and \(B\) in the dual model space \R{4*} is defined by
\begin{equation}
A\vee B = \J^{-1}(\J(A)\vee\J(B))).
\end{equation}

The wedge and the join possess the following properties in \T{3}:\\
\(\tb{a}\wedge\tb{b}\) is a line at the intersection of the planes \(\tb{a}\) and \(\tb{b}\)\\
\(\tb{a}\wedge\tb{b}\) is a line at infinity if \(\tb{a}\) and \(\tb{b}\) belong to the same stack but do not coincide\\
\(\tb{a}\wedge\tb{b}=0\) if \(\tb{a}\) and \(\tb{b}\) represent the same plane\\
\(\tb{a}\vee\tb{b}=0\) \\
\(\tb{a}\wedge\mb{\Lambda}\) is a point where the plane \(\tb{a}\) and the line \(\mb{\Lambda}\) intersect \\
\(\tb{a}\wedge\mb{\Lambda}=0\) if the line \(\mb{\Lambda}\) is in the plane \(\tb{a}\)\\
\(\tb{a}\vee\mb{\Lambda}=0\) \\
\(\tb{a}\wedge\tb{b}\wedge\tb{c}\) is a point at the intersection of the three planes \(\tb{a}\), \(\tb{b}\), and \(\tb{c}\)\\ 
\(\mb{\Lambda}\vee\tb{P}\) is a plane that passes through the line \(\mb{\Lambda}\) and the point \(\tb{P}\)\\
\(\mb{\Lambda}\vee\tb{P}=0\) if the point \(\tb{P}\) lies on the line \(\mb{\Lambda}\)\\
\(\mb{\Lambda}\wedge\tb{P}=0\)\\
\(\tb{P}\vee\tb{Q}\) is a line that passes through the points \(\tb{P}\) and \(\tb{Q}\)\\
\(\tb{P}\vee\tb{Q}=0\) if \(\tb{P}\) and \(\tb{Q}\) represent the same point\\
\(\tb{P}\wedge\tb{Q}=0\)\\
\(\tb{P}\vee\tb{Q}\vee\tb{R}\) is a plane that passes through the points \(\tb{P}\), \(\tb{Q}\), and \(\tb{R}\)

The following equalities hold:
\begin{equation}
\begin{aligned}
\mb{\Lambda}\wedge\mb{\Phi}=(\mb{\Lambda}\vee\mb{\Phi})\I,\\
\tb{a}\wedge\tb{P}=-(\tb{a}\vee\tb{P})\I,
\end{aligned}
\end{equation}
where
\begin{equation}
\mb{\Lambda}\vee\mb{\Phi}=-( p_{10}q_{23}+p_{20}q_{31}+p_{30}q_{12}+p_{23}q_{10}+p_{31}q_{20}+p_{12}q_{30})
\end{equation}
for \(\mb{\Lambda}=p_{10}\e_{10}+p_{20}\e_{20}+p_{30}\e_{30}+p_{23}\e_{23}+p_{31}\e_{31}+p_{12}\e_{12}\)
and \(\mb{\Phi}=q_{10}\e_{10}+q_{20}\e_{20}+q_{30}\e_{30}+q_{23}\e_{23}+q_{31}\e_{31}+q_{12}\e_{12}\),
and
\begin{equation}
\tb{a}\vee\tb{P}=-(dw+ax+by+cz)
\end{equation}
for \(\tb{a}=d\e_0+a\e_1+b\e_2+c\e_3\) and \(\tb{P}=w\e_{123}+x\e_{320}+y\e_{130}+z\e_{210}\).

\subsection{Orientation}
\begin{figure}[t!]
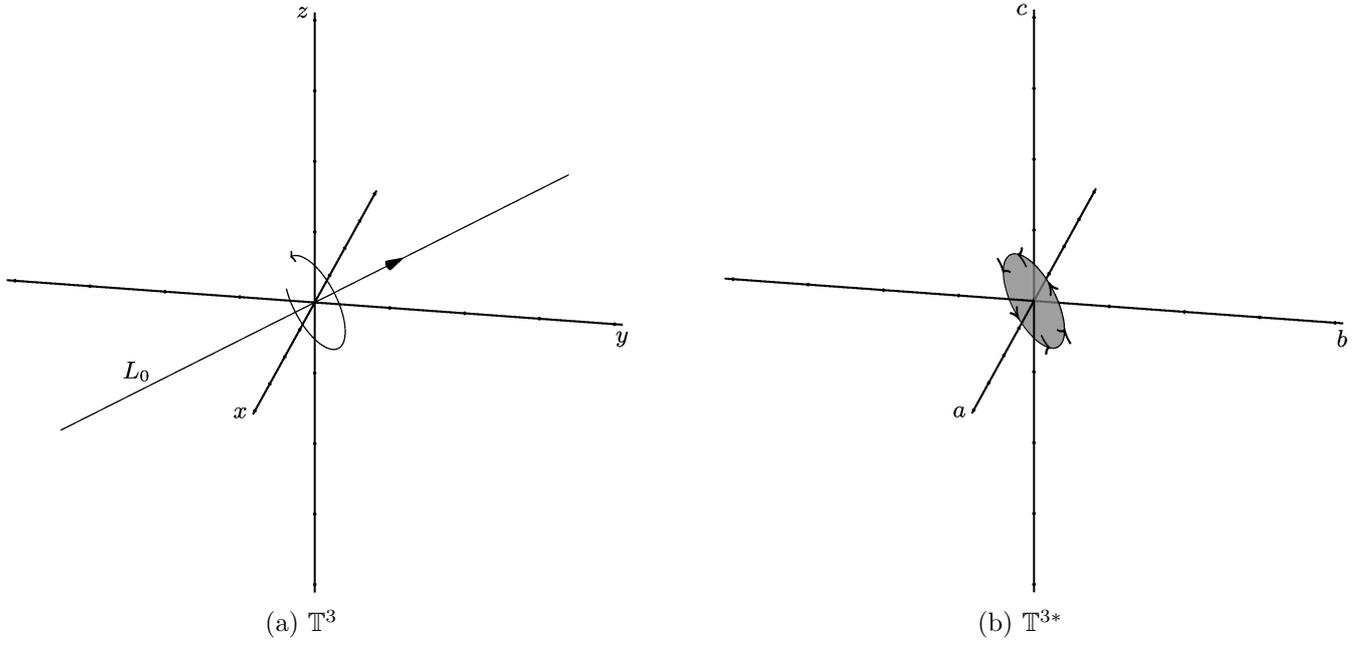
\hspace{-1cm}
\begin{subfloatenv}{\T{3}}
\begin{asy}
import Figure3D;
Figure3D f = Figure3D();

var L = Line(0,0,0,1,-1,-1/3);
f.line(L, "$L_0$", position=0.8, align=(0,-2,0.5), orientation_theta=+50);
vector((-1,1,1/3));

\end{asy}
\end{subfloatenv}\hfill%
\begin{subfloatenv}{\T{3*}}
\begin{asy}
import Figure3D;
Figure3D f = Figure3D(space=DUAL);

var L = Line(0,0,0,1,-1,-1/3);
bivector((1,-1,-1/3));

\end{asy}
\end{subfloatenv}
\caption{The top-down and bottom-up orientation of a line}
\label{orienting line through origin in 3D}
\end{figure}

Since a sheaf of planes in \T{3} is an extension of a sheaf of lines in \T{2},
it can be oriented in a similar fashion.
Thus, the top-down orientation of a finite line in \T{3} is defined by selecting a specific sense of rotation around the line.
A bivector \(\mb{\Lambda}_0=p_{23}\e_{23}+p_{31}\e_{31}+p_{12}\e_{12}\) represents a plane in \R{4*}
that lies entirely in the hyperplane \(d=0\) and does not intersect \T{3*} at any finite point.
It represents a line at infinity in \T{3*}, which is dual to a sheaf of planes in \T{3} whose axis
passes through the origin of \T{3}.
So, \(\mb{\Lambda}_0\) dually represents a line \(L_0\) in \T{3} passing through the origin.
For instance, a line dually represented by  \(\mb{\Lambda}_0=\e_{23}-\e_{31}-\tfrac{1}{3}\e_{12}\) is shown in
 Figure~\ref{orienting line through origin in 3D}(a).
\(\mb{\Lambda}_0\) is defined in \R{4*} but it lies entirely in the hyperplane \(d=0\),
so it can be identified with a bivector in \T{3*}, which is shown in Figure~\ref{orienting line through origin in 3D}(b).
Similarly, the bivector \(\Id(\mb{\Lambda}_0)\) is defined in \R{4} but lies entirely in the hyperplane \(w=0\),
so it can be identified with a bivector in \T{3}, whose orientation is shown with an arc 
in Figure~\ref{orienting line through origin in 3D}(a).
This orientation determines the top-down orientation of \(L_0\).
In other words, the top-down orientation of \(L_0\) is determined by the orientation of 
\(\Id(\mb{\Lambda}_0)=\e^{23}-\e^{31}-\tfrac{1}{3}\e^{12}\)
seen as a bivector in \T{3}.
The relationship between the orientation of \(\mb{\Lambda}_0\) and \(L_0\) is illustrated in Figure~\ref{orienting line through origin in 3D}.

The line \(L_0\) shown in Figure~\ref{orienting line through origin in 3D}(a) is directly represented by 
\(\J(\mb{\Lambda}_0)=\e^{10}-\e^{20}-\tfrac{1}{3}\e^{30}=\e^0\vee\tb{x}\), where \(\tb{x}=-\e^1+\e^2+\tfrac{1}{3}\e^3\)
can be identified with a vector in \T{3}, which provides the bottom-up orientation of the line \(L_0\)
and is shown in Figure~\ref{orienting line through origin in 3D}(a).

\begin{figure}[ht!]
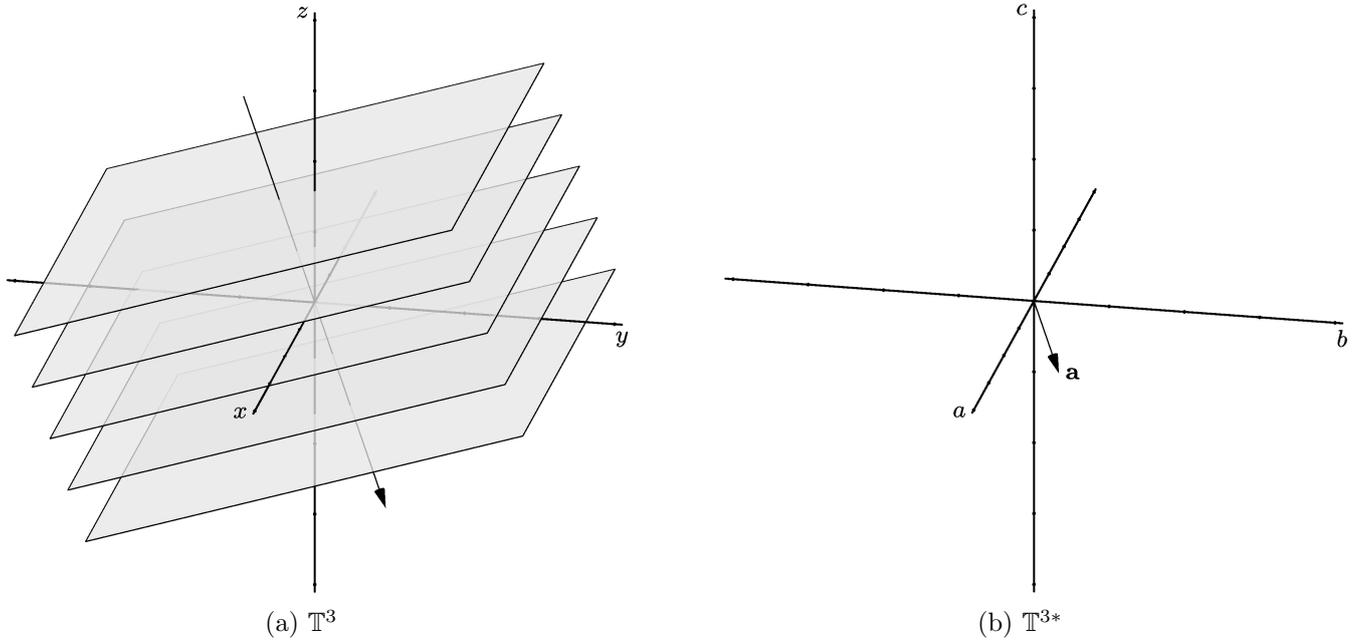
\hspace{-1cm}
\begin{subfloatenv}{\T{3}}
\begin{asy}
import Figure3D;
Figure3D f = Figure3D();

triple n = (0,1/3,-1);
n = normalise(n);

MV[] ps = {};
ps.push(Plane(0,n.x,n.y,n.z));
ps.push(Plane(0.75,n.x,n.y,n.z));
ps.push(Plane(-0.75,n.x,n.y,n.z));
ps.push(Plane(1.5,n.x,n.y,n.z));
ps.push(Plane(-1.5,n.x,n.y,n.z));

real osize = 3;
draw((-osize*n)--(osize*n), Arrow3);

for(var p: ps) { f.plane(p, draw_orientation=false); }

\end{asy}
\end{subfloatenv}\hfill%
\begin{subfloatenv}{\T{3*}}
\begin{asy}
import Figure3D;
Figure3D f = Figure3D(space=DUAL);

triple n = (0,1/3,-1);
draw(Label("$\textbf{a}$", 1, align=E), (0,0,0)--n, Arrow3);

\end{asy}
\end{subfloatenv}
\caption{Orienting a line at infinity in \T{3}}
\label{orienting line at infinity in T3}
\end{figure}

On the other hand, \(\mb{\Lambda}_\infty=p_{10}\e_{10}+p_{20}\e_{20}+p_{30}\e_{30}\) 
represents a line passing through the origin of \T{3*} and, therefore,
dually represents a line at infinity in \T{3}, which corresponds to a stack of planes.
The bivector \(\mb{\Lambda}_\infty\) can be written as \(\mb{\Lambda}_\infty=\e_0\wedge\tb{a}\), where \(\tb{a}\)
can be assumed to lie entirely in the hyperplane \(d=0\), i.e.\ \(\tb{a}=a\e_1+b\e_2+c\e_3\), 
since any non-zero component along \(\e_0\) cancels out by the outer product with \(\e_0\).
For example,  \(\mb{\Lambda}_\infty=-\tfrac{1}{3}\e_{20}+\e_{30}=\e_0\wedge(\tfrac{1}{3}\e_2-\e_3)\)
 dually represents a stack of planes in \T{3}, i.e.\ a line at infinity,
shown in Figure~\ref{orienting line at infinity in T3}(a).
The orientation of \(\Id(\tb{a})=\tfrac{1}{3}\e^2-\e^3\) seen as a vector in \T{3} determines the top-down orientation of the stack of planes,
as illustrated in Figure~\ref{orienting line at infinity in T3}.
The same line at infinity is directly represented by the bivector \(\J(\mb{\Lambda}_\infty)=-\tfrac{1}{3}\e^{31}+\e^{12}\), 
which can be seen as a bivector in \T{3} since it lies entirely in the hyperplane \(w=0\).
It provides the bottom-up orientation of the line at infinity, which can be defined as a particular sense of rotation in
the planes of the stack.

Note that  \(\e_0\wedge\tb{a}\) can be viewed as the intersection of a plane \(\tb{a}\) and the plane at infinity,
which is dually represented by \(\e_0\).

\begin{figure}[t]
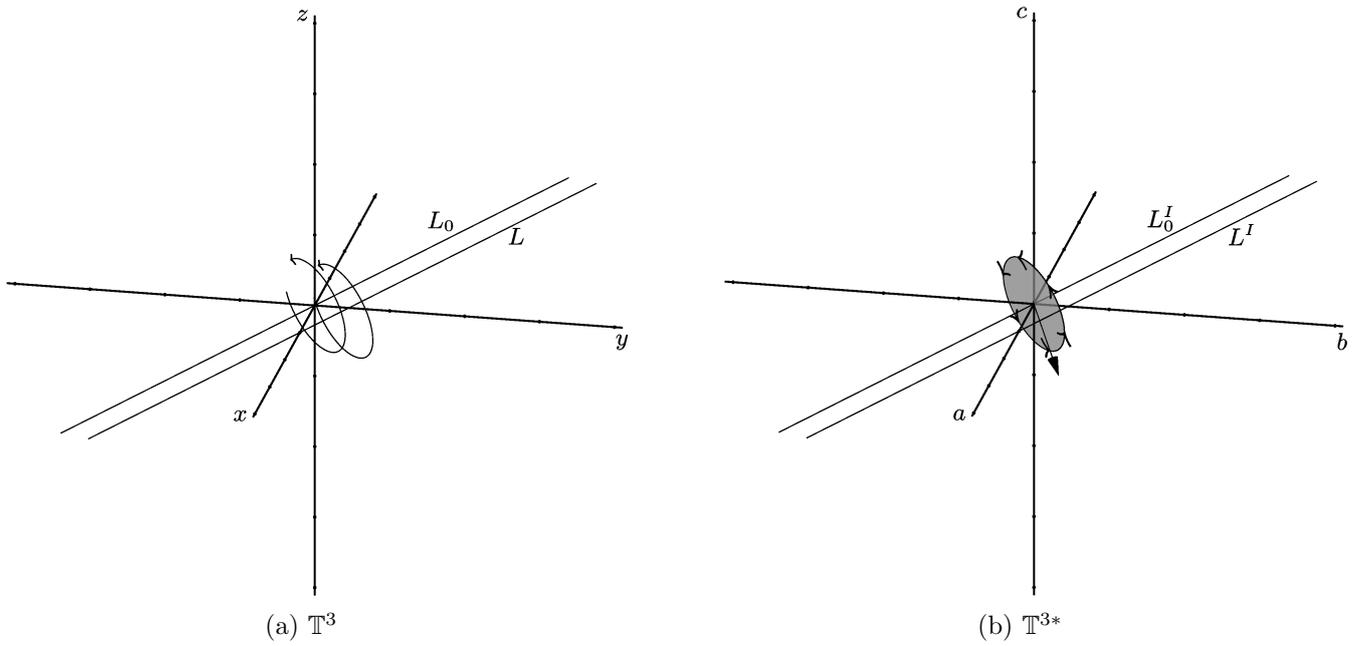
\hspace{-1cm}
\begin{subfloatenv}{\T{3}}
\begin{asy}
import Figure3D;
Figure3D f = Figure3D();

var L0 = Line(0,0,0,1,-1,-1/3);
var L = Line(0,-1/3,1,1,-1,-1/3);
f.line(L0, "$L_0$", position=0.2, align=(0,-2,0.5), orientation_theta=+50);
f.line(L, "$L$", position=0.2, align=(0,2,0), orientation_theta=+50);

\end{asy}
\end{subfloatenv}\hfill%
\begin{subfloatenv}{\T{3*}}
\begin{asy}
import Figure3D;
Figure3D f = Figure3D(space=DUAL);

var L0 = Line(0,0,0,1,-1,-1/3);
var L = Line(0,-1/3,1,1,-1,-1/3);
f.line(L0, "$L_0^I$", position=0.2, align=(0,-2,0.5), draw_orientation=false);
f.line(L, "$L^I$", position=0.2, align=(0,2,0), draw_orientation=false);
bivector((1,-1,-1/3));
vector(-(0,-1/3,1));

\end{asy}
\end{subfloatenv}
\caption{A finite line in \T{3} shifted from the origin}
\label{orienting lines in 3D}
\end{figure}

The bivector \(\mb{\Lambda}=\mb{\Lambda}_0+\mb{\Lambda}_\infty\), 
where \(\mb{\Lambda}_0=\e_{23}-\e_{31}-\tfrac{1}{3}\e_{12}\)  
and \(\mb{\Lambda}_\infty=-\tfrac{1}{3}\e_{20}+\e_{30}=\e_0\wedge\tb{a}\),
is simple and, therefore, dually represents a line in \T{3}.
The line is finite since \(\mb{\Lambda}_0\ne0\),
but unlike \(\mb{\Lambda}_0\), the line  \(\mb{\Lambda}\) is shifted from the origin since \(\mb{\Lambda}_\infty\ne0\).
However, its top-down orientation is still determined by the orientation of the bivector \(\Id(\mb{\Lambda}_0)\)
as illustrated in Figure~\ref{orienting lines in 3D}.
Similarly, its bottom-up orientation is the same as that of the line \(\mb{\Lambda}_0\).
The addition of \(\mb{\Lambda}_\infty\) to a line passing through the origin causes it to shift from the origin
but it has no effect on the orientation.
The direction of the shift depends on the direction of \(\Id(\tb{a})\) seen as a vector in \T{3}
and the amount of the shift depends on the ratio of the weights of  \(\mb{\Lambda}_0\) and  \(\mb{\Lambda}_\infty\).
In Euclidean space, the shift is perpendicular to \(\Id(\tb{a})\) seen as a vector in \T{3}.
 The line can be visualised easily by computing its central point.

\begin{figure}[h!]
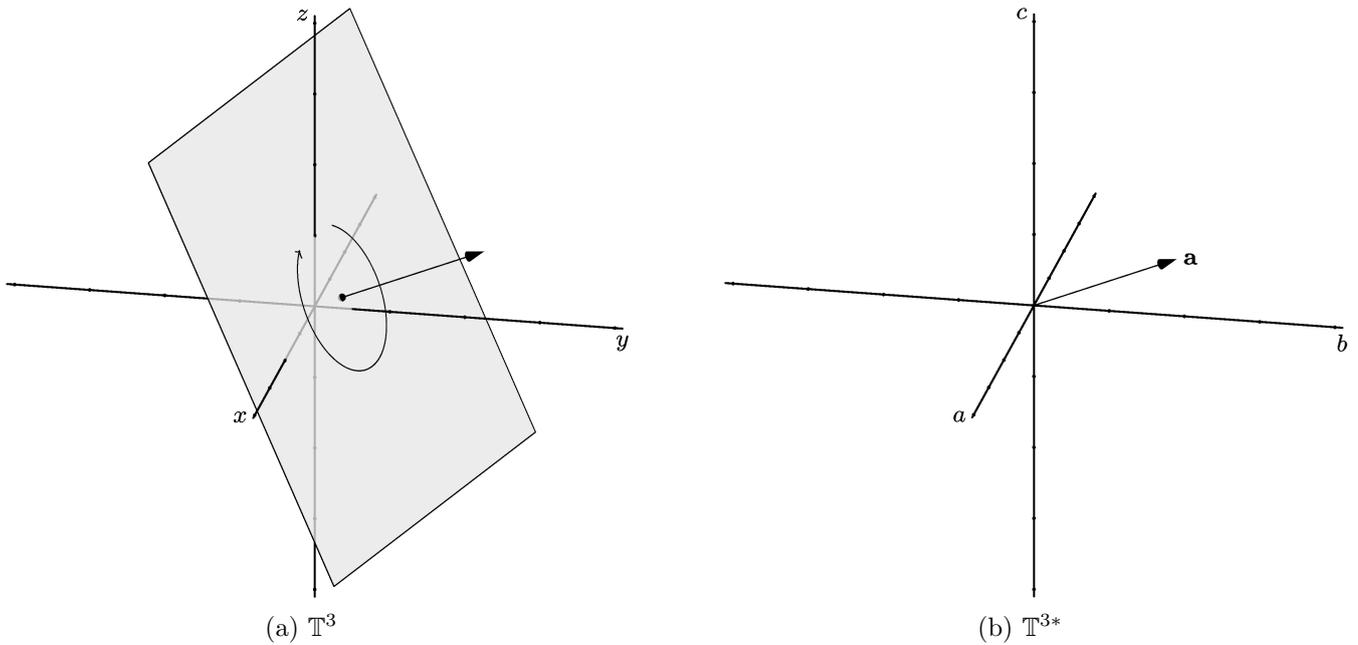
\hspace{-1cm}
\begin{subfloatenv}{\T{3}}
\begin{asy}
import Figure3D;
Figure3D f = Figure3D();

var p = Plane(-1,1/2,2,1);
f.plane(p,draw_bottom_up_orientation_experimental=true);
//bivector((-1/2,-2,-1));

\end{asy}
\end{subfloatenv}\hfill%
\begin{subfloatenv}{\T{3*}}
\begin{asy}
import Figure3D;
Figure3D f = Figure3D(space=DUAL);

triple n = (1/2,2,1);
draw(Label("$\textbf{a}$", 1, align=E), (0,0,0)--n, Arrow3);

\end{asy}
\end{subfloatenv}
\caption{Orienting a plane in \T{3}}
\label{orienting plane in T3}
\end{figure}

The trivector \(\e_{123}\) dually represents the origin of \T{3}.
I identify its top-down orientation with the orientation of  \(\e_{123}\),
which is referred to as right-handed.
On the other hand, \(-\e_{123}\) is the origin with the left-handed orientation.
A finite point with the right-handed orientation and located at \((x,y,z)\) in \T{3} is dually represented by
\(\tb{P}_{\!rh}=\e_{123}+x\e_{320}+y\e_{130}+z\e_{210}\).
The same point with the left-handed orientation is dually represented by 
 \(\tb{P}_{\!lh}=-\tb{P}_{\!rh}=-\e_{123}-x\e_{320}-y\e_{130}-z\e_{210}\).
The point in \T{3} dually represented by \(\tb{P}_{\!rh}\) is directly represented by 
\(\J(\tb{P}_{\!rh})=\e^0+x\e^1+y\e^2+z\e^3\) and is assumed to have the positive bottom-up orientation, whereas
\(\J(\tb{P}_{\!lh})=-\e^0-x\e^1-y\e^2-z\e^3\) represents the same point with the negative bottom-up orientation.
Hence, I will also use the notation \(\tb{P}_{\!+}=  \tb{P}_{\!rh}\) and \(\tb{P}_{\!-}=  \tb{P}_{\!lh}\).
 
Any trivector that can be written as \(\e_0\wedge\mb{\Lambda}\), 
where \(\mb{\Lambda}\) is a line, 
dually represents a point at infinity in \T{3}.
It can be seen as the intersection of a finite line and the plane at infinity.
Note that it is sufficient to consider lines passing through the origin, since any shift from the origin is canceled out by the outer product with \(\e_0\).
Points at infinity inherit their top-down and bottom-up orientation from lines.
The top-down orientation of \(\e_0\wedge\mb{\Lambda}\) is determined by the orientation of \(\Id(\mb{\Lambda}_0)\)
seen as a bivector in \T{3} (\(\mb{\Lambda}_0\) refers to the component of the line \(\mb{\Lambda}\) which passes through the origin).
The bottom-up orientation is determined by the orientation of the vector \(\tb{x}=\J(\e_0\wedge\mb{\Lambda})\),
which lies entirely in the hyperplane \(w=0\) and, hence, can be seen as a vector in \T{3}.
The other way to get the same bottom-up orientation is to express \(\J(\mb{\Lambda}_0)\) as \(\e^0\vee\tb{x}\).
In both cases, if \(\mb{\Lambda}_0=p_{23}\e_{23}+p_{31}\e_{31}+p_{12}\e_{12}\), then 
\(\tb{x}=-p_{23}\e^1-p_{31}\e^2-p_{12}\e^3\) and the bottom-up orientation is provided by 
the vector \((-p_{23},-p_{31},-p_{12})\) in \T{3}.

The top-down orientation of a plane in \T{3} is closely related to the top-down orientation of the line at infinity that lies in the plane.
If a plane is dually represented by  \(\tb{a}=d\e_0+a\e_1+b\e_2+c\e_3\), the top down orientation is found
by taking \(\Id(\tb{a})\), removing the component along \(\e^0\), and considering the result as a vector in \T{3}.
So, the top-down orientation of a plane dually represented by 
\(\tb{a}=d\e_0+a\e_1+b\e_2+c\e_3\) is given by the vector \((a,b,c)\) in \T{3}
(see Figure~\ref{orienting plane in T3}, where \(\tb{a}=-\e_0+\tfrac{1}{2}\e_1+2\e_2+\e_3\)).
The same plane is directly represented by the trivector \(\J(\tb{a})=d\e^{123}+a\e^{320}+b\e^{130}+c\e^{210}\).
Since the component along \(\e^{123}\) is responsible for the shift from the origin only and has no effect on the orientation,
it can be ignored. 
This gives the trivector \(T=a\e^{320}+b\e^{130}+c\e^{210}=\e^0\vee B\),
where the bivector \(B=-a\e^{23}-b\e^{31}-c\e^{12}\) lies entirely in the hyperplane \(w=0\) and can be seen as a bivector in \T{3}.
The bivector \(B\) provides the bottom-up orientation of the plane.
The bottom-up orientation of a plane selects a sense of rotation in the plane.

In the top-down view, the plane at infinity \(d\e_0\) is oriented towards the origin if \(d>0\) and away from the origin if \(d<0\).
In the bottom-up view, it is either right-handed if  \(d>0\)  or left-handed if \(d<0\).


Summary (from \(\bigwedge\R{4*}\) to \T{3})


Trivectors \(\rightarrow\) Points in \T{3}: \\
\(\tb{P}=\e_{123}+x\e_{320}+y\e_{130}+z\e_{210} \rightarrow\) a right-handed (positive) point at \((x,y,z)\). \\
\(\alpha\tb{P}\), \(\alpha>0\) \(\rightarrow\) same as above but with a different weight.\\
\(-\tb{P}\) \(\rightarrow\) same as \(\tb{P}\)   but with the opposite orientation (left-handed or negative).\\
\(\tb{N}=\e_0\wedge\mb{\Lambda}_0\) where \(\mb{\Lambda}_0=p_{23}\e_{23}+p_{31}\e_{31}+p_{12}\e_{12}\) \(\rightarrow\) a point at infinity 
that lies on a line \(\mb{\Lambda}_0\);
its top-down orientation is given by \(\Id(\mb{\Lambda}_0)=p_{23}\e^{23}+p_{31}\e^{31}+p_{12}\e^{12}\),
while the bottom-up orientation by \(\J(\e_0\wedge\mb{\Lambda}_0)=-p_{23}\e^{1}-p_{31}\e^{2}-p_{12}\e^{3}\), both seen as blades in \T{3}.\\
\(\alpha\tb{N}\), \(\alpha>0\) \(\rightarrow\) same as above but with a different weight.\\
\(-\tb{N}\) \(\rightarrow\) same as \(\tb{N}\)  but with the opposite orientation.

Simple bivectors \(\rightarrow\) Lines in \T{3}: \\
\(\mb{\Lambda}_0=p_{23}\e_{23}+p_{31}\e_{31}+p_{12}\e_{12}\) \(\rightarrow\)
a line passing through the origin of \T{3} and coinciding with the linear subspace of the vector \((p_{23},p_{31},p_{12})\)
with the top-down orientation given by \(\Id(\mb{\Lambda}_0)=p_{23}\e^{23}+p_{31}\e^{31}+p_{12}\e^{12}\)
and bottom-up orientation by \(\J(\e_0\wedge\mb{\Lambda}_0)=-p_{23}\e^{1}-p_{31}\e^{2}-p_{12}\e^{3}\).\\
\(\mb{\Phi}=\e_0\wedge\tb{a}_0\) where \(\tb{a}_0=a\e_1+b\e_2+c\e_3\) \(\rightarrow\) a line at infinity in \T{3} that lies in the plane \(ax+by+cz=0\);
its top-down orientation is given by \(\Id(\tb{a}_0)=a\e^1+b\e^2+c\e^3\), or \((a,b,c)\),
and bottom-up orientation by \(\J(\e_0\wedge\tb{a}_0)=-a\e^{23}-b\e^{31}-c\e^{12}\).\\
\(\mb{\Lambda}=p_{10}\e_{10}+p_{20}\e_{20}+ p_{30}\e_{30}+ p_{23}\e_{23}+ p_{31}\e_{31}+ p_{12}\e_{12}\), \(p_{10}p_{23}+p_{20}p_{31}+p_{30}p_{12}=0\) \(\rightarrow\) 
a finite line in \T{3} consisting of \((x,y,z)=\pm(p_{23},p_{31},p_{12})\tau+(x_c,y_c,z_c)\), where \(\tau\in\R{}\) and \((x_c,y_c,z_c)\)
is the centre of the line;
the orientation of \(\mb{\Lambda}\) is the same as that of \(\mb{\Lambda}_0\). 

Vectors \(\rightarrow\) finite planes in \T{3} and the plane at infinity:\\
\(\tb{a}=d\e_0+a\e_1+b\e_2+c\e_3\) \(\rightarrow\) a plane in \T{3} defined by \(ax+by+cz+d=0\);
the top-down orientation by \(\Id(\tb{a}_0)=a\e^1+b\e^2+c\e^3\), where \(\tb{a}_0=a\e_1+b\e_2+c\e_3\), 
and the bottom-up orientation by \(\J(\e_0\wedge\tb{a}_0)=-a\e^{23}-b\e^{31}-c\e^{12}\).\\
\(\alpha\tb{a}\), \(\alpha>0\) \(\rightarrow\) same plane with a different weight.\\
\(-\tb{a}\) \(\rightarrow\) same as \(\tb{a}\) but with the opposite orientation.\\ 
\(d\e_0\) \(\rightarrow\) the plane at infinity in \T{3} (in the top-down view, towards the origin if \(d>0\) or away from the origin if \(d<0\)).


\subsection{The metric and Clifford algebra}
The standard basis vectors are assumed to be orthogonal.
For vectors \(\tb{a}_1\) and \(\tb{a}_2\) in \R{4*}, the inner product is given by
\begin{equation}
\tb{a}_1\cdot\tb{a}_2=d_1d_2\e_0\cdot\e_0+a_1a_2\e_1\cdot\e_1+b_1b_2\e_2\cdot\e_2 + c_1c_2\e_3\cdot\e_3,
\end{equation}
where \(\tb{a}_1=d_1\e_0+a_1\e_1+b_1\e_2+c_1\e_3\) and  \(\tb{a}_2=d_2\e_0+a_2\e_1+b_2\e_2+c_2\e_3\),
with the following list of options for the inner product of the standard basis vectors:
\begin{center}
\begin{tabular}{cccccl}
\(\e_0\cdot\e_0\) & \(\e_1\cdot\e_1\) & \(\e_2\cdot\e_2\) & \(\e_3\cdot\e_3\) & \T{3} & \\
0 & 1 & 1 & 1 &\E{3}& Euclidean \\ 
1& 1& 1 & 1& \El{3} & Elliptic \\
\m1 & 1& 1 & 1& \Hy{3} & Hyperbolic \\
0 & 1  & 1 & \m1& \M{3} &Minkowski (pseudo-Euclidean)  \\
 1 & 1 & 1 & \m1& \(dS_{3}\) & de-Sitter \\
 \m1 & 1 & 1 & \m1& \(AdS_{3}\) & Anti de-Sitter \\
\end{tabular}
\end{center}
For instance, \(\tb{a}_1\cdot\tb{a}_2=-d_1d_2+a_1a_2+b_1b_2+c_1c_2\) if the metric is hyperbolic.
Adopting one of these options yields a specific metric geometry in \T{3}.

In \(\bigwedge\R{4*}\), the geometric product between a \(k\)-vector \(A_k\) and an \(l\)-vector \(B_l\)  includes at most three multivectors
of grades \(|k-l|\), \(|k-l|+2\), and \(k+l\).
It exhibits the following properties:
\begin{equation}
\begin{split}
&\tb{a}\tb{b}=\tb{a}\cdot\tb{b}+\tb{a}\wedge\tb{b} \textrm{ (scalar + bivector)},\quad
\tb{a}\cdot\tb{b}=\tb{b}\cdot\tb{a},\quad
\tb{a}\wedge\tb{b}=-\tb{b}\wedge\tb{a},\\
&\tb{a}\mb{\Lambda}=\tb{a}\cdot\mb{\Lambda}+\tb{a}\wedge\mb{\Lambda} \textrm{ (vector + trivector)},\quad
\tb{a}\cdot\mb{\Lambda}=-\mb{\Lambda}\cdot\tb{a},\quad
\tb{a}\wedge\mb{\Lambda}=\mb{\Lambda}\wedge\tb{a},\\
&\tb{a}\tb{P}=\tb{a}\cdot\tb{P}+\tb{a}\wedge\tb{P} \textrm{ (bivector + pseudoscalar)},\quad
\tb{a}\cdot\tb{P}=\tb{P}\cdot\tb{a},\quad
\tb{a}\wedge\tb{P}=-\tb{P}\wedge\tb{a},\\
&\tb{a}\I=\tb{a}\cdot\I \textrm{ (trivector)},\quad
\tb{a}\cdot\I=-\I\cdot\tb{a},\\
&\mb{\Lambda}\mb{\Phi}=\mb{\Lambda}\cdot\mb{\Phi}+\mb{\Lambda}\times\mb{\Phi}+\mb{\Lambda}\wedge\mb{\Phi} \textrm{ (scalar + bivector + pseudoscalar)},\\
&\quad\quad\quad\mb{\Lambda}\cdot\mb{\Phi}=\mb{\Phi}\cdot\mb{\Lambda},\quad
\mb{\Lambda}\wedge\mb{\Phi}=\mb{\Phi}\wedge\mb{\Lambda},\\
&\mb{\Lambda}\tb{P}=\mb{\Lambda}\cdot\tb{P}+\mb{\Lambda}\times\tb{P} \textrm{ (vector + trivector)},\quad
\mb{\Lambda}\cdot\tb{P}=\tb{P}\cdot\mb{\Lambda},\\
&\mb{\Lambda}\I=\mb{\Lambda}\cdot\I \textrm{ (bivector)},\quad
\mb{\Lambda}\cdot\I = \I\cdot\mb{\Lambda},\\
&\tb{P}\tb{Q}=\tb{P}\cdot\tb{Q}+\tb{P}\times\tb{Q} \textrm{ (scalar + bivector)},\quad
\tb{P}\cdot\tb{Q}=\tb{Q}\cdot\tb{P},\\
&\tb{P}\I=\tb{P}\cdot\I \textrm{ (vector)},\quad
\tb{P}\cdot\I=-\I\cdot\tb{P}.
\end{split}
\end{equation}
For completeness, note the following relations: \\
\(\tb{a}\vee\tb{P}=-\tb{P}\vee\tb{a}\),\;
  \(\mb{\Lambda}\vee\tb{P}=\tb{P}\vee\mb{\Lambda}\),\;
\(\tb{P}\vee\tb{Q}=-\tb{Q}\vee\tb{P}\),\;
\(\mb{\Lambda}\vee\mb{\Phi}=\mb{\Phi}\vee\mb{\Lambda}\).

It follows that
\begin{equation}
\begin{split}
&\tb{a}\cdot\tb{b}=\tfrac{1}{2}(\tb{a}\tb{b}+\tb{b}\tb{a}),\quad 
\tb{a}\wedge\tb{b}=\tfrac{1}{2}(\tb{a}\tb{b}-\tb{b}\tb{a}),\quad \tb{a}\times\tb{b}=\tb{a}\wedge\tb{b},  \\
&\tb{a}\cdot\mb{\Lambda}=\tfrac{1}{2}(\tb{a}\mb{\Lambda}-\mb{\Lambda}\tb{a}),\quad
\tb{a}\wedge\mb{\Lambda}=\tfrac{1}{2}(\tb{a}\mb{\Lambda}+\mb{\Lambda}\tb{a}),\quad  
\tb{a}\times\mb{\Lambda}=\tb{a}\cdot\mb{\Lambda},  \\
&\tb{a}\cdot\tb{P}=\tfrac{1}{2}(\tb{a}\tb{P}+\tb{P}\tb{a}),\quad
\tb{a}\wedge\tb{P}=\tfrac{1}{2}(\tb{a}\tb{P}-\tb{P}\tb{a}),\quad  \tb{a}\times\tb{P}=\tb{a}\wedge\tb{P},  \\
&\mb{\Lambda}\cdot\mb{\Phi}+\mb{\Lambda}\wedge\mb{\Phi}
=\tfrac{1}{2}(\mb{\Lambda}\mb{\Phi}+\mb{\Phi}\mb{\Lambda}),\quad
\mb{\Lambda}\times\mb{\Phi}=\tfrac{1}{2}(\mb{\Lambda}\mb{\Phi}-\mb{\Phi}\mb{\Lambda}),\\
&\mb{\Lambda}\cdot\tb{P}=\tfrac{1}{2}(\mb{\Lambda}\tb{P}+\tb{P}\mb{\Lambda}),\quad
\mb{\Lambda}\times\tb{P}=\tfrac{1}{2}(\mb{\Lambda}\tb{P}-\tb{P}\mb{\Lambda}),\\
&\tb{P}\cdot\tb{Q}=\tfrac{1}{2}(\tb{P}\tb{Q}+\tb{Q}\tb{P}), \quad
\tb{P}\times\tb{Q}=\tfrac{1}{2}(\tb{P}\tb{Q}-\tb{Q}\tb{P}).\\
\end{split}
\end{equation}

The following identities are useful:
\(\mb{\Lambda}\times\tb{P}=-(\mb{\Lambda}\vee\tb{P})\I,\quad
\tb{P}\times\mb{\Lambda}=(\tb{P}\vee\mb{\Lambda})\I,\quad
\tb{P}\times\tb{Q}=-(\tb{P}\vee\tb{Q})\I\).
The reverse and grade involution of multivector \(M=s+\tb{a}+\mb{\Lambda}+\tb{P}+\alpha\I\) are given by
\(\reverse{M}=s+\tb{a}-\mb{\Lambda}-\tb{P}+\alpha\I\) and \(\involute{M}=s-\tb{a}+\mb{\Lambda}-\tb{P}+\alpha\I\),
respectively.
As in \R{3*}, \(\reverse{AB}=\reverse{B}\reverse{A}\) and \(\involute{AB}=\involute{A}\involute{B}\) for any multivectors \(A\) and \(B\).
The same expressions apply to the inner and outer products and the commutator.
I also define \((M)_{14}=s-\tb{a}+\mb{\Lambda}+\tb{P}-\alpha\I\) for reversing the sign of grades 1 and 4.
In general, \((M)_G\) reverses the sign of all grades in \(M\) that are in the list of grades \(G\),
so that \(\reverse{M}=(M)_{23}\) and \(\involute{M}=(M)_{13}\).

The norm of multivector \(M\) is defined by
\begin{equation}
\norm{M}=|M\reverse{M}(M\reverse{M})_{14}|^{\tfrac{1}{4}}.
\end{equation}
If the norm of \(M\) is not zero, its inverse is given by
\begin{equation}
M^{-1}=\frac{\reverse{M}(M\reverse{M})_{14}}{M\reverse{M}(M\reverse{M})_{14}}.
\end{equation}
For blades, \(A_k\reverse{A}_k\) is a scalar.
So, the expressions for the norm and the inverse simplify:
\begin{equation}
\norm{A_k}=|A_k\reverse{A}_k|^{\tfrac{1}{2}}, \quad A_k^{-1}=\frac{\reverse{A}_k}{A_k\reverse{A}_k}.
\end{equation}

The algebra of even multivectors in \(\bigwedge\R{4*}\) is 8-dimensional and consists of multivectors
\begin{equation}
E=u+ p_{10}\e_{10}+p_{20}\e_{20}+p_{30}\e_{30}+p_{23}\e_{23}+p_{31}\e_{31}+p_{12}\e_{12} +v\I,
\end{equation}
where \(u,v\in\R{}\).
This algebra is associative and, therefore, it cannot be isomorphic to the algebra of octonians which is 8-dimensional too but  not associative.
In fact, the algebra of even multivectors in \(\bigwedge\R{4*}\) is isomorphic to the algebra of biquaternions.

\subsection{Euclidean space \E{3}}

For a plane \(\tb{a}=d\e_0+a\e_1+b\e_2+c\e_3\), 
a line \(\mb{\Lambda}=p_{10}\e_{10}+p_{20}\e_{20}+p_{30}\e_{30}+p_{23}\e_{23}+p_{31}\e_{31}+p_{12}\e_{12}\),
where \(p_{10}p_{23}+p_{20}p_{31}+ p_{30}p_{12}=0\), and 
a point \(\tb{P}=w\e_{123}+x\e_{320}+y\e_{130}+z\e_{210}\), the norm is given by 
\begin{equation}
\norm{\tb{a}}=\sqrt{a^2+b^2+c^2},\quad
\norm{\mb{\Lambda}}=\sqrt{p_{23}^2+p_{31}^2+p_{12}^2},\quad
\norm{\tb{P}}=|w|.
\end{equation}
Furthermore, \(\tb{a}^2=a^2+b^2+c^2\), \(\mb{\Lambda}^2=-(p_{23}^2+p_{31}^2+p_{12}^2)\), and \(\tb{P}^2=-w^2\).
A normalised point can be right-handed (positive) or left-handed (negative).
Points, lines, and the plane at infinity have zero norm.

\begin{figure}[h!]
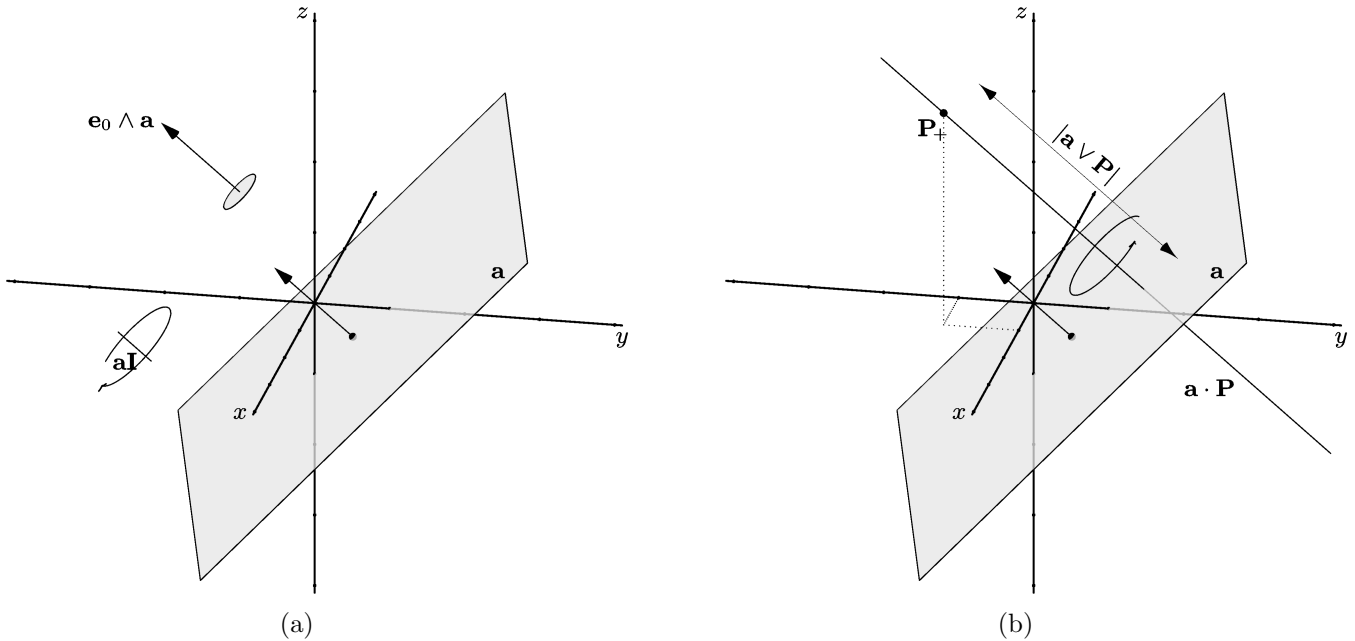
\hspace{-1cm}
\begin{subfloatenv}{ }
\begin{asy}
import Figure3D;
Figure3D f = Figure3D();

var a = Plane(1,1/4,-1,1);
f.plane(a, "$\textbf{a}$", align=(1,-3,-1));

var L = wedge(e_0,a);
f.line_at_infinity(L, (0,-1,1.5), "$\textbf{e}_0\wedge\textbf{a}$", align=(0,-1,0));

var N = a*I;
f.point_at_infinity(N, (2,-2,0), "$\textbf{a}\textbf{I}$", align=(0,-0.5,-1));

\end{asy}
\end{subfloatenv}\hfill%
\begin{subfloatenv}{ }
\begin{asy}
import Figure3D;
Figure3D f = Figure3D();

var a = Plane(1,1/4,-1,1);
f.plane(a, "$\textbf{a}$", align=(1,-3,-1));

var P = Point(1,1,-1,3);
f.point(P, "$\textbf{P}$", align=(0,-0.5,-1));

var L = dot(a,P);
f.line(L, "$\textbf{a}\cdot\textbf{P}$", align=(0,-1,-1));

var Q = wedge(a,dot(a,P));
f.stretch(P, Q, "$|\textbf{a}\vee\textbf{P}|$", label_angle=-40,align=(0.4,0.6), shift_distance=0.65);

\end{asy}
\end{subfloatenv}
\caption{Basic properties of points and planes  in \E{3}}
\label{basic  E3}
\end{figure}

The line at infinity \(\e_0\wedge\tb{a}\) provides the orientation of a plane \(\tb{a}\).
I will depict it as an arrow attached to a small circular base, which can be seen
as a small segment of one of the planes in the stack that represents the line at infinity.
The arrow indicates the orientation of the stack.
On the other hand, \(\tb{a}\I\) is a point at infinity which lies in the direction
perpendicular to the plane \(\tb{a}\) (it is called the polar point of \(\tb{a}\)).
I will depict it as a short line segment indicating the direction towards the point at infinity
and an arc around it indicating its orientation.
Note that even though this depiction is similar to that of a finite line,
the set of planes in \T{3} representing a point at infinity is different from a sheaf of planes,
which represents a finite line.
The plane \(\tb{a}=\e_0+\tfrac{1}{4}\e_1-\e_2+\e_3\), the line at infinity \(\e_0\wedge\tb{a}\) associated to it,
and its polar point \(\tb{a}\I\) are shown on Figure~\ref{basic  E3}(a).

The line that is perpendicular to a plane \(\tb{a}\) and passes through a point \(\tb{P}\)
is obtained by taking the inner product \(\tb{a}\cdot\tb{P}\)
(see Figure~\ref{basic  E3}(b) where \(\tb{P}=\e_{123}+\e_{320}-\e_{130}+3\e_{210}\)).
Recall that the top-down view of a point is a bundle of planes which consists of all planes passing through the point.
The dot product picks out the bundle's planes that are perpendicular to \(\tb{a}\), which
results in a sheaf of planes whose line of intersection passes through \(\tb{P}\) and is perpendicular to \(\tb{a}\).
For a positively oriented point \(\tb{P}\) indicated in the figures by the subscript \(+\),
the orientation of the line  \(\tb{a}\cdot\tb{P}\)  and the plane \(\tb{a}\) conform to the right-hand rule.
The distance  between the point and the plane is given by \(|\tb{a}\vee\tb{P}|\)
if  both \(\tb{a}\) and \(\tb{P}\) are normalised.
Moreover, \(\tb{a}\vee\tb{P}\) is negative if \(\tb{a}\) is oriented towards a positively oriented point \(\tb{P}\).
In general, \(|d|\) gives the distance from the origin to a normalised plane \(\tb{a}=d\e_0+a\e_1+b\e_2+c\e_3\) and
\(d\) is positive if the plane is oriented towards the origin.

For a finite point \(\tb{P}=w\e_{123}+x\e_{320}+y\e_{130}+z\e_{210}\), I get \(\tb{P}\I=w\e_0\), which is the plane at infinity weighted by \(w\).

For two points \(\tb{P}=\e_{123}+\e_{320}\) and \(\tb{Q}=\tfrac{1}{2}\e_{123}+\tfrac{1}{2}\e_{130}+\tfrac{1}{6}\e_{210}\),
the join 
\(\tb{P}\vee\tb{Q}
=
-\tfrac{1}{6}\e_{20}+\tfrac{1}{2}\e_{30}
+\tfrac{1}{2}\e_{23}-\tfrac{1}{2}\e_{31}-\tfrac{1}{6}\e_{12}\) 
is a line that passes through the points (see Figure~\ref{basic  E3 2}(a)).
The distance \(r\) between \(\tb{P}\) and \(\tb{Q}\) is given by 
\begin{equation}
r=\norm{\tb{P}\vee\tb{Q}}
\end{equation}
if the points are normalised.
Taking normalisation into account, the distance between the points evaluates to  \(\tfrac{\sqrt{19}}{3}\).
The commutator \(\tb{P}\times\tb{Q}\) is a line at infinity whose stack consists of the  planes perpendicular to the line \(\tb{P}\vee\tb{Q}\).
The orientation of \(\tb{P}\times\tb{Q}\) and \(\tb{P}\vee\tb{Q}\) conform to the right-hand rule
if both points are right-handed (or left-handed).
Similar results obtain for the join \(\tb{P}\vee\tb{N}\) of a finite point \(\tb{P}\) and a points at infinity  \(\tb{N}\) as
shown in Figure~\ref{basic  E3 2}(b) where \(\tb{P}=\e_{123}+\e_{320}\) 
and \(\tb{N}=\e_0\wedge (\tfrac{1}{2}\e_{23}-\tfrac{1}{2}\e_{31}-\tfrac{1}{6}\e_{12})\).

\begin{figure}[t!]\hspace{-1cm}
\begin{subfloatenv}{ }
\begin{asy}
import Figure3D;
Figure3D f = Figure3D();

var P = Point(1,1,0,0);
var Q = Point(1/2,0,1/2,1/6);
f.point(P, "$\textbf{P}$", align=(0,-1,0.5));
f.point(Q, "$\textbf{Q}$", align=(0,1,-0.5));

var L = join(P,Q);
f.line(L, "$\textbf{P}\vee\textbf{Q}$", align=(0,-1,0.5));

var L0 = cross(P,Q);
f.line_at_infinity(L0, (0,-1,1.5), "$\textbf{P}\times\textbf{Q}$", align=(0,0,-1));

\end{asy}
\end{subfloatenv}\hfill%
\begin{subfloatenv}{ }
\begin{asy}
import Figure3D;
Figure3D f = Figure3D();

var P = Point(1,1,0,0);
var N = wedge(e_0, Line(0,0,0,1/2,-1/2,-1/6));
f.point(P, "$\textbf{P}$", align=(0,-1,0.5));
f.point_at_infinity(N, (2,2,0), "$\textbf{N}$", align=(0,-2.5,-1.5));

var L = join(P,N);
f.line(L, "$\textbf{P}\vee\textbf{N}$", align=(0,-1,0.5));

var L0 = cross(P,N);
f.line_at_infinity(L0, (0,-1,1.5), "$\textbf{P}\times\textbf{N}$", align=(0,0,-1));

\end{asy}
\end{subfloatenv}
\caption{The join of two points in \E{3}}
\label{basic  E3 2}
\hspace{-1cm}
\begin{subfloatenv}{ }
\begin{asy}
import Figure3D;
Figure3D f = Figure3D();

var a = Plane(1,0,-1/6,1/2);
var b = Plane(-1/2,0,-1/3,1);

f.plane(a,"$\textbf{a}$",align=(0,-1,1));
f.plane(b,"$\textbf{b}$",align=(0,-1,1));

var L = wedge(a,b);
f.line_at_infinity(L,(0,-0.75,2.25),"$\mathbf{\Phi}$");

\end{asy}
\end{subfloatenv}
\begin{subfloatenv}{ }
\begin{asy}
import Figure3D;
Figure3D f = Figure3D();

MV a = Plane(-1,1,4/3,-1); //  sqrt(34)/3
MV b = Plane(-1,1, 5/6,1/2); //  sqrt(70)/6

a = a/norm(a);
b = b/norm(b);

var L = wedge(a,b);
var u = (L.p23,L.p31,L.p12);

var centre = totriple(toline(L).centre());

f.plane(a, u, "$\textbf{a}$", align=(0,2,-0.25), draw_orientation=true, O=centre);
f.plane(b, u, "$\textbf{b}$", align=(0,2,-0.25), draw_orientation=true, O=centre);

var v1 = (a.a,a.b,a.c);
var v2 = (b.a,b.b,b.c);

f.arc(centre, v1, v2,"$\alpha$",  (0,1,0.5));
f.line(L,"$\mathbf{\Lambda}$",0.07, align=(0,0,-1));

\end{asy}
\end{subfloatenv}
\caption{Parallel and intersecting planes in \E{3}}
\label{intersecting planes in E3}
\end{figure}

For normalised parallel (and anti-parallel) planes, \(\tb{a}\tb{b}=\pm1+\tb{a}\wedge\tb{b}\), where
\(\mb{\Phi}=\tb{a}\wedge\tb{b}\) is a line at infinity that lies in both planes.
The sign of the scalar component is positive if both planes are oriented in the same direction (parallel planes) and negative otherwise
(anti-parallel planes.)
The orientation of  \(\mb{\Phi}\) is from \(\tb{a}\) to \(\tb{b}\) if both planes are oriented in the same direction.
Otherwise, it is from \(\tb{b}\) to \(\tb{a}\).
The line at infinity \(\mb{\Phi}=\tb{a}\wedge\tb{b}\) for \(\tb{a}=\e_0-\tfrac{1}{6}\e_2+\tfrac{1}{2}\e_3\) and 
\(\tb{b}=-\tfrac{1}{2}\e_0-\tfrac{1}{3}\e_2+\e_3\) 
is shown in Figure~\ref{intersecting planes in E3}(a).
The distance between the planes is given by \(\norm{\e_{123}\vee(\tb{a}\wedge\tb{b})}\) if both planes are normalised.
For normalised intersecting planes \(\tb{a}\) and \(\tb{b}\), the geometric product yields
\begin{equation}
\tb{a}\tb{b}=\cos{\alpha}+\mb{\Lambda}\sin{\alpha},
\end{equation}
where \(\mb{\Lambda}=\tb{a}\wedge\tb{b}/\sin{\alpha}\) and
\begin{equation}
\cos{\alpha}=\tb{a}\cdot\tb{b}
\end{equation}
defines the angle \(\alpha\) between the planes, or more precisely the angle between the orientation directions of the planes.
The normalised planes \(\tb{a}=\tfrac{3}{\sqrt{34}}(-\e_0+\e_1 +\tfrac{4}{3}\e_2-\e_3)\) and
\(\tb{b}=\tfrac{6}{\sqrt{70}}(-\e_0+\e_1 +\tfrac{5}{6}\e_2+\tfrac{1}{2}\e_3)\) shown in Figure~\ref{intersecting planes in E3}(b)
 intersect along the line 
\(\mb{\Lambda}=\tfrac{1}{\sqrt{19}}(-\e_{20}+3\e_{30}+3\e_{23}-3\e_{31}-\e_{12})\),
which is also normalised by construction.
For the angle between the planes, I get \(\alpha\approx54^\circ\).

\begin{figure}[t]
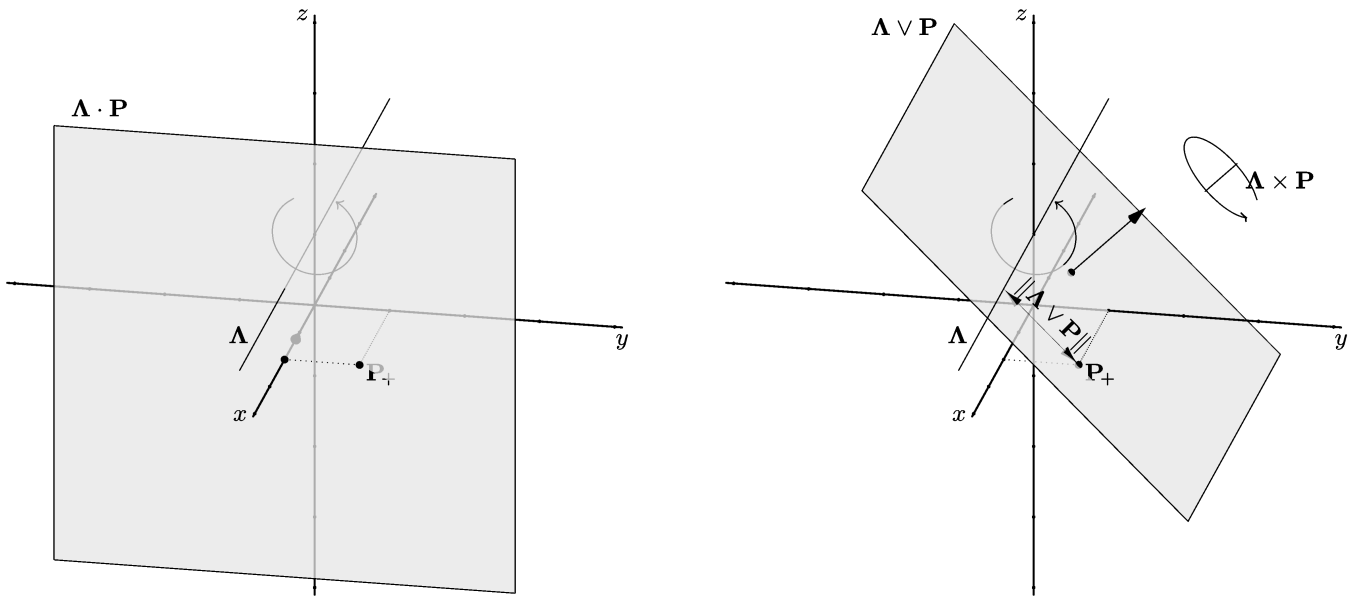
\hspace{-1cm}
\begin{subfloatenv}{}
\begin{asy}
import Figure3D;
Figure3D f = Figure3D();

var R = Point(1,1,0,1);
var Q = Point(1,0,0,1);
var L = join(R,Q);

var P = Point(1,2,1,0);

f.line(L, size=5,"$\mathbf{\Lambda}$", position=0.9, align=(0,-1,0.5));
f.point(P, "$\textbf{P}$",align=(0,1,-0.5));

f.plane(dot(L,P),"$\mathbf{\Lambda}\cdot\textbf{P}$", align=(0,2.5,2));

//var K = join(P,wedge(L,dot(P,L)));
//triple n = (K.p23,K.p31,K.p12);
//f.plane(dot(P,L),u=n,"$\textbf{P}\cdot\mathbf{\Lambda}$", align=(0,-2.5,-1));

\end{asy}
\end{subfloatenv}\hfill%
\begin{subfloatenv}{}
\begin{asy}
import Figure3D;
Figure3D f = Figure3D();

var R = Point(1,1,0,1);
var Q = Point(1,0,0,1);
var L = join(R,Q);

var P = Point(1,2,1,0);

f.line(L, size=5,"$\mathbf{\Lambda}$", position=0.9, align=(0,-1,0.5));
f.point(P, "$\textbf{P}$",align=(0,1,-0.5));

var K = join(P,wedge(L,dot(P,L)));
triple n = (K.p23,K.p31,K.p12);

f.plane(join(L,P),u=n,"$\mathbf{\Lambda}\vee\textbf{P}$", align=(0,-2.5,-0.5));
f.point_at_infinity(cross(L,P),centre=(0,2.5,2),"$\mathbf{\Lambda}\times\textbf{P}$", align=(0,3.5,0));

var Q = dot(P,L)/L;
f.stretch(P, Q, "$\lVert\mathbf{\Lambda}\vee\textbf{P}\rVert$", label_angle=-45,align=(0.4,0.4), shift_distance=0);

\end{asy}
\end{subfloatenv}
\caption{Points and lines in \E{3}}
\label{basic lines and points in E3}
\end{figure}

The dot product \(\mb{\Lambda}\cdot\tb{P}\) of a point \(\tb{P}\) and a line \(\mb{\Lambda}\)
is a plane that passes through \(\tb{P}\) and is perpendicular to \(\mb{\Lambda}\).
Moreover, if both the point and the line are normalised, then the resulting plane is also normalised.
For instance, the dot product of \(\tb{P}=\e_{123}+2\e_{320}+\e_{130}\)
and \(\mb{\Lambda}=-\e_{20}+\e_{23}\) yields \(\mb{\Lambda}\cdot\tb{P}=2\e_0-\e_1\) (see Figure~\ref{basic lines and points in E3}(a)).
The join \(\tb{P}\vee\mb{\Lambda}\) of a point \(\tb{P}\) and a line \(\mb{\Lambda}\) is a plane that passes through both the point and the line.
The distance from \(\tb{P}\) to \(\mb{\Lambda}\) is given by \(\norm{\mb{\Lambda}\vee\tb{P}}\)
if both \(\tb{P}\) and  \(\mb{\Lambda}\) are normalised.
This is illustrated in Figure~\ref{basic lines and points in E3}(b), where  \(\mb{\Lambda}\vee\tb{P}=-\e_0+\e_2+\e_3\). 
Since both the point and the line are normalised, 
the distance between them equals \(\norm{\mb{\Lambda}\vee\tb{P}}=\sqrt{1^2+1^2}=\sqrt{2}\).
Recall that \(\mb{\Lambda}\times\tb{P}=-(\mb{\Lambda}\vee\tb{P})\I\) for any bivector \(\mb{\Lambda}\).
So, the commutator \(\mb{\Lambda}\times\tb{P}\) is a point at infinity
 in the direction perpendicular to the plane \(\mb{\Lambda}\vee\tb{P}\).
The planes \(\mb{\Lambda}\cdot\tb{P}\) and \(\mb{\Lambda}\vee\tb{P}\) are perpendicular to each other
and their intersection \(\mb{\Lambda}_P=(\mb{\Lambda}\cdot\tb{P})\wedge(\mb{\Lambda}\vee\tb{P})\)
gives a line passing through \(\tb{P}\) and perpendicular to \(\mb{\Lambda}\).
For the point and the line used in Figure~\ref{basic lines and points in E3}, 
I get \(\mb{\Lambda}_P=\e_{10}-2\e_{20}-2\e_{30}+\e_{31}-\e_{12}\).

The join \(\tb{P}\vee\tb{Q}\vee\tb{R}\) of three points \(\tb{P}\), \(\tb{Q}\), and \(\tb{R}\) is a plane that passes through them. 
For instance, the join of the normalised right-handed points \(\tb{P}=\e_{123}\), 
\(\tb{Q}=\e_{123}+\e_{320}\), and \(\tb{R}=\e_{123}+\e_{130}\) located 
at \((0,0,0)\), \((1,0,0)\), and \((0,1,0)\), respectively,
yields \(-\e_3\), i.e.\ the plane \(z=0\) with the downward orientation.
The area of the triangle defined by the points is given by \(\tfrac{1}{2}\norm{\tb{P}\vee\tb{Q}\vee\tb{R}}\) if the points are normalised.
The join of four points is a scalar and \(\tfrac{1}{3!}|\tb{P}\vee\tb{Q}\vee\tb{R}\vee\tb{S}|\) gives 
the volume of a 3-simplex\footnote{A 3-simplex is a tetrahedron. In general, \(n\)-simplex is a generalisation of the concept of triangle
to \(n\) dimensions.}
defined by points \(\tb{P}\), \(\tb{Q}\), \(\tb{R}\), and \(\tb{S}\) if they are normalised.

On the other hand,  \(\tb{a}\wedge\tb{b}\wedge\tb{c}\) is a point where the planes 
\(\tb{a}\), \(\tb{b}\), and \(\tb{c}\) intersect.
For instance, \(\e_1\wedge\e_2\wedge\e_3=\e_{123}\) is a normalised point at the origin, with the right-handed orientation.
\(\tb{a}\wedge\mb{\Lambda}\) yields a point where  \(\mb{\Lambda}\) intersects \(\tb{a}\).
For instance, \(\e_1\wedge\e_{23}\) is a point where the \(x\)-axis, dually represented by \(\e_{23}\),
intersects the plane \(x=0\), dually represented by \(\e_1\).
If the line is parallel to the plane, a point at infinity is obtained.

\begin{figure}[ht!]
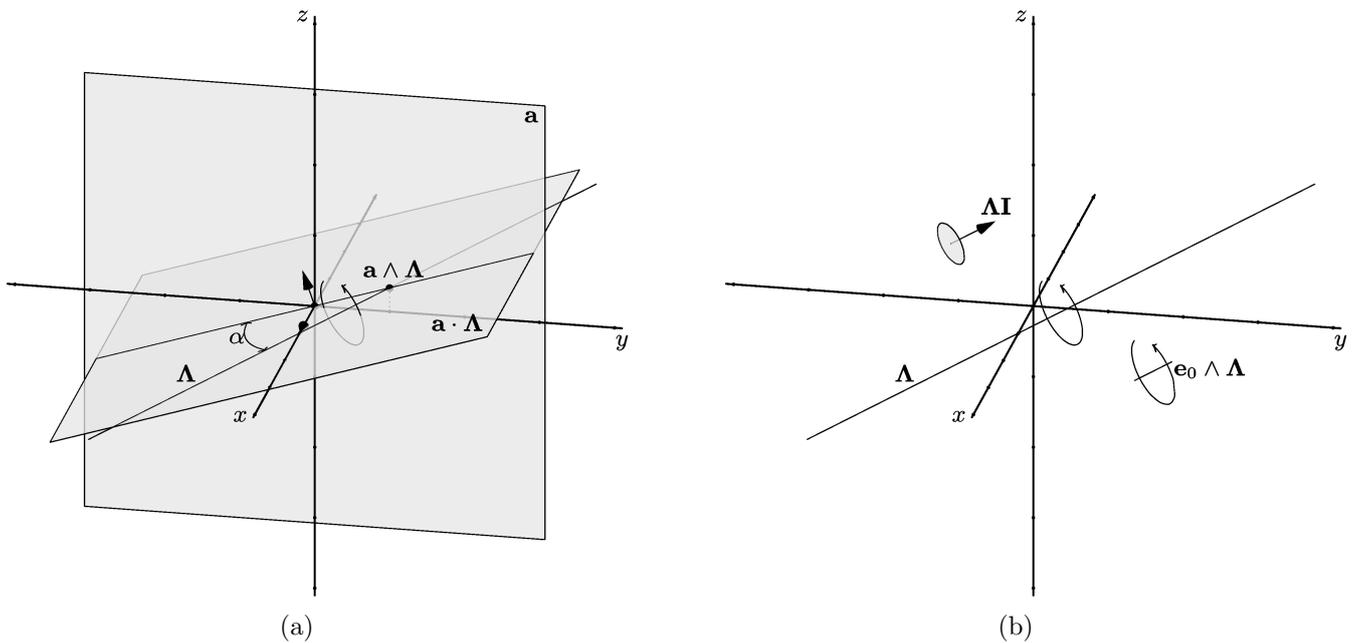
\hspace{-1cm}
\begin{subfloatenv}{ }
\begin{asy}
import Figure3D;
Figure3D f = Figure3D();

var P = Point(1,1,0,0);
var Q = Point(1,0,1,1/3);
var L = join(P,Q)/2;

var a = Plane(0,1,0,0);

f.plane(a, "$\textbf{a}$",align=(0,-1,-1));
f.line(L, "$\mathbf{\Lambda}$", position=0.8, align=(0,-0.25,1));

f.plane(dot(a,L),"$\textbf{a}\cdot\mathbf{\Lambda}$", align=(0.5,-0.7,1));
f.point(wedge(a,L),"$\textbf{a}\wedge\mathbf{\Lambda}$", align=(1.5,0.5,2.5),draw_orientation=false);

var c = totriple(wedge(a,L));
var v1 = (L.p23,L.p31,L.p12)/norm(L)*2;
var K = dot(L,a)/a;
var v2 = (K.p23,K.p31,K.p12)/norm(K)*2;
f.arc(c,v1,v2,"$\alpha$",align=(1,-0.5,0.5));

\end{asy}
\end{subfloatenv}
\begin{subfloatenv}{ }
\begin{asy}
import Figure3D;
Figure3D f = Figure3D();

var P = Point(1,1,0,0);
var Q = Point(1,0,1,1/3);
var L = join(P,Q)/2;

f.line(L, "$\mathbf{\Lambda}$", position=0.8, align=(0,-0.25,1));

f.point_at_infinity(wedge(e_0,L),"$\textbf{e}_0\wedge\mathbf{\Lambda}$", centre=(2,2,0), align=(0,3,0.5));
f.line_at_infinity(L*I,"$\mathbf{\Lambda}\textbf{I}$", centre=(-2,-1.5,0),align=(0,0,1));

\end{asy}
\end{subfloatenv}
\caption{Lines and planes in \E{3}}
\label{basic lines and planes in E3}
\end{figure}

The dot product \(\tb{a}\cdot\mb{\Lambda}\) of a plane \(\tb{a}\) and a line \(\mb{\Lambda}\)
is a plane that contains  \(\mb{\Lambda}\) and is perpendicular to \(\tb{a}\).
If both the plane and the line are normalised, then \(\norm{\tb{a}\cdot\mb{\Lambda}}=\cos{\alpha}\),
where \(\alpha\) is the angle between the line and the plane (\(\tb{a}\cdot\mb{\Lambda}=0\) if \(\mb{\Lambda}\) is perpendicular to  \(\tb{a}\)).
The top-down view of a line is a sheaf of planes.
The dot product picks out one of the sheaf's planes, the one that is perpendicular to  \(\tb{a}\).
This is illustrated in Figure~\ref{basic lines and planes in E3}(a), where \(\tb{a}=\e_1\) and 
\(\mb{\Lambda}=\tfrac{1}{2}(-\tfrac{1}{3}\e_{20}+\e_{30}+\e_{23}-\e_{31}-\tfrac{1}{3}\e_{12})\).
The trivector \(\e_0\wedge\mb{\Lambda}\) is a point at infinity where \(\mb{\Lambda}\) intersects  the plane at infinity
(see Figure~\ref{basic lines and planes in E3}(b)).
\(\mb{\Lambda}\I\) is a line at infinity whose stack consists of the planes perpendicular to \(\mb{\Lambda}\).

\begin{figure}[t]
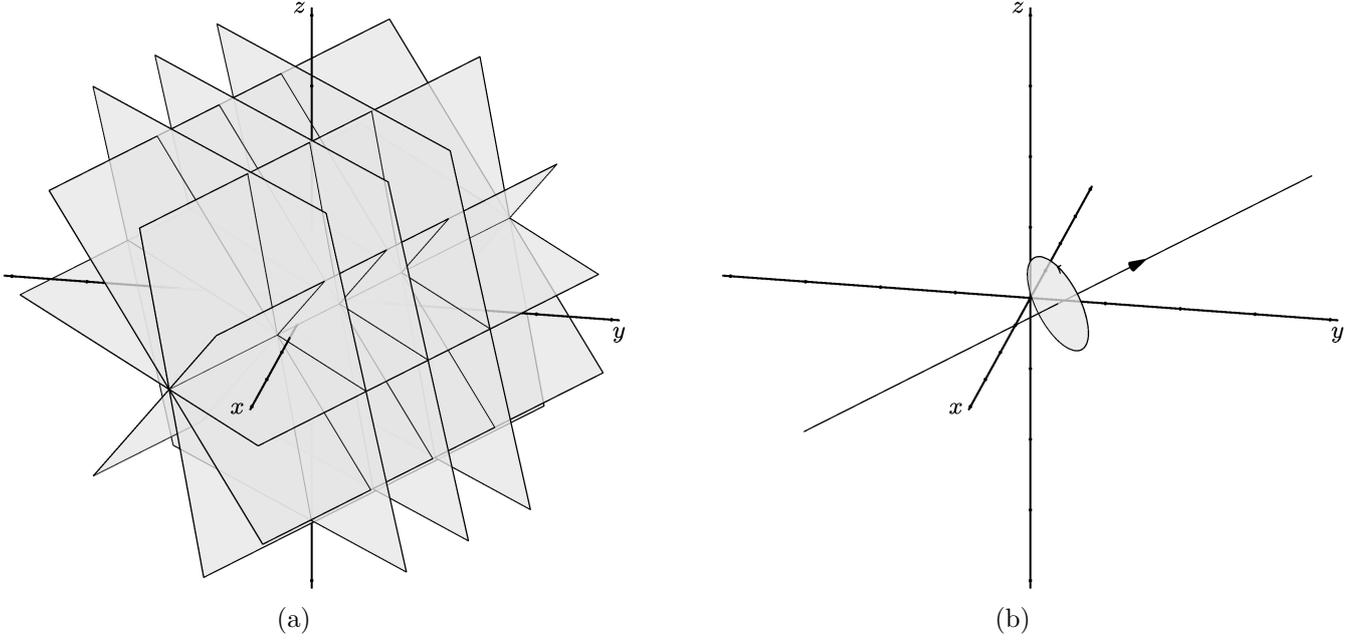
\hspace{-1cm}
\begin{subfloatenv}{ }
\begin{asy}
import Figure3D;
Figure3D f = Figure3D();

var P = Point(1,1,0,0);
var Q = Point(1,0,1,1/3);
var L = join(P,Q);

var C = L + I*L;

var L = (1 - 0.5*join(C,C)/dot(C,C)*I)*C;
var Lp = 0.5*join(C,C)/dot(C,C)*I*C;

L = L/norm(L);
triple n = (L.p23,L.p31,L.p12);
var a = join(L,Point(1,0,0,1));
int N = 4;
for(int i: sequence(N)) { real alpha = i/N*pi; var Rot = cos(alpha/2) - sin(alpha/2)*L; f.plane(Rot*a/Rot,u=n,draw_orientation=false,size=2.75);}

var b = join(O,L*I);
int M = 1;
for(int i: sequence(-M,M)) { real alpha = i/M*pi; var Tr = 1 - i/2*Lp/1.45; write(Tr); f.plane(Tr*b/Tr,draw_orientation=false,size=2.75);}

\end{asy}
\end{subfloatenv}
\begin{subfloatenv}{ }
\begin{asy}
import Figure3D;
Figure3D f = Figure3D();

var P = Point(1,1,0,0);
var Q = Point(1,0,1,1/3);
var L = join(P,Q);

var cL = L + I*L;

f.compound_line(cL);

\end{asy}
\end{subfloatenv}
\caption{Visualising a non-simple bivector in \E{3}}
\label{non-simple bivector in E3}
\end{figure}

Any non-simple bivector \(\mb{\Lambda}\) can be uniquely decomposed into a line \(\mb{\Lambda}_0\)
that passes through the origin and a line at infinity \(\mb{\Lambda}_\infty=\e_0\wedge\tb{a}\).
Shifting the focus to the dual model space for a moment, 
I note that the bivector \(\mb{\Lambda}_0\) and the vector \(\tb{a}\) lie in the 3-dimensional subspace of \R{4*} where \(d=0\).
Using Euclidean metric I can decompose \(\tb{a}\) into two components, \(\tb{a}=\tb{a}_\parallel+\tb{a}_\perp\), where
\(\tb{a}_\parallel\) lies in the plane of the bivector \(\mb{\Lambda}_0\) 
and \(\tb{a}_\perp\) is perpendicular to \(\mb{\Lambda}_0\). 
This defines two lines at infinity in \E{3}, 
\(\mb{\Lambda}_{\infty\parallel}=\e_0\wedge\tb{a}_\parallel\) and \(\mb{\Lambda}_{\infty\perp}=\e_0\wedge\tb{a}_\perp\).
Then, \(\mb{\Lambda}_1=\mb{\Lambda}_0+\mb{\Lambda}_{\infty\parallel}\) is a finite line parallel to \(\mb{\Lambda}_0\)
and shifted from the origin in accord with  \(\tb{a}_\parallel\), whereas
\(\mb{\Lambda}_2=\mb{\Lambda}_{\infty\perp}\) is a component of \(\mb{\Lambda}_\infty\) 
that cannot be combined with \(\mb{\Lambda}_0\) to produce a line.
Planes that belong to the stack dually represented by  \(\mb{\Lambda}_{\infty\perp}\) are perpendicular to the line \(\mb{\Lambda}_1\).
Thus, in Euclidean space, a non-simple bivector \(\mb{\Lambda}=\mb{\Lambda}_1+\mb{\Lambda}_2\) 
can be visualised as a compound object consisting of a finite line
and a line at infinity perpendicular to it in the sense explained above (see Figure~\ref{non-simple bivector in E3}(a)).
The lines \(\mb{\Lambda}_1\) and \(\mb{\Lambda}_2\) are called the axes of a non-simple bivector \(\mb{\Lambda}\),
the finite axis and the axis at infinity.

Since \(\I\mb{\Lambda}_0\) is a line at infinity whose stack consists of the planes perpendicular to \(\mb{\Lambda}_0\),
I can assume \(\mb{\Lambda}_{\infty\perp}=a\I\mb{\Lambda}_0\) for some scalar factor \(a\).
To find \(a\), I compute 
\(\mb{\Lambda}\vee\mb{\Lambda}
=
(\mb{\Lambda}_0+\mb{\Lambda}_{\infty\parallel}+\mb{\Lambda}_{\infty\perp})\vee(\mb{\Lambda}_0
+\mb{\Lambda}_{\infty\parallel}+\mb{\Lambda}_{\infty\perp})
=2\mb{\Lambda}_0\vee\mb{\Lambda}_{\infty\perp}+2\mb{\Lambda}_0\vee\mb{\Lambda}_{\infty\parallel}
+2\mb{\Lambda}_{\infty\parallel}\vee\mb{\Lambda}_{\infty\perp}\),
where the last two terms vanish 
since \(\mb{\Lambda}_0+\mb{\Lambda}_{\infty\parallel}\) and \(\mb{\Lambda}_{\infty\parallel}+\mb{\Lambda}_{\infty\perp}\) 
are simple bivectors.
So, 
\(\mb{\Lambda}\vee\mb{\Lambda}=2\mb{\Lambda}_0\vee\mb{\Lambda}_{\infty\perp}=
2\mb{\Lambda}_0\vee(a\I\mb{\Lambda}_0)=2a\mb{\Lambda}_0\cdot\mb{\Lambda}_0=
2a\mb{\Lambda}\cdot\mb{\Lambda}\) and
\begin{equation}
\mb{\Lambda}_1=(1-a\I)\mb{\Lambda}, \quad \mb{\Lambda}_2=a\I\mb{\Lambda},
\quad\textrm{and } 
a =\frac{\mb{\Lambda}\vee\mb{\Lambda}}{2\mb{\Lambda}\cdot\mb{\Lambda}}.
\end{equation}
A non-simple bivector is visualised in Figure~\ref{non-simple bivector in E3}(b) by displaying its axes.

\begin{figure}[t]
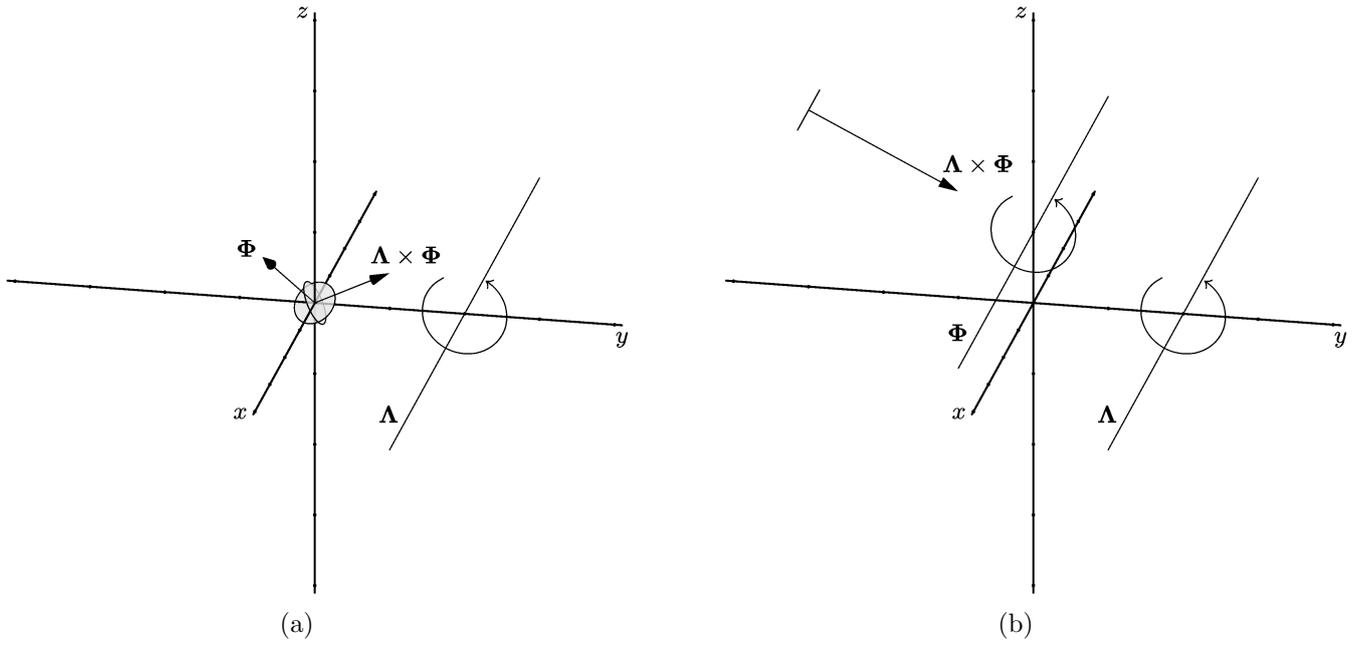
\hspace{-1cm}
\begin{subfloatenv}{ }
\begin{asy}
import Figure3D;
Figure3D f = Figure3D();

var P = Point(1,0,2,0);
var Q = Point(1,-2,2,0);
var L = join(P,Q);
L /= norm(L);

var K = e_0*(2e_1-e_2+2*e_3)/2;

f.line(L, size=5,"$\mathbf{\Lambda}$", position=0.9, align=(0,-1,0.5));
f.line_at_infinity(K,"$\mathbf{\Phi}$", align=(0,-1,0.5));

f.line_at_infinity(cross(L,K),centre=(0,0,0),"$\mathbf{\Lambda}\times\mathbf{\Phi}$",align=(0,0.5,1.5));

\end{asy}
\end{subfloatenv}
\begin{subfloatenv}{ }
\begin{asy}
import Figure3D;
Figure3D f = Figure3D();

var P = Point(1,0,2,0);
var Q = Point(1,-2,2,0);
var L = join(P,Q);
L /= norm(L);

var R = Point(1,0,0,1);
var Q = Point(1,1,0,1);
var K = -join(R,Q)/4;
K = K/norm(K);

f.line(L, size=5,"$\mathbf{\Lambda}$", position=0.9, align=(0,-1,0.5));
f.line(K, size=5,"$\mathbf{\Phi}$", position=0.9, align=(0,-1,0.5));

f.line_at_infinity(cross(L,K),centre=(0,-3,2.5),"$\mathbf{\Lambda}\times\mathbf{\Phi}$",align=(0,1,3));

\end{asy}
\end{subfloatenv}
\caption{Commutator of lines in \E{3}}
\label{basic lines 1 in E3}
\end{figure}

\begin{figure}[h!]
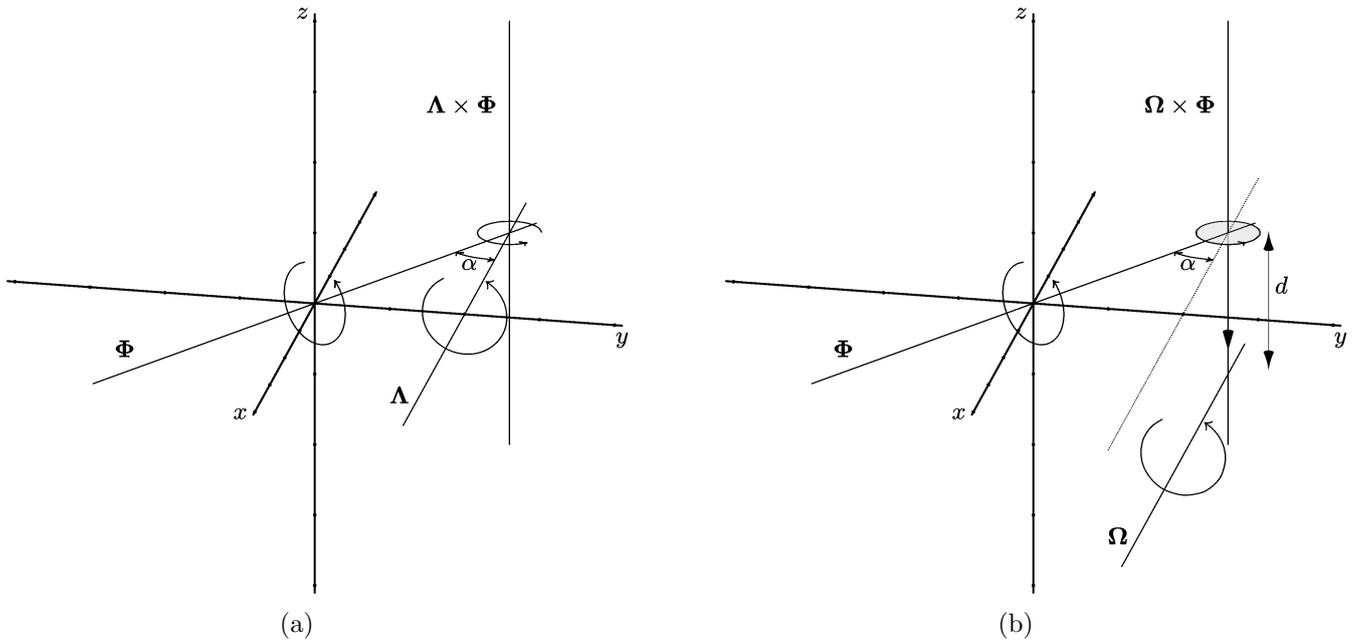
\hspace{-1cm}
\begin{subfloatenv}{ }
\begin{asy}
import Figure3D;
Figure3D f = Figure3D();

var P = Point(1,0,2,0);
var Q = Point(1,-2,2,0);
var L = join(P,Q);
L /= norm(L);

var R = Point(1,0,0,0);
var Q = Point(1,-3,2,0);
var K = join(R,Q)/4;
K = K/norm(K);

f.line(L,"$\mathbf{\Lambda}$", position=0.9, align=(0,-1,0.5));
f.line(K,"$\mathbf{\Phi}$", position=0.9, align=(0,-1,1.5),orientation_theta=-10);

f.line(cross(L,K),size=3,"$\mathbf{\Lambda}\times\mathbf{\Phi}$",align=(0,-2,0),position=0.8,orientation_theta=-15);

triple centre = (-3,2,0);
var v1 = 1*(L.p23,L.p31,L.p12);
var v2 = 1*(K.p23,K.p31,K.p12);
f.arc(centre, v1, v2,"$\alpha$",  (2,0,0));

\end{asy}
\end{subfloatenv}
\begin{subfloatenv}{ }
\begin{asy}
import Figure3D;
Figure3D f = Figure3D();

var P = Point(1,0,2,-2);
var Q = Point(1,-2,2,-2);
var L = join(P,Q);
L /= norm(L);

var R = Point(1,0,0,0);
var Q = Point(1,-3,2,0);
var K = join(R,Q)/4;
K = K/norm(K);

f.line(L,"$\mathbf{\Omega}$", position=0.9, align=(0,-1,0.5));
f.line(K,"$\mathbf{\Phi}$", position=0.9, align=(0,-1,1.5),orientation_theta=-10);

f.compound_line(cross(L,K),size=3,"$\mathbf{\Omega}\times\mathbf{\Phi}$",align=(0,-2,0),position=0.8,orientation_theta=-15);

var T = exp(-0.5*2*e_0*e_3);
var M = T*L/T;
f.line(M, size=5, draw_orientation=false, pen=dotted);

triple centre = (-3,2,0);
var v1 = 1*(M.p23,M.p31,M.p12);
var v2 = 1*(K.p23,K.p31,K.p12);
f.arc(centre, v1, v2,"$\alpha$",  (2,0,0));

f.stretch(Point(1,-3,2,0), Point(1,-3,2,-2), "$d$", label_angle=0,align=(1,1), shift_distance=0.1,shift_angle=-90);

\end{asy}
\end{subfloatenv}
\caption{Commutator of lines in \E{3} (continued)}
\label{basic lines in E3}
\end{figure}

The configuration arising from the cross product of a finite line and a line at infinity
is shown in Figure~\ref{basic lines 1 in E3}(a) for \(\mb{\Lambda}=2\e_{20}+\e_{23}\)
and \(\mb{\Phi}=-\e_{10} +\tfrac{1}{2}\e_{20} -\e_{30}\), which gives a line at infinity
\(\mb{\Lambda}\times\mb{\Phi}=-\e_{20}-\tfrac{1}{2}\e_{30}\).

If  \(\mb{\Lambda}\) and \(\mb{\Phi}\) are parallel as shown in Figure~\ref{basic lines 1 in E3}(b),
the commutator is a line at infinity, with the orientation pointing from \(\mb{\Phi}\) 
to \(\mb{\Lambda}\) (or \(\mb{\Lambda}\) to \(\mb{\Phi}\) 
if the lines are anti-parallel).
For \(\mb{\Lambda}=2\e_{20}+\e_{23}\) and \(\mb{\Phi}=-\e_{20}+\e_{23}\),
I get \(\mb{\Lambda}\times\mb{\Phi}=-2\e_{20}+\e_{30}\).
The distance between the lines is given by \(\norm{\e_{123}\vee(\mb{\Lambda}\times\mb{\Phi})}\) if the lines are normalised.
Note that the join of parallel or anti-parallel lines is zero.

For normalised lines \(\mb{\Lambda}\) and 
\(\mb{\Phi}\),
the geometric product is as follows:
\begin{equation}
\mb{\Lambda}\mb{\Phi}
=
-\cos{\alpha}+\mb{\Lambda}\times\mb{\Phi} \pm \I \,r\sin{\alpha},
\end{equation}
where \(\alpha=\arccos{(-\mb{\Lambda}\cdot\mb{\Phi})}\) is the angle between the lines (more precisely, the angle between their bottom-up orientation vectors), 
\(r\) is the distance between the lines (the sign depends on the relative orientation of the lines), 
and the commutator \(\mb{\Lambda}\times\mb{\Phi}\) is in general a non-simple bivector
whose finite axis passes through both lines and is perpendicular to both of them.
Since  \(\mb{\Lambda}\wedge\mb{\Phi}=(\mb{\Lambda}\vee\mb{\Phi})\I\),
the distance \(r\) can be computed from \(|\mb{\Lambda}\vee\mb{\Phi}|=r\sin{\alpha}\)
if \(\alpha\ne0,\pi\).

If \(\mb{\Lambda}\) and \(\mb{\Phi}\) intersect, 
the commutator is a finite line that is perpendicular to both \(\mb{\Lambda}\) and \(\mb{\Phi}\), 
as shown in Figure~\ref{basic lines in E3}(a), where 
 \(\mb{\Lambda}=2\e_{30}+\e_{23}\),
\(\mb{\Phi}=\tfrac{1}{\sqrt{13}}(3\e_{23}-2\e_{31})\) and, therefore,
\(\alpha=\arccos{\tfrac{3}{\sqrt{13}}}\approx34^\circ\), \(d=0\), and 
\(\mb{\Lambda}\times\mb{\Phi}
=\tfrac{1}{\sqrt{13}}(-4\e_{10}-6\e_{20}+2\e_{12})\) is a simple bivector.

The commutator of two finite non-intersecting lines is a non-simple bivector.
For example, \(\mb{\Omega}=2\e_{20}+2\e_{30}+\e_{23}\)  and \(\mb{\Phi}=\tfrac{1}{\sqrt{13}}(3\e_{23}-2\e_{31})\)
give \(\mb{\Omega}\times\mb{\Phi}=\tfrac{1}{\sqrt{13}}(-4\e_{10}-6\e_{20}+6\e_{30}+2\e_{12})\), which
is not simple (see Figure~\ref{basic lines in E3}(b)). 
\(\mb{\Omega}\) can be obtained  from  \(\mb{\Lambda}=2\e_{30}+\e_{23}\), shown in Figure~\ref{basic lines in E3}(b)
with a dotted line,
by shifting it downward by 2.
So, the angle \(\alpha\) between  \(\mb{\Omega}\) and \(\mb{\Phi}\) is the same as the angle between \(\mb{\Lambda}\) and \(\mb{\Phi}\).
I get \(r=2\) for the distance as expected. 
The commutator splits into a finite line
and a line at infinity perpendicular to it:
\(\mb{\Omega}\times\mb{\Phi}=\tfrac{1}{\sqrt{13}}(-4\e_{10}-6\e_{20}+2\e_{12})+
\tfrac{6}{\sqrt{13}}\e_{30}\).

\begin{figure}[t!]
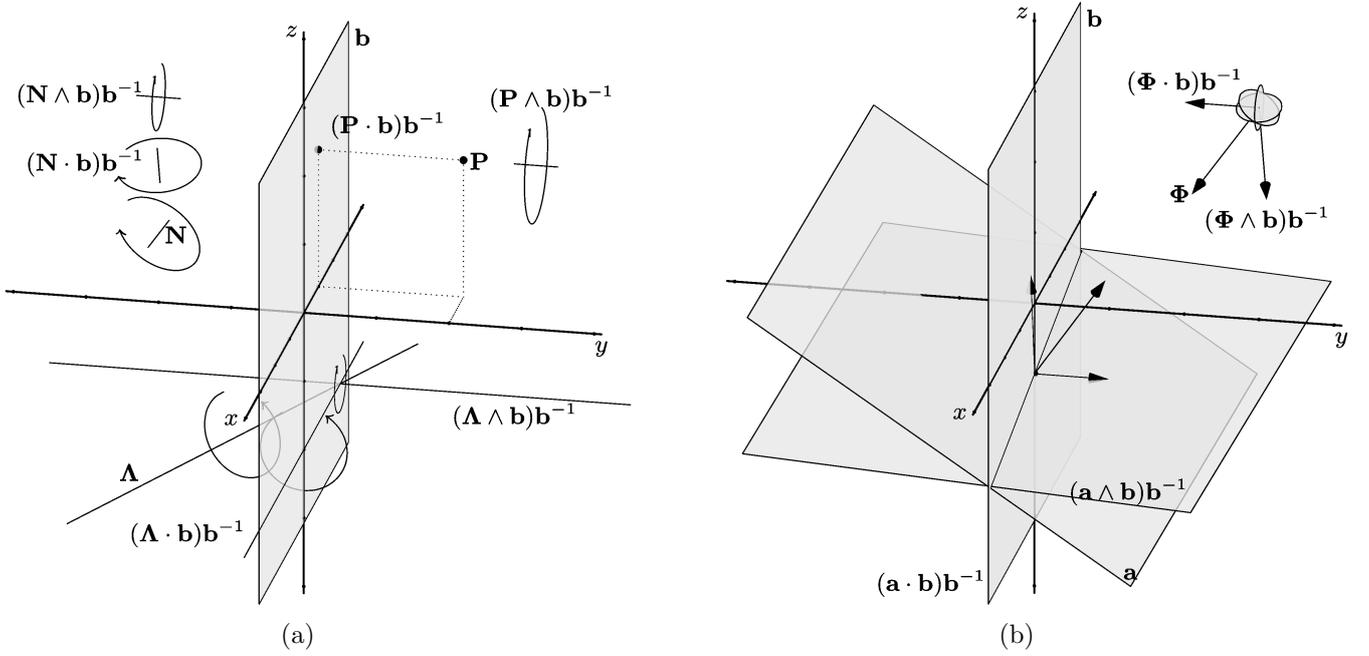
\hspace{-1cm}
\begin{subfloatenv}{ }
\begin{asy}
import Figure3D;
Figure3D f = Figure3D();

var a = join(O, e_0*e_2);
f.plane(a,"$\textbf{b}$", draw_orientation=false, align=(0,1,-1));

var P = Point(1,-1,2,2);
f.point(P,"$\textbf{P}$", draw_helper_lines=true,draw_orientation=false);
f.point(dot(P,a)/a, "$(\textbf{P}\cdot\textbf{b})\textbf{b}^{-1}$", align=(0,2,2),draw_orientation=false);
draw((-1,0,2)--(-1,2,2),dotted);
f.point_at_infinity(wedge(P,a)/a, (-1,3,2), "$(\textbf{P}\wedge\textbf{b})\textbf{b}^{-1}$", align=(0,1,8));

var N = Point(0,1,2,3)/3;
f.point_at_infinity(N,(0,-2,1),"$\textbf{N}$",orientation_theta=-10);
f.point_at_infinity(dot(N,a)/a,(0,-2,2),"$(\textbf{N}\cdot\textbf{b})\textbf{b}^{-1}$",align=(0,-2.5,0));
f.point_at_infinity(wedge(N,a)/a,(0,-2,3),"$(\textbf{N}\wedge\textbf{b})\textbf{b}^{-1}$",align=(0,-2.5,0));

var L = join(Point(1,2,-2,-2),Point(1,-2.5,0,-2))/4;

f.line(L,"$\mathbf{\Lambda}$",position=0.8,align=(0,-0.5,1));
f.line(dot(L,a)/a, "$(\mathbf{\Lambda}\cdot\textbf{b})\textbf{b}^{-1}$",position=0.9,align=(0,-2,0));
f.line(wedge(L,a)/a, size=4,"$(\mathbf{\Lambda}\wedge\textbf{b})\textbf{b}^{-1}$",align=(0,0,-1));

\end{asy}
\end{subfloatenv}
\begin{subfloatenv}{ }
\begin{asy}
import Figure3D;
Figure3D f = Figure3D();

var a = join(O, e_0*e_2);
f.plane(a,"$\textbf{b}$", draw_orientation=false, align=(0,1,-1));

var K = Line(1,2,3,0,0,0)/2;
f.line_at_infinity(K,(0,3,3),"$\mathbf{\Phi}$",align=(0,-0.5,0));
f.line_at_infinity(dot(K,a)/a,(0,3,3), "$(\mathbf{\Phi}\cdot\textbf{b})\textbf{b}^{-1}$",align=(0,0,1));
f.line_at_infinity(wedge(K,a)/a,(0,3,3), "$(\mathbf{\Phi}\wedge\textbf{b})\textbf{b}^{-1}$",align=(0,0,-1));

var b = join(Point(1,0,0,-1),e_0*(0.25*e_1+e_2+1.5*e_3));
triple u = cross((b.a,b.b,b.c),(1,0,0));
f.plane(b,u,"$\textbf{a}$",align=(0,0,1),draw_orientation=true,O=(0,0,-1));
f.plane(dot(b,a)/a,u,"$(\textbf{a}\cdot\textbf{b})\textbf{b}^{-1}$",align=(0,-1,1),draw_orientation=true,O=(0,0,-1),intersect=false);
f.plane(wedge(b,a)/a,u,"$(\textbf{a}\wedge\textbf{b})\textbf{b}^{-1}$",align=(0,-1,1),draw_orientation=true,O=(0,0,-1));

\end{asy}
\end{subfloatenv}
\caption{Projection on and rejection by a plane in \E{3}}
\label{projection on plane in E3}
\end{figure}

\begin{figure}[h!]
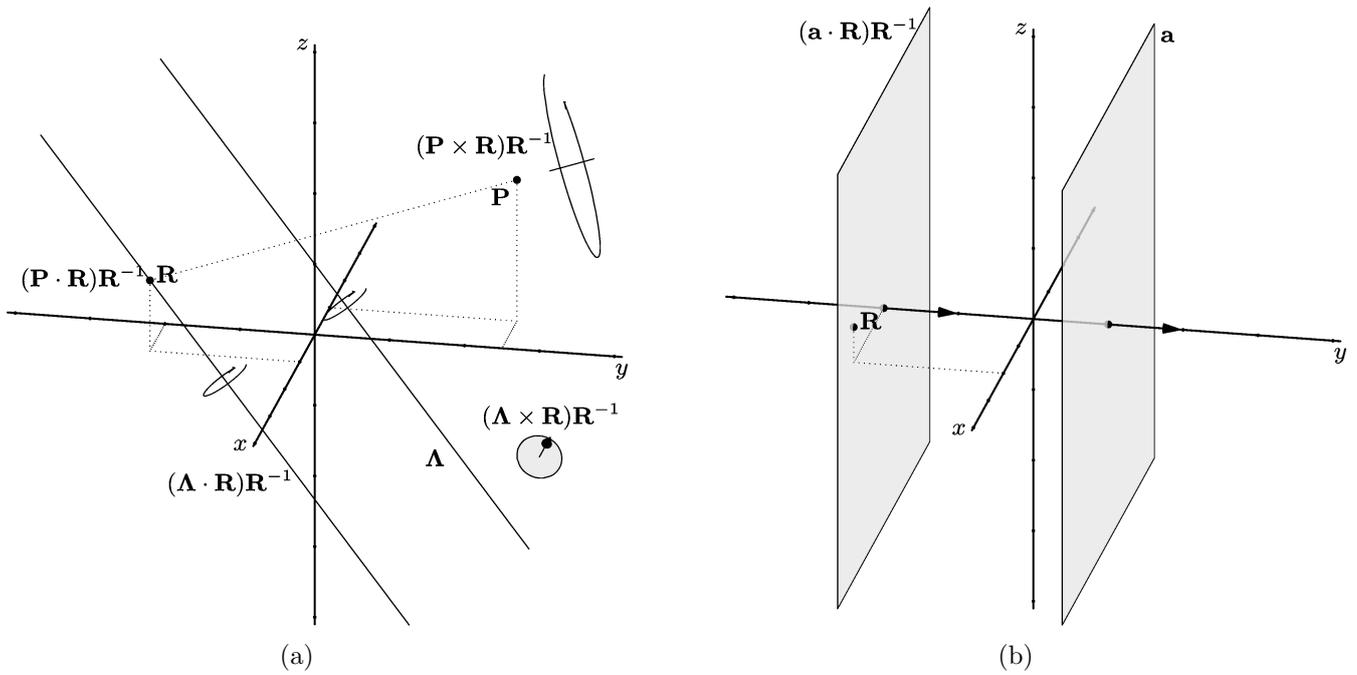
\hspace{-1cm}
\begin{subfloatenv}{ }
\begin{asy}
import Figure3D;
Figure3D f = Figure3D();

var R = Point(1,1,-2,1);
f.point(R,"$\textbf{R}$", draw_helper_lines=true,draw_orientation=false,align=(0,1,0.5));

var P = Point(1,-1,2.5,2);
f.point(P,"$\textbf{P}$", draw_helper_lines=true,draw_orientation=false,align=(1,-1,-1));
f.point(dot(P,R)/R, "$(\textbf{P}\cdot\textbf{R})\textbf{R}^{-1}$", align=(0,-1,0),draw_orientation=false, draw_helper_lines=false);
draw(totriple(P)--totriple(dot(P,R)/R),dotted);
triple s = totriple(P + 0.15*(P-R));
f.point_at_infinity(cross(P,R)/R,s,"$(\textbf{P}\times\textbf{R})\textbf{R}^{-1}$",align=(0,-3,1.5));

var M = join(Point(1,0,0,1),Point(1,0.5,1,0))/4;

f.line(M,'$\mathbf{\Lambda}$',align=(0.5,-1.5,-0.5));
f.line(dot(M,R)/R,"$(\mathbf{\Lambda}\cdot\textbf{R})\textbf{R}^{-1}$",position=0.3,align=(1,-1,0));
f.line_at_infinity(cross(M,R)/R,(0,3,-1.5),"$(\mathbf{\Lambda}\times\textbf{R})\textbf{R}^{-1}$",align=(0,0,1));

\end{asy}
\end{subfloatenv}
\begin{subfloatenv}{ }
\begin{asy}
import Figure3D;
Figure3D f = Figure3D();

var a = -e_0+e_2;
f.plane(a,"$\textbf{a}$", align=(0,1,-1));

var R = Point(1,2,-2,1/2);
f.point(R,"$\textbf{R}$", draw_helper_lines=true,draw_orientation=false,align=(0,1,0.5));

f.plane(dot(a,R)/R, "$(\textbf{a}\cdot\textbf{R})\textbf{R}^{-1}$", align=(0,-2,-2));

\end{asy}
\end{subfloatenv}
\caption{Projection on and rejection by a point in \E{3}}
\label{projection on point in E3}
\end{figure}

\begin{figure}[t!]
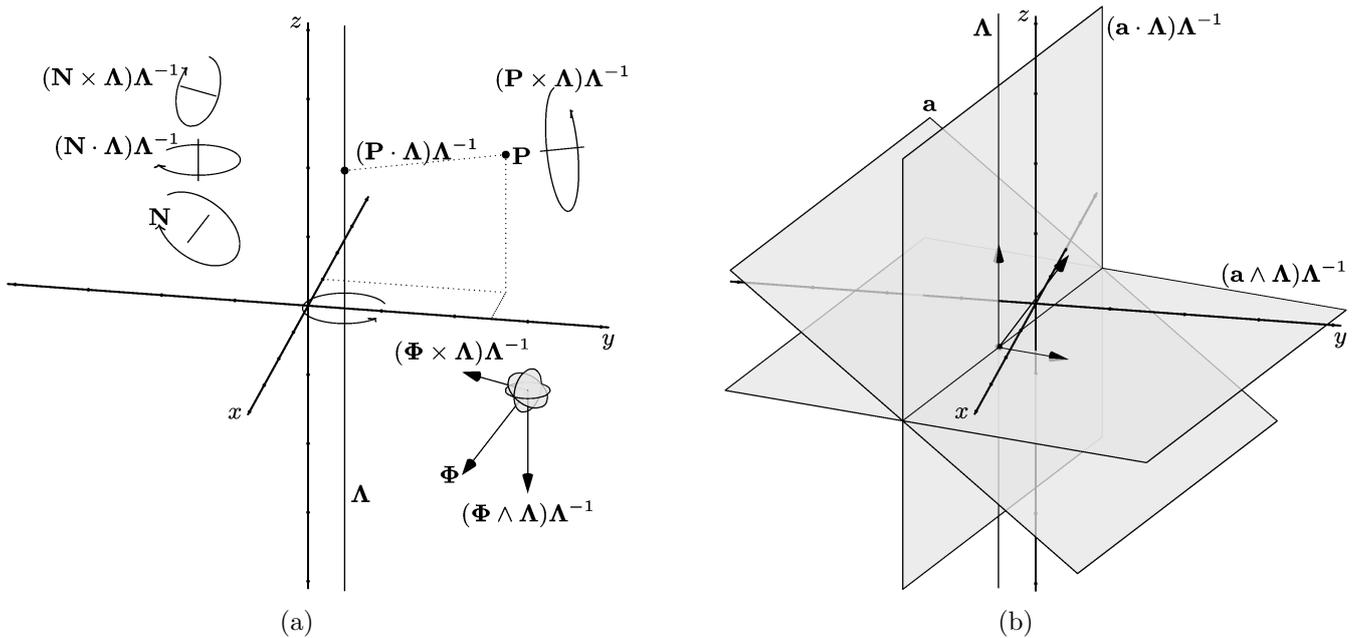
\hspace{-1cm}
\begin{subfloatenv}{ }
\begin{asy}
import Figure3D;
Figure3D f = Figure3D();

var L = (0.5*e_0-e_2)*e_1;
f.line(L,"$\mathbf{\Lambda}$", position=0.2, align=(0,1,-1));

var P = Point(1,-1,2.5,2);
f.point(P,"$\textbf{P}$", draw_helper_lines=true,draw_orientation=false);
f.point(dot(P,L)/L, "$(\textbf{P}\cdot\mathbf{\Lambda})\mathbf{\Lambda}^{-1}$", align=(0,1.8,1.5),draw_orientation=false, draw_helper_lines=false);
draw(totriple(P)--totriple(dot(P,L)/L),dotted);
triple s = (-1,2.5,2) + 0.35*((-1,2.5,2)-(0,0.5,2));
f.point_at_infinity(cross(P,L)/L,s,"$(\textbf{P}\times\mathbf{\Lambda})\mathbf{\Lambda}^{-1}$",align=(0,0,8.5));

var N = Point(0,1,2,3)/3;
f.point_at_infinity(N,(0,-1.5,1),"$\textbf{N}$",orientation_theta=-10,align=(0,-4,1.1));
f.point_at_infinity(dot(N,L)/L,(0,-1.5,2),"$(\textbf{N}\cdot\mathbf{\Lambda})\mathbf{\Lambda}^{-1}$",align=(0,-3,1));
f.point_at_infinity(cross(N,L)/L,(0,-1.5,3),"$(\textbf{N}\times\mathbf{\Lambda})\mathbf{\Lambda}^{-1}$",align=(0,-3,1));

var K = Line(1,2,3,0,0,0)/2;
f.line_at_infinity(K,(0,3,-1),"$\mathbf{\Phi}$",align=(0,-0.5,0));
f.line_at_infinity(cross(K,L)/L,(0,3,-1), "$(\mathbf{\Phi}\times\mathbf{\Lambda})\mathbf{\Lambda}^{-1}$",align=(0,0,1));
f.line_at_infinity(wedge(K,L)/L,(0,3,-1), "$(\mathbf{\Phi}\wedge\mathbf{\Lambda})\mathbf{\Lambda}^{-1}$",align=(0,0,-1));

\end{asy}
\end{subfloatenv}
\begin{subfloatenv}{ }
\begin{asy}
import Figure3D;
Figure3D f = Figure3D();

var L = (-0.5*e_0-e_2)*e_1;
f.line(L,"$\mathbf{\Lambda}$", position=0.95, draw_orientation=false, align=(0,-1,1));

var b = join(Point(1,0,0,-1),e_0*(0.25*e_1+e_2+1.5*e_3));
triple c = totriple(wedge(b,L));
f.plane(b,(0,0,1),"$\textbf{a}$",align=(0,0,1),draw_orientation=true,O=c);
f.plane(dot(b,L)/L,(0,0,1),"$(\textbf{a}\cdot\mathbf{\Lambda})\mathbf{\Lambda}^{-1}$",align=(0,1,-1),draw_orientation=true,O=c);

var K = wedge(b,dot(b,L)/L);
triple n = (K.p23,K.p31,K.p12);

f.plane(wedge(b,L)/L,n,"$(\textbf{a}\wedge\mathbf{\Lambda})\mathbf{\Lambda}^{-1}$",align=(2,-2,4),draw_orientation=true,O=c);

\end{asy}
\end{subfloatenv}
\caption{Projection on and rejection by a line in \E{3}}
\label{projection on line in E3}
\end{figure}

\begin{figure}[t!]
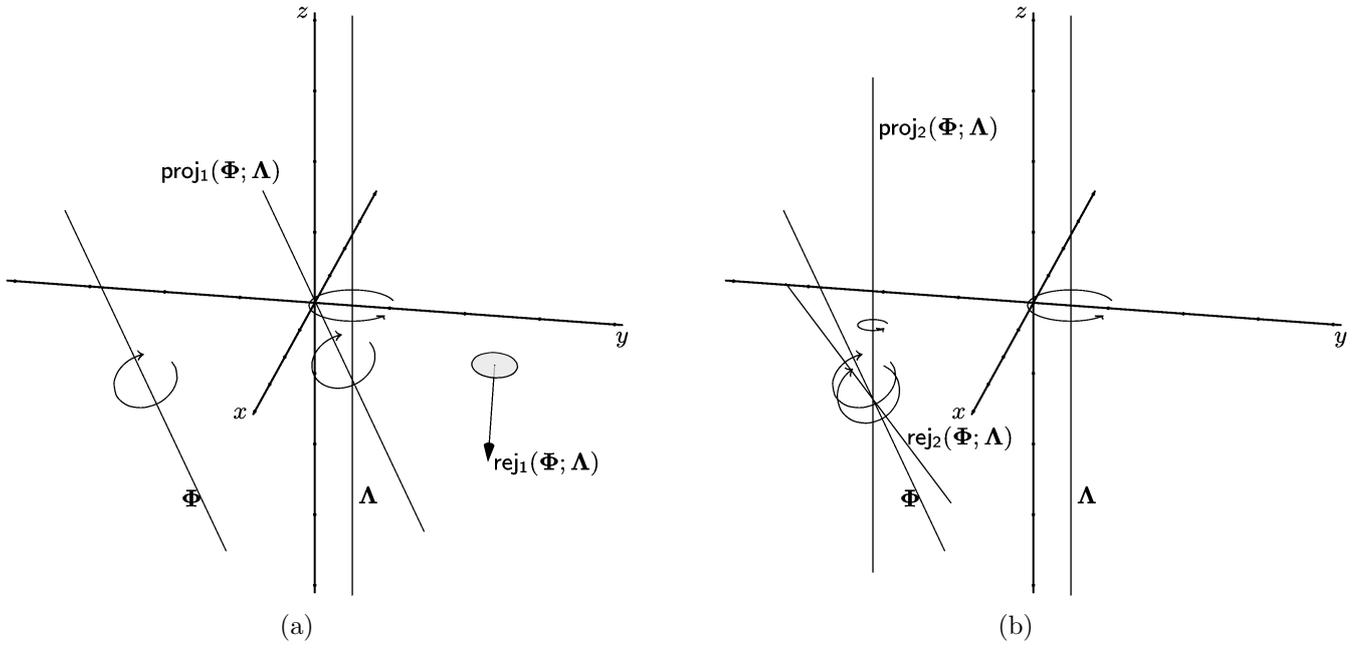
\hspace{-1cm}
\begin{subfloatenv}{ }
\begin{asy}
import Figure3D;
Figure3D f = Figure3D();

var L = (0.5*e_0-e_2)*e_1;
f.line(L,"$\mathbf{\Lambda}$", position=0.2, align=(0,1,-1));

var M = join(Point(1,1,-2,-1),Point(0,2,1,-0.5))/4;

f.line(M,'$\mathbf{\Phi}$',align=(0.75,0,-0.8));
f.line(dot(M,L)/L+axis(cross(M,L))/L,position=1,align=(0,-0.65,1),"$\mathsf{proj_1}(\mathbf{\Phi};\mathbf{\Lambda})$");
f.line_at_infinity((cross(M,L)-axis(cross(M,L))+wedge(M,L))/L,(0.5,2.5,-0.5),"$\mathsf{rej_1}(\mathbf{\Phi};\mathbf{\Lambda})$");

\end{asy}
\end{subfloatenv}
\begin{subfloatenv}{ }
\begin{asy}
import Figure3D;
Figure3D f = Figure3D();

var L = (0.5*e_0-e_2)*e_1;
f.line(L,"$\mathbf{\Lambda}$", position=0.2, align=(0,1,-1));

var M = join(Point(1,1,-2,-1),Point(0,2,1,-0.5))/4;

f.line(M,'$\mathbf{\Phi}$',align=(0.75,0,-0.8));
f.line(dot(M,L)/L+(cross(M,L)-axis(cross(M,L)))/L,size=3.5,position=0.9,align=(0,1,0),"$\mathsf{proj_2}(\mathbf{\Phi};\mathbf{\Lambda})$");
f.line((axis(cross(M,L))+wedge(M,L))/L,position=0.3,"$\mathsf{rej_2}(\mathbf{\Phi};\mathbf{\Lambda})$");

//import Figure3D;
//Figure3D f = Figure3D();
//var L = (e_0-e_2)*e_1;
//f.line(L,"$\mathbf{\Lambda}$", position=0.9, align=(0,1,1));

//var Q = Point(1,3,-2,1/2)/4;
//f.point(Q,"$\textbf{Q}$", draw_helper_lines=true,draw_orientation=false,align=(0,0,0.5));
//f.plane(join(Q,L)/L, "$(\textbf{Q}\vee\mathbf{\Lambda})\mathbf{\Lambda}^{-1}$", align=(0,1.8,-0.5),draw_orientation=false);
//f.line(join(L,Q)/Q, "$(\mathbf{\Lambda}\vee\textbf{Q})\textbf{Q}^{-1}$", position=0.1, align=(0,4,0), shift=-3.5);

\end{asy}
\end{subfloatenv}
\caption{Projection and rejection of a line on/by another line}
\label{components of rejection in E3}
\end{figure}

Projection and rejection in \E{3} can be defined by splitting the geometric product between two blades into two components
in the same way it was done in \E{2}.
Projection and rejection are illustrated in Figures~\ref{projection on plane in E3},~\ref{projection on point in E3},~\ref{projection on line in E3}.
In general, the geometric product of two finite lines consists of three components rather than two:
the inner product, the commutator, and the outer product.
If the two lines intersect, the outer product is zero and, therefore, the geometric product splits naturally into two components,
which give the projection and rejection.
On the other hand, if the lines do not intersect, I can assemble the three elements of the geometric product into
two components as follows.
The commutator can be split into the two axes, which can then be combined with the inner and outer products in two different ways
resulting in two distinct projections and two distinct rejections.
Firstly,
\begin{equation}
\begin{aligned}
&\ts{proj_1}(\mb{\Phi};\mb{\Lambda}) 
=(\mb{\Phi}\cdot\mb{\Lambda}+(\mb{\Phi}\times\mb{\Lambda})_{fa})\mb{\Lambda}^{-1},\\
&\ts{rej_1}(\mb{\Phi};\mb{\Lambda}) 
=((\mb{\Phi}\times\mb{\Lambda})_{ia}+\mb{\Phi}\wedge\mb{\Lambda})\mb{\Lambda}^{-1},
\end{aligned}
\end{equation}
where \((\mb{\Phi}\times\mb{\Lambda})_{fa}\) and \((\mb{\Phi}\times\mb{\Lambda})_{ia}\)
are the finite and the infinite axes of the commutator, respectively.
The line \(\ts{proj_1}(\mb{\Phi};\mb{\Lambda}) \) is parallel to \(\mb{\Phi}\) and intersects  \(\mb{\Lambda}\)
at the point where the finite axis \((\mb{\Phi}\times\mb{\Lambda})_{fa}\) intersects \(\mb{\Lambda}\)
(see Figure~\ref{components of rejection in E3}{a}).


Secondly,
\begin{equation}
\begin{aligned}
&\ts{proj_2}(\mb{\Phi};\mb{\Lambda}) 
=(\mb{\Phi}\cdot\mb{\Lambda}+(\mb{\Phi}\times\mb{\Lambda})_{ia})\mb{\Lambda}^{-1},\\
&\ts{rej_2}(\mb{\Phi};\mb{\Lambda}) 
=((\mb{\Phi}\times\mb{\Lambda})_{fa}+\mb{\Phi}\wedge\mb{\Lambda})\mb{\Lambda}^{-1},
\end{aligned}
\end{equation}
where I swapped the finite axis and the axis at infinity.
Projection and rejection of the second kind are illustrated in Figure~\ref{components of rejection in E3}(b).
They can be described as rotational projection and rejection since they corresponds
to a rotation of \(\mb{\Phi}\) without translation,
the projection is parallel to \(\mb{\Lambda}\) and the rejection is perpendicular to it.
On the other hand, projection and rejection of the first kind can be called translational projection and rejection
since they corresponds to a translation of \(\mb{\Phi}\) without rotation,
where both the projection and rejection are parallel to \(\mb{\Phi}\) and the projection passes through \(\mb{\Lambda}\),
while the rejection is at infinity.

Note that the projection of any geometric object at infinity on a finite point  is zero.
Also note that the rejection of a plane by a finite point is a plane at infinity.

\begin{figure}[t!]
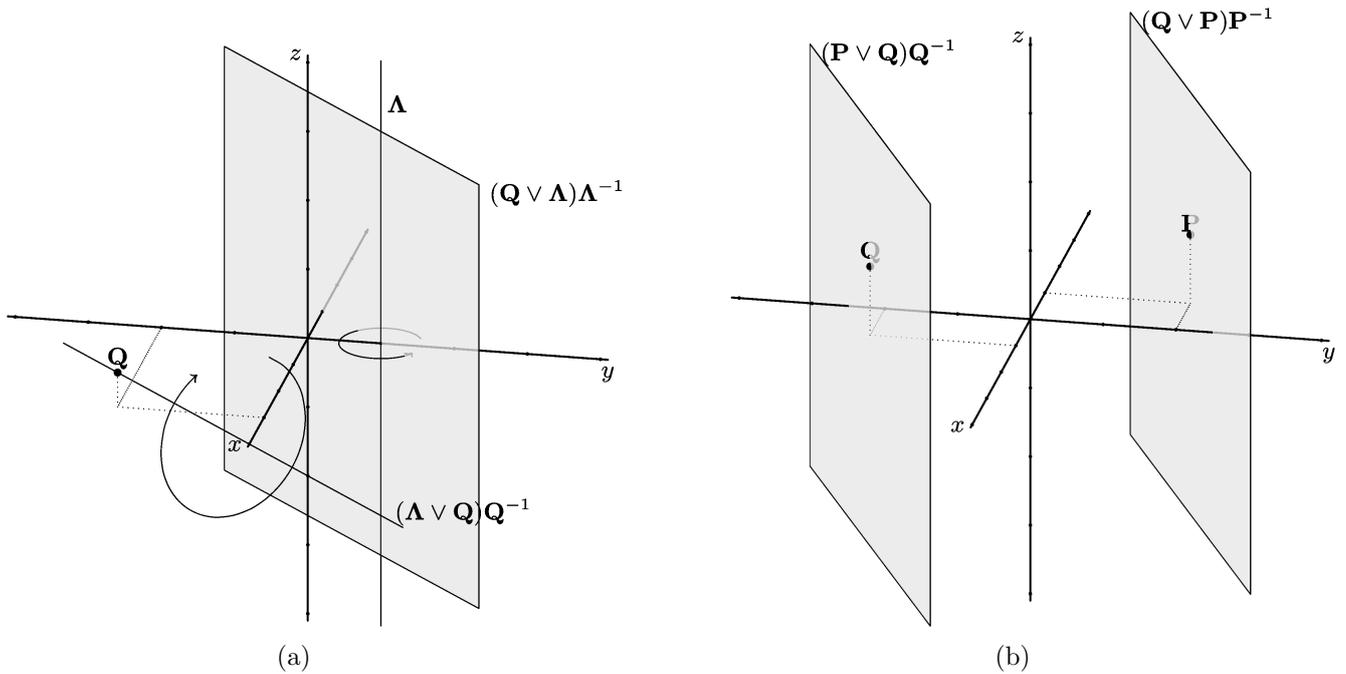
\hspace{-1cm}
\begin{subfloatenv}{ }
\begin{asy}
import Figure3D;
Figure3D f = Figure3D();
var L = (e_0-e_2)*e_1;
f.line(L,"$\mathbf{\Lambda}$", position=0.9, align=(0,1,1));

var Q = Point(1,3,-2,1/2)/4;
f.point(Q,"$\textbf{Q}$", draw_helper_lines=true,draw_orientation=false,align=(0,0,0.5));
f.plane(join(Q,L)/L, "$(\textbf{Q}\vee\mathbf{\Lambda})\mathbf{\Lambda}^{-1}$", align=(0,1.8,-0.5),draw_orientation=false);
f.line(join(L,Q)/Q, "$(\mathbf{\Lambda}\vee\textbf{Q})\textbf{Q}^{-1}$", position=0.1, align=(0,4,0), shift=-3.5);

\end{asy}
\end{subfloatenv}
\begin{subfloatenv}{ }
\begin{asy}
import Figure3D;
Figure3D f = Figure3D();

var Q = Point(1,1,-2,1)/4;
var P = Point(1,-1,2,1)/4;

f.point(Q,"$\textbf{Q}$", draw_helper_lines=true,draw_orientation=false,align=(0,0,0.5));
f.point(P,"$\textbf{P}$", draw_helper_lines=true,draw_orientation=false,align=(0,0,0.5));

f.plane(join(P,Q)/Q, "$(\textbf{P}\vee\textbf{Q})\textbf{Q}^{-1}$", align=(0,1.8,-0.5),draw_orientation=false);
f.plane(join(Q,P)/P, "$(\textbf{Q}\vee\textbf{P})\textbf{P}^{-1}$", label_position=3, align=(0,1.8,-0.5),draw_orientation=false);

\end{asy}
\end{subfloatenv}
\caption{Auxiliary objects}
\label{auxiliary objects in E3}
\end{figure}

The join of a line and a point can be used to define two auxiliary geometric objects of interest,
as illustrated in Figure~\ref{auxiliary objects in E3}(a).
The join of two points can be used to define planes shown in Figure ~\ref{auxiliary objects in E3}(b).

\begin{figure}[h!]
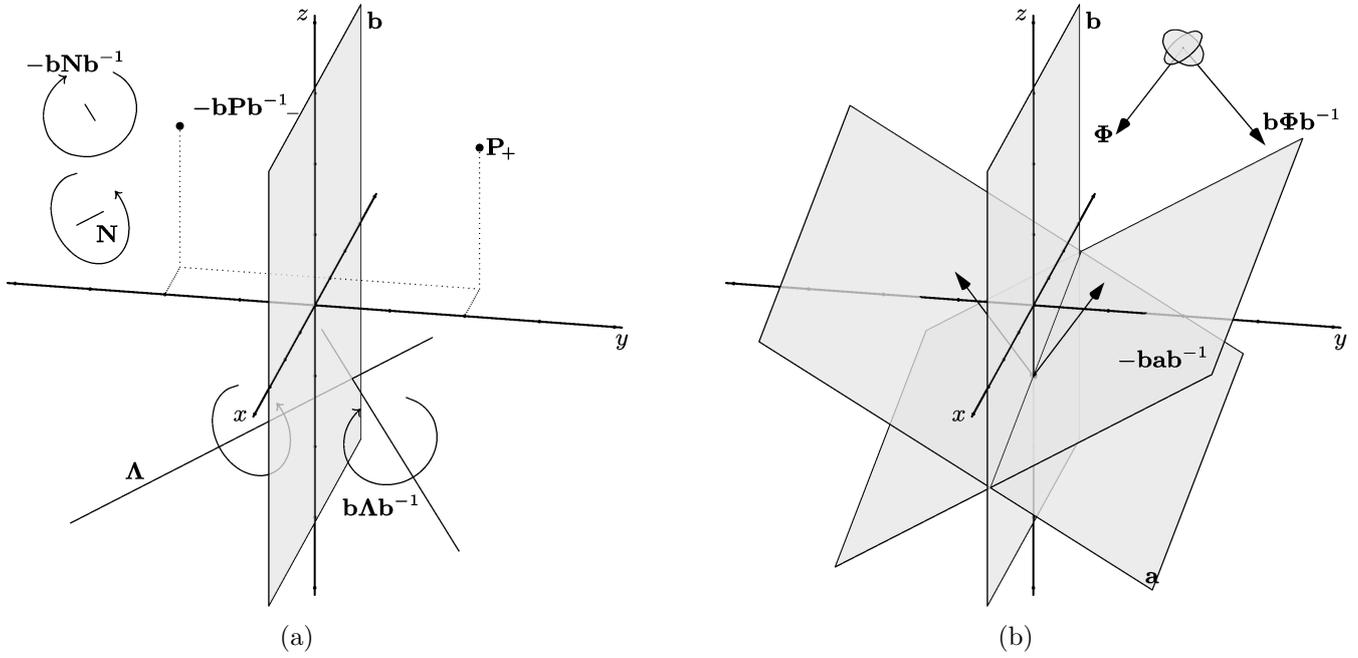

\hspace{-1cm}
\begin{subfloatenv}{ }
\begin{asy}
import Figure3D;
Figure3D f = Figure3D();

var a = join(O, e_0*e_2);
f.plane(a,"$\textbf{b}$", draw_orientation=false, align=(0,1,-1));

var P = Point(1,-1,2,2);
f.point(P,"$\textbf{P}$", draw_helper_lines=true);
f.point(-a*P/a, "$-\textbf{b}\textbf{P}\textbf{b}^{-1}$", align=(0,2,2));

var L = join(Point(1,2,-2,-2),Point(1,-2.5,0,-2))/4;

f.line(L,"$\mathbf{\Lambda}$",position=0.8,align=(0,-0.5,1),orientation_theta=-10);
f.line(a*L/a, "$\textbf{b}\mathbf{\Lambda}\textbf{b}^{-1}$",position=0.2,align=(0,-2,0),orientation_theta=-10);

var N = wedge(e_0,L); //Point(0,1,2,3)/3;
f.point_at_infinity(N,(0,-3,1),"$\textbf{N}$",align=(0,1,-1),orientation_theta=-10);
f.point_at_infinity(-a*N/a,(0,-3,2.5),"$-\textbf{b}\textbf{N}\textbf{b}^{-1}$",align=(0,-1,6));

\end{asy}
\end{subfloatenv}
\begin{subfloatenv}{ }
\begin{asy}
import Figure3D;
Figure3D f = Figure3D();

var a = join(O, e_0*e_2);
f.plane(a,"$\textbf{b}$", draw_orientation=false, align=(0,1,-1));

var K = Line(1,2,3,0,0,0)/2;
f.line_at_infinity(K,(0,2,3.8),"$\mathbf{\Phi}$",align=(0,-0.5,0));
f.line_at_infinity(a*K/a,(0,2,3.8), "$\textbf{b}\mathbf{\Phi}\textbf{b}^{-1}$",align=(0,1.5,2.5));

var b = join(Point(1,0,0,-1),e_0*(0.25*e_1+e_2+1.5*e_3));
f.plane(b,u=(0,1,0),"$\textbf{a}$",align=(0,0,1),draw_orientation=true,O=(0,0,-1));
f.plane(-a*b/a,u=(0,1,0),"$-\textbf{b}\textbf{a}\textbf{b}^{-1}$",align=(0,-1,1),draw_orientation=true,O=(0,0,-1));

\end{asy}
\end{subfloatenv}
\caption{Top-down reflection in a plane in \E{3}}
\label{reflection in a plane in E3}
\end{figure}

\begin{figure}[t!]
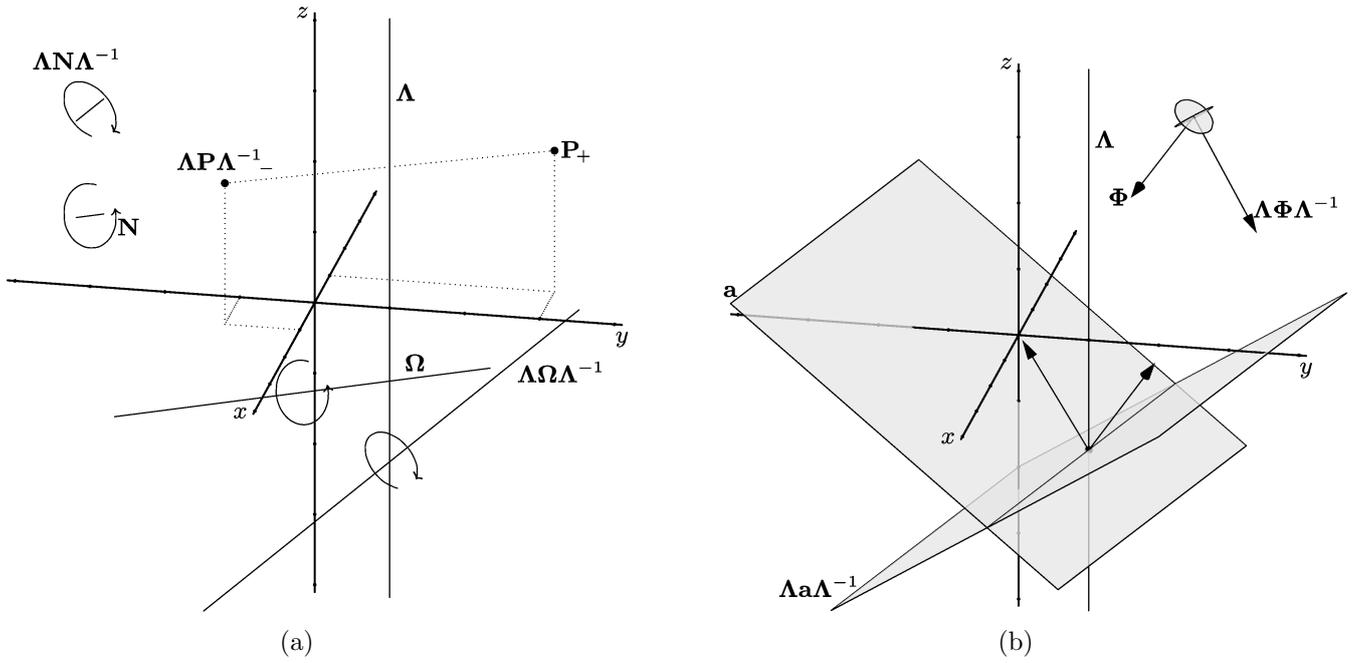

\hspace{-1cm}
\begin{subfloatenv}{ }
\begin{asy}
import Figure3D;
Figure3D f = Figure3D();

var L = e_1*e_2 + e_0*e_1;
f.line(L,"$\mathbf{\Lambda}$", position=0.9, draw_orientation=false, align=(0,1,-1));

var P = Point(1,-1,3,2);
f.point(P,"$\textbf{P}$", draw_helper_lines=true);
f.point(-L*P/L, "$\mathbf{\Lambda}\textbf{P}\mathbf{\Lambda}^{-1}$", align=(0,0,2));
draw(totriple(P)--totriple(L*P/L),dotted);

var K = join(Point(1,1,-0.5,-1),Point(1,-1,0.5,-1.5))/4;

f.line(K,"$\mathbf{\Omega}$",position=0.2,align=(0,0,1),orientation_theta=-10);
f.line(L*K/L, "$\mathbf{\Lambda}\mathbf{\Omega}\mathbf{\Lambda}^{-1}$",position=0.8,align=(0,2,0),orientation_theta=-10,shift=-1);

var N = wedge(e_0,K); //Point(0,1,2,3)/3;
f.point_at_infinity(N,(0,-3,1),"$\textbf{N}$",align=(0,4,-1),orientation_theta=-10);
f.point_at_infinity(L*N/L,(0,-3,2.5),"$\mathbf{\Lambda}\textbf{N}\mathbf{\Lambda}^{-1}$",align=(0,-1,6));

\end{asy}
\end{subfloatenv}
\begin{subfloatenv}{ }
\begin{asy}
import Figure3D;
Figure3D f = Figure3D();

var L = e_1*e_2 + e_0*e_1;
f.line(L,"$\mathbf{\Lambda}$", position=0.9, draw_orientation=false, align=(0,1,-1));

var K = Line(1,2,3,0,0,0)/2;
f.line_at_infinity(K,(0,2.5,3.5),"$\mathbf{\Phi}$",align=(0,-0.5,0));
f.line_at_infinity(L*K/L,(0,2.5,3.5), "$\mathbf{\Lambda}\mathbf{\Phi}\mathbf{\Lambda}^{-1}$",align=(0,1.5,2.5));

var b = join(Point(1,0,0,-1),e_0*(0.25*e_1+e_2+1.5*e_3));

var M = wedge(b,L*b/L);
triple u = (M.p23,M.p31,M.p12);
triple c = totriple(wedge(L,b));
f.plane(b,u,"$\textbf{a}$",align=(0,0,1),draw_orientation=true,O=c);
f.plane(L*b/L,u,"$\mathbf{\Lambda}\textbf{a}\mathbf{\Lambda}^{-1}$",align=(0,-0.5,2),draw_orientation=true,O=c);

\end{asy}
\end{subfloatenv}
\caption{Top-down reflection in a line in \E{3}}
\label{reflection in a line in E3}
\end{figure}

\begin{figure}[h!]
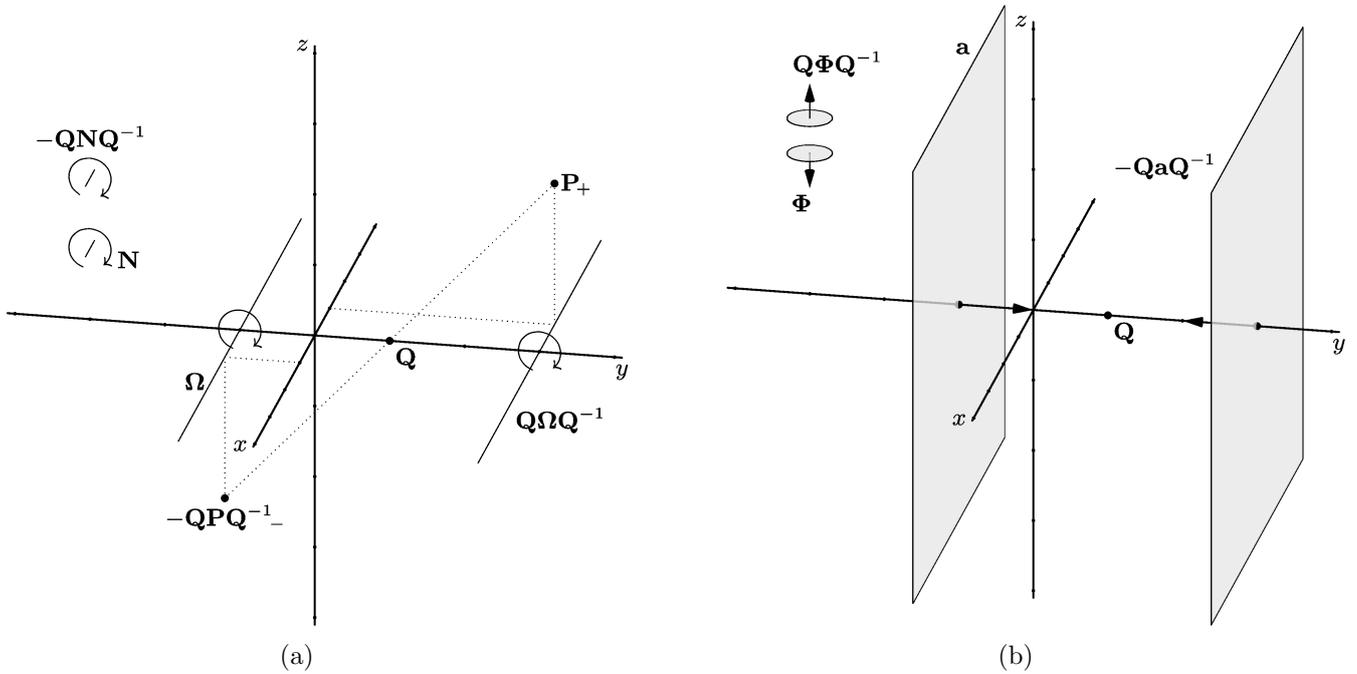
\hspace{-1cm}
\begin{subfloatenv}{ }
\begin{asy}
import Figure3D;
Figure3D f = Figure3D();

var Q = Point(1,0,1,0);
f.point(Q,"$\textbf{Q}$", draw_orientation=false, align=(0,1,-1));

var P = Point(1,-1,3,2);
f.point(P,"$\textbf{P}$", draw_helper_lines=true);
f.point(-Q*P/Q, "$-\textbf{Q}\textbf{P}\textbf{Q}^{-1}$", align=(0,0,-1));
draw(totriple(P)--totriple(Q*P/Q),dotted);

var K = join(Point(1,0,-1,0),Point(1,1,-1,0))/4;

f.line(K,"$\mathbf{\Omega}$",position=0.2,align=(0,-0.5,1),orientation_theta=-10);
f.line(Q*K/Q, "$\textbf{Q}\mathbf{\Omega}\textbf{Q}^{-1}$",position=0.2,align=(0,2,0),orientation_theta=-10,shift=0);

var N = wedge(e_0,K); //Point(0,1,2,3)/3;
f.point_at_infinity(N,(0,-3,1),"$\textbf{N}$",align=(0,4,-1),orientation_theta=-10);
f.point_at_infinity(-Q*N/Q,(0,-3,2),"$-\textbf{Q}\textbf{N}\textbf{Q}^{-1}$",align=(0,0,4));

\end{asy}
\end{subfloatenv}
\begin{subfloatenv}{ }
\begin{asy}
import Figure3D;
Figure3D f = Figure3D();

var Q = Point(1,0,1,0);
f.point(Q,"$\textbf{Q}$", draw_orientation=false, align=(0,1,-1));

var K = Line(0,0,1,0,0,0)/2;
f.line_at_infinity(K,(0,-3,2),"$\mathbf{\Phi}$",align=(0,-0.5,-1));
f.line_at_infinity(Q*K/Q,(0,-3,2.5), "$\textbf{Q}\mathbf{\Phi}\textbf{Q}^{-1}$",align=(0,0.5,1));

var b = e_2+e_0;
f.plane(b,"$\textbf{a}$",align=(0,-5,-6));
f.plane(-Q*b/Q,"$-\textbf{Q}\textbf{a}\textbf{Q}^{-1}$",align=(0,-1.5,2));

\end{asy}
\end{subfloatenv}
\caption{Top-down reflection in a point in \E{3}}
\label{reflection in a point in E3}
\end{figure}

The top-down reflection of a plane \(\tb{a}\) in an invertible plane \(\tb{b}\) is given by \(-\tb{b}\tb{a}\tb{b}^{-1}\).
The top-down reflection of lines and points in an invertible plane can be constructed from the suitable plane reflections.
For instance, if a line \(\mb{\Lambda}\) is written as the intersection of two suitable planes, e.g.\ \(\mb{\Lambda}=\tb{a}_1\wedge\tb{a}_2\), then
the top-down reflection of \(\mb{\Lambda}\) in \(\tb{b}\) is given by 
\(\tb{b}\mb{\Lambda}\tb{b}^{-1}=\tb{b}(\tb{a}_1\wedge\tb{a}_2)\tb{b}^{-1}=(-\tb{b}\tb{a}_1\tb{b}^{-1})\wedge(-\tb{b}\tb{a}_2\tb{b}^{-1})\).
Similarly, for a point \(\tb{P}=\tb{a}_1\wedge\tb{a}_2\wedge\tb{a}_3\), I get 
\(-\tb{b}\tb{P}\tb{b}^{-1}=-\tb{b}(\tb{a}_1\wedge\tb{a}_2\wedge\tb{a}_3)\tb{b}^{-1}
=(-\tb{b}\tb{a}_1\tb{b}^{-1})\wedge(-\tb{b}\tb{a}_2\tb{b}^{-1})\wedge(-\tb{b}\tb{a}_3\tb{b}^{-1})\).
Examples are given in Figure~\ref{reflection in a plane in E3}.
The top-down reflection of a plane \(\tb{a}\) in an invertible line \(\mb{\Lambda}\) is given by \(\mb{\Lambda}\tb{a}\mb{\Lambda}^{-1}\).
It can be thought of as two consecutive reflections in two perpendicluar planes intersecting at the line \(\mb{\Lambda}\).
The top-down reflection of lines and points in the line \(\mb{\Lambda}\) can be constructed from suitable plane reflections as before, e.g.\ 
\(\mb{\Lambda}\tb{P}\mb{\Lambda}^{-1}=(\mb{\Lambda}\tb{a}_1\mb{\Lambda}^{-1})\wedge(\mb{\Lambda}\tb{a}_2\mb{\Lambda}^{-1})\wedge(\mb{\Lambda}\tb{a}_3\mb{\Lambda}^{-1})\).
See Figure~\ref{reflection in a line in E3} for examples.
The top-down reflection of a plane \(\tb{a}\) in an invertible point \(\tb{Q}\) is given by \(-\tb{Q}\tb{a}\tb{Q}^{-1}\).
This can be thought of as three consecutive reflection in three perpendicular planes intersecting at \(\tb{Q}\).
As each of the three reflections carries the minus sign, it is retained in the final expression.
The other reflections in an invertible point are constructed from suitable plane reflections (see Figure~\ref{reflection in a point in E3}).
In general, the top-down reflection of a blade \(B_l\) in an invertible blade \(A_k\) is given by \((-1)^{kl}A_k B_lA_k^{-1}\).

Observe the similarity between the reflection of points at infinity and finite lines.
This is related to the fact that a point at infinity is the intersection point of a suitable finite line with the plane at infinity.
Moreover, a point at infinity shares its orientation with the orientation of the finite line.
Lines at infinity and finite planes exhibit the same similarity.
Indeed, a line at infinity is the intersection of a suitable plane with the plane at infinity, and they share the orientation as well.
The reflection of a line in a finite point does not change its top-down orientation and, therefore, it is equivalent to a translation.
The reflection of lines in \E{2} and planes in \E{3} exhibits similar correspondence.
Note that the reflection of a line in a plane which is parallel to it results in a line with the opposite top-down orientation.
This is consistent with the fact that a line is a sheaf of planes in the top-down view.
Therefore, in order to obtain its reflection in a given plane, one needs to reflect each plane in the sheaf,
which obviously produces the opposite top-down orientation.
The same analysis can be applied to other reflections.

The bottom-up reflection of a point \(\tb{P}\) in an invertible point \(\tb{Q}\) is assumed to be given by \(\tb{Q}\tb{P}\tb{Q}^{-1}\).
If a line \(\mb{\Lambda}\) can be written as the join of two points, e.g. \(\mb{\Lambda}=\tb{P}_1\vee\tb{P}_2\), then
the reflection of \(\mb{\Lambda}\) in \(\tb{Q}\) is given by the join of the corresponding point reflections.
Furthermore, I get \((\tb{Q}\tb{P}_1\tb{Q}^{-1})\vee(\tb{Q}\tb{P}_2\tb{Q}^{-1})=-\tb{Q}(\tb{P}_1\vee\tb{P}_2)\tb{Q}^{-1}\) and,
therefore, the bottom-up reflection of \(\mb{\Omega}\) in \(\tb{Q}\) is given by \(-\tb{Q}\mb{\Omega}\tb{Q}^{-1}\).
Similar analysis for a plane \(\tb{a}\) gives \(\tb{Q}\tb{a}\tb{Q}^{-1}\) for the bottom-up reflection of \(\tb{a}\) in \(\tb{Q}\).
I further assume that the bottom-up reflection of \(\tb{P}\) in an invertible line \(\mb{\Lambda}\) and in an invertible plane \(\tb{b}\)
are given by \(\mb{\Lambda}\tb{P}\mb{\Lambda}^{-1}\) and \(\tb{b}\tb{P}\tb{b}^{-1}\) respectively.
The bottom-up reflection of lines and planes in invertible lines and planes is constructed by means of suitable point reflections and the join.
For instance, I get \((\tb{b}\tb{P}_1\tb{b}^{-1})\vee(\tb{b}\tb{P}_2\tb{b}^{-1})=-\tb{b}(\tb{P}_1\vee\tb{P}_2)\tb{b}^{-1}\) and, therefore,
the bottom-up reflection of \(\mb{\Omega}\) in \(\tb{b}\) is given by \(-\tb{b}\mb{\Omega}\tb{b}^{-1}\).
I also get \((\tb{b}\tb{P}_1\tb{b}^{-1})\vee(\tb{b}\tb{P}_2\tb{b}^{-1})\vee(\tb{b}\tb{P}_3\tb{b}^{-1})=\tb{b}(\tb{P}_1\vee\tb{P}_2\vee\tb{P}_3)\tb{b}^{-1}\),
so the bottom-up reflection of \(\tb{a}\) in \(\tb{b}\) is given by \(\tb{b}\tb{a}\tb{b}^{-1}\).
In general, the bottom-up reflection of a blade \(B_l\) in an invertible blade \(A_k\) is given by \((-1)^{nk(l-1)}A_kB_lA_k^{-1}\), where
\(n\) is the dimension of the space, e.g. \(n=3\) in \E{3}.


\begin{figure}[t!]
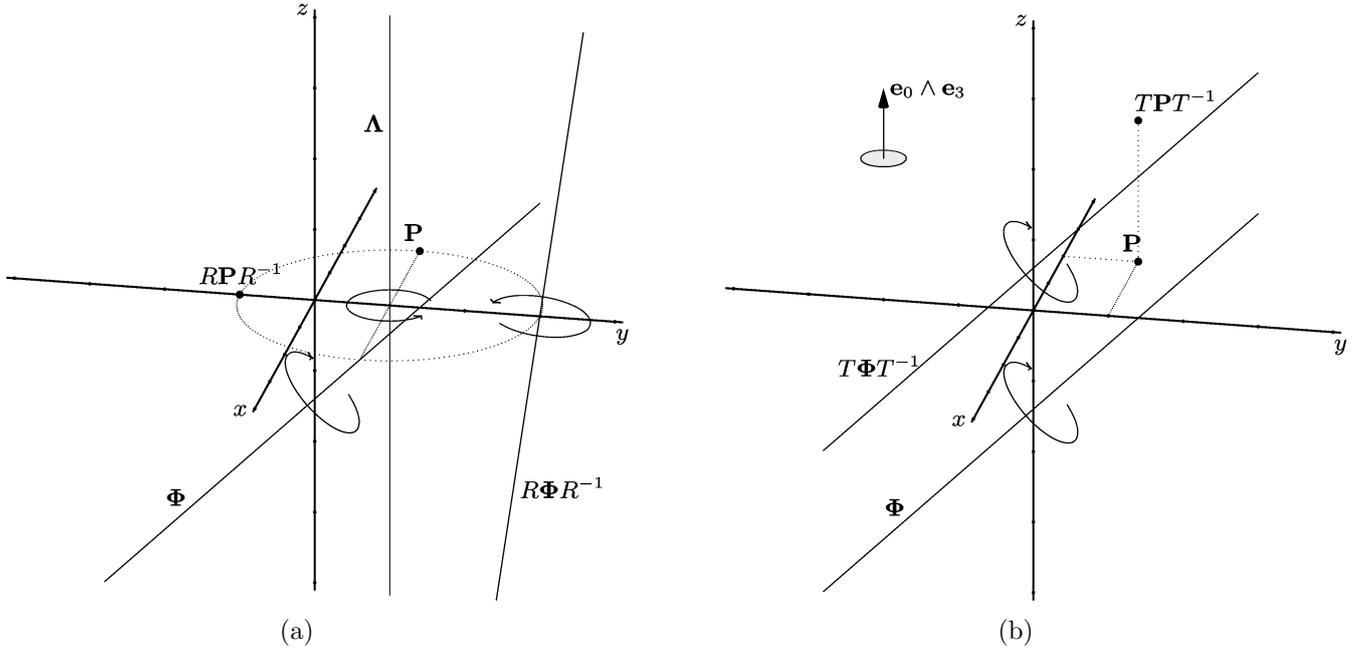

\hspace{-1cm}
\begin{subfloatenv}{ }
\begin{asy}
import Figure3D;
Figure3D f = Figure3D();

real alpha = 90*pi/180;
var L = -join(Point(1,0,1,0),Point(1,0,1,1));
L/=norm(L);
var R = exp(-1/2*alpha*L);

f.line(L,"$\mathbf{\Lambda}$", position=0.8, align=(0,-1,0.5));

var P = Point(1,-2,1,0);
var M = join(Point(1,2,1,0),Point(1,2,2,1));

f.line(M,"$\mathbf{\Phi}$", position=0.8, align=(0,-1,0.5),orientation_theta=90);
f.point(P,"$\textbf{P}$", align=(0,-0.5,1.5),draw_orientation=false,draw_helper_lines=false);

//f.line(R*M/R,size=3,"$R\mathbf{\Phi}R^{-1}$",align=(0,1,0),position=0.8);
f.line(R*M/R,"$R\mathbf{\Phi}R^{-1}$",align=(0,1,0),position=0.8);
f.point(R*P/R,"$R\textbf{P}R^{-1}$", align=(0,0,1.5),draw_orientation=false,draw_helper_lines=false);

draw(circle((0,1,0),2,(0,0,1)), dotted);
draw((-2,1,0)--(2,1,0),dotted);

real t = 2;
var S = exp(-1/2*(alpha-t*I)*L); 

\end{asy}
\end{subfloatenv}
\begin{subfloatenv}{ }
\begin{asy}
import Figure3D;
Figure3D f = Figure3D();

real t = 2;
var T = exp(-1/2*t*e_0*e_3);

var P = Point(1,-2,1,0);
var M = join(Point(1,2,1,0),Point(1,2,2,1));

f.line(M,"$\mathbf{\Phi}$", position=0.8, align=(0,-1,0.5),orientation_theta=90);
f.point(P,"$\textbf{P}$", align=(0,-0.5,1.5),draw_orientation=false);

f.line(T*M/T,"$T\mathbf{\Phi}T^{-1}$",align=(0,-2,-1),position=0.75,shift=-1.41,orientation_theta=90);
f.point(T*P/T,"$T\textbf{P}T^{-1}$", align=(0,1,1.5),draw_orientation=false);

//draw((2,1,0)--(2,1,2), dotted);

f.line_at_infinity(e_0*e_3,(0,-2,2),"$\textbf{e}_0\wedge\textbf{e}_3$");

\end{asy}
\end{subfloatenv}
\caption{Rotation and translation in \E{3}}
\label{rotation in E3}
\end{figure}

If a blade \(A_k\) can be split into projection and rejection with respect to an invertible blade \(B_l\), i.e.\ 
\(A_k  = \ts{proj}(A_k;B_l) +\ts{rej}(A_k; B_l)\),
then scaling  \(A_k\) by a factor of \(\gamma\) with respect to \(B_l\) is given by
\begin{equation}
\ts{scale}(A_k;B_l,\gamma)  = \ts{proj}(A_k;B_l) +\gamma\ts{rej}(A_k; B_l).
\label{scaling 3D}
\end{equation}
The right-hand side can also be written as \(A_k +(\gamma-1)\ts{rej}(A_k;B_l)\).
For lines, two kinds of scaling can be defined depending on whether one uses
rotational or translational projection and rejection.

Any proper motion is \E{3} can be obtained via an even number of consecutive reflections in normalised planes which square to unity.
The definition and properties of the Spin group in \E{2} equally apply to \E{3} if I replace lines in \E{2} with planes in \E{3}, e.g.\ 
a spinor in \E{3} can be written as the product of an even number of planes, all of which square to unity.
Moreover, any proper motion can be generated by a spinor that has the form \(S=e^A\) for some bivector \(A\).
Three distinct kinds of proper motion are possible in \E{3}: rotation, translation, and simultaneous rotation and translation along the axis of rotation.
The latter is obtained when \(A\) is a non-simple bivector.

In \E{3}, a spinor can be written as
\begin{equation}
S=e^{-\tfrac{1}{2}(\alpha-\lambda\I)\mb{\Lambda}},
\end{equation}
where \(\alpha,\lambda\in\R{}\) and \(\mb{\Lambda}\) is a normalised line, i.e.\ a simple bivector that satisfies \(\norm{\mb{\Lambda}}=1\).
When it acts on a blade \(A_k\) as
\begin{equation}
SA_kS^{-1},
\end{equation}
it results in 1) a rotation of \(A_k\) by the angle \(\alpha\) around the line \(\mb{\Lambda}\) 
(the sense of rotation is the same as the orientation of \(\mb{\Lambda}\) for \(\alpha>0\))
and 2) a translation by \(\lambda\) in the direction of the top-down orientation of \(\I\mb{\Lambda}\) if \(\lambda>0\)
(note that in \E{3} the top-down orientation of \(\I\mb{\Lambda}\) coincides with the bottom-up orientation of \(\mb{\Lambda}\)).
A pure rotation results if \(\lambda=0\), while \(\alpha=0\) gives a pure translation.
A translation can also be performed with  a spinor
\begin{equation}
T=e^{-\tfrac{1}{2}\lambda\e_0\wedge\tb{a}},
\end{equation}
where \(\e_0\wedge\tb{a}\) is a line at infinity and \(\tb{a}\) is a normalised plane.

For the spinor \(R=e^{-\tfrac{1}{2}\alpha\mb{\Lambda}}\), where
\(\alpha=\tfrac{\pi}{2}\) and \(\mb{\Lambda}=-\e_{10}+\e_{12}\),
its action on  \(\tb{P}=\e_{123}-2\e_{320}+\e_{130}\)
and  \(\mb{\Phi}=\e_{10}-2\e_{20}+2\e_{30}-\e_{31}-\e_{12}\)
gives
\(R\tb{P}R^{-1}=\e_{123}-\e_{130}\) and 
\(R\mb{\Phi}R^{-1}=3\e_{10}+3\e_{30}+\e_{23}-\e_{12}\)
shown in Figure~\ref{rotation in E3}(a).
Applying \(T=e^{-\tfrac{1}{2}\lambda\e_0\wedge\tb{a}}\), where
\(\lambda=2\) and \(\tb{a}=\e_3\), to the same point and line gives
 \(T\tb{P}T^{-1}=\e_{123}-2\e_{320}+\e_{130}+2\e_{210}\)
and  \(T\mb{\Phi}T^{-1}=-\e_{10}-2\e_{20}+2\e_{30}-\e_{31}-\e_{12}\)
shown in Figure~\ref{rotation in E3}(b).
Since \(\tb{a}\) in the spinor \(T\) is perpendicular to \(\mb{\Lambda}\), 
bivectors \(\e_0\wedge\tb{a}\) and \(\mb{\Lambda}\) commute and, therefore,
their action can be combined in a single spinor \(S=e^{-\tfrac{1}{2}(\alpha-\lambda\I)\mb{\Lambda}}\).
Indeed, \(-\I\mb{\Lambda}=\e_0\wedge\tb{a}\) for \(\mb{\Lambda}=-\e_{10}+\e_{12}\)
and  \(\tb{a}=\e_3\).

The same applies to rotation and translation of planes, e.g.\ \(R\tb{b}R^{-1}\)
gives a rotation of \(\tb{b}\) by the spinor \(R\).
Translation of points and lines at infinity does not change them.
In general, combined rotational and translational action of a spinor \(S\) on multivector \(M\) 
is given by \(SMS^{-1}\).

\newpage
\section{4-dimensional geometry}

\subsection{Introduction}

In this section on the 4-dimensional geometry, I will make use the intuitions developed in the previous sections.
I will rely on Grassmann and Clifford algebras to develop the insight into the properties of geometric objects in four dimensions,
since effective visualisation is difficult.
However, the basic setup of the model is similar to that in the 2- and 3-dimensional cases.
It is briefly summarised below.

In the target space \T{4}, a 3-dimensional hyperplane which does not pass through the origin can be defined by 
\begin{equation}\label{plane in 4D}
ax+by+cz+ht+1=0,
\end{equation}
where \(a\), \(b\), \(c\), and \(h\) are fixed and \(x\), \(y\), \(z\) and \(t\) are variable.
The dual space \T{4*} is a linear space of the coefficients \((a,b,c,h)\),
where every point defines a hyperplane\footnote{In this section, a hyperplane refers to a 3-dimensional hyperplane unless stated otherwise.}
 in the target space \T{4} via Equation~(\ref{plane in 4D}). 
For instance,  \(Q=(a,b,c,h)=(0,0,0,2)\) in \T{4*} is dual to 
the hyperplane \(Q^*\) in \T{4} given by \(t=-\tfrac{1}{2}\).
On the other hand, \(P=(x,y,z,t)\) in \T{4} is dual to the hyperplane \(P^*\) in \T{4*} defined by Equation~(\ref{plane in 4D}),
where \(x\), \(y\), \(z\), and \(t\) are seen as fixed coefficients and \(a\), \(b\), \(c\), \(h\) are variable.
For instance, \(P=(x,y,z,t)=(0,-3,0,0)\) in \T{4} is dual to the hyperplane \(P^*\) in \T{4*} given by \(b=\tfrac{1}{3}\).
The dual hyperplane \(P^*\) enables the top-down view of the point \(P\).
Every point in the hyperplane \(P^*\) corresponds to a hyperplane in \T{4} passing through \(P\).
So, \(P\) can be identified with a set of all hyperplanes passing through the point, a bundle of hyperplanes.
The bundle has two possible orientations and thus provides the top-down orientation of the point \(P\).
At least four hyperplanes are needed to define a point as their intersection.

A line in \T{4*} is dual to a plane in \T{4}.
Let \(L\) be a line in \T{4*} and \(L^*\) is its dual plane in \T{4}. 
Points which lie on \(L\) correspond to a set of hyperplanes in \T{4}, intersecting along the plane \(L^*\).
This set of hyperplanes is a simple extension of a 2-dimensional sheaf of lines
into two extra dimensions
(recall that a sheaf of planes in \T{3} is an extension of a 2-dimensional sheaf into one extra dimension).
So, in the top-down view, a plane in \T{4} is represented by a sheaf of hyperplanes intersecting along the plane.
The plane's top-down orientation is given by the sense of rotation of the hyperplanes in the sheaf.
At least two hyperplanes are needed to define a plane as their intersection.

A plane in \T{4*} is dual to a line in \T{4}.
Points which lie on the plane correspond to the hyperplanes in \T{4} intersecting along the dual line.
In the top-down view, a line in \T{4} is represented by a bundle of hyperplanes intersecting along the line.
This bundle has two possible orientations and provides the top-down orientation of the line.
Note that the bundle of hyperplanes that defines a line is distinct from the bundle that defines a point, even though I use the same term for both.
At least three hyperplanes are required to define a line as their intersection.

The origin of \T{4*} corresponds to the hyperplane at infinity in \T{4},
while the origin of \T{4} corresponds to the hyperplane at infinity in \T{4*}.
In the top-down view, the hyperplane at infinity can be oriented either towards the origin or away from the origin.

A stack of hyperplanes consists of all hyperplanes
\begin{equation}
\label{stack of hyperplanes in 4D}
ax+by+cz+ht+d=0,
\end{equation}
where \((a,b,c,h)\) is fixed and \(d\) spans all possible values in \(\R{}\).
A stack of hyperplanes represents a plane at infinity in \T{4}, which can be thought of as the intersection of the hyperplanes in the stack.
The orientation of the stack provides the top-down orientation of the plane at infinity it represents.
A plane at infinity in \T{4}, i.e.\ a stack of hyperplanes, is dual to a line in \T{4*} passing through the origin.
A line at infinity in \T{4} is dual to a plane in \T{4*} passing through the origin.
Similarly, lines and planes in \T{4} passing through the origin are dual to planes and lines at infinity in \T{4*}, respectively.
A point at infinity in \T{4} is dual to a hyperplane in \T{4*} passing through the origin,
while a hyperplane in \T{4} passing through the origin is dual to a point at infinity in \T{4*}.

\subsection{Embedding and Grassmann algebra}

I will use \((w,x,y,z,t)\) to refer to vectors in the target model space \R{5} and \((d,a,b,c,h)\) for vectors in
the dual model space \R{5*}.
The target space \T{4} is embedded in the target model space \R{5} and the dual space \T{4*} is embedded in the dual model space \R{5*}
in the same way as in the lower-dimensional cases.
Namely, \((1,x,y,z,t)\) in \R{5} is identified with \((x,y,z,t)\) in \T{4}
and \((1,a,b,c,h)\) in \R{5*} is identified with \((a,b,c,h)\) in \T{4*}.
So, \T{4} is a 4-dimensional hyperplane embedded in the target model space \R{5} at \(w=1\),
and \T{4*} is a 4-dimensional hyperplane at \(d=1\) in \R{5*}.
The intersection of linear subspaces of \R{5} with \T{4} gives points, lines, planes, and hyperplanes in \T{4}.
The same applies to linear subspaces in \R{5*}.
The duality and identity transformations, defined in the same fashion as in the lower-dimensional case,
are extended to linear subspaces of the model spaces \R{5} and \R{5*}.

I employ the following notation for the standard basis of \R{5}:\\
\(\e^0=(1,0,0,0,0)\), \(\e^1=(0,1,0,0,0)\), \(\e^2=(0,0,1,0,0)\), \(\e^3=(0,0,0,1,0)\), \(\e^4=(0,0,0,0,1)\),\\
so that \((w,x,y,z,t)=w\e^0+x\e^1+y\e^2+z\e^3+t\e^4\).\\
And in \R{5*}:\\
\(\e_0=(1,0,0,0,0)\), \(\e_1=(0,1,0,0,0)\), \(\e_2=(0,0,1,0,0)\), \(\e_3=(0,0,0,1,0)\), \(\e_4=(0,0,0,0,1)\), \\
so that \((d,a,b,c,h)=d\e_0+a\e_1+b\e_2+c\e_3+h\e_4\).

The basis of Grassmann algebra \(\bigwedge\R{5*}\) consists of \(2^5=32\) multivectors:\\
1, \\
\(\e_0, \e_1, \e_2, \e_3, \e_4\), \\
\(\e_{10}, \e_{20}, \e_{30}, \e_{40}, \e_{23}, \e_{31}, \e_{12}, \e_{41}, \e_{42}, \e_{43}\),\\
\(\e_{234}, \e_{314}, \e_{124}, \e_{321}, \e_{410}, \e_{420}, \e_{430}, \e_{230}, \e_{310}, \e_{120}\), \\
\(\e_{1234}, \e_{2340}, \e_{3140}, \e_{1240}, \e_{3210}\),\\
\(\I_4=\e_{01234}\).

Note that \(\e_{3210}=\e_{0123}\) but \(\e_{321}=-\e_{123}\). As usual, I will use \(\I=\I_4\).

The duality transformation \(\J:\bigwedge\R{5*}\to\bigvee\R{5}\)  is defined on the basis by the following table

\begin{tabular*}{\textwidth}{lcccccccccccccccc}
\(X\)      & 1   & \(\e_{0}\) & \(\e_{1}\) & \(\e_{2}\) & \(\e_{3}\) &  \(\e_4\) &    \(\e_{1234}\) & \(\e_{2340}\) & \(\e_{3140}\) & \(\e_{1240}\) & \(\e_{3210}\)    & \(\e_{01234}\) \\
\cline{1-13} \\[-10pt]
\(\J(X)\) & \(\e^{01234}\) & \(\e^{1234}\) & \(\e^{2340}\) & \(\e^{3140}\) & \(\e^{1240}\) & \(\e^{3210}\) & \(\e^{0}\) & \(\e^{1}\) & \(\e^{2}\) & \(\e^{3}\) &  \(\e^4\) & 1 \\
\end{tabular*}

\begin{tabular*}{\textwidth}{lcccccccccccccccc}
\(X\)   &  \(\e_{10}\) & \(\e_{20}\) & \(\e_{30}\) & \(\e_{40}\) & \(\e_{23}\) & \(\e_{31}\) & \(\e_{12}\)  & \(\e_{41}\) & \(\e_{42}\) & \(\e_{43}\)  \\
\cline{1-11} \\[-10pt]
\(\J(X)\)  & \(\e^{234}\) & \(\e^{314}\) & \(\e^{124}\) & \(\e^{321}\) & \(\e^{410}\) & \(\e^{420}\) & \(\e^{430}\)  & \(\e^{230}\) & \(\e^{310}\) & \(\e^{120}\)  \\
\end{tabular*}

\begin{tabular*}{\textwidth}{lcccccccccccccccc}
\(X\)  & \(\e_{234}\) & \(\e_{314}\) & \(\e_{124}\) & \(\e_{321}\) & \(\e_{410}\) & \(\e_{420}\) & \(\e_{430}\) & \(\e_{230}\) & \(\e_{310}\) & \(\e_{120}\)  \\
\cline{1-11} \\[-10pt]
\(J(X)\)   &  \(\e^{10}\) & \(\e^{20}\) & \(\e^{30}\) & \(\e^{40}\) & \(\e^{23}\) & \(\e^{31}\) & \(\e^{12}\) & \(\e^{41}\) & \(\e^{42}\) & \(\e^{43}\)  \\
\end{tabular*}

It extends to general multivectors by linearity.
The inverse transformation is obtained by lowering and raising the indices.

The join of multivectors \(A\) and \(B\) in  \(\bigwedge\R{5*}\) is defined by
\begin{equation}
A\vee B = \J^{-1}(\J(A)\vee\J(B))).
\end{equation}

I will use capital Roman letters in bold font for four-vectors in \(\bigwedge\R{5*}\),
which dually represent points in \T{4},
capital Greek letters for trivectors, which dually represent lines (if the trivector is simple),
and small Greek letters for bivectors, which dually represent planes (if the bivector is simple).
So, a general multivector in \R{5*} can be written as
\begin{equation}
M=s+\tb{a}+\mb{\sigma}+\mb{\Lambda}+\tb{P}+\alpha\I,
\end{equation}
where \(s,\alpha\in\R{}\).

In a 5-dimensional space, most bivectors and trivectors are not simple,
i.e.\ they cannot be written as the outer product of two vectors and three vectors, respectively.
For instance, trivector \(\e_{104}+\e_{234}\) is not simple.
A bivector \(\mb{\sigma}\) is simple if and only if \(\mb{\sigma}\wedge\mb{\sigma}=0\).
A trivector \(\mb{\Lambda}\) is simple if and only if \(\mb{\Lambda}\vee\mb{\Lambda}=0\)
(note that \(\mb{\Lambda}\wedge\mb{\Lambda}=0\) for any trivector in \R{5*}).
In other words, a trivector is simple if and only if its dual bivector is simple.
In general, \(\mb{\sigma}\wedge\mb{\sigma}\) is a point and \(\mb{\Lambda}\vee\mb{\Lambda}\) is a hyperplane.

Computing the outer product of a bivector
\begin{equation}
\mb{\sigma} = p_{10}\e_{10}+p_{20}\e_{20}+ p_{30}\e_{30}+p_{40}\e_{40}+ p_{23}\e_{23}+ p_{31}\e_{31}+ p_{12}\e_{12}
+ p_{41}\e_{41}+ p_{42}\e_{42}+p_{43}\e_{43}
\end{equation}
with itself yields 
\[
\begin{aligned}
\mb{\sigma}\wedge\mb{\sigma} &=
-2(p_{23}p_{41}+p_{31}p_{42}+p_{12}p_{43})\e_{1234}+\\
&\quad+2(p_{23}p_{40}+p_{30}p_{42}-p_{20}p_{43})\e_{2340}
+2(p_{31}p_{40}+p_{10}p_{43}-p_{30}p_{41})\e_{3140}\\
&\quad+2(p_{12}p_{40}+p_{20}p_{41}-p_{10}p_{42})\e_{1240}-2(p_{10}p_{23}+p_{12}p_{30}+p_{20}p_{31})\e_{3210}.
\end{aligned}
\]
So, \(\mb{\sigma}\) is simple if 
\begin{equation}
\left\{
\begin{aligned}
p_{23}p_{41}+p_{31}p_{42}+p_{12}p_{43}=0,\\
p_{20}p_{43}-p_{30}p_{42}-p_{23}p_{40}=0,\\
p_{30}p_{41}-p_{10}p_{43}-p_{31}p_{40}=0,\\
p_{10}p_{42}-p_{20}p_{41}-p_{12}p_{40}=0,\\
p_{10}p_{23}+p_{20}p_{31}+ p_{30}p_{12}=0.
\end{aligned}
\right.
\label{bivector simple constraints}
\end{equation}
If \(\mb{\sigma}\) is simple, it represents a linear subspace in \R{5*}.
It can be obtained by solving
\begin{equation}
\tb{a}\wedge\mb{\sigma}=0
\end{equation}
for \(\tb{a}=d\e_0+a\e_1+b\e_2+c\e_3+h\e_4\), which gives
\begin{equation}
\begin{aligned}
\left\{
\begin{aligned}
a p_{23} +b p_{31}+c p_{12}=0,\\
d p_{23} +b p_{30}-c p_{20}=0,\\
d p_{31} +c p_{10}-a p_{30}=0,\\
d p_{12} +a p_{20}-b p_{10}=0,\\
\end{aligned}
\right.\quad\quad
\left\{
\begin{aligned}
h p_{23}=b p_{43}-c p_{42},\\
h p_{31}=c p_{41}-a p_{43},\\
h p_{12}=a p_{42}-b p_{41},
\end{aligned}
\right.\quad\quad
\left\{
\begin{aligned}
a p_{40}=d p_{41}+h p_{10},\\
b p_{40}=d p_{42}+h p_{20},\\
c p_{40}=d p_{43}+h p_{30}.\\
\end{aligned}
\right.
\end{aligned}
\label{bivector plane in T4*}
\end{equation}
Of these equations, the last three are linearly independent.
The other equations follow from the last three and the constraints~(\ref{bivector simple constraints}).
Given that the last three equations in~(\ref{bivector plane in T4*}) involve five variables, 
they define a plane in \R{5*}.
Its intersection with \T{4*} is found by setting \(d=1\),
which gives a line in \T{4*}.
Let \(L\) be a line at the intersection of \T{4*} and the 2-dimensional subspace represented by \(\mb{\sigma}\).
The same simple bivector \(\mb{\sigma}\) dually represents a linear subspace in \R{5} defined by
\begin{equation}
\tb{x}\vee\J(\mb{\sigma})=0,
\end{equation}
where \(\tb{x}=w\e^0+x\e^1+ y\e^2 +z\e^3 +t\e^4\).
Using the duality transformation gives
\[
\begin{aligned}
\tb{x}
\vee &
(p_{10}\e^{234}+p_{20}\e^{314}+ p_{30}\e^{124}+p_{40}\e^{321} +      p_{23}\e^{410}+ p_{31}\e^{420}+
p_{12}\e^{430}
+ p_{41}\e^{230}+ p_{42}\e^{310}+p_{43}\e^{120})=\\
&=
+(xp_{10}+p_{20}y+zp_{30}+tp_{40})\e^{1234}-\\
&\quad-(wp_{10}+yp_{12}-zp_{31}-tp_{41})\e^{2340}
-(wp_{20}+zp_{23}-xp_{12}-tp_{42})\e^{3140}-\\
&\quad-(wp_{30}+xp_{31}-yp_{23}-tp_{43})\e^{1240}
-(wp_{40}+xp_{41}+yp_{42}+zp_{43})\e^{3210}=0
\end{aligned}
\]
and, therefore,  
\begin{equation}
\left\{
\begin{aligned}
xp_{10}+yp_{20}+zp_{30}+tp_{40}=0,\\
wp_{10}+yp_{12}-zp_{31}-tp_{41}=0,\\
wp_{20}+zp_{23}-xp_{12}-tp_{42}=0,\\
wp_{30}+xp_{31}-yp_{23}-tp_{43}=0,\\
wp_{40}+xp_{41}+yp_{42}+zp_{43}=0.\\
\end{aligned}
\right.
\label{L* in T4}
\end{equation}
In this system, only the first and the last equations are linearly independent.
They define a 3-dimensional subspace in \R{5} whose intersection with \T{4},
obtained by setting \(w=1\), gives \(L^*\), which is a plane dual to the line \(L\).
The plane \(L^*\) in \T{4} is represented dually by the bivector \(\mb{\sigma}\in\bigwedge\R{5*}\)
and directly by the trivector \(\J(\mb{\sigma})\in\bigvee\R{5}\).

The central point of the plane \(L^*\) is defined as the intersection point
of \(L^*\) and \(p_0(L^I)\), which refers to the plane passing through the line \(L^I\) and the origin of \T{4}.
A vector \(\tb{x}\) lies on the plane \(p_0(L^I)\) if \(\tb{x}\vee(\e^0\vee\Id(\mb{\sigma}))=0\) and, therefore,
its coordinates in the standard basis satisfy
\begin{equation}
\left\{
\begin{aligned}
yp_{43}-zp_{42}-tp_{23}=0,\\
zp_{41}-xp_{43} -tp_{31}=0,\\
xp_{42}-yp_{41} -tp_{12} =0,\\
xp_{23}+yp_{31}+zp_{12}=0.\\
\end{aligned}
\right.
\label{p0(L^I) in T4}
\end{equation}
Only two of these equations are linearly independent, so the system defines a plane in \T{4}.
Solving Equations~(\ref{L* in T4}) and~(\ref{p0(L^I) in T4}) yields
the coordinates of the central point:
\begin{equation}
\begin{aligned}
x_c=\frac{p_{20}p_{12}-p_{30}p_{31}-p_{40}p_{41}}{p_{23}^2+p_{31}^2+p_{12}^2+p_{41}^2+p_{42}^2+p_{43}^2},\quad
y_c=\frac{p_{30}p_{23}-p_{10}p_{12}-p_{40}p_{42}}{p_{23}^2+p_{31}^2+p_{12}^2+p_{41}^2+p_{42}^2+p_{43}^2},\\
z_c=\frac{p_{10}p_{31}-p_{20}p_{23}-p_{40}p_{43}}{p_{23}^2+p_{31}^2+p_{12}^2+p_{41}^2+p_{42}^2+p_{43}^2},\quad
t_c=\frac{p_{10}p_{41}+p_{20}p_{42}+p_{30}p_{43}}{p_{23}^2+p_{31}^2+p_{12}^2+p_{41}^2+p_{42}^2+p_{43}^2}.\\
\end{aligned}
\end{equation}
Note that the plane \(L^*\) shares a common line at infinity with the attitude of
\begin{equation}
B_\parallel(L^*)
=
p_{41}\e^{23}+p_{42}\e^{31}+p_{43}\e^{12}
+p_{23}\e^{41}+p_{31}\e^{42}+p_{12}\e^{43},
\end{equation}
which is viewed as a bivector in \T{4}, rather than \R{5}.
The orientation of \(B_\parallel(L^*)\) provides the bottom-up orientation of the plane \(L^*\).
The top-down orientation of the plane \(L^*\) is provided by
\begin{equation}
B_\perp(L^*)
=
p_{23}\e^{23}+p_{31}\e^{31}+p_{12}\e^{12}
+p_{41}\e^{41}+p_{42}\e^{42}+p_{43}\e^{43}.
\end{equation}
The bivector \(\e_0\wedge\tb{p}\), where \(\tb{p}=-p_{10}\e_1-p_{20}\e_2-p_{30}\e_3-p_{40}\e_4\),
defines a plane at infinity in \T{4} with the top-down orientation provided by
\(\Id(\tb{p})=-p_{10}\e^1-p_{20}\e^2-p_{30}\e^3-p_{40}\e^4\), which is seen as a vector in \T{4} rather than \R{5}.
Since \(\mb{\sigma}\) is simple and Equalities~(\ref{bivector simple constraints}) are satisfied,
\(\Id(\tb{p})\) lies in the attitude of the bivector \(B_\perp(L^*)\).
It determines the direction in which the attitude of \(B_\parallel(L^*)\)
must be shifted from the origin to arrive at the plane \(L^*\).
This configuration is similar to the one that defines a line shifted from the origin in \T{3}.

The central point of the line \(L\) in \T{4*} is defined as a point where
 \(L\) intersect \(q_0((L^*)^I)\), which refers to the hyperplane
passing through the plane \((L^*)^I\) and the origin of \T{4*} (\((L^*)^I\) is the identical counterpart of the plane \(L^*\)).
A point in \T{4*} represented by a vector \(\tb{a}=d\e_0+a\e_1+b\e_2+c\e_3+h\e_4\) 
lies in the hyperplane \(q_0((L^*)^I)\) if \(\tb{a}\wedge(\e_0\wedge\Id^{-1}(\J(\mb{\sigma})))=0\),
which gives
\begin{equation}
ap_{10}+bp_{20}+cp_{30}+hp_{40}=0.
\end{equation}
The central point of the line \(L\) is obtained by intersecting \(L\) with the hyperplane \(q_0((L^*)^I)\), which gives
the coordinates
\begin{equation}
\begin{aligned}
a_c=-\frac{p_{20}p_{12}-p_{30}p_{31}-p_{40}p_{41}}{p_{10}^2+p_{20}^2+p_{30}^2+p_{40}^2},\quad
b_c=-\frac{p_{30}p_{23}-p_{10}p_{12}-p_{40}p_{42}}{p_{10}^2+p_{20}^2+p_{30}^2+p_{40}^2},\\
c_c=-\frac{p_{10}p_{31}-p_{20}p_{23}-p_{40}p_{43}}{p_{10}^2+p_{20}^2+p_{30}^2+p_{40}^2},\quad
h_c=-\frac{p_{10}p_{41}+p_{20}p_{42}+p_{30}p_{43}}{p_{10}^2+p_{20}^2+p_{30}^2+p_{40}^2}.\\
\end{aligned}
\end{equation}
The line \(L\) shares a common point at infinity with the attitude of \(\tb{p}=-p_{10}\e_1-p_{20}\e_2-p_{30\e_3}-p_{40}\e_4\)
seen as a vector in \T{4*}.
If the space is Euclidean, I get \(d_{C(L^*)}d_{C(L)^I}=1\) for
the distance \(d_{C(L^*)}\) from the central point of the plane \(L^*\) to the origin of \T{4}
and the distance \(d_{C(L)^I}\) from the identical counterpart of the central point \(C(L)\) to the origin.

For a  trivector
\begin{multline}
\mb{\Lambda}
=
s_{234}\e_{234} + s_{314}\e_{314} + s_{124}\e_{124} + s_{321}\e_{321}+ \\
 + s_{420}\e_{420} + s_{430}\e_{430}+s_{410}\e_{410} + s_{230}\e_{230}
  + s_{310}\e_{310} + s_{120}\e_{120}
\end{multline}
the join with itself yields
\begin{equation}
\begin{aligned}
\mb{\Lambda}\vee\mb{\Lambda}=
-2(s_{410}s_{230}+s_{420}s_{310}+s_{430}s_{120})\e_0
-2(s_{120}s_{314}-s_{124}s_{310}-s_{410}s_{321})\e_1\\
-2(s_{124}s_{230}-s_{120}s_{234}-s_{420}s_{321})\e_2
-2(s_{234}s_{310}-s_{230}s_{314}-s_{430}s_{321})\e_3\\
-2(s_{410}s_{234}+s_{420}s_{314} +s_{430}s_{124} )\e_4,\\
\end{aligned}
\end{equation}
and so \(\mb{\Lambda}\) is simple if
\begin{equation}
\left\{
\begin{aligned}
s_{410}s_{230}+s_{420}s_{310}+s_{430}s_{120}=0,\\
s_{120}s_{314}-s_{124}s_{310}-s_{410}s_{321}=0,\\
s_{124}s_{230}-s_{120}s_{234}-s_{420}s_{321}=0,\\
s_{234}s_{310}-s_{230}s_{314}-s_{430}s_{321}=0,\\
s_{410}s_{234}+s_{420}s_{314} +s_{430}s_{124}=0.
\end{aligned}
\right.
\end{equation}
A simple trivector corresponds to a 3-dimensional linear
subspace in \R{5*}, whose vectors \(\tb{a}\)
satisfy \(\tb{a}\wedge\mb{\Lambda}=0\), which gives
\begin{equation}
\left\{
\begin{aligned}
as_{234}+bs_{314}+cs_{124}+hs_{321}=0,\\
ds_{234}+bs_{430}-cs_{420}-hs_{230}=0,\\
ds_{314}+cs_{410}-as_{430}-hs_{310}=0,\\
ds_{124}+as_{420}-bs_{410}-hs_{120}=0,\\
ds_{321}+as_{230}+bs_{310}+cs_{120}=0,\\
\end{aligned}
\right.
\end{equation}
and, therefore,
the trivector represents a plane in \T{4*}, obtained by setting \(d=1\).
Let \(K\) be the plane at the intersection of \(\mb{\Lambda}\) and \T{4*}.
The same trivector \(\mb{\Lambda}\) dually represents a linear subspace in \R{5},
whose vectors \(\tb{x}\) satisfy \(\tb{x}\vee\J(\mb{\Lambda})=0\),
which gives
\begin{equation}
\left\{
\begin{aligned}
xs_{410}+ys_{420}+zs_{430}=0,\\
ws_{410}+ys_{124}-zs_{314}=0,\\
ws_{420}+zs_{234}-xs_{124}=0,\\
ws_{430}+xs_{314}-ys_{234}=0,\\
\end{aligned}
\right.\quad\quad
\left\{
\begin{aligned}
ts_{410} = ys_{120}-zs_{310},\\
ts_{420} = zs_{230}-xs_{120},\\
ts_{430} = xs_{310}-ys_{230},\\
\end{aligned}
\right.\quad\quad
\left\{
\begin{aligned}
xs_{321}=ws_{230}+ts_{234},\\
ys_{321}=ws_{310}+ts_{314},\\
zs_{321}=ws_{120}+ts_{124}.\\
\end{aligned}
\right.
\end{equation}
Setting \(w=1\) gives the intersection of the subspace with \T{4},
i.e.\ the line \(K^*\) in \T{4} dual to the plane \(K\) and dually represented by the  trivector \(\mb{\Lambda}\).
In the standard basis, the coordinates of the central point of the line \(K^*\) are given by
\begin{equation}
\begin{aligned}
x_c=-\frac{s_{314}s_{430}-s_{124}s_{420}-s_{321}s_{230}}{s_{234}^2+s_{314}^2+s_{124}^2+s_{321}^2},\quad
y_c=-\frac{s_{124}s_{410}-s_{234}s_{430}-s_{321}s_{310}}{s_{234}^2+s_{314}^2+s_{124}^2+s_{321}^2},\\
z_c=-\frac{s_{234}s_{420}-s_{314}s_{410}-s_{321}s_{120}}{s_{234}^2+s_{314}^2+s_{124}^2+s_{321}^2},\quad
t_c=-\frac{s_{234}s_{230}+s_{314}s_{310}+s_{124}s_{120}}{s_{234}^2+s_{314}^2+s_{124}^2+s_{321}^2},\\
\end{aligned}
\end{equation}
and the coordinates of the central point of  \(K\) are given by
\begin{equation}
\begin{aligned}
a_c=\frac{s_{314}s_{430}-s_{124}s_{420}-s_{321}s_{230}}{s_{410}^2+s_{420}^2+s_{430}^2+s_{230}^2+s_{310}^2+s_{120}^2},\quad
b_c=\frac{s_{124}s_{410}-s_{234}s_{430}-s_{321}s_{310}}{s_{410}^2+s_{420}^2+s_{430}^2+s_{230}^2+s_{310}^2+s_{120}^2},\\
c_c=\frac{s_{234}s_{420}-s_{314}s_{410}-s_{321}s_{120}}{s_{410}^2+s_{420}^2+s_{430}^2+s_{230}^2+s_{310}^2+s_{120}^2},\quad
h_c=\frac{s_{234}s_{230}+s_{314}s_{310}+s_{124}s_{120}}{s_{410}^2+s_{420}^2+s_{430}^2+s_{230}^2+s_{310}^2+s_{120}^2}.\\
\end{aligned}
\end{equation}
In Euclidean space, the distances of the central points of the line \(K^*\) and the plane \(K^I\)
are related in the usual way.
Note that the line \(K^*\) shares a common point at infinity with the attitude of the vector 
\((s_{234},s_{314},s_{124},s_{321})\) in \T{4}.
So, points that lie on \(K^*\) satisfy \((x,y,z,t)=\pm(s_{234},s_{314},s_{124},s_{321})\tau+(x_c,y_c,z_c,t_c)\),
where \(\tau\in\R{}\) and \((x_c,y_c,z_c,t_c)\) is the central point of the line.
The top-down orientation of the line is provided by
\(\Id(\mb{\Lambda}_0)=s_{234}\e^{234} +s_{314}\e^{314}+ s_{124}\e^{124} +s_{321}\e^{321}\)
seen as a trivector in \T{4},
while the bottom-up orientation is provided by \(\tb{x}=-s_{234}\e^1 -s_{314}\e^2- s_{124}\e^3 -s_{321}\e^4\)
seen as a vector in \T{4} (\(\tb{x}\) satisfies \(\J(\mb{\Lambda}_0)=\e^0\vee\tb{x}\) as in the lower-dimensional cases).

Any trivector that can be written as \(\e_0\wedge\mb{\sigma}\), where \(\mb{\sigma}\) is a plane, 
dually represents a line at infinity in \T{4}.
This line can be thought of as the intersection of the plane \(\mb{\sigma}\) and the hyperplane at infinity.
The top-down orientation of the line at infinity coincides with the orientation of  \(\mb{\sigma}\).
A line at infinity in \T{4} could be depicted as a small planar segment, indicating all possible directions towards the line at infinity,
and a 2-dimensional arc complementary to the planar segment to indicate the line's top-down orientation.
This depiction would be similar to that of a line at infinity in \T{3}, except that in \T{4}
the top-down orientation is given by a complementary arc rather than a complementary arrow.

The subspace represented by a vector \(\tb{s}\) in \R{5*} is found by solving \(\tb{a}\wedge\tb{s}=0\) for
\(\tb{a}=d\e_0+a\e_1+b\e_2+c\e_3+h\e_4\).
For example, \(\tb{s}=3\e_0-\e_1+2\e_4\) represents a 1-dimensional subspace in \R{5*} given by
\begin{multline*}
(d\e_0+a\e_1+b\e_2+c\e_3+h\e_4)\wedge (3\e_0-\e_1+2\e_4)=
(3a+d)\e_{10} + (3h -2d)\e_{40}
+b\e_{12}-c\e_{31}\\ 
-(2a+h)\e_{41}
-2b\e_{42}-2c\e_{43}=0,
\end{multline*}
which implies \(3a=-d\), \(3h=2d\), \(2a=-h\), \(b=0\), \(c=0\) and for \(d=1\) yields a point in \T{4*} at \((-\tfrac{1}{3},0,0,\tfrac{2}{3})\).
So, \(\tb{s}=3\e_0-\e_1+2\e_4\) dually represents the hyperplane \(-\tfrac{1}{3}x +\tfrac{2}{3}t+1=0\), which is the same as
\(-x +2t+3=0\).
The top-down orientation of this hyperplane is provided by \(\tb{x}=-\e^1+2\e^4 \), 
viewed as a vector in \T{4}.
\(d\e_0\) directly represents the origin of \T{4*} and, therefore, dually represents the hyperplane at infinity in \T{4}.
I assume it to be oriented towards the origin if \(d>0\) and away from the origin if \(d<0\).

Let \(Q\) be a point in \T{4*} directly represented by  \(\tb{s}=d\e_0+a\e_1+b\e_2+c\e_3+h\e_4\).
It is dual to the hyperplane \(Q^*\) defined by \(ax+by+cz+ht+d=0\).
The central point of \(Q^*\) is defined as a point in \T{4}
where \(Q^*\) intersects the line passing through \(Q^I\)  and the origin of \T{4}.
It is given by
\begin{equation}
C(Q^*)=-\frac{(a,b,c,h)}{a^2+b^2+c^2+h^2}.
\end{equation}
In Euclidean space, \(d_{Q_I}d_{C(Q^*)}=1\) for the distance \(d_{Q_I}\) between \(Q^I\) and the origin
and the distance \(d_{C(Q^*)}\) between the central point of the hyperplane \(Q^*\) and the origin.

A four-vector \(\tb{P}\) represents a 4-dimensional hyperplane in \R{5*}, which can be obtained by solving \(\tb{a}\wedge\tb{P}=0.\)
For example, for 
\(\tb{P}=2\e_{1234}+\e_{2340}+3\e_{3140}+4\e_{1240}-\e_{3210}\), I get
\[
(d\e_0+a\e_1+b\e_2+c\e_3+h\e_4)\wedge(2\e_{1234}+\e_{2340}+3\e_{3140}+4\e_{1240}-\e_{3210})=(2d+a+3b+4c-h)\e_{01234}=0,
\]
which implies \(2d+a+3b+4c-h=0\).
The 4-dimensional hyperplane's intersection with \T{4*} is obtained by setting \(d=1\), which gives a 3-dimensional hyperplane
defined by \(2+a+3b+4c-h=0\) or \(\tfrac{1}{2}a+\tfrac{3}{2}b+2c-\tfrac{1}{2}h+1=0\).
This hyperplane corresponds to the point \((\tfrac{1}{2},\tfrac{3}{2},2,-\tfrac{1}{2})\) in the target space \T{4}.
The point is positively oriented.
In general, a point \(\tb{P}=w\e_{1234}+x\e_{2340}+y\e_{3140}+z\e_{1240}+t\e_{3210}\) is assumed to be
positively oriented if \(w>0\) and negatively oriented otherwise.
I will use the same terminology for the top-down and bottom-up orientation, since the positive top-down orientation
corresponds to the positive bottom-up orientation (the same for the negative orientation).

Any four-vector that can be written as \(\tb{N}=\e_0\wedge\mb{\Lambda}\), where \(\mb{\Lambda}\) is a line 
(it is sufficient to consider lines passing through the origin),
dually represents a point at infinity. 
It is located at the intersection of the line \(\mb{\Lambda}\) and the hyperplane at infinity and its orientation
coincides with the orientation of \(\mb{\Lambda}\).
It could be depicted as a small linear segment indicating the direction towards the point at infinity 
and a complementary trivector, an oriented sphere, to indicate the point's top-down orientation.
This depiction would be similar to that of a point at infinity in \T{3},
except that in \T{4} a point at infinity is oriented with a complementary oriented sphere rather than an arc.

The wedge and the join possess the following properties:\\
\(\tb{a}\wedge\tb{b}\) is a plane at the intersection of the hyperplanes \(\tb{a}\) and \(\tb{b}\)\\
\(\tb{a}\wedge\tb{b}\) is a plane at infinity if \(\tb{a}\) and \(\tb{b}\) belong to the same stack but do not coincide\\
\(\tb{a}\wedge\tb{b}=0\) if \(\tb{a}\) and \(\tb{b}\) represent the same hyperplane\\
\(\tb{a}\vee\tb{b}=0\) \\
\(\tb{a}\wedge\mb{\sigma}\) is a line where the hyperplane \(\tb{a}\) and the plane \(\mb{\sigma}\) intersect \\
\(\tb{a}\wedge\mb{\sigma}=0\) if the plane \(\mb{\sigma}\) is in the hyperplane \(\tb{a}\)\\
\(\tb{a}\vee\mb{\sigma}=0\) \\
\(\tb{a}\wedge\mb{\Lambda}\) is a point where the hyperplane \(\tb{a}\) and the line \(\mb{\Lambda}\) intersect\\
\(\tb{a}\wedge\mb{\Lambda}=0\) if the line \(\mb{\Lambda}\) is in the hyperplane \(\tb{a}\)\\
\(\tb{a}\vee\mb{\Lambda}=0\) \\
\(\tb{a}\wedge\tb{b}\wedge\tb{c}\) is a line at the intersection of the three hyperplanes \(\tb{a}\), \(\tb{b}\), and \(\tb{c}\)\\
\(\tb{a}\wedge\tb{b}\wedge\tb{c}\wedge\tb{d}\) 
is a point at the intersection of the four hyperplanes \(\tb{a}\), \(\tb{b}\), \(\tb{c}\), and \(\tb{d}\) \\
\(\tb{a}\wedge\tb{P}=0\) if the point \(\tb{P}\) lies in the hyperplane \(\tb{a}\)\\
\(\mb{\sigma}\wedge\mb{\rho}\) is a point where the planes  \(\mb{\sigma}\) and  \(\mb{\rho}\) intersect\\
\(\mb{\sigma}\wedge\mb{\rho}=0\) if the planes \(\mb{\sigma}\) and  \(\mb{\rho}\) have a common line or coincide\\
\(\mb{\sigma}\wedge\mb{\Lambda}=0\) if the plane  \(\mb{\sigma}\) and the line \(\mb{\Lambda}\) intersect\\
\(\mb{\sigma}\wedge\tb{P}=0\)\\
\(\mb{\sigma}\vee\tb{P}\) is a hyperplane that contains both the plane \(\mb{\sigma}\) and the point \(\tb{P}\)\\
\(\mb{\sigma}\vee\tb{P}=0\) if the point \(\tb{P}\) lies on the  plane  \(\mb{\sigma}\)\\
\(\mb{\Lambda}\wedge\mb{\Phi}=0\)\\
\(\mb{\Lambda}\vee\mb{\Phi}\) is a hyperplane that contains both  lines \(\mb{\Lambda}\) and \(\mb{\Phi}\)\\
\(\mb{\Lambda}\vee\mb{\Phi}=0\) if the lines \(\mb{\Lambda}\) and \(\mb{\Phi}\) lie in the same plane\\
\(\mb{\Lambda}\wedge\tb{P}=0\)\\
\(\mb{\Lambda}\vee\tb{P}\) is a plane that contains both the point \(\tb{P}\) and the line \(\mb{\Lambda}\)\\
\(\mb{\Lambda}\vee\tb{P}=0\) if the point \(\tb{P}\) lies on the line \(\mb{\Lambda}\)\\
\(\tb{P}\wedge\tb{Q}=0\) \\
\(\tb{P}\vee\tb{Q}\) is a line that passes through the points \(\tb{P}\) and \(\tb{Q}\)\\
\(\tb{P}\vee\tb{Q}=0\) if \(\tb{P}\) and \(\tb{Q}\) represent the same point\\
\(\tb{P}\vee\tb{Q}\vee\tb{R}\) is a plane that passes through the  points \(\tb{P}\), \(\tb{Q}\), and \(\tb{R}\)\\
\(\tb{P}\vee\tb{Q}\vee\tb{R}\vee\tb{S}\) is a hyperplane that passes through the points \(\tb{P}\), \(\tb{Q}\), \(\tb{R}\), and \(\tb{S}\)

The following equalities hold:
\begin{equation}
\begin{aligned}
&\mb{\sigma}\wedge\mb{\Lambda}=(\mb{\sigma}\vee\mb{\Lambda})\I,\\
&\tb{a}\wedge\tb{P}=(\tb{a}\vee\tb{P})\I,
\end{aligned}
\end{equation}
where
 \begin{equation}
\mb{\sigma}\vee\mb{\Lambda}=
-(p_{10} s_{234}+p_{20} s_{314}+p_{30} s_{124}+p_{40} s_{321}
+p_{23} s_{410}+p_{31} s_{420}+p_{12} s_{430}
+p_{41} s_{230}+p_{42} s_{310}+p_{43} s_{120})
\end{equation}
and
\begin{equation}
\tb{a}\vee\tb{P}=dw+ax+by+cz+ht.
\end{equation} 

Summary (from \(\bigwedge\R{5*}\) to \T{4})


Four-vectors \(\rightarrow\) Points in \T{4}: \\
\(\tb{P}=\e_{1234}+x\e_{2340}+y\e_{3140}+z\e_{1240} + t\e_{3210} \rightarrow\) a positive point at \((x,y,z,t)\). \\
\(\alpha\tb{P}\), \(\alpha>0\) \(\rightarrow\) same as above but with a different weight.\\
\(-\tb{P}\) \(\rightarrow\) same as \(\tb{P}\)  but with the opposite orientation (negative).\\
\(\tb{N}=\e_0\wedge\mb{\Lambda}\) \(\rightarrow\) a point at infinity 
that lies on the line \(\mb{\Lambda}\);
the top-down  orientation is given by \(\Id(\mb{\Lambda}_0)=s_{234}\e^{234}+s_{314}\e^{314}+s_{124}\e^{124}+s_{321}\e^{321}\)
and the bottom-up orientation by \(\J(\e_0\wedge\mb{\Lambda}_0)=-s_{234}\e^1-s_{314}\e^2-s_{124}\e^3-s_{321}\e^4\).\\
\(\alpha\tb{N}\), \(\alpha>0\) \(\rightarrow\) same as above but with a different weight.\\
\(-\tb{N}\) \(\rightarrow\) same as \(\tb{N}\) but with the opposite orientation.

Simple trivectors \(\rightarrow\) Lines in \T{4}: \\
\(\mb{\Lambda}_0=s_{234}\e_{234}+s_{314}\e_{314}+s_{124}\e_{124}+s_{321}\e_{321}\) \(\rightarrow\)
a line passing through the origin of \T{4} and coinciding with the linear subspace of the vector \((s_{234},s_{314},s_{124},s_{321})\);
the top-down  orientation is given by \(\Id(\mb{\Lambda}_0)=s_{234}\e^{234}+s_{314}\e^{314}+s_{124}\e^{124}+s_{321}\e^{321}\)
and the bottom-up orientation by \(\J(\e_0\wedge\mb{\Lambda}_0)=-s_{234}\e^1-s_{314}\e^2-s_{124}\e^3-s_{321}\e^4\).\\
\(\mb{\Lambda}\), with \(\mb{\Lambda}\vee\mb{\Lambda}=0\) \(\rightarrow\) 
a finite line in \T{4} whose points are given by \((x,y,z,t)=\pm(s_{234},s_{314},s_{124},s_{321})\tau+(x_c,y_c,z_c,t_c)\), 
where \(\tau\in\R{}\) and \((x_c,y_c,z_c,t_c)\)
is the centre of the line;
the orientation of \(\mb{\Lambda}\) is the same as that of  \(\mb{\Lambda}_0\).\\
\(\mb{\Phi}=\e_0\wedge\mb{\sigma}_0\) \(\rightarrow\) 
a line at infinity in \T{4} which lies in the plane \(\mb{\sigma}_0\); the orientation of \(\mb{\Phi}\) is the same as that of \(\mb{\sigma}_0\).

Simple bivectors \(\rightarrow\) Planes in \T{4}: \\
\(\mb{\sigma}_0=p_{23}\e_{23}+p_{31}\e_{31}+p_{12}\e_{12}+p_{41}\e_{41}+p_{42}\e_{42}+p_{43}\e_{43}\), with 
\(\mb{\sigma}_0\wedge\mb{\sigma}_0=0\) \(\rightarrow\)
a plane passing through the origin of \T{4} whose attitude coincides with that of the bivector 
\(B_\parallel=p_{41}\e^{23}+p_{42}\e^{31}+p_{43}\e^{12}+p_{23}\e^{41}+p_{31}\e^{42}+p_{12}\e^{43}\);
the top-down orientation is given by the bivector \(\Id(\mb{\sigma}_0)=B_\perp=p_{23}\e^{23}+p_{31}\e^{31}+p_{12}\e^{12}+p_{41}\e^{41}+p_{42}\e^{42}+p_{43}\e^{43} \)
and the bottom-up orientation by \(\J(\e_0\wedge\mb{\sigma}_0)=B_\parallel\).\\
\(\mb{\sigma}\), \(\mb{\sigma}\wedge\mb{\sigma}=0\) \(\rightarrow\) 
a finite plane in \T{4} whose points are given by \((x,y,z,t)=\tau\tb{x}+\rho\tb{y}+(x_c,y_c,z_c,t_c)\), 
where \(\tb{x}\wedge\tb{y}=B_\parallel\), \(\tau,\rho\in\R{}\) and \((x_c,y_c,z_c,t_c)\)
is the centre of the plane (vectors \(\tb{x}\) and \(\tb{y}\) are  in \T{4} rather than \R{5});
the orientation of \(\mb{\sigma}\) is the same as that of \(\mb{\sigma}_0\).\\
\(\mb{\eta}=\e_0\wedge\tb{a}_0\) where \(\tb{a}_0=a\e_1+b\e_2+c\e_3+h\e_4\) \(\rightarrow\) a plane at infinity in \T{4} which lies in the hyperplane \(ax+by+cz+ht=0\);
its top-down orientation is given by \(\Id(\tb{a}_0)=a\e^1+b\e^2+c\e^3+h\e^4\) and bottom-up orientation by \(\J(\e_0\wedge\tb{a}_0)=-a\e^{234}-b\e^{314}-c\e^{124}-h\e^{321}\).

Vectors \(\rightarrow\) finite hyperplanes in \T{4} and the hyperplane at infinity:\\
\(\tb{a}=d\e_0+a\e_1+b\e_2+c\e_3+h\e_4\) \(\rightarrow\) the hyperplane \(ax+by+cz+ht+d=0\);
its top-down orientation is given by \(\Id(\tb{a}_0)=a\e^1+b\e^2+c\e^3+h\e^4\), where \(\tb{a}_0=a\e_1+b\e_2+c\e_3+h\e_4\),
and the bottom-up orientation by \(\J(\e_0\wedge\tb{a}_0)=-a\e^{234}-b\e^{314}-c\e^{124}-h\e^{321}\).\\
\(\alpha\tb{a}\), \(\alpha>0\) \(\rightarrow\) same hyperplane with a different weight.\\
\(-\tb{a}\) \(\rightarrow\) same hyperplane with the opposite orientation.\\ 
\(d\e_0\) \(\rightarrow\) the hyperplane at infinity in \T{4} (in the top-down view, towards the origin if \(d>0\) or away from the origin if \(d<0\)).


\subsection{The metric and Clifford algebra}

For vectors \(\tb{a}_1\) and \(\tb{a}_2\) in \R{5*}, the inner product is given by
\begin{equation}
\tb{a}_1\cdot\tb{a}_2=d_1d_2\e_0\cdot\e_0+a_1a_2\e_1\cdot\e_1+b_1b_2\e_2\cdot\e_2 + c_1c_2\e_3\cdot\e_3+h_1h_2\e_4\cdot\e_4,
\end{equation}
where \(\tb{a}_1=d_1\e_0+a_1\e_1+b_1\e_2+c_1\e_3+h_1\e_4\) and  \(\tb{a}_2=d_2\e_0+a_2\e_1+b_2\e_2+c_2\e_3+h_2\e_4\),
with the following list of options for the inner product of the standard basis vectors:
\begin{center}
\begin{tabular}{ccccccl}
\(\e_0\cdot\e_0\) & \(\e_1\cdot\e_1\) & \(\e_2\cdot\e_2\) & \(\e_3\cdot\e_3\) & \(\e_4\cdot\e_4\) & \T{4} & \\
    0 & 1 & 1 & 1 &     1 &\E{4}& Euclidean \\ 
    1 & 1 & 1 & 1 &     1 & \El{4} & Elliptic \\
\m1 & 1 & 1 & 1 &     1 & \Hy{4} & Hyperbolic \\
    0 & 1 & 1 & 1 & \m1 & \M{4} & Minkowski (pseudo-Euclidean)  \\
    1 & 1 & 1 & 1 & \m1 & \dS{4} & de-Sitter \\
\m1 & 1 & 1 & 1 & \m1 & \AdS{4} & Anti de-Sitter \\
\end{tabular}
\end{center}
For instance, \(\tb{a}_1\cdot\tb{a}_2=-d_1d_2+a_1a_2+b_1b_2+c_1c_2+h_1h_2\) if the metric is hyperbolic.
Adopting one of these options results in a specific metric geometry in \T{4}.

In \(\bigwedge\R{5*}\), the geometric product of a \(k\)-vector \(A_k\) and an \(l\)-vector \(B_l\)  includes at most three multivectors
of grades \(|k-l|\), \(|k-l|+2\), and \(k+l\).
It exhibits the following properties:
\[
\begin{split}
&\tb{a}\tb{b}=\tb{a}\cdot\tb{b}+\tb{a}\wedge\tb{b} \textrm{ (scalar + bivector)},\quad
\tb{a}\cdot\tb{b}=\tb{b}\cdot\tb{a},\;
\tb{a}\wedge\tb{b}=-\tb{b}\wedge\tb{a},\\
&\tb{a}\mb{\sigma}=\tb{a}\cdot\mb{\sigma}+\tb{a}\wedge\mb{\sigma} \textrm{ (vector + trivector)},\quad
\tb{a}\cdot\mb{\sigma}=-\mb{\sigma}\cdot\tb{a},\;
\tb{a}\wedge\mb{\sigma}=\mb{\sigma}\wedge\tb{a},\\
&\tb{a}\mb{\Lambda}=\tb{a}\cdot\mb{\Lambda}+\tb{a}\wedge\mb{\Lambda} \textrm{ (bivector + four-vector)},\quad
\tb{a}\cdot\mb{\Lambda}=\mb{\Lambda}\cdot\tb{a},\;
\tb{a}\wedge\mb{\Lambda}=-\mb{\Lambda}\wedge\tb{a},\\
&\tb{a}\tb{P}=\tb{a}\cdot\tb{P}+\tb{a}\wedge\tb{P} \textrm{ (trivector + pseudoscalar)},\quad
\tb{a}\cdot\tb{P}=-\tb{P}\cdot\tb{a},\;
\tb{a}\wedge\tb{P}=\tb{P}\wedge\tb{a},\\
&\mb{\sigma}\mb{\rho}=\mb{\sigma}\cdot\mb{\rho}+\mb{\sigma}\times\mb{\rho}+\mb{\sigma}\wedge\mb{\rho} 
\textrm{ (scalar + bivector + four-vector)},
\quad\mb{\sigma}\cdot\mb{\rho}=\mb{\rho}\cdot\mb{\sigma},\;
\mb{\sigma}\wedge\mb{\rho}=\mb{\rho}\wedge\mb{\sigma},\\
&\mb{\sigma}\mb{\Lambda}=\mb{\sigma}\cdot\mb{\Lambda}+\mb{\sigma}\times\mb{\Lambda} + \mb{\sigma}\wedge\mb{\Lambda}
\textrm{ (vector + trivector + pseudoscalar)},\quad
\mb{\sigma}\cdot\mb{\Lambda}=\mb{\Lambda}\cdot\mb{\sigma},\;
\mb{\sigma}\wedge\mb{\Lambda}=\mb{\Lambda}\wedge\mb{\sigma},\\
&\mb{\sigma}\tb{P}=\mb{\sigma}\cdot\tb{P}+\mb{\sigma}\times\tb{P} \textrm{ (bivector + four-vector)},\quad
\mb{\sigma}\cdot\tb{P}=\tb{P}\cdot\mb{\sigma},\\
&\mb{\Lambda}\mb{\Phi}=\mb{\Lambda}\cdot\mb{\Phi}+\mb{\Lambda}\times\mb{\Phi}+(\mb{\Lambda}\vee\mb{\Phi})\I 
\textrm{ (scalar + bivector + four-vector)},
\quad\mb{\Lambda}\cdot\mb{\Phi}=\mb{\Phi}\cdot\mb{\Lambda},\;
\mb{\Lambda}\vee\mb{\Phi}=\mb{\Phi}\vee\mb{\Lambda},\\
&\mb{\Lambda}\tb{P}=\mb{\Lambda}\cdot\tb{P}+(\mb{\Lambda}\vee\tb{P})\I \textrm{ (vector + trivector)},
\quad\mb{\Lambda}\cdot\tb{P}=-\tb{P}\cdot\mb{\Lambda},\;
\mb{\Lambda}\vee\tb{P}=\tb{P}\vee\mb{\Lambda},\\
&\tb{P}\tb{Q}=\tb{P}\cdot\tb{Q}+\tb{P}\times\tb{Q} \textrm{ (scalar + bivector)},
\quad\tb{P}\cdot\tb{Q}=\tb{Q}\cdot\tb{P}.
\end{split}
\]
\(M\I=M\cdot\I\) and  \(M\I=\I M\) for any multivector \(M\).

Note the following relations:
\(\tb{a}\vee\tb{P}=\tb{P}\vee\tb{a}\),\;
\(\mb{\sigma}\vee\tb{P}=-\tb{P}\vee\mb{\sigma}\),\;
\(\mb{\Lambda}\vee\tb{P}=\tb{P}\vee\mb{\Lambda}\),\;
\(\tb{P}\vee\tb{Q}=-\tb{Q}\vee\tb{P}\),\;
\(\mb{\sigma}\vee\mb{\Lambda}=\mb{\Lambda}\vee\mb{\sigma}\),\;
\(\mb{\sigma}\vee\tb{P}=-\tb{P}\vee\mb{\sigma}\),\;
\(\mb{\Lambda}\vee\mb{\Phi}=\mb{\Phi}\vee\mb{\Lambda}\).

It follows that
\[
\begin{split}
&\tb{a}\cdot\tb{b}=\tfrac{1}{2}(\tb{a}\tb{b}+\tb{b}\tb{a}),\quad 
\tb{a}\wedge\tb{b}=\tfrac{1}{2}(\tb{a}\tb{b}-\tb{b}\tb{a}),\quad \tb{a}\times\tb{b}=\tb{a}\wedge\tb{b},  \\
&\tb{a}\cdot\mb{\sigma}=\tfrac{1}{2}(\tb{a}\mb{\sigma}-\mb{\sigma}\tb{a}),\quad
\tb{a}\wedge\mb{\sigma}=\tfrac{1}{2}(\tb{a}\mb{\sigma}+\mb{\sigma}\tb{a}),\quad  
\tb{a}\times\mb{\sigma}=\tb{a}\cdot\mb{\sigma}  \\
&\tb{a}\cdot\mb{\Lambda}=\tfrac{1}{2}(\tb{a}\mb{\Lambda}+\mb{\Lambda}\tb{a}),\quad
\tb{a}\wedge\mb{\Lambda}=\tfrac{1}{2}(\tb{a}\mb{\Lambda}-\mb{\Lambda}\tb{a}),\quad  
\tb{a}\times\mb{\Lambda}=\tb{a}\wedge\mb{\Lambda}  \\
&\tb{a}\cdot\tb{P}=\tfrac{1}{2}(\tb{a}\tb{P}-\tb{P}\tb{a}),\quad
\tb{a}\wedge\tb{P}=\tfrac{1}{2}(\tb{a}\tb{P}+\tb{P}\tb{a}),\quad  \tb{a}\times\tb{P}=\tb{a}\cdot\tb{P}  \\
&\mb{\sigma}\cdot\mb{\rho}+\mb{\sigma}\wedge\mb{\rho}
=\tfrac{1}{2}(\mb{\sigma}\mb{\rho}+\mb{\rho}\mb{\sigma}),\quad
\mb{\sigma}\times\mb{\rho}=\tfrac{1}{2}(\mb{\sigma}\mb{\rho}-\mb{\rho}\mb{\sigma}),\\
&\mb{\sigma}\cdot\mb{\Lambda}+\mb{\sigma}\wedge\mb{\Lambda}
=\tfrac{1}{2}(\mb{\sigma}\mb{\Lambda}+\mb{\Lambda}\mb{\sigma}),\quad
\mb{\sigma}\times\mb{\Lambda}=\tfrac{1}{2}(\mb{\sigma}\mb{\Lambda}-\mb{\Lambda}\mb{\sigma}),\\
&\mb{\sigma}\cdot\tb{P}=\tfrac{1}{2}(\mb{\sigma}\tb{P}+\tb{P}\mb{\sigma}),\quad
\mb{\sigma}\times\tb{P}=\tfrac{1}{2}(\mb{\sigma}\tb{P}-\tb{P}\mb{\sigma}),\\
&\mb{\Lambda}\cdot\mb{\Phi}+(\mb{\Lambda}\vee\mb{\Phi})\I
=\tfrac{1}{2}(\mb{\Lambda}\mb{\Phi}+\mb{\Phi}\mb{\Lambda}),\quad
\mb{\Lambda}\times\mb{\Phi}=\tfrac{1}{2}(\mb{\Lambda}\mb{\Phi}-\mb{\Phi}\mb{\Lambda}),\\
&\mb{\Lambda}\cdot\tb{P}=\tfrac{1}{2}(\mb{\Lambda}\tb{P}-\tb{P}\mb{\Lambda}),\quad
(\mb{\Lambda}\vee\tb{P})\I=\tfrac{1}{2}(\mb{\Lambda}\tb{P}+\tb{P}\mb{\Lambda}),\quad
\mb{\Lambda}\times\tb{P}=\mb{\Lambda}\cdot\tb{P}\\
&\tb{P}\cdot\tb{Q}=\tfrac{1}{2}(\tb{P}\tb{Q}+\tb{Q}\tb{P}), \quad
\tb{P}\times\tb{Q}=\tfrac{1}{2}(\tb{P}\tb{Q}-\tb{Q}\tb{P}).\\
\end{split}
\]

The following identities are useful:
\(\mb{\sigma}\times\tb{P}=(\mb{\sigma}\vee\tb{P})\I,\;
\tb{P}\times\tb{Q}=(\tb{P}\vee\tb{Q})\I\).

The reverse and grade involution of multivector \(M=s+\tb{a}+\mb{\sigma}+\mb{\Lambda}+\tb{P}+\alpha\I\) are given by
\(\reverse{M}=s+\tb{a}-\mb{\sigma}-\mb{\Lambda}+\tb{P}+\alpha\I\) and \(\involute{M}=s-\tb{a}+\mb{\sigma}-\mb{\Lambda}+\tb{P}-\alpha\I\),
respectively. 
As in lower dimensions, \(\reverse{AB}=\reverse{B}\reverse{A}\) and \(\involute{AB}=\involute{A}\involute{B}\) 
(the same applies to inner and outer products and the commutator).
I will also make use of \((M)_{14}=s-\tb{a}+\mb{\sigma}+\mb{\Lambda}-\tb{P}+\alpha\I\).
In the subscript notation, I have \(\reverse{M}=(M)_{23}\) and \(\involute{M}=(M)_{135}\).
For a Study number \(N=\alpha+\beta\I\), the conjugate is defined by \(\overbar{N}=\alpha-\beta\I\).

The norm of multivector \(M\) is defined by
\begin{equation}
\norm{M}=|N\overbar{N}|^{\tfrac{1}{8}}, \quad\textrm{where } N=M\reverse{M}(M\reverse{M})_{14}.
\end{equation}
If the norm  is not zero, the inverse of \(M\) is given by
\begin{equation}
M^{-1}=\frac{\reverse{M}(M\reverse{M})_{14}\overbar{N}}{N\overbar{N}}.
\end{equation}
For blades, this simplifies to
\begin{equation}
\norm{A_k}=|A_k\reverse{A}_k|^{\tfrac{1}{2}}, \quad A_k^{-1}=\frac{\reverse{A}_k}{A_k\reverse{A}_k}.
\end{equation}

The algebra of even multivectors in \(\bigwedge\R{5*}\) is 16-dimensional and consists of multivectors
\(
E=u+ \mb{\sigma} +\tb{P},
\)
where \(u\in\R{}\).

\subsection{Euclidean hyperspace \E{4}}

For a hyperplane \(\tb{a}\), a plane \(\mb{\sigma}\) (simple bivector), a line \(\mb{\Lambda}\) (simple trivector), and a point \(\tb{P}\) the norm is given by
\begin{equation}
\begin{split}
&\norm{\tb{a}}=\sqrt{a^2+b^2+c^2+h^2},\quad
\norm{\mb{\sigma}}=\sqrt{p_{23}^2+p_{31}^2+p_{12}^2+p_{41}^2+p_{42}^2+p_{43}^2},\quad\\
&\norm{\mb{\Lambda}}=\sqrt{s_{234}^2+s_{314}^2+s_{124}^2+s_{321}^2},\quad
\norm{\tb{P}}=|w|.
\end{split}
\end{equation}
Unlike the 2- and 3-dimensional cases, \(\tb{P}^2=w^2\) in \E{4}.
Furthermore, \(\tb{a}^2=\norm{\tb{a}}^2\), \(\mb{\sigma}^2=-\norm{\mb{\sigma}}^2\), \(\mb{\Lambda}^2=-\norm{\mb{\Lambda}}^2\).
Note that the above expressions for \(\mb{\sigma}\) and \(\mb{\Lambda}\) do not apply to non-simple bivectors and trivectors.

Possible geometric configurations in \E{4} are similar to those observed in \E{3}.
A notable exception is a configuration involving two planes intersecting at a point if \(\mb{\sigma}\wedge\mb{\rho}\ne0\).
I will call such planes complementary if \(\mb{\sigma}\cdot\mb{\rho}=0\) and \(\mb{\sigma}\times\mb{\rho}=0\).
On the other hand, if two planes intersect along a line, there is a hyperplane that contains both planes
and the configuration is 3-dimensional in that hyperplane.
Two planes intersecting along a line can be perpendicular in the 3-dimensional sense if \(\mb{\sigma}\cdot\mb{\rho}=0\),
but they are not complementary.

The plane at infinity \(\e_0\wedge\tb{a}\) has the same top-down orientation as a hyperplane \(\tb{a}\),
while the point \(\tb{a}\I\) is at infinity in the direction perpendicular to \(\tb{a}\).
The point \(\tb{a}\I\) is called the polar point of the hyperplane \(\tb{a}\).
The line \(\tb{a}\cdot\tb{P}\) passes through  \(\tb{P}\) and is perpendicular to  \(\tb{a}\), so that
\(\tb{a}\I\) lies on the line \(\tb{a}\cdot\tb{P}\) for any point \(\tb{P}\).
The distance between \(\tb{a}\) and \(\tb{P}\) is given by \(|\tb{a}\vee\tb{P}|\) if both \(\tb{a}\) and \(\tb{P}\) are normalised.
The plane \(\tb{a}\cdot\mb{\Lambda}\) passes through  \(\mb{\Lambda}\) and is perpendicular to  \(\tb{a}\).
The angle \(\alpha\) between a normalised line \(\mb{\Lambda}\) and a normalised hyperplane \(\tb{a}\) satisfies
\(\norm{\tb{a}\cdot\mb{\Lambda}}=\cos{\alpha}\).
The hyperplane \(\tb{a}\cdot\mb{\sigma}\) passes through \(\mb{\sigma}\) and is perpendicular to \(\tb{a}\), while
\(\norm{\tb{a}\cdot\mb{\sigma}}=\cos{\alpha}\) gives the angle \(\alpha\) between them if they are normalised.
The dot product \(\tb{a}\cdot\tb{b}=\cos{\alpha}\) gives the angle \(\alpha\) between \(\tb{a}\) and \(\tb{b}\) if they are normalised.

The hyperplane \(\mb{\Lambda}\cdot\tb{P}\) passes through  \(\tb{P}\) and is perpendicular to  \(\mb{\Lambda}\).
The plane \(\mb{\Lambda}\vee\tb{P}\) passes through both \(\mb{\Lambda}\) and \(\tb{P}\),
while \(\norm{\mb{\Lambda}\vee\tb{P}}\) gives the distance between the line and the point if they are normalised.
The line \(\tb{P}\vee\tb{Q}\) passes through both points and \(\norm{\tb{P}\vee\tb{Q}}\) gives the distance between them if they are normalised.
The plane \(\mb{\sigma}\cdot\tb{P}\) passes through \(\tb{P}\) and is complementary to \(\mb{\sigma}\).
The planes \(\mb{\sigma}\) and \(\mb{\sigma}\cdot\tb{P}\) intersect at the projection of \(\tb{P}\) on \(\mb{\sigma}\).
The hyperplane \(\mb{\sigma}\vee\tb{P}\) contains both \(\mb{\sigma}\) and \(\tb{P}\),
while the distance between \(\mb{\sigma}\) and \(\tb{P}\) is given by \(\norm{\mb{\sigma}\vee\tb{P}}\), provided they are normalised.

The hyperplane \(\mb{\Lambda}\vee\mb{\Phi}\) contains both  lines \(\mb{\Lambda}\) and \(\mb{\Phi}\),
but \(\mb{\Lambda}\vee\mb{\Phi}=0\) if the lines intersect.
In \E{4}, the configuration involving two lines is 3-dimensional in the hyperplane \(\mb{\Lambda}\vee\mb{\Phi}\).
It is essentially identical to that observed in \E{3}, with only minor modifications.
For instance, \(\mb{\Lambda}\cdot\mb{\Phi}=-\cos{\alpha}\) gives the angle \(\alpha\) between the lines.
It is equal to the angle between \(\mb{\Lambda}\) and the translational projection of \(\mb{\Phi}\) on \(\mb{\Lambda}\).
Moreover, \(\norm{\mb{\Lambda}\vee\mb{\Phi}}=r\sin{\alpha}\), where \(r\) is the distance between the lines.
The commutator \(\mb{\Lambda}\times\mb{\Phi}\) is a non-simple bivector unless the lines intersect.
Recall that in \E{3} the commutator of two lines can be written as the sum of a finite line
and a line at infinity perpendicular to it.
This structure is completely determined in the hyperplane \(\mb{\Lambda}\vee\mb{\Phi}\), which is essentially a 3-dimensional Euclidean space.
By extending this structure into the forth dimension perpendicular to \(\mb{\Lambda}\vee\mb{\Phi}\),
one obtains a finite plane derived from the finite axis and a plane at infinity derived from the axis at infinity.
This 4-dimensional structure is captured by the commutator \(\mb{\Lambda}\times\mb{\Phi}\),
which can be decomposed into the two complementary planes.

The hyperplane \(\mb{\sigma}\cdot\mb{\Lambda}\) passes through \(\mb{\Lambda}\) and is perpendicular to \(\mb{\sigma}\),
so that  \(\mb{\Lambda}\wedge(\mb{\sigma}\cdot\mb{\Lambda})=0\) and \(\mb{\sigma}\cdot(\mb{\sigma}\cdot\mb{\Lambda})=0\).
On the other hand, \(\mb{\sigma}\wedge(\mb{\sigma}\cdot\mb{\Lambda})\) is a line in the hyperplane \(\mb{\sigma}\cdot\mb{\Lambda}\).
Both lines \(\mb{\Lambda}\)  and \(\mb{\sigma}\wedge(\mb{\sigma}\cdot\mb{\Lambda})\) lie in \(\mb{\sigma}\cdot\mb{\Lambda}\).
In fact, the configuration involving \(\mb{\sigma}\) and \(\mb{\Lambda}\)
is determined by the relationship between the lines \(\mb{\Lambda}\)  and \(\mb{\sigma}\wedge(\mb{\sigma}\cdot\mb{\Lambda})\),
which is 3-dimensional in character.
For instance,
the angle \(\alpha\) between \(\mb{\sigma}\) and  \(\mb{\Lambda}\), which satisfies 
\(\norm{\mb{\sigma}\cdot\mb{\Lambda}}=\cos{\alpha}\) for the normalised objects,
coincides with the angle between \(\mb{\Lambda}\) and the line \(\mb{\sigma}\wedge(\mb{\sigma}\cdot\mb{\Lambda})\).
Moreover, \(|\mb{\sigma}\vee\mb{\Lambda}|=r\sin{\alpha}\) gives the distance \(r\) between the lines and, therefore, the distance
between \(\mb{\Lambda}\) and \(\mb{\sigma}\).
So, there is a unique point on  \(\mb{\Lambda}\) whose distance to  \(\mb{\sigma}\) is minimal.
The commutator \(\mb{\sigma}\times\mb{\Lambda}\)
is a non-simple trivector in general. 

If a trivector \(\mb{\Phi}\) is not simple and finite, \(\mb{\Phi}\vee\mb{\Phi}\ne0\) and \(\mb{\Phi}\cdot\mb{\Phi}\ne0\), 
it can be uniquely decomposed into two simple trivectors, 
\(\mb{\Phi}=\mb{\Phi}_{fc}+\mb{\Phi}_{ic}\),
representing a finite line and a line at infinity, where
\begin{equation}
\mb{\Phi}_{fc}=\left(1-\frac{(\mb{\Phi}\vee\mb{\Phi})\I}{2\mb{\Phi}\cdot\mb{\Phi}}\right)\mb{\Phi}, \textrm{ and }
\mb{\Phi}_{ic}=\frac{(\mb{\Phi}\vee\mb{\Phi})\I}{2\mb{\Phi}\cdot\mb{\Phi}}\mb{\Phi}.
\end{equation}
\(\mb{\Phi}_{fc}\) is the finite component of \(\mb{\Phi}\) and \(\mb{\Phi}_{ic}\) the infinite component.
Both \(\mb{\Phi}_{fc}\) and \(\mb{\Phi}_{ic}\) lie in the 3-dimensional hyperplane \(\mb{\Phi}\vee\mb{\Phi}\).
If \(\mb{\Phi}\) is not simple but \(\mb{\Phi}\cdot\mb{\Phi}=0\), e.g.\ \(\mb{\Phi}=\e_{120}+2\e_{430}\), 
it can be written as \(\mb{\Phi}=\e_0\wedge\mb{\pi}\),
where \(\mb{\pi}\) is a non-simple bivector (\(\mb{\pi}\) can be assumed to pass through the origin).
Then the decomposition of \(\mb{\Phi}\) depends on whether \(\mb{\pi}\) can be decomposed into two complementary planes.
Note that \(\mb{\Phi}\vee\mb{\Phi}\) is a hyperplane at infinity in this case.

Let \(\mb{\pi}\) be a non-simple bivector.
If \((\mb{\pi}\cdot\mb{\pi})^2\ne(\mb{\pi}\wedge\mb{\pi})^2\),
 \(\mb{\pi}\) can be split uniquely into two complementary planes, \(\mb{\pi}=\mb{\pi}_1+\mb{\pi}_2\).
Otherwise the decomposition is not unique.
Squaring \(\mb{\pi}=\mb{\pi}_1+\mb{\pi}_2\) and using the complementarity of \(\mb{\pi}_1\) and \(\mb{\pi}_2\) yields
\begin{equation}
\begin{split}
&\mb{\pi}_1\mb{\pi}_2=\tfrac{1}{2}\mb{\pi}\wedge\mb{\pi},\\
&\mb{\pi}_1^2+\mb{\pi}_2^2=\mb{\pi}\cdot\mb{\pi}.
\end{split}
\end{equation}
Note also that \(\mb{\pi}_1\) and \(\mb{\pi}_2\) commute with \(\mb{\pi}_1\wedge\mb{\pi}_2\) and
\(\mb{\pi}\) commutes with \(\mb{\pi}\wedge\mb{\pi}\).
Squaring the first equation above and combining it with the second yields
the following quadratic equation 
\begin{equation}
\mb{\pi}_{1,2}^4
-(\mb{\pi}\cdot\mb{\pi})\mb{\pi}_{1,2}^2
+\tfrac{1}{4}(\mb{\pi}\wedge\mb{\pi})^2=0
\label{p1 p2}
\end{equation}
for \(\mb{\pi}_1^2\) and \(\mb{\pi}_2^2\).
In Euclidean space \E{4}, \((\mb{\pi}\cdot\mb{\pi})^2 \ge (\mb{\pi}\wedge\mb{\pi})^2\ge0\) is satisfied.
If \((\mb{\pi}\cdot\mb{\pi})^2>(\mb{\pi}\wedge\mb{\pi})^2\) and \((\mb{\pi}\wedge\mb{\pi})^2\ne0\),
Equation~(\ref{p1 p2}) has two solutions
\begin{equation}
\mb{\pi}_{1,2}^2=
\tfrac{1}{2}\left(
\mb{\pi}\cdot\mb{\pi}
\pm
\sqrt{({\mb{\pi}\cdot\mb{\pi}})^2 - (\mb{\pi}\wedge\mb{\pi})^2}
\right).
\end{equation}
Both \(\mb{\pi}_1^2\) and \(\mb{\pi}_2^2\) are non-zero
and, therefore,
\[
\mb{\pi}_1=\frac{1}{2}(\mb{\pi}\wedge\mb{\pi}) \mb{\pi}_2/\mb{\pi}_2^2,
\quad 
\mb{\pi}_2=\frac{1}{2}(\mb{\pi}\wedge\mb{\pi}) \mb{\pi}_1/\mb{\pi}_1^2,
\]
which gives
\begin{equation}
\mb{\pi}_1
=
\mb{\pi}/ (1+\tfrac{1}{2}\mb{\pi}\wedge\mb{\pi} / \mb{\pi}_1^2),
\quad
\mb{\pi}_2
=
\mb{\pi}/(1+\tfrac{1}{2}\mb{\pi}\wedge\mb{\pi} / \mb{\pi}_2^2)
\label{s1 and s2}
\end{equation}
and
\begin{equation}
\mb{\pi}_1
=
\frac{(1-\frac{1}{2}\mb{\pi}\wedge\mb{\pi} / \mb{\pi}_1^2)\mb{\pi}}
{1-\frac{1}{4}(\mb{\pi}\wedge\mb{\pi})^2 / \mb{\pi}_1^4},
\quad
\mb{\pi}_2
=
\frac{(1-\frac{1}{2}\mb{\pi}\wedge\mb{\pi} / \mb{\pi}_2^2)\mb{\pi}}
{1-\frac{1}{4}(\mb{\pi}\wedge\mb{\pi})^2 / \mb{\pi}_2^4}.
\end{equation}

If \((\mb{\pi}\cdot\mb{\pi})^2=(\mb{\pi}\wedge\mb{\pi})^2\), then 
\(\mb{\pi}_1^2=\mb{\pi}_2^2=\tfrac{1}{2}\mb{\pi}\cdot\mb{\pi}\) and 
\(\mb{\pi} = \mb{\pi}_{1,2} (1+\mb{\pi}\wedge\mb{\pi} / \mb{\pi}\cdot\mb{\pi})\)
but multivector \((1+\mb{\pi}\wedge\mb{\pi} / \mb{\pi}\cdot\mb{\pi})\) is not invertible
and, therefore, there is no unique solution (note that in this case \(\mb{\pi}\) is not invertible either).
The decomposition into complementary planes is not unique.
For instance,  \(\mb{\pi}=\e_{12}+\e_{34}\)
can be decomposed into \(\mb{\pi}_1=\e_{12}\) and \(\mb{\pi}_2=\e_{34}\)
or into \(\mb{\pi}'_1=\tfrac{1}{2}(\e_{12}+\e_{34}+\e_{23}+\e_{41})\) 
and \(\mb{\pi}'_2=\tfrac{1}{2}(\e_{12}+\e_{34}-\e_{23}-\e_{41})\).

If \((\mb{\pi}\wedge\mb{\pi})^2=0\) and \(\mb{\pi}\cdot\mb{\pi}\ne0\), then
\(\mb{\pi}\wedge\mb{\pi}\) is a point at infinity and
\(\mb{\pi}_1^2=\mb{\pi}\cdot\mb{\pi}\), \(\mb{\pi}_2^2=0\).
\(\mb{\pi}\) can be decomposed into a finite plane and and a plane at infinity as follows:
\begin{equation}
\displaystyle\mb{\pi}_1=\left(1-\frac{\mb{\pi}\wedge\mb{\pi}}{2\mb{\pi}\cdot\mb{\pi}}\right)\mb{\pi},
\quad
\mb{\pi}_2=\frac{\mb{\pi}\wedge\mb{\pi}}{2\mb{\pi}\cdot\mb{\pi}}\mb{\pi}.
\end{equation}

Finally,  \(\mb{\pi}\cdot\mb{\pi}=0\) implies that \(\mb{\pi}\) is simple (a plane at infinity).

Projection and rejection in \E{4} are defined in the usual way.
For instance, the projection of a hyperplane \(\tb{a}\) on a hyperplane \(\tb{b}\)
is given by \((\tb{a}\cdot\tb{b})\tb{b}^{-1}\), while the rejection is given by
\((\tb{a}\wedge\tb{b})\tb{b}^{-1}\).
If the geometric product consists of three terms,
as is the case for \(\mb{\sigma}\mb{\rho}\),
\(\mb{\sigma}\mb{\Lambda}\), and  \(\mb{\Phi}\mb{\Lambda}\),
one can split the commutator into two components and define two kinds of projection and rejection.
Projection and rejection can be used to define scaling in the same way as in \E{2} and \E{3}.

The reflection of a blade \(B_l\) in an invertible blade \(A_k\)
is given by \((-1)^{kl}A_k B_l A_k^{-1}\) for the top-down reflection and \((-1)^{nk(l-1)} A_k B_l A_k^{-1}\) for the bottom-up reflection 
(\(n\) is the space dimension). 
For instance, the top-down reflection of a hyperplane \(\tb{a}\) in a hyperplane \(\tb{b}\) is given by \(-\tb{b}\tb{a}\tb{b}^{-1}\).
Since \(n=4\) in \E{4}, the bottom-up reflection formula simplifies to \(A_k B_l A_k^{-1}\) in \E{4}.
Similar to other dimensions, all top-down reflections can be understood in terms of reflections in suitable hyperplanes,
On the other hand, the bottom-up reflection are constructed from point reflections.

Any proper motion in \E{4} can be obtained by performing an even number of consecutive reflections in hyperplanes, which square to unity.
A spinor in \E{4} is defined  as the product of an even number of hyperplanes, which all square to unity.
So, any proper motion can be generated by a spinor.
Moreover, any proper motion can be generated by a spinor that has the form \(S=e^A\) where \(A\) is a bivector.
Depending on the properties of the bivector, different kinds of proper motion are obtained.

For multivector \(M\) in \E{4}, the proper motion induced by a spinor \(S\) is obtained by computing
\(SMS^{-1}\).
The spinor \(S\)
assumes one of the following forms:\\
\(S=e^{-\tfrac{1}{2}\lambda\e_0\wedge\tb{a}}\), where \(\norm{\tb{a}}=1\),
yields a translation by \(\lambda\) along lines  perpendicular to the hyperplane \(\tb{a}\) in the direction of the top-down orientation of \(\tb{a}\).\\
\(S=e^{-\tfrac{1}{2}\alpha\mb{\sigma}}\), where \(\mb{\sigma}\) is a finite plane and \(\norm{\mb{\sigma}}=1\),
yields a rotation by the angle \(\alpha\) around the plane \(\mb{\sigma}\),
the axis of rotation,
in the direction of the top-down orientation of \(\mb{\sigma}\).
Note that points and lines that lie on the axis are not affected by the rotation,
so the axis is invariant.\\
\(S=e^{-\tfrac{1}{2}(\alpha-\lambda\tb{a}\I)\mb{\sigma}}\),
where \(\mb{\sigma}\) is a finite plane, \(\norm{\tb{a}}=1\), \(\norm{\mb{\sigma}}=1\), and \(\tb{a}\cdot\mb{\sigma}=0\),
yields a combination of a rotation as above and a translation along lines perpendicular to \(\tb{a}\) in the direction that
depends on the orientation of \(\tb{a}\) and the axis of rotation;
\(\tb{a}\cdot\mb{\sigma}=0\) implies \((\tb{a}\I)\vee\mb{\sigma}=0\),
i.e.\ the point at infinity \(\tb{a}\I\)
lies on the axis of rotation \(\mb{\sigma}\) and, therefore,
 ensures that the translation
occurs in one of the directions parallel to the axis.\\
\(S=e^{-\tfrac{1}{2}(\alpha\mb{\sigma}_1+\beta\mb{\sigma}_2)}\),
where \(\alpha\ne\beta\) and planes \(\mb{\sigma}_1\) and \(\mb{\sigma}_2\)
are normalised and complementary, 
yields a combination of a rotation by the angle \(\alpha\)
around \(\mb{\sigma}_1\)
and a rotation by the angle \(\beta\) around \(\mb{\sigma}_2\).
In other words, this spinor yields a double rotation around two complementary axes.
The point \(\mb{\sigma}_1\wedge\mb{\sigma}_2\) where the axes intersect
is not affected by the double rotation.
Other points and lines that lie on the axes undergo a simple rotation,
that is they do not leave the axes.
For instance, points in \(\mb{\sigma}_1\) are rotated by the angle \(\beta\)
around \(\mb{\sigma}_1\wedge\mb{\sigma}_2\),
and points in \(\mb{\sigma}_2\) are rotated by the angle \(\alpha\).
In this sense the axes of a double rotation are invariant.\\
\(S=e^{-\tfrac{1}{2}\alpha(\mb{\sigma}_1+\mb{\sigma}_2)}\),
where the planes \(\mb{\sigma}_1\) and \(\mb{\sigma}_2\)
are normalised and complementary as above but \(\alpha=\beta\),
yields an isoclinic rotation by the angle \(\alpha\).
It is characterised by many invariant axes since 
\(\mb{\sigma}=\mb{\sigma}_1+\mb{\sigma}_2\) cannot be uniquely decomposed into two complementary planes.
The point \(\mb{\sigma}_1\wedge\mb{\sigma}_2\)
is invariant under the isoclinic rotation.
Note, however, that not every plane that passes through the point \(\mb{\sigma}_1\wedge\mb{\sigma}_2\)
is invariant.\\
So, there are five distinct kinds of proper motions in \E{4}:
translation, simple rotation, rotation and translation in the axis,
double rotation, and isoclinic rotation.

\newpage

\section{1-dimensional geometry}
In a 1-dimensional space, points in the target space correspond to
points in the dual space.
These spaces are embedded in the model  spaces in the standard way.
The standard basis of \R{2} consists of
\(\e^0=(1,0)\) and \(\e^1=(0,1)\), so that \(\tb{x}=w\e^0+x\e^1\),
and the standard basis of \R{2*} consists of
\(\e_0=(1,0)\) and \(\e_1=(0,1)\), so that \(\tb{a}=d\e_0+a\e_1\).
The basis of Grassmann algebra of \(\bigwedge\R{2*}\) consists of:
1,
\(\e_0\), \(\e_1\),
\(\e_{01}\).
The duality transformation and its inverse are defined on the basis by

\begin{tabular*}{0.4\textwidth}{lcccc}
\(X\)   &  \(1\) & \(\e_{0}\) & \(\e_{1}\) & \(\e_{01}\) \\
\cline{1-5} \\[-10pt]
\(\J(X)\)  & \(-\e^{01}\) & \(-\e^{1}\) & \(\e^{0}\) & \(1\)  \\
\end{tabular*}%
\begin{tabular*}{0.4\textwidth}{lcccc}
\(Y\)      & 1   & \(\e^{0}\) & \(\e^{1}\) & \(\e^{01}\)  \\
\cline{1-5} \\[-10pt]
\(\J^{-1}(Y)\) & \(\e_{01}\) & \(\e_{1}\) & \(-\e_{0}\) & \(-1\)  \\
\end{tabular*}

The join of two multivectors in \(\bigwedge\R{2*}\) is defined in the usual way: 
\(\tb{a}\vee\tb{b}=\J^{-1}(\J(\tb{a})\vee\J(\tb{b}))\), which for
\(\tb{a}_1=d_1\e_0+a_1\e_1\) and \(\tb{a}_2=d_2\e_0+a_2\e_1\) yields \(\tb{a}_1\vee\tb{a}_2=-(d_1a_2-a_1d_2)\).
On the other hand, \(\tb{a}_1\wedge\tb{a}_2=(d_1a_2-a_1d_2)\e_{01}\) and, therefore, 
\(\tb{a}\wedge\tb{b}=-(\tb{a}\vee\tb{b})\e_{01}\) for any vectors \(\tb{a}\) and \(\tb{b}\).

\begin{figure}[h]
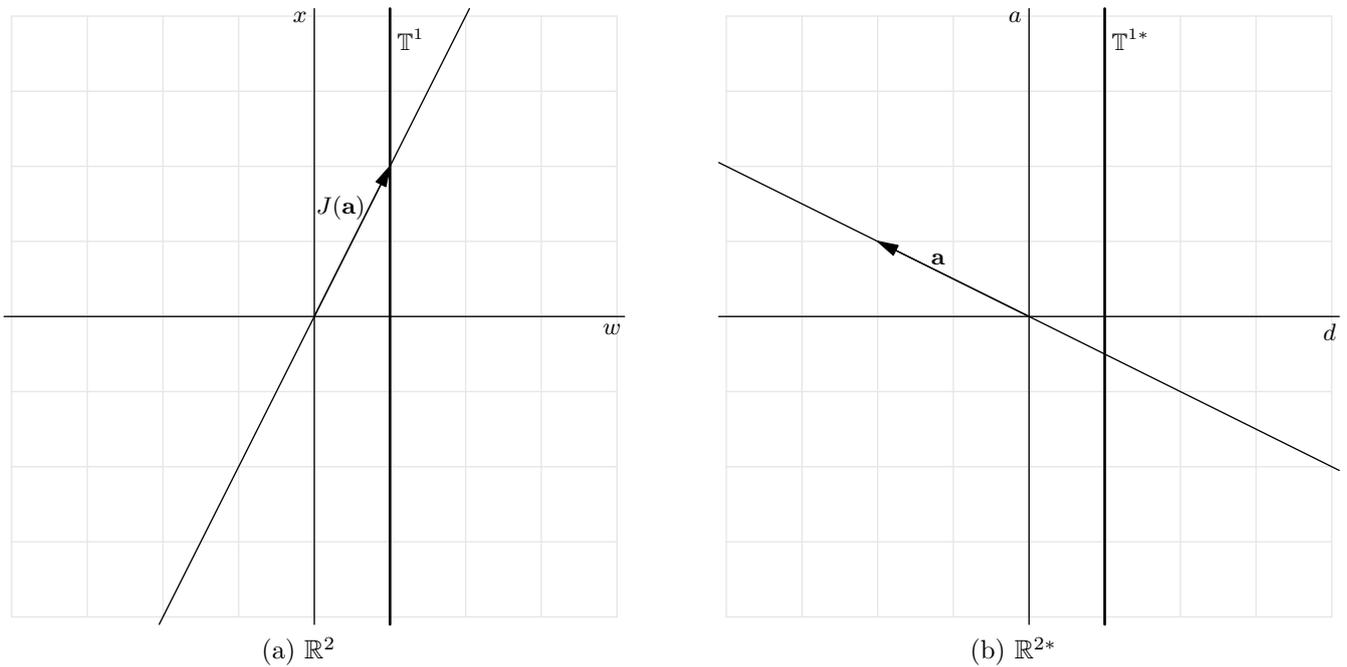

\begin{subfloatenv}{\R{2}}
\begin{asy}
import Figure2D;
Figure f = Figure(xaxis_name="$w$",yaxis_name="$x$");

f.line(e_1-e_0,draw_orientation=false,label="$\mathbb{T}^1$",position=0.05,align=(1,0),pen=currentpen+1);

MV a = Line(0,-2,1);
f.line(a,draw_orientation=false);
draw(Label("$J(\textbf{a})$",0.6,align=(-0.6,1)), (0,0)--(1,2), Arrow);

\end{asy}
\end{subfloatenv}\hfill%
\begin{subfloatenv}{\R{2*}}
\begin{asy}
import Figure2D;
Figure f = Figure(xaxis_name="$d$",yaxis_name="$a$");

f.line(e_1-e_0,draw_orientation=false,label="$\mathbb{T}^{1*}$",position=0.05,align=(1,0),pen=currentpen+1);

MV a = Line(0,1,2);
f.line(a,draw_orientation=false);
draw(Label("$\textbf{a}$",0.6,align=(0,1)), (0,0)--(-2,1), Arrow);

\end{asy}
\end{subfloatenv}
\caption{The duality relationship between points in \T{1} and points in \T{1*}}
\label{target and dual spaces in 1d}
\end{figure}

If \(x\ne0\), then \(\tb{a}=-x\e_0+\e_1\) directly represents a point in the dual space \T{1*} at \(a=-1/x\), 
and dually represents a point in the target space \T{1} at \(x\).
These points satisfy the usual duality relationship 
\begin{equation}
ax+1=0.
\end{equation}
The same point in \T{1} is directly represented by \(\J(\tb{a})=\e^0+x\e^1\)
(see Figure~\ref{target and dual spaces in 1d} where \(a=-\tfrac{1}{2}\) and \(x=2\)).
A finite point dually represented by \(\tb{a}=-x\e_0+\e_1\) will be considered positively oriented
(\(x\) could be positive or negative),
while the same point represented by \(-\tb{a}=x\e_0-\e_1\) will be negatively oriented.
The point shown in Figure~\ref{target and dual spaces in 1d}(a) is positively oriented.

\(\alpha\tb{a}\) dually represents the same point but with a different weight
and \(-\tb{a}\) dually represents the same point but with the opposite orientation. 
\(\alpha\e_{01}\) dually represents \(\alpha\) in \T{1} and a scalar in
\(\bigwedge\R{2*}\) dually represents the whole \T{1}.
The join of two points is a scalar and, therefore, gives the whole space \T{1}.
\(\e_0\) dually represents the point at infinity in \T{1}
and \(\e_1\) dually represents the origin.


Euclidean, elliptic, or hyperbolic space is selected by choosing the appropriate metric:
\begin{center}
\begin{tabular}{cccl}
\(\e_0\cdot\e_0\) & \(\e_1\cdot\e_1\) & \T{1}&\\
 0 & 1 & \E{1}& Euclidean \\ 
 1 &1 & \(E_{1}\) & Elliptic \\
\m1 &1 & \Hy{1} & Hyperbolic \\
\end{tabular}
\end{center}
Since the spaces are 1-dimensional, no kinematic structure can be defined.

Resulting Clifford algebra exhibits the following properties:\\
\(\tb{a}\tb{b}=\tb{a}\cdot\tb{b}+\tb{a}\wedge\tb{b}\) (scalar + pseudoscalar),\(\quad\)
\(\tb{a}\cdot\tb{b}=\tb{b}\cdot\tb{a}\),\;
\(\tb{a}\wedge\tb{b}=-\tb{b}\wedge\tb{a}\),\\
\(\tb{a}\e_{01}=\tb{a}\cdot\e_{01}\) (vector),\(\quad\)
\(\tb{a}\cdot\e_{01}=-\e_{01}\cdot\tb{a}\).\\
The following identities apply:\\
\(\tb{a}\cdot\tb{b}=\tfrac{1}{2}(\tb{a}\tb{b}+\tb{b}\tb{a}), \quad
\tb{a}\wedge\tb{b}=\tfrac{1}{2}(\tb{a}\tb{b}-\tb{b}\tb{a}),\quad
\tb{a}\wedge\tb{b}=-(\tb{a}\vee\tb{b})\e_{01}.
\)

A multivector in \(\bigwedge\R{2*}\) can be written as
\begin{equation}
M=s+d\e_0+a\e_1+\alpha\e_{01}.
\end{equation}
Its norm is defined by
\begin{equation}
\norm{M}=|M\reverse{M}(M\reverse{M})_{1}|^{\tfrac{1}{4}},
\end{equation}
where \(\reverse{M}=s+d\e_0+a\e_1-\alpha\e_{01}\) and \((M)_1=s-d\e_0-a\e_1+\alpha\e_{01}\).
For multivectors with a non-zero norm, the inverse is given by
\begin{equation}
M^{-1}=\frac{\reverse{M}(M\reverse{M})_{1}}{M\reverse{M}(M\reverse{M})_{1}}.
\end{equation}
This simplifies to \(\tb{a}^{-1}=\tb{a}/\tb{a}^2\) for vectors
and \(E^{-1}=\reverse{E}/(E\reverse{E})\) for even multivectors.
The bivector \(\e_{01}\)  is not invertible in \E{1}, but
\(\e_{01}^{-1}=-\e_{01}\) in \El{1}, and \(\e_{01}^{-1}=\e_{01}\) in \Hy{1}.
Note also that \(\e_{01}^2=0\) in \E{1}, \(\e_{01}^2=-1\) in \El{1}, and \(\e_{01}^2=1\) in \Hy{1}.
Since \(\alpha\tb{a}=\tb{a}\alpha\) and \(\beta\e_{01}\tb{a}=-\tb{a}\beta\e_{01}\),
I get \(\tb{a}E=\reverse{E}\tb{a}\) for any vector \(\tb{a}\) and any even multivector \(E=\alpha+\beta\e_{01}\).
If \(E\reverse{E}=1\), I get \(\tb{a}E=E^{-1}\tb{a}\).

The algebra of even multivectors
consists of multivectors \(E=\alpha+\beta\e_{01}\) and is, therefore, 2-dimensional. 
It is isomorphic to complex numbers if the metric is elliptic, which implies \(\e_{01}^2=-1\),
so that the bivector \(\e_{01}\) can play the role of the imaginary unit.

\subsection{Euclidean line \E{1}}

In Euclidean space \E{1}, the norm of \(\tb{a}=d\e_0+a\e_1\) is given by
\(\norm{\tb{a}}=|a|\).
The distance \(r\) between two points dually represented by normalised vectors \(\tb{a}=-x\e_0+\e_1\) and \(\tb{b}=-y\e_0+\e_1\) is given by
\begin{equation}
r=|\tb{a}\vee\tb{b}|=|x-y|.
\end{equation}

The polar point \(\tb{a}\e_{01}\) of a finite point \(\tb{a}\) is the point at infinity, weighted by the weight of \(\tb{a}\).

A point is proper if the square of the representing vector is positive, \(\tb{a}^2>0\).
Any finite point is proper in \E{1}.
Spinors in \E{1} are multivectors that can be written as the product of an even number of normalised proper points.
Furthermore, any spinor can be written as either \(S=e^{\alpha\e_{01}}=1+\alpha\e_{01}\) or \(-e^{\alpha\e_{01}}\).
Since \(\e_1 e^{x\e_{01}}=\e_1(1+x\e_{01})=\e_1-x\e_0\), I can replace \(\tb{a}=-x\e_0+\e_1\) with
\begin{equation}\label{point via spinor in E1}
\tb{a}=\e_1 e^{x\e_{01}},
\end{equation}
which is a useful representation as it helps simplify many algebraic manipulations.
Since spinors are even multivectors and each spinor satisfies \(S\reverse{S}=1\), I have 
\begin{equation}\label{e1 commute with S in E1}
\e_1S=S^{-1}\e_1.
\end{equation}
The translation by \(\lambda\), in the direction of positive \(x\) if \(\lambda>0\), 
of a point dually represented by \(\tb{a}\)
is given by \(T\tb{a}T^{-1}\),
where \(T=e^{-\tfrac{1}{2}\lambda\e_{01}}=1-\tfrac{1}{2}\lambda\e_{01}\).
\(T\) is a spinor and, therefore, \(T\tb{a}T^{-1}=\tb{a}T^{-2}=\tb{a}e^{\lambda\e_{01}}\),
so Equation~(\ref{point via spinor in E1}) can be interpreted as the translation of \(\e_1\), the origin of \E{1},
by \(x\).
There is only one kind of proper motion in \E{1}, translation,
and any translation can be generated by a suitable spinor.

The top-down reflection of \(\tb{a}\) in \(\tb{b}\) is given by \(-\tb{b}\tb{a}\tb{b}^{-1}\). It reverses the orientation of \(\tb{a}\).
The bottom-up reflection is given by \(\tb{b}\tb{a}\tb{b}^{-1}\). It preserves the orientation \(\tb{a}\).
Scaling \(\tb{a}\) by \(\gamma\) with respect to \(\tb{b}\) is 
the usual combination of projection and rejection,
\(\ts{scale}(\tb{a};\tb{b},\gamma)=(\tb{a}\cdot\tb{b})\tb{b}^{-1}+\gamma(\tb{a}\wedge\tb{b})\tb{b}^{-1}\).
Using (\ref{point via spinor in E1}) and (\ref{e1 commute with S in E1}), the following equalities can be readily obtained:
\begin{gather}
x[T\tb{a}T^{-1}]=x[\tb{a}]+\lambda,\\
x[-\tb{b}\tb{a}\tb{b}^{-1}]=x[\tb{b}]-(x[\tb{a}]-x[\tb{b}]),\\
x[\ts{scale}(\tb{a};\tb{b},\gamma)]
=x[\tb{b}]+\gamma(x[\tb{a}]-x[\tb{b}]),
\end{gather}
where \(x[d\e_0+a\e_1]=-d/a\) gives the point in \T{1} dually represented by \(\tb{a}=d\e_0+a\e_1\).

\newpage
\appendix
\small
\section{Practical Considerations}
In general, it is not practical to perform symbolic computations with Clifford algebra by hand.
To give an extreme example, the geometric product of two multivectors in \R{5*} includes 1024 terms.
It is too onerous and error prone to perform symbolic manipulations manually.
For this tutorial, I performed most symbolic computations with the Geometric algebra module for \verb|sympy|, 
which is described at\\
\hspace*{1cm}\verb|http://docs.sympy.org/dev/modules/galgebra/GA/GAsympy.html|\\
Its usage to compute the geometric product of two multivectors \verb|m1| and \verb|m2| in \R{3*} is illustrated below
\begin{verbatim}
from sympy import *
from sympy.galgebra.GA import *
set_main(sys.modules[__name__])
e0,e1,e2 = MV.setup('e0 e1 e2', '# 0 0, 0 # 0, 0 0 #')
MV.set_str_format(str_mode=2)
s1,d1,a1,b1,w1,x1,y1,p1 = symbols('s1 d1 a1 b1 w1 x1 y1 p1')
s2,d2,a2,b2,w2,x2,y2,p2 = symbols('s2 d2 a2 b2 w2 x2 y2 p2')

m1 = s1 + d1*e0+a1*e1+b1*e2 + w1*e1*e2+x1*e2*e0+y1*e0*e1 + p1*e0*e1*e2
m2 = s2 + d2*e0+a2*e1+b2*e2 + w2*e1*e2+x2*e2*e0+y2*e0*e1 + p2*e0*e1*e2
print m1*m2
\end{verbatim}
Computation with Clifford algebra is not routinely available in popular programming languages.
For numerical computations and visualisation, I made use of \verb|Asymptote| available at\\
\hspace*{1cm}\verb|http://asymptote.sourceforge.net/|\\
\verb|Asymptote| includes a high-level graphics and general programming language based on C++
and is designed to provide 2D and 3D plotting capabilities for LaTeX.

An elementary example of using \verb|Asymptote|
to compute the join of two points and display the resulting line
 is given below

\begin{verbatim}
struct MV {
  real s, d, a, b, w, x, y, p;
  void operator init(real s, real d, real a, real b, real w, real x, real y, real p) {
    this.s = s;
    this.d = d; this.a = a; this.b = b;
    this.w = w; this.x = x; this.y = y;
    this.p = p;
  }
}

void write(MV m) { unravel m; write(s,d,a,b,w,x,y,p); }

MV join(MV m1, MV m2) {
  from m1 unravel s as s1, d as d1, a as a1, b as b1, w as w1, x as x1, y as y1, p as p1;
  from m2 unravel s as s2, d as d2, a as a2, b as b2, w as w2, x as x2, y as y2, p as p2;
  MV m = MV(0,0,0,0,0,0,0,0);
  unravel m;
  s = -(-a1*x2-a2*x1-b1*y2-b2*y1-d1*w2-d2*w1-p1*s2-p2*s1);
  d = -(-d1*p2-d2*p1+x1*y2-x2*y1);
  a = -(-a1*p2-a2*p1-w1*y2+w2*y1);
  b = -(-b1*p2-b2*p1+w1*x2-w2*x1);
  w = p1*w2+p2*w1;
  x = p1*x2+p2*x1;
  y = p1*y2+p2*y1;
  p = p1*p2;
  return m;
}

MV P = MV(0,0,0,0,1,2,1,0);
MV Q = MV(0,0,0,0,1,-1,-1,0);
MV v = join(P,Q);
write(v);

import graph;
size(5cm,5cm);
axes("$x$", "$y$", min=-(5,5), max=(5,5));

unravel P; dot((x,y));
unravel Q; dot((x,y));

unravel v;
pair n = (b,-a);
pair centre = -(a,b)*d/(a^2+b^2);
path line = (centre-5*n/length(n))--(centre+5*n/length(n));
draw(line);
\end{verbatim}

\bibliographystyle{plain}
\bibliography{g.bib}

\end{document}